\documentclass[14pt]{article}
\usepackage{extsizes}
\usepackage{amsfonts,amssymb,amsthm}
\usepackage[intlimits]{amsmath}
\usepackage[T1]{fontenc}
\usepackage[utf8]{inputenc}
\usepackage[english, ukrainian]{babel}
\usepackage{tasks}
\usepackage{enumitem}
\usepackage{tikz}
\usepackage{bigints}
\usepackage{tocloft}
\usepackage{hyperref}
\hypersetup{
    colorlinks=true,
    linkcolor=black,
    filecolor=black,
    urlcolor=black,
    citecolor=blue
}
\usepackage[top=2cm,bottom=2cm,left=2cm,right=2cm]{geometry}

\newcommand{\Int}{\int\limits}

\DeclareMathOperator{\tgh}{th}

\settasks{
before-skip = {-3mm},
after-skip={-2mm},
after-item-skip={-1mm}
}

\begin{document}

\pagenumbering{gobble}
\begin{center}

\vspace{2cm}

{\large \textbf{I.P. Blazhievska, R. Riba Garcia}}

\vspace{7cm}

{\Large \textbf{INTEGRAL CALCULUS \\[0.3cm] OF ONE-DIMENSIONAL FUNCTIONS}}
\vspace{1cm}

{\Large \textbf{Personal tasks and Samples}}

\vspace{4cm}
Interactive textbook is developed for students\\
who study on  technical specialities

\vspace{7cm}
2020
\end{center}

\newpage

\parbox[b][4cm][t]{35mm}{\includegraphics[scale=0.3]{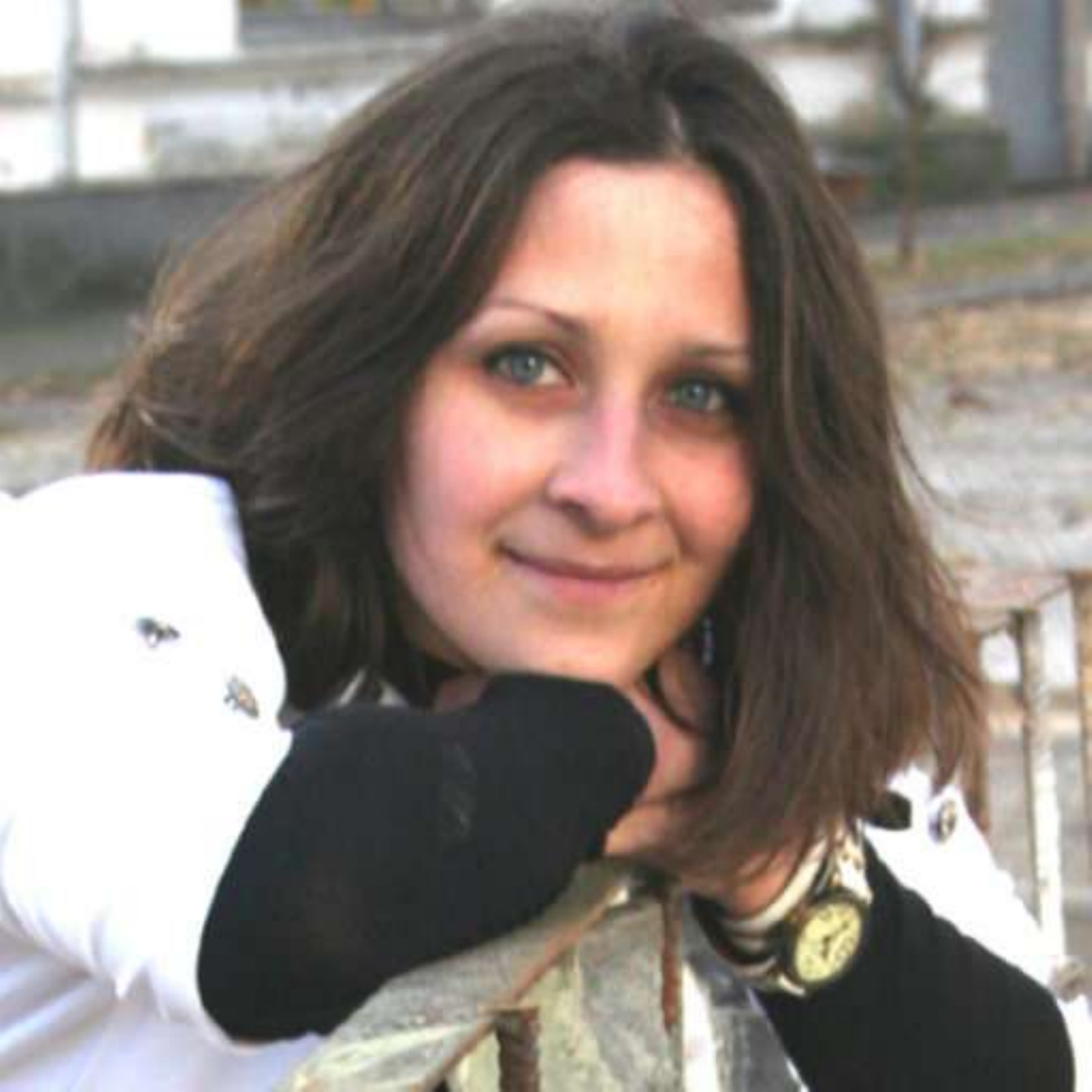}}
\hfill
\parbox[b][4cm][t]{120mm}{\vspace{4mm} \textit{Irina Blazhievska}, PhD in Physics and Mathematics, Senior Lecturer, Department of Mathematical Analysis and Probability Theory, NTUU ``Igor Sikorsky Kyiv Polytechnic Institute,'' Av. Peremogy, 37, 03056 Kyiv, Ukraine\\ e-mail: i.blazhievska@gmail.com}
\vspace{20mm}

\parbox[b][4cm][t]{35mm}{\includegraphics[scale=0.5]{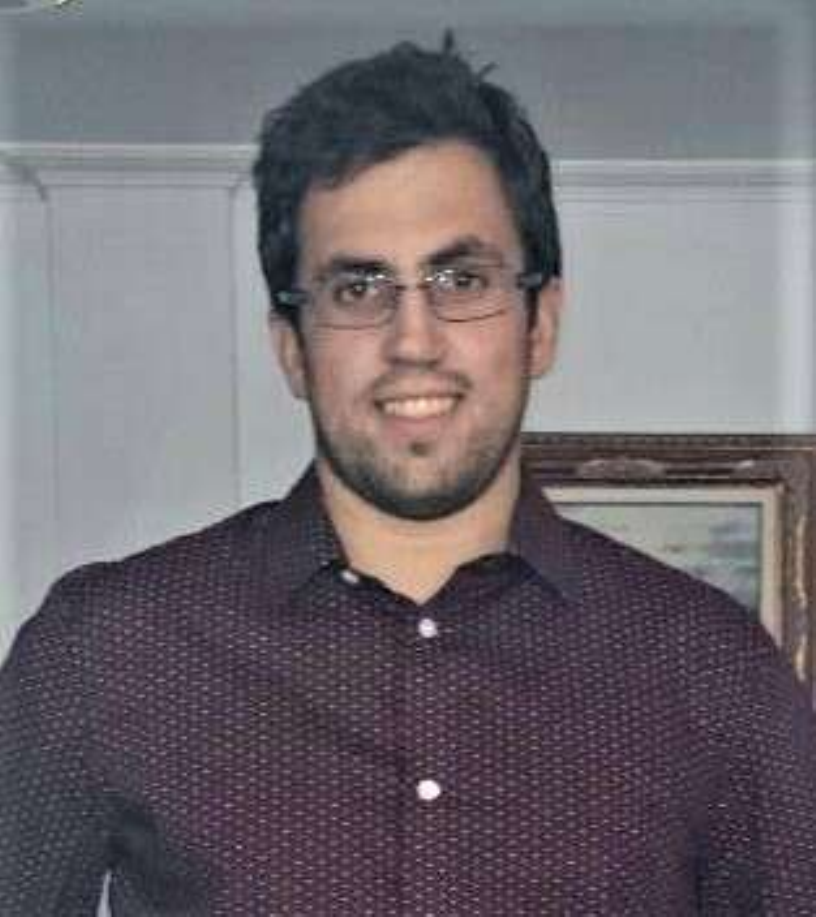}}
\hfill
\parbox[b][4cm][t]{120mm}{\vspace{4mm} \textit{Ricard Riba Garcia}, PhD in Mathematics, Associate Professor, Departament de Matem\`{a}tiques,
Universitat Aut\`{o}noma de Barcelona, Edifici Cc, 08193 Bellaterra (Barcelona), Spain\\ e-mail: ricard.riba@gmail.com}

\vspace{0.5cm}
\begin{center}
\section*{INTEGRAL CALCULUS \\ OF ONE-DIMENSIONAL FUNCTIONS\\
Personal tasks and Samples}
{\small{Electronic network tutorial edition}}
\end{center}

\vspace{1cm}
The interactive textbook is developed for English-speaking students whose study on Mathematical Calculus is based on classic programs of Ukrainian higher educational institutions. Its structure is as follows: the first part consists of 30 personal tasks with problems in one-dimensional integration (indefinite and definite integrals, geometric applications); the second part proposes the algorithmic solutions of the sample task with graphics built in Wolfram Mathematica 11.1, and 20 video-lessons guided by authors.
\vspace{2cm}
\begin{flushright}
$\copyright$ I.P. Blazhievska, R. Riba Garcia, 2020\\
\end{flushright}

\newpage

\pagenumbering{arabic}
\setcounter{page}{3}
\selectlanguage{english}
\tableofcontents

\newpage
\setcounter{secnumdepth}{0}

\section{Introduction}

The changes in higher education during the last decade have increased the amount of individual work of students. As result, there appears the necessity of clear-developed guidelines and resources for learning (on-line learning) of fundamental sciences.

This textbook is created for English-speaking students whose study in Mathematical Calculus is based on modern programs of Ukrainian higher educational institutions. It includes personal tasks and samples on ``Integral Calculus of One-Dimensional Functions'', one of the classic parts of Calculus course of the first year degree on technical specialities. All problems are the new ones; they were generated/tested by applying Wolfram Mathematica technologies. Guided by authors, 20 video-lessons are available for on-line learning on our educational YouTube channel.

Let us be more specific. The textbook contains 30 personal tasks and a solved sample task with guidelines (both written and interactive versions). The mathematical skeleton of each personal task consists of 4 parts:
 \begin{itemize}
   \item methods of the indefinite integration;
   \item evaluation of definite integrals;
   \item geometric applications of definite integrals, which are related to finding the metric characteristics of the curves, regions and solids of revolution;
   \item integration of improper integrals.
 \end{itemize}
 Authors illustrate on-line how to solve 20 problems from the sample task (an access from the textbook employs links and QR-codes).

In more detail, each personal task consists of 40 problems divided on 12 categories referring to fixed techniques (\cite[Chapter 6]{Efimov_Demidovich}, \cite[Chapters 7-8]{Pismennyi}, \cite[Chapters 8-10]{Fikhtengolts}, \cite[Chapter 4]{Sagan} and \cite[Part 3]{Spivak}):
\begin{enumerate}[label=\textbf{\arabic*.}]
  \item to be solved, the proposed 4 problems use the reduction to the table of integration combined with general properties of indefinite integrals, and the method of substitution under the differential (see, \cite[p. 267-273]{Efimov_Demidovich}, \cite[p. 193-202]{Pismennyi} and \cite[Chapter 13]{Spivak});
  \item to be solved, the proposed 4 problems use the integration of fractions with quadratic functions (see, \cite[p. 286-288]{Efimov_Demidovich} and \cite[p. 214-216]{Pismennyi});
  \item to be solved, the proposed 6 problems use the integration by parts (3 common classes) or suitable substitutions (see, \cite[p. 272-276]{Efimov_Demidovich} and \cite[p. 200-203]{Pismennyi});
  \item to be solved, the proposed 3 problems use the integration techniques for polynomial fractions (see, \cite[p. 276-281]{Efimov_Demidovich} and \cite[p. 203-211]{Pismennyi});
  \item to be solved, the proposed 3 problems use the integration of trigonometric functions: products and rational functions of sines and cosines with equal arguments, products of mentioned functions with different arguments (see, \cite[p.281-286]{Efimov_Demidovich}, \cite[p. 212-214]{Pismennyi} and \cite[Chapter 15]{Spivak});
  \item to be solved, the proposed 2 problems use the integration techniques for fractions with radicals by applying the trigonometric and hyperbolic substitutions (see, \cite[p. 286-289]{Efimov_Demidovich}, \cite[p. 214-219]{Pismennyi} and \cite[Chapters 15, 17]{Spivak});
  \item to be solved, the proposed 4 problems use the integration techniques for definite Riemann integrals with applying the Newton-Leibnitz formula. Here, the general evaluation, integration by parts and applying of suitable substitutions are covered (see, \cite[p. 290-300]{Efimov_Demidovich} and \cite[p. 221-233]{Pismennyi});
  \item to be solved, the proposed 3 problems require the building of correct regions in Cartesian, Parametric and Polar coordinates, applying of suitable formulas and valid evaluation of definite integrals (see, \cite[p. 306-311]{Efimov_Demidovich}, \cite[p. 237-242]{Pismennyi} and \cite[p. 233-257]{Spivak});
  \item to be solved, the proposed 3 problems require the building of correct arc segments in Cartesian, Parametric and Polar coordinates, applying of suitable formulas and valid evaluation of definite integrals (see, \cite[p. 311-314]{Efimov_Demidovich}, \cite[p. 242-245]{Pismennyi} and \cite[p. 257-261]{Spivak});
  \item to be solved, the proposed 3 problems require the building of correct generatrices in Cartesian, Parametric and Polar coordinates, analysis of surfaces of revolution, applying of suitable formulas and valid evaluation of definite integrals (see, \cite[p. 314-316]{Efimov_Demidovich}, \cite[p. 247-248]{Pismennyi} and \cite[p. 267]{Spivak});
  \item to be solved, the proposed 3 problems require the building of correct rotating regions in Cartesian and Polar coordinates, analysis of solids of revolution, applying of suitable formulas and valid evaluation of definite integrals (see, \cite[p. 317-319]{Efimov_Demidovich}, \cite[p. 245-247]{Pismennyi} and \cite[p. 265-267]{Spivak});
  \item to be solved, the proposed 2 problems use the integration techniques for improper Riemann integrals of the 1st and 2nd kinds (see, \cite[p. 300-305]{Efimov_Demidovich} and \cite[p. 233-237]{Pismennyi}).
\end{enumerate}

\noindent The following classic list of the curves are employed in the textbook:
   \begin{itemize}
     \item  in Cartesian coordinates: ellipses, exponentials, circles, catenaries, logarithmic functions, squared and cubic parabolas (power functions), straight lines and general curves;
     \item  in Parametric coordinates: astroids, ellipses, circles, cycloids, involutes of classic curves and spirals;
     \item  in Polar coordinates: Archimedean and logarithmic spirals, cardioids, lemniscates (8-shaped and $\infty$-shaped) and petaled roses.
   \end{itemize}
If a parametric range is omitted in the problem, it is supposed to be the whole region of the curve's well-definiteness.

In order to visualize the shapes of the curves, regions, surfaces and solids, we strongly recommend students to use the following links:
 \begin{itemize}
     \item https://www.mathcurve.com;
     \item http://mathworld.wolfram.com;
     \item http://old.nationalcurvebank.org/volrev/volrev.htm
 \end{itemize}
The suitable formulas for evaluation of metric characteristics are available as well.

\vspace{0.5cm}
The sample task consists on algorithmic solutions of the proposed 40 problems. Here, we point out the applied techniques for indefinite and definite integrals (categories 1--7), develop  step-by-step algorithms with 2D/3D-simulations via Wolfram Mathematica 11.1 for geometric problems (categories 8-11) and solve the problems with improper integrals (category 12). Based on problems from the sample task, we created  20 illustrative video-lessons on techniques of integration. These lessons were prepared with KPI TV$^\copyright$ Team and nowadays are available on our educational YouTube channel:
\begin{itemize}
     \item \url{https://www.youtube.com/channel/UCHrLMGhP7cM664BCF_e32Qg}
\end{itemize}
The list of problems covered in video-lessons is attached below (see, pages 70-73); the fast access is possible to achieve within links and QR-codes.

This interactive textbook finishes with references on the best books in our framework complemented by useful external links.
\vspace{0.5cm}

To remove possible misunderstandings, the next page exhibits the differences between some notations of Post-Soviet countries and the foreign ones.

\newpage
\section{On notions used in the textbook}
Following by mathematical traditions stated in Ukraine (Post-Soviet countries), we have used the following notions for some trigonometric, hyperbolic functions and their inverses:
$$\arraycolsep=0.4cm
\begin{array}{llll}
\bullet \ \textrm{tg}\ x & \text{instead of} & \tan{x}, & \textrm{tg}\ x=\dfrac{\sin x}{\cos x}; \\[0.7cm]
\bullet \ \textrm{ctg}\ x & \text{instead of} & \cot{x}, & \textrm{ctg}\ x=\dfrac{\cos x}{\sin x}; \\[0.7cm]
\bullet \ \textrm{arctg}\ x & \text{instead of} & \arctan{x}; \\[0.7cm]
\bullet \ \textrm{arcctg}\ x & \text{instead of} & \textrm{arccot}\ x, & \textrm{arcctg}\ x=\dfrac{\pi}{2}-\textrm{arctg}\ x; \\[0.7cm]
\bullet \ \textrm{sh}\ x & \text{instead of} & \sinh{x}, & \textrm{sh}\ x=\dfrac{e^x-e^{-x}}{2}; \\[0.7cm]
\bullet \ \textrm{ch}\ x & \text{instead of} & \cosh{x}, & \textrm{ch}\ x=\dfrac{e^x+e^{-x}}{2}; \\[0.7cm]
\bullet \ \textrm{th}\ x & \text{instead of} & \tanh{x}, & \textrm{th}\ x=\dfrac{\textrm{sh}\ x}{\textrm{ch}\ x}; \\[0.7cm]
\bullet \ \textrm{cth}\ x & \text{instead of} & \textrm{coth}\ {x}, & \textrm{cth}\ x=\dfrac{\textrm{ch}\ x}{\textrm{sh}\ x}; \\[0.7cm]
\bullet \ \textrm{arcsh}\ x & \text{instead of} & \textrm{arcsinh}\ {x}, & \textrm{arcsh}\ x=\ln\left(x+\sqrt{x^2+1}\right); \\[0.7cm]
\bullet \ \textrm{arcch}\ x & \text{instead of} & \textrm{arccosh}\ {x}, & \textrm{arcch}\ x=\ln\left(x+\sqrt{x^2-1}\right), \ |x|\geq1; \\[0.7cm]
\bullet \ \textrm{arcth}\ x & \text{instead of} & \textrm{arctanh}\ {x}, & \textrm{arcth}\ x=\dfrac{1}{2}\ln\left(\dfrac{1+x}{1-x}\right), \ |x|<1; \\[0.7cm]
\bullet \ \textrm{arccth}\ x & \text{instead of} & \textrm{arccoth}\ {x}, & \textrm{arccth}\ x=\dfrac{1}{2}\ln\left(\dfrac{x+1}{x-1}\right), \ |x|>1;
\end{array}$$
\vspace{0.6cm}

For more information about hyperbolic functions, inverse hyperbolic functions and their transforms see \cite[p. 285-289]{Efimov_Demidovich} and \cite[Chapter 17]{Spivak}.

\newpage

\section{Personal task 1}
   \begin{minipage}[t]{90mm}
\textbf{1.}\hspace{1mm}Integrate using the table and \\ substitution under differential:
\begin{tasks}[ label-align=left, label-offset={1mm}, label-width={3mm}, item-indent={0mm} ](2)
\task $\Int x(3x^2-2)^3 dx$;
\task $\Int \textrm{tg}^2x dx$;
\task $\Int \sin xe^{2\cos x-5}dx$;
\task $\Int\cos^2(1-x)dx$.
\end{tasks}
\textbf{2.}\hspace{1mm}Integrate the quadratic fractions:
\begin{tasks}[ label-align=left, label-offset={1mm}, label-width={3mm},column-sep={-40pt}, item-indent={0mm} , after-item-skip={0mm}](2)
\task $\int\dfrac{2-3x}{x^2-4} dx$;
\task $\Int\dfrac{dx}{\sqrt{20+24x-9x^2}}$;
\task $\Int\dfrac{(x+4)dx}{x^2+6x+5}$
\task $\Int\dfrac{(2x-5)dx}{\sqrt{x^2-2x-15}}$.
\end{tasks}
\textbf{3.}\hspace{1mm}Integrate by parts or using the \\ suitable substitutions:
\begin{tasks}[ label-align=left, label-offset={1mm}, label-width={3mm}, column-sep={-30pt}, item-indent={0mm} ,after-item-skip={-1mm}](2)
\task $\Int\dfrac{dx}{\sqrt{e^x+1}}$;
\task $\Int\dfrac{dx}{x\sqrt{x^2-1}}$;
\task $\Int\ln(x-3)dx$;
\task $\Int (x^2+5) \cos2x dx$;
\task $\Int\arcsin2x dx$;
\task $\Int\dfrac{\ln(\cos x)}{\cos^2x}dx$.
\end{tasks}
\textbf{4.}\hspace{1mm}Integrate the polynomial fractions:
\begin{tasks}[ label-align=left, label-offset={1mm}, label-width={3mm}, item-indent={0mm} ,after-item-skip={-1mm}](1)
\task $\Int\dfrac{3x^2+14x+19}{(x^2+4x+3)(x+5)} dx$;
\task $\Int\dfrac{x^3+1}{x^3-2x^2+x}dx$;
\task $\Int\dfrac{10x+6}{(x^2+2x+5)(x-1)}dx$.
\end{tasks}
\textbf{5.}\hspace{1mm}Integrate trigonometric expressions:
\begin{tasks}[ label-align=left, label-offset={1mm}, label-width={3mm}, item-indent={0mm} ,column-sep={-30pt}, after-item-skip={-1mm}](2)
\task $\Int\sin3x\cos xdx$;
\task $\Int\cos^4 3x\sin^2 3x dx$;
\task $\Int\dfrac{dx}{5+2\sin x+3\cos x}$.
\end{tasks}
\textbf{6.}\hspace{1mm}Integrate the fractions with radicals:
\begin{tasks}[ label-align=left, label-offset={1mm}, label-width={3mm}, item-indent={0mm},column-sep={-50pt} ](1)
\task $\Int\dfrac{\sqrt{9-x^2}}{x^4}dx$;
\task $\Int\dfrac{(1-\sqrt{x+1})dx}{(1+\sqrt[3]{x+1})\sqrt{x+1}}$.
\end{tasks}
\end{minipage}
\vline\;
      \begin{minipage}[t]{90mm}

\textbf{7.}\hspace{1mm}Solve the definite integrals:
\begin{tasks}[ label-align=left, label-offset={1mm}, label-width={3mm}, item-indent={0mm},after-item-skip=-1mm,  , column-sep={-20pt}](2)
\task $\Int_{2}^{e+1}x\ln(x-1)dx$;
\task $\Int_{0}^{\pi/4}\dfrac{\sqrt{\textrm{tg}^3x}}{\cos^2x}dx$;
\task $\Int_{0}^{\sqrt{3}} x\sqrt{1+x^2}dx$;
\task $\Int_{3}^{29}\dfrac{\sqrt[3]{(x-2)^2} \ dx}{3+\sqrt[3]{(x-2)^2}}.$
\end{tasks}
\textbf{8.}\hspace{1mm}Find the area of the figure bounded \\ by the curves:
\begin{tasks}[ label-align=left, label-offset={2mm}, label-width={3mm}, item-indent={0mm}, after-item-skip=-1mm,     ](1)
\task $y=2x^2-8x+6,\ y=x^2-3x$;
\task $\left\lbrace\begin{array}{l} x=4\cos^3t, \\ y=4\sin^3t;\end{array}\right. \ x=2 \ (x\geq 2)$;
\task $\rho=6\cos3\phi.$
\end{tasks}
\textbf{9.}\hspace{1mm}Find the arc-length of the curve:
\begin{tasks}[ label-align=left, label-offset={2mm}, label-width={3mm}, item-indent={0mm}, after-item-skip=-1mm,     ](1)
\task $y=\ln x,\ \sqrt{3}\leq x\leq \sqrt{8}$;
\task $\left\lbrace\begin{array}{l}
x=2(t-\sin t), \\ y=2(1-\cos t);
\end{array}\right. \ 0\leq t\leq2\pi$;
\task $\rho=e^{-\phi}, \ 0\leq\phi\leq\pi$.
\end{tasks}
\textbf{10.}\hspace{1mm}Find the area of the surface formed \\ by rotating the curves around the $l$-axis:
\begin{tasks}[ label-align=left, label-offset={2mm}, label-width={3mm}, item-indent={0mm}, after-item-skip=-1mm,    ](1)
\task $y=\frac{1}{3}x^3, \ -1\leq x\leq 1, \ l=OX$;
\task $\left\lbrace\begin{array}{l} x=1+2\cos t, \\ y=3+2\sin t;\end{array}\right.  \  (x\geq 0) \ l=OY$;
\task $\rho=2\sqrt{\sin2\big(\phi-\frac{\pi}{4}\big)}, \ l=o\rho$.
\end{tasks}
\textbf{11.}\hspace{1mm}Find the volume of the body formed \\ by rotating the curves around the $l$-axis:
\begin{tasks}[ label-align=left, label-offset={2mm}, label-width={3mm}, item-indent={0mm}, after-item-skip=-1mm,    ](1)
\task $y=x^2-2x+1,\ y=x+1, \ l=OX$;
\task $x=-y^2+5y-6,\ x=0, \ l=OY$;
\task $\rho=2(1+\cos\phi), \ l=o\rho$.
\end{tasks}
\textbf{12.}\hspace{1mm}Solve the improper integrals:
\begin{tasks}[ label-align=left, label-offset={1mm}, label-width={3mm}, item-indent={0mm},   ,   ,column-sep={-40pt}](2)
\task $\Int_{0}^{\infty}\dfrac{x}{x^4+16}dx$;
\task $\Int_{0}^{1/2}\dfrac{dx}{\sqrt[3]{2-4x}}$.
\end{tasks}
\end{minipage}

\section{Personal task 2}

\begin{minipage}[t]{90mm}
\textbf{1.}\hspace{1mm}Integrate using the table and \\ substitution under differential:
\begin{tasks}[ label-align=left, label-offset={1mm}, label-width={3mm}, column-sep={-25pt}, item-indent={0mm} , ](2)
\task $\Int \dfrac{(3-\sqrt[3]{x})^2 dx}{\sqrt[3]{x}}$;
\task $\Int \dfrac{(\cos 2x)dx}{\cos^2x \sin^2x}$;
\task $\Int x3^{5x^2-2}dx$;
\task $\Int\sin^2(1-2x)dx$.
\end{tasks}
\textbf{2.}\hspace{1mm}Integrate the quadratic fractions:
\begin{tasks}[ label-align=left, label-offset={1mm}, label-width={3mm},column-sep={-18pt}, item-indent={0mm},   ,, ](2)
\task $\Int\dfrac{4x+5}{x^2-1} dx$;
\task $\Int\dfrac{dx}{\sqrt{5-8x-4x^2}}$;
\task $\Int\dfrac{(2x-5)dx}{x^2+6x+13}$;
\task $\Int\dfrac{(4x+1)dx}{\sqrt{x^2+4x-12}}$.
\end{tasks}
\textbf{3.}\hspace{1mm}Integrate by parts or using the \\ suitable substitutions:
\begin{tasks}[ label-align=left, label-offset={1mm}, label-width={3mm}, column-sep={-40pt}, item-indent={0mm},  ](2)
\task $\Int\dfrac{(4x+3)dx}{(x-2)^3}$;
\task $\Int\dfrac{dx}{x\sqrt{12x^2+4x+1}}$;
\task! $\Int (1-5x+x^2) \sin 3x dx$;
\task $\Int x^3\ln xdx$;
\task $\Int\textrm{arctg} x\ dx$;
\task $\Int\cos (\ln x)dx$.
\end{tasks}
\textbf{4.}\hspace{1mm}Integrate the polynomial fractions:
\begin{tasks}[ label-align=left, label-offset={1mm}, label-width={3mm},after-item-skip=-0.5mm, item-indent={0mm} ,, ](1)
\task $\Int\dfrac{15x^2+15x-54}{(x^2+x-2)(x-2)} dx$;
\task $\Int\dfrac{x^3-2x^2-2x+1}{x^3-x^2} dx$;
\task $\Int\dfrac{2x^2+4x-26}{(x^2-4x+8)(x-1)}dx$.
\end{tasks}
\textbf{5.}\hspace{1mm}Integrate trigonometric expressions:
\begin{tasks}[ label-align=left, label-offset={1mm}, label-width={3mm}, column-sep={-20pt}, item-indent={0mm},  ](2)
\task $\Int\sin3x\sin 9xdx$;
\task $\Int\cos^3 x\sqrt[7]{\sin^4 x} dx$;
\task $\Int\dfrac{dx}{4+5\sin^2 x-3\cos^2 x}$.
\end{tasks}
\textbf{6.}\hspace{1mm}Integrate the functions with radicals:
\begin{tasks}[ label-align=left, label-offset={1mm}, label-width={3mm}, item-indent={0mm},column-sep={-30pt} ](2)
\task $\Int x^4\sqrt{4-x^2}dx$;
\task $\Int\sqrt{\dfrac{x+4}{x-4}} \dfrac{dx}{(x+4)^2}$.
\end{tasks}
\end{minipage}
\vline\;
      \begin{minipage}[t]{100mm}
\textbf{7.}\hspace{1mm}Solve the definite integrals:
\begin{tasks}[ label-align=left, label-offset={1mm}, label-width={3mm},column-sep={-60pt}, item-indent={0mm},after-item-skip=-1mm,  ](2)
\task $\Int_{0}^{1}(x-1)^2 e^xdx$;
\task $\Int_{0}^{2\pi}\cos^3(\frac{x}{4})\sin^3(\frac{x}{4})dx$;
\task $\Int_{0}^{\sqrt[6]{2e-1}}\dfrac{12x^5}{x^6+1}dx$;
\task $\Int_{0}^{\ln2}\dfrac{dx}{e^x(3+e^{-x})}.$
\end{tasks}

\textbf{8.}\hspace{1mm}Find the area of the figure bounded \\
by the curves:
\begin{tasks}[ label-align=left, label-offset={2mm}, label-width={3mm}, item-indent={0mm}, after-item-skip=-1mm,     ](1)
\task $y=(x-2)^2,\ y=4x-8$;
\task $\left\lbrace\begin{array}{l} x=3\cos t, \\ y=4\sin t; \end{array}\right. \ y=2 \ (y\geq2)$;
\task  $\rho=\sin 2\phi.$
\end{tasks}

\textbf{9.}\hspace{1mm}Find the arc-length of the curve:
\begin{tasks}[ label-align=left, label-offset={2mm}, label-width={3mm}, item-indent={0mm}, after-item-skip=-1mm,     ](1)
\task  $y=\sqrt{1-x^2},\ -\frac{\sqrt{2}}{2}\leq x\leq 1$;
\task  $\left\lbrace\begin{array}{l} x=4e^{t}\cos t, \\ y=4e^{t}\sin t; \end{array}\right. \ 0\leq t\leq2\pi$;
\task  $\rho=2(1-\cos\phi)$.
\end{tasks}

\textbf{10.}\hspace{1mm}Find the area of the surface formed \\ by rotating the curves around the $l$-axis:
\begin{tasks}[ label-align=left, label-offset={2mm}, label-width={3mm}, item-indent={0mm}, after-item-skip=-1mm,    ](1)
\task $y=e^{x}, \ 0\leq x\leq \frac{1}{2}\ln8, \ l=OX$;
\task $\left\lbrace\begin{array}{l} x=2(t-\sin t), \\ y=2(1-\cos t);\end{array}\right. \begin{array}{l} 0\leq t\leq2\pi, \\ l=OX;\end{array}$
\task $\rho=3\sqrt{\cos2\phi}, \ l=o\rho$.
\end{tasks}

\textbf{11.}\hspace{1mm}Find the volume of the body formed \\
by rotating the curves around the $l$-axis:
\begin{tasks}[ label-align=left, label-offset={2mm}, label-width={3mm}, item-indent={0mm}, after-item-skip=-1mm,     ](1)
\task  $y=x^2+2x+1,\ y=1-x, \ l=OX$;
\task  $x=-y^2+2y,\ x=4y-2y^2, \ l=OY$;
\task  $\rho=4\sin\phi, \ l=o\rho$.
\end{tasks}

\textbf{12.}\hspace{1mm}Solve the improper integrals:
\begin{tasks}[ label-align=left, label-offset={1mm}, label-width={3mm}, item-indent={0mm},   ,   ,column-sep={-40pt}](2)
\task $\Int_{1}^{\infty}\dfrac{16x}{16x^4-1}dx$;
\task $\Int_{\pi/4}^{\pi/2}\dfrac{\sin x}{\cos^3x}dx$.
\end{tasks}
\end{minipage}

\section{Personal task 3}

   \begin{minipage}[t]{90mm}
\textbf{1.}\hspace{1mm}Integrate using the table and \\ substitution under differential:
\begin{tasks}[ label-align=left, label-offset={1mm}, label-width={3mm}, column-sep={5pt}, item-indent={0mm} ](2)
\task $\Int \sqrt{x}(1-2x^4)^2 dx$;
\task $\Int \textrm{ctg}^2x dx$;
\task $\Int \dfrac{\ln^3(2x-5)}{2x-5}dx$;
\task $\Int\cos^2(x+3)dx$.
\end{tasks}
\textbf{2.}\hspace{1mm}Integrate the quadratic fractions:
\begin{tasks}[ label-align=left, label-offset={1mm}, label-width={3mm},column-sep={-45pt}, item-indent={0mm} ,, ](2)
\task $\Int\dfrac{6x-5}{x^2-9} dx$;
\task $\Int\dfrac{dx}{\sqrt{-55+16x-x^2}}$;
\task $\Int\dfrac{(3x+4)dx}{x^2+2x+5}$;
\task $\Int\dfrac{(2x+3)dx}{\sqrt{x^2+2x-15}}$.
\end{tasks}
\textbf{3.}\hspace{1mm}Integrate by parts or using the \\ suitable substitutions:
\begin{tasks}[ label-align=left, label-offset={1mm}, label-width={3mm}, column-sep={-20pt}, item-indent={0mm},  ](2)
\task $\Int\dfrac{e^x}{\sqrt{e^{2x}+1}}dx$;
\task $\Int\dfrac{\sqrt{4+x^2}}{x}dx$;
\task $\Int 8x^2 \sin 2x dx$;
\task $\Int \dfrac{3x^2+6x-3}{e^{3x}}dx$;
\task $\Int\arccos 4x dx$;
\task $\Int\dfrac{x}{\cos^2x}dx$.
\end{tasks}
\textbf{4.}\hspace{1mm}Integrate the polynomial fractions:
\begin{tasks}[ label-align=left, label-offset={1mm}, label-width={3mm},after-item-skip=-0.5mm, item-indent={0mm} ,, ](1)
\task $\Int\dfrac{4x^2-10x}{(x^2-4x+3)(x-2)} dx$;
\task  $\Int\dfrac{2x^3+3x^2-x+12}{x^3+4x^2} dx$;
\task  $\Int\dfrac{-3x^2+24x-63}{(x^2-4x+13)(x+1)}dx$.
\end{tasks}
\textbf{5.}\hspace{1mm}Integrate trigonometric expressions:
\begin{tasks}[ label-align=left, label-offset={1mm}, label-width={3mm},after-item-skip=-1mm, item-indent={0mm} , ](2)
\task $\Int\cos7x\cos3xdx$;
\task $\Int\dfrac{\sin^3 x}{\cos^7 x} dx$;
\task $\Int\dfrac{\cos x dx}{(\sin x+1)^2-\cos^2x}$.
\end{tasks}
\textbf{6.}\hspace{1mm}Integrate the fractions with radicals:
\begin{tasks}[ label-align=left, label-offset={1mm}, label-width={3mm}, item-indent={0mm},column-sep={-1in} ](2)
\task $\Int\dfrac{dx}{x^2\sqrt{x^2-1}}$;
\task $\Int\dfrac{\sqrt[3]{(x+1)^2}+\sqrt[6]{x+1}}{\sqrt{x+1}+\sqrt[3]{(x+1)^2}} dx$.
\end{tasks}
\end{minipage}
\vline\;
      \begin{minipage}[t]{100mm}
\textbf{7.}\hspace{1mm}Solve the definite integrals:
\begin{tasks}[ label-align=left, label-offset={1mm}, label-width={3mm},column-sep={-40pt}, item-indent={0mm},after-item-skip=-1mm,  ](2)
\task  $\Int_{0}^{2e-1}\dfrac{\ln^2(x+1)}{x+1}dx$;
\task $\Int_{0}^{\pi} \sin^8x dx$;
\task $\Int_{0}^{\sqrt{3}} \dfrac{x^2}{1+x^2}dx$;
\task $\Int_{0}^{5}\dfrac{dx}{2x+\sqrt{3x+1}}.$
\end{tasks}

\textbf{8.}\hspace{1mm}Find the area of the figure bounded \\
by the curves:
\begin{tasks}[ label-align=left, label-offset={2mm}, label-width={3mm}, item-indent={0mm}, after-item-skip=-1mm,     ](1)
\task $x+y=7,\ xy=6$;
\task $\left\lbrace\begin{array}{l} x=2\cos^3t,\\ y=4\sin^3t; \end{array}\right.\ y=0 \;(y\geq 0)$;
\task $\rho=\sqrt{\sin2\phi}.$
\end{tasks}

\textbf{9.}\hspace{1mm}Find the arc-length of the curve:
\begin{tasks}[ label-align=left, label-offset={2mm}, label-width={3mm}, item-indent={0mm}, after-item-skip=-1mm,     ](1)
\task $y=\dfrac{1}{3}x^3,\ -\frac{1}{2}\leq x\leq \frac{1}{2}$;
\task $\left\lbrace\begin{array}{l} x=4(t-\sin t),\\ y=4(1-\cos t); \end{array}\right. \ 0\leq t\leq\pi$;
\task $\rho=e^{\frac{3\phi}{4}}, \ 0\leq\phi\leq 2\pi$.
\end{tasks}

\textbf{10.}\hspace{1mm}Find the area of the surface formed \\ by rotating the curves around the $l$-axis:
\begin{tasks}[ label-align=left, label-offset={2mm}, label-width={3mm}, item-indent={0mm}, after-item-skip=-1mm,    ](1)
\task $y=\frac{1}{4}\textrm{ch}4x, \ 0\leq x\leq \ln8, \ l=OX$;
\task $\left\lbrace\begin{array}{l} x=2+\cos t,\\ y=\sin t; \end{array}\right. \ l=OY$;
\task $\rho=6\cos\phi, \ l=o\rho$.
\end{tasks}

\textbf{11.}\hspace{1mm}Find the volume of the body formed \\
by rotating the curves around the $l$-axis:
\begin{tasks}[ label-align=left, label-offset={2mm}, label-width={3mm}, item-indent={0mm}, after-item-skip=-1mm,     ](1)
\task $y=x^2-2x+5,\ y=x+5, \ l=OX$;
\task $y=\ln x, \ x=0, \ y=0, \ y=\ln4, \ l=OY$;
\task $\rho=3(1-\cos\phi), \ l=o\rho$.
\end{tasks}
\textbf{12.}\hspace{1mm}Solve the improper integrals:
\begin{tasks}[ label-align=left, label-offset={1mm}, label-width={3mm}, item-indent={0mm},   ,   ](2)
\task $\Int_{0}^{\infty}\dfrac{x^6}{\sqrt[5]{(x^7+1)^3}}dx$;
\task $\Int_{0}^{\pi/2}\dfrac{\cos^2x}{\sin^4x}dx$.
\end{tasks}
\end{minipage}

\section{Personal task 4}

   \begin{minipage}[t]{90mm}
\textbf{1.}\hspace{1mm}Integrate using the table and \\ substitution under differential:
\begin{tasks}[ label-align=left, label-offset={1mm}, label-width={3mm}, column-sep={5pt}, item-indent={0mm} , ,after-item-skip=-1mm](2)
\task  $\Int \sqrt{x}(1-x^2)^3 dx$;
\task $\Int \dfrac{(\sin 2x)dx}{(\cos x-1)^2-1}$;
\task $\Int x^3e^{2x^4-1}dx$;
\task $\Int\textrm{th}^2xdx$.
\end{tasks}
\textbf{2.}\hspace{1mm}Integrate the quadratic fractions:
\begin{tasks}[ label-align=left, label-offset={1mm}, label-width={3mm},column-sep={-20pt}, item-indent={0mm} ,, ](2)
\task $\Int\dfrac{-2x-5}{x^2-4} dx$;
\task $\Int\dfrac{dx}{\sqrt{12+12x-9x^2}}$;
\task $\Int\dfrac{(-2x+3)dx}{x^2-6x+10}$;
\task $\Int\dfrac{(4x-5)dx}{\sqrt{x^2+2x-8}}$.
\end{tasks}
\textbf{3.}\hspace{1mm}Integrate by parts or using the \\ suitable substitutions:
\begin{tasks}[ label-align=left, label-offset={1mm}, label-width={3mm}, column-sep={-35pt}, item-indent={0mm},  ,after-item-skip=-1mm,](2)
\task $\Int\dfrac{2x+5}{(x+3)^7}dx$;
\task $\Int \dfrac{dx}{(x+1)^2\sqrt{-2x-x^2}}$;
\task! $\Int x^2(\cos x+\sin x) dx$;
\task $\Int \dfrac{\ln(x+1)}{x+1}dx$;
\task $\Int\sin (\ln x)dx$;
\task $\Int\textrm{arctg} 2x\ dx$.
\end{tasks}
\textbf{4.}\hspace{1mm}Integrate the polynomial fractions:
\begin{tasks}[ label-align=left, label-offset={1mm}, label-width={3mm},after-item-skip=-0.5mm, item-indent={0mm} , ](1)
\task $\Int\dfrac{3x^2+42}{(x^2+5x+4)(x-2)} dx$;
\task $\Int\dfrac{4x^3-2x^2-9}{x^3+2x^2+x} dx$;
\task $\Int\dfrac{5x^2+22x+27}{(x^2+4x+5)(x+2)}dx$.
\end{tasks}
\textbf{5.}\hspace{1mm}Integrate trigonometric expressions:
\begin{tasks}[ label-align=left, label-offset={1mm}, label-width={3mm},after-item-skip=-1mm, item-indent={0mm} ,column-sep={5pt}, ](2)
\task $\Int\cos4x\sin 7xdx$;
\task $\Int\dfrac{\cos^5 x}{\sin^7 x} dx$;
\task $\Int\dfrac{dx}{5-4\sin x+2\cos x}$.
\end{tasks}
\textbf{6.}\hspace{1mm}Integrate the fractions with radicals:
\begin{tasks}[ label-align=left, label-offset={1mm}, label-width={3mm}, item-indent={0mm},column-sep={-30pt} ,, ](2)
\task $\Int\dfrac{\sqrt{9-x^2}}{x^2}dx$;
\task $\Int\dfrac{\sqrt[4]{x}+\sqrt{x}}{1+\sqrt{x}} dx$.
\end{tasks}
\end{minipage}
\vline\;
      \begin{minipage}[t]{100mm}
\textbf{7.}\hspace{1mm}Solve the definite integrals:
\begin{tasks}[ label-align=left, label-offset={1mm}, label-width={3mm},column-sep={-40pt}, item-indent={0mm},after-item-skip=-1mm,  ](2)
\task $\Int_{1/\ln5}^{1}\dfrac{e^{-\frac{1}{x}}}{x^2}dx$;
\task $\Int_{0}^{\pi/8}\cos^4 2x\sin^4 2xdx$;
\task $\Int_{0}^{1}\dfrac{12x+5}{x^2-4}dx$;
\task $\Int_{-2}^{2}x^3\sqrt{4-x^2}dx.$
\end{tasks}
\textbf{8.}\hspace{1mm}Find the area of the figure bounded \\
by the curves:
\begin{tasks}[ label-align=left, label-offset={2mm}, label-width={3mm}, item-indent={0mm}, after-item-skip=-1mm,     ](1)
\task $y=(x+3)^2,\ y=2x+6$;
\task  $\left\lbrace\begin{array}{l} x=4\cos t,\\ y=4\sin t; \end{array}\right. \ x=2 \ (x\geq2)$;
\task  $\rho=\sin 3\phi.$
\end{tasks}

\textbf{9.}\hspace{1mm}Find the arc-length of the curve:
\begin{tasks}[ label-align=left, label-offset={2mm}, label-width={3mm}, item-indent={0mm}, after-item-skip=-1mm,     ](1)
\task $y=\ln(\sin x),\ \frac{\pi}{3}\leq x\leq \frac{\pi}{2}$;
\task $\left\lbrace\begin{array}{l} x=2\cos t-\cos2t,\\ y=2\sin t-\sin2t; \end{array}\right. \ 0\leq t \leq2\pi$;
\task $\rho=\phi, \ 0\leq \phi\leq\pi$.
\end{tasks}

\textbf{10.}\hspace{1mm}Find the area of the surface formed \\ by rotating the curves around the $l$-axis:
\begin{tasks}[ label-align=left, label-offset={2mm}, label-width={3mm}, item-indent={0mm}, after-item-skip=-1mm,    ](1)
\task $y=\sqrt{2x+6}, \ 0\leq x\leq 3, \ l=OX$;
\task $\left\lbrace\begin{array}{l} x=\cos^3 t, \\ y=\sin^3 t; \end{array}\right. \ l=OY$;
\task $\rho=3\sqrt{\cos2\phi}, \ l=o\rho$.
\end{tasks}

\textbf{11.}\hspace{1mm}Find the volume of the body formed \\
by rotating the curves around the $l$-axis:
\begin{tasks}[ label-align=left, label-offset={2mm}, label-width={3mm}, item-indent={0mm}, after-item-skip=-1mm,     ](1)
\task $y=x^2-4x+4,\ y=x, \ l=OX$;
\task $x=9-(y-2)^2,\ x=0, \ l=OY$;
\task $\rho=6\sin\phi, \ l=o\rho$.
\end{tasks}

\textbf{12.}\hspace{1mm}Solve the improper integrals:
\begin{tasks}[ label-align=left, label-offset={1mm}, label-width={3mm}, item-indent={0mm},   ,   ](2)
\task  $\Int_{2}^{\infty}\dfrac{4x-2}{4x^2-4x-3}dx$;
\task $\Int_{0}^{1}\dfrac{dx}{\sqrt{3-2x-x^2}}$.
\end{tasks}
\end{minipage}

\section{Personal task 5}

 \begin{minipage}[t]{90mm}
\textbf{1.}\hspace{1mm}Integrate using the table and \\ substitution under differential:
\begin{tasks}[ label-align=left, label-offset={1mm}, label-width={3mm}, column-sep={-18pt}, item-indent={0mm} , ,after-item-skip=-1mm](2)
\task! $\Int \sqrt{x}(3+2x\sqrt{x})^2 dx$;
\task $\Int \dfrac{dx}{\cos2x+\sin^2x}$;
\task $\Int \dfrac{x^3 dx}{\sqrt[3]{x^4+9}}$;
\task $\Int(1-\sin\frac{x}{2})^2dx$.
\end{tasks}
\textbf{2.}\hspace{1mm}Integrate the quadratic fractions:
\begin{tasks}[ label-align=left, label-offset={1mm}, label-width={3mm},column-sep={-20pt}, item-indent={0mm} ,, ](2)
\task $\Int\dfrac{-4x+7}{x^2-16} dx$;
\task $\Int\dfrac{dx}{\sqrt{15+6x-9x^2}}$;
\task $\Int\dfrac{(-2x-5)dx}{x^2-8x+20}$;
\task $\Int\dfrac{(6x-3)dx}{\sqrt{x^2+4x-12}} $.
\end{tasks}
\textbf{3.}\hspace{1mm}Integrate by parts or using the \\ suitable substitutions:
\begin{tasks}[ label-align=left, label-offset={1mm}, label-width={3mm}, column-sep={-20pt}, item-indent={0mm},  ,after-item-skip=-1mm,](2)
\task $\Int\dfrac{dx}{\sqrt{9+e^x}}$;
\task $\Int\dfrac{dx}{x^2\sqrt{x^2+4}}$;
\task $\Int \ln^2 x dx$;
\task $\Int (1-6x^2) \cos 3x dx$;
\task $\Int e^{\sqrt{x+3}}dx$;
\task $\Int x\cdot\textrm{arctg}x dx$.
\end{tasks}
\textbf{4.}\hspace{1mm}Integrate the polynomial fractions:
\begin{tasks}[ label-align=left, label-offset={1mm}, label-width={3mm},after-item-skip=-0.5mm, item-indent={0mm} , ](1)
\task $\Int\dfrac{-5x^2-4x+9}{(x^2+4x+3)(x+2)} dx$;
\task $\Int\dfrac{x^3+3x^2+3x+4}{x^3+4x^2+4x} dx$;
\task $\Int\dfrac{5x^2-24x+24}{(x^2-6x+10)(x-1)}dx$.
\end{tasks}
\textbf{5.}\hspace{1mm}Integrate trigonometric expressions:
\begin{tasks}[ label-align=left, label-offset={1mm}, label-width={3mm},after-item-skip=-1mm, item-indent={0mm} ,, ](2)
\task $\Int\sin2x\sin 4xdx$;
\task $\Int\sqrt[3]{\dfrac{\cos^2 x}{\sin^8 x}} dx$;
\task $\Int\dfrac{\cos x dx}{6-\sin^2 x+3\cos^2 x}$.
\end{tasks}
\textbf{6.}\hspace{1mm}Integrate the functions with radicals:
\begin{tasks}[ label-align=left, label-offset={1mm}, label-width={3mm}, item-indent={0mm},column-sep={-20pt} ,, ](2)
\task $\Int x^2\sqrt{1-x^2}dx$;
\task $\Int\sqrt[3]{\dfrac{x+1}{x-1}}\dfrac{dx}{(x-1)^3}$.
\end{tasks}
\end{minipage}
\vline\;
      \begin{minipage}[t]{100mm}
\textbf{7.}\hspace{1mm}Solve the definite integrals:
\begin{tasks}[ label-align=left, label-offset={1mm}, label-width={3mm},column-sep={-30pt}, item-indent={0mm},after-item-skip=-1mm,  ](2)
\task $\Int_{0}^{2}(2-x)^2 e^{2x}dx$;
\task $\Int_{-\pi/4}^{\pi/4}\dfrac{\sin^4x}{\cos^6x}dx$;
\task $\Int_{-1/2}^{1}\dfrac{dx}{\sqrt{8+2x-x^2}}$;
\task $\Int_{0}^{7}\dfrac{\sqrt[3]{(x+1)^2} \ dx}{4+\sqrt[3]{(x+1)^2}}.$\end{tasks}

\textbf{8.}\hspace{1mm}Find the area of the figure bounded \\
by the curves:
\begin{tasks}[ label-align=left, label-offset={2mm}, label-width={3mm}, item-indent={0mm}, after-item-skip=-1mm,     ](1)
\task $x+y=-3,\ xy=2$;
\task $\left\lbrace\begin{array}{l} x=5\cos t,\\ y=2\sin t; \end{array}\right. \ y=0\  (y\geq0)$;
\task $\rho=4(1+\sin\phi)$.
\end{tasks}

\textbf{9.}\hspace{1mm}Find the arc-length of the curve:
\begin{tasks}[ label-align=left, label-offset={2mm}, label-width={3mm}, item-indent={0mm}, after-item-skip=-1mm,     ](1)
\task $y=3+\textrm{ch}x,\ 0\leq x\leq \ln2$;
\task $\left\lbrace\begin{array}{l} x=4\cos^3 t,\\ y=4\sin^3 t; \end{array}\right. \ 0\leq t\leq\frac{\pi}{2}$;
\task $\rho=\sqrt{2}e^{\phi}, \ 0\leq\phi\leq \frac{\pi}{2}$.
\end{tasks}

\textbf{10.}\hspace{1mm}Find the area of the surface formed \\ by rotating the curves around the $l$-axis:
\begin{tasks}[ label-align=left, label-offset={2mm}, label-width={3mm}, item-indent={0mm}, after-item-skip=-1mm,    ](1)
\task $y=\frac{1}{2}x^2, \ 0\leq x\leq 2\sqrt{6}, \ l=OY$;
\task $\left\lbrace\begin{array}{l} x=3(t-\sin t), \\ y=3(1-\cos t); \end{array}\right. \begin{array}{l}
 t\leq2\pi, \\ l=OX
\end{array}$;
\task $\rho=8\sin\phi, \ l=o\rho$.
\end{tasks}

\textbf{11.}\hspace{1mm}Find the volume of the body formed \\
by rotating the curves around the $l$-axis:
\begin{tasks}[ label-align=left, label-offset={2mm}, label-width={3mm}, item-indent={0mm}, after-item-skip=-1mm,     ](1)
\task $y=x^2-2x+2,\ y=x+2, \ l=OX$;
\task $x=y^2+4,\ x=0, \ |y|=2, \ l=OY$;
\task $\rho=2\sin2\phi, \ l=o\rho$.
\end{tasks}

\textbf{12.}\hspace{1mm}Solve the improper integrals:
\begin{tasks}[ label-align=left, label-offset={1mm}, label-width={3mm}, item-indent={0mm},   ,   ,column-sep={-40pt}](2)
\task $\Int_{1/\sqrt[3]{3}}^{\infty}\dfrac{4x^2}{9x^6+1}dx$;
\task  $\Int_{1}^{e}\dfrac{dx}{x\sqrt{1-\ln^2x}}$.
\end{tasks}
\end{minipage}

\section{Personal task 6}

 \begin{minipage}[t]{90mm}
\textbf{1.}\hspace{1mm}Integrate using the table and \\ substitution under differential:
\begin{tasks}[ label-align=left, label-offset={1mm}, label-width={3mm}, column-sep={0pt}, item-indent={0mm} , ,after-item-skip=-1mm](2)
\task! $\Int x(6\sqrt[4]{x}-\sqrt{x})^2 dx$
\task $\Int \dfrac{(\sin2x)dx}{(2-\sin x)^2-4} $;
\task $\Int \dfrac{\ln^7(3x+1)}{3x+1}dx$;
\task $\Int\cos^2(2x-5)dx$.
\end{tasks}
\textbf{2.}\hspace{1mm}Integrate the quadratic fractions:
\begin{tasks}[ label-align=left, label-offset={1mm}, label-width={3mm},column-sep={-30pt}, item-indent={0mm} ,, ](2)
\task $\Int\dfrac{12x+7}{x^2-4} dx$;
\task  $\Int\dfrac{dx}{\sqrt{4x^2-20x+24}}$;
\task  $\Int\dfrac{(-3x+1)dx}{x^2+6x+13}$;
\task  $\Int\dfrac{-4x+3}{\sqrt{-2x-x^2}} dx$.
\end{tasks}
\textbf{3.}\hspace{1mm}Integrate by parts or using the \\ suitable substitutions:
\begin{tasks}[ label-align=left, label-offset={1mm}, label-width={3mm}, column-sep={5pt}, item-indent={0mm},  ,after-item-skip=-1mm,](2)
\task  $\Int\dfrac{e^{4x} dx}{\sqrt{e^{8x}+16}}$;
\task $\Int\dfrac{dx}{x^2\sqrt{3-4x+x^2}}$;
\task $\Int (3-4x) \sin 4x dx$;
\task $\Int \dfrac{5x^2-1}{e^{5x}}dx$;
\task $\Int(\arcsin x)^2 dx$;
\task $\Int\dfrac{\ln x}{x^3}dx$.
\end{tasks}
\textbf{4.}\hspace{1mm}Integrate the polynomial fractions:
\begin{tasks}[ label-align=left, label-offset={1mm}, label-width={3mm},after-item-skip=-0.5mm, item-indent={0mm} , ](1)
\task $\Int\dfrac{x^2+14}{(x^2-5x+4)(x-2)} dx$;
\task$\Int\dfrac{3x^3+7x^2+9x+9}{x^3+3x^2} dx$;
\task $\Int\dfrac{9x^2-24x+46}{(x^2-6x+10)(x+2)}dx$.
\end{tasks}
\textbf{5.}\hspace{1mm}Integrate trigonometric fractions:
\begin{tasks}[ label-align=left, label-offset={1mm}, label-width={3mm},after-item-skip=-1mm, item-indent={0mm} ,, ](2)
\task $\Int\sin7x\cos2xdx$
\task $\Int\dfrac{\sin^2 2x}{\cos^4 2x} dx$;
\task $\Int\dfrac{5\sin x dx}{14-12\cos x-9\sin^2x}$.
\end{tasks}
\textbf{6.}\hspace{1mm}Integrate the fractions with radicals:
\begin{tasks}[ label-align=left, label-offset={1mm}, label-width={3mm}, item-indent={0mm},column-sep={-30pt}  ](2)
\task $\Int\dfrac{\sqrt{x^2-4}}{x^5}dx$;
\task $\Int\sqrt{\dfrac{x+2}{x-2}}\dfrac{dx}{(x+2)^2}$.
\end{tasks}
\end{minipage}
\vline\;
   \begin{minipage}[t]{100mm}
\textbf{7.}\hspace{1mm}Solve the definite integrals:
\begin{tasks}[ label-align=left, label-offset={1mm}, label-width={3mm},column-sep={-40pt}, item-indent={0mm},after-item-skip=-1mm,  ](2)
\task $\Int_{3}^{e+2}\dfrac{\ln^5(x-2)}{x-2}dx$;
\task $\Int_{-\pi/2}^{\pi/2} \cos^3x\sin^8x dx$;
\task $\Int_{3}^{8} \dfrac{x^2}{4-x^2}dx$;
\task $\Int_{-2}^{6}\dfrac{\sqrt{2x+4}}{4+\sqrt{2x+4}}dx.$
\end{tasks}

\textbf{8.}\hspace{1mm}Find the area of the figure bounded \\
by the curves:
\begin{tasks}[ label-align=left, label-offset={2mm}, label-width={3mm}, item-indent={0mm}, after-item-skip=-1mm,     ](1)
\task $x+2y=5,\ xy=2$;
\task  $\left\lbrace\begin{array}{l} x=3\cos^3t,\\ y=3\sin^3t; \end{array}\right. \ x=0 \ (x\geq 0)$;
\task  $\rho=2\cos4\phi.$
\end{tasks}

\textbf{9.}\hspace{1mm}Find the arc-length of the curve:
\begin{tasks}[ label-align=left, label-offset={2mm}, label-width={3mm}, item-indent={0mm}, after-item-skip=-1mm,     ](1)
\task $y=\ln(x+\sqrt{x^2-4}),\ 2\leq x\leq 2\sqrt{5}$;
\task $\left\lbrace\begin{array}{l} x=2(t-\sin t),\\ y=2(1-\cos t); \end{array}\right. \ 0\leq t\leq 2\pi$;
\task $\rho=6\cos\phi-6\sin\phi$.
\end{tasks}

\textbf{10.}\hspace{1mm}Find the area of the surface formed \\ by rotating the curves around the $l$-axis:
\begin{tasks}[ label-align=left, label-offset={2mm}, label-width={3mm}, item-indent={0mm}, after-item-skip=-1mm,    ](1)
\task $y=1+\frac{1}{3}\textrm{ch}3x, \ -1\leq x\leq 1, \ l=OX$;
\task $\left\lbrace\begin{array}{l}x=1+2\cos t,\\ y=2\sin t; \end{array}\right. \ (x\geq 0), \ l=OY$;
\task $\rho=1+\cos\phi, \ l=o\rho$.
\end{tasks}

\textbf{11.}\hspace{1mm}Find the volume of the body formed \\
by rotating the curves around the $l$-axis:
\begin{tasks}[ label-align=left, label-offset={2mm}, label-width={3mm}, item-indent={0mm}, after-item-skip=-1mm,     ](1)
\task $y=x^2-2x+3,\ y=x+3, \ l=OX$;
\task  $y=\sqrt{1-x^2}, \ y=x, \ x=0, \ l=OY$;
\task  $\rho=4\phi, \ 0\leq\phi\leq\pi, \ l=o\rho$.
\end{tasks}

\textbf{12.}\hspace{1mm}Solve the improper integrals:
\begin{tasks}[ label-align=left, label-offset={1mm}, label-width={3mm}, item-indent={0mm},   ,   ,column-sep={-40pt}](2)
\task  $\Int_{0}^{\infty}\dfrac{x^2}{x^2+1}dx$;
\task $\Int_{0}^{\pi/2}e^{-\textrm{ctg}x}\dfrac{dx}{\sin^2x}$.
\end{tasks}
\end{minipage}

\section{Personal task 7}

   \begin{minipage}[t]{90mm}
\textbf{1.}\hspace{1mm}Integrate using the table and \\ substitution under differential:
\begin{tasks}[ label-align=left, label-offset={1mm}, label-width={3mm}, column-sep={-10pt}, item-indent={0mm} , ,after-item-skip=-1mm](2)
\task $\Int \dfrac{(x-1)^2}{\sqrt[3]{x}} dx$;
\task $\Int\dfrac{\textrm{tg}^3x}{\cos^2x}dx$;
\task $\Int x\sqrt{x^2-4}dx$;
\task $\Int 2\sin^2\Big(\frac{x}{2}\Big)dx$;
\end{tasks}
\textbf{2.}\hspace{1mm}Integrate the quadratic fractions:
\begin{tasks}[ label-align=left, label-offset={1mm}, label-width={3mm},column-sep={-20pt}, item-indent={0mm} ,, ](2)
\task $\Int\dfrac{2x-7}{x^2-9} dx$;
\task $\Int\dfrac{dx}{\sqrt{15-6x-9x^2}}$;
\task $\Int\dfrac{(-2x+9)dx}{x^2+4x+20}$;
\task $\Int\dfrac{(6x+3)dx}{\sqrt{x^2-2x-8}}$.
\end{tasks}
\textbf{3.}\hspace{1mm}Integrate by parts or using the \\ suitable substitutions:
\begin{tasks}[ label-align=left, label-offset={1mm}, label-width={3mm}, column-sep={0pt}, item-indent={0mm},  ,after-item-skip=-1mm,](2)
\task  $\Int\dfrac{dx}{x\sqrt{x+4}}$;
\task  $\Int\dfrac{\sqrt{1+\ln x}}{x\ln x}dx$;
\task!  $\Int (2-3x^2) \sin 2x dx$;
\task  $\Int x4^xdx$;
\task  $\Int\textrm{arcctg}\sqrt{x} dx$;
\task  $\Int\ln(x^2+1)dx$.
\end{tasks}
\textbf{4.}\hspace{1mm}Integrate the polynomial fractions:
\begin{tasks}[ label-align=left, label-offset={1mm}, label-width={3mm},after-item-skip=-0.5mm, item-indent={0mm} , ](1)
\task $\Int\dfrac{3x^2+5x-6}{(x^2+4x+3)(x+2)} dx$;
\task $\Int\dfrac{2x^3+11x^2+17x+9}{x^3+5x^2+8x+4} dx$;
\task $\Int\dfrac{x^2+x-2}{x^3+2x^2+2x}dx$.
\end{tasks}
\textbf{5.}\hspace{1mm}Integrate trigonometric expressions:
\begin{tasks}[ label-align=left, label-offset={1mm}, label-width={3mm},after-item-skip=-1mm, item-indent={0mm} ,, ](2)
\task! $\Int\cos x\sin x\cos3x dx$;
\task $\Int\cos^5 x\sqrt[3]{\sin x} dx$;
\task $\Int\dfrac{dx}{5-3\cos x}$.
\end{tasks}
\textbf{6.}\hspace{1mm}Integrate the functions with radicals:
\begin{tasks}[ label-align=left, label-offset={1mm}, label-width={3mm}, item-indent={0mm},column-sep={0in} ,, ](2)
\task $\Int\sqrt{1-4x-x^2}dx$;
\task $\Int\dfrac{dx}{x(\sqrt{x}+\sqrt[3]{x^2})}$.
\end{tasks}
\end{minipage}
\vline\;
      \begin{minipage}[t]{100mm}
\textbf{7.}\hspace{1mm}Solve the definite integrals:
\begin{tasks}[ label-align=left, label-offset={1mm}, label-width={3mm},column-sep={-40pt}, item-indent={0mm},after-item-skip=-1mm,  ](2)
\task $\Int_{-3}^{\ln2}(x+3)^2 e^xdx$;
\task $\Int_{0}^{2\pi}\cos^2(\frac{x}{2})\sin^4(\frac{x}{2})dx$;
\task $\Int_{0}^{1}\dfrac{4x^3}{x^8+1}dx$;
\task $\Int_{0}^{-\ln2}\sqrt{1-e^{2x}}dx.$
\end{tasks}

\textbf{8.}\hspace{1mm}Find the area of the figure bounded \\
by the curves:
\begin{tasks}[ label-align=left, label-offset={2mm}, label-width={3mm}, item-indent={0mm}, after-item-skip=-1mm,     ](1)
\task $y=(x+1)^2,\ y=6x+6$;
\task  $\left\lbrace\begin{array}{l} x=12\cos t,\\ y=5\sin t;\end{array}\right.
 \ x=6 \ (x\geq6)$;
\task  $\rho=2\sqrt{\cos2\phi}.$
\end{tasks}

\textbf{9.}\hspace{1mm}Find the arc-length of the curve:
\begin{tasks}[ label-align=left, label-offset={2mm}, label-width={3mm}, item-indent={0mm}, after-item-skip=-1mm,     ](1)
\task  $y=\frac{1}{2}\ln(\cos2x),\ 0\leq x\leq \frac{\pi}{12}$;
\task  $\left\lbrace\begin{array}{l} x=4\cos t,\\ y=4\sin t; \end{array}\right. \ y=2 \ (y\geq 2)$;
\task  $\rho=6\sin^2\Big(\frac{\phi}{2}\Big)$.
\end{tasks}

\textbf{10.}\hspace{1mm}Find the area of the surface formed \\ by rotating the curves around the $l$-axis:
\begin{tasks}[ label-align=left, label-offset={2mm}, label-width={3mm}, item-indent={0mm}, after-item-skip=-1mm,    ](1)
\task $y=\sqrt{3}x+7, \ 0\leq x\leq 2\sqrt{3}, \ l=OX$;
\task $\left\lbrace\begin{array}{l} x=3\cos^3 t, \\ y=3\sin^3 t; \end{array}\right. \ l=OY$;
\task $\rho=2e^{-\phi}, 0\leq\phi\leq\pi, \ l=o\rho$.
\end{tasks}

\textbf{11.}\hspace{1mm}Find the volume of the body formed \\
by rotating the curves around the $l$-axis:
\begin{tasks}[ label-align=left, label-offset={2mm}, label-width={3mm}, item-indent={0mm}, after-item-skip=-1mm,     ](1)
\task $y=(x-2)^2,\ y=x+4, \ l=OX$;
\task $y=9-x^2,\ y=9-3x^2, y=0, \ l=OY$;
\task $\rho=4\sin2\phi, \ l=o\rho$.
\end{tasks}

\textbf{12.}\hspace{1mm}Solve the improper integrals:
\begin{tasks}[ label-align=left, label-offset={1mm}, label-width={3mm}, item-indent={0mm},   ,   ,column-sep={-40pt}](2)
\task $\Int_{1}^{\infty}\dfrac{dx}{x^2(x+1)}$;
\task  $\Int_{0}^{2/\pi}\sin\left(\dfrac{1}{x}\right)\dfrac{dx}{x^2}$.
\end{tasks}
\end{minipage}

\section{Personal task 8}

   \begin{minipage}[t]{90mm}
\textbf{1.}\hspace{1mm}Integrate using the table and \\ substitution under differential:
\begin{tasks}[ label-align=left, label-offset={1mm}, label-width={3mm}, column-sep={10pt}, item-indent={0mm} , ,after-item-skip=-1mm](2)
\task! $\Int\dfrac{x^2-2x^3e^x+4}{x^3} dx$;
\task $\Int \dfrac{4+\cos^2x}{1+\cos2x} dx$;
\task $\Int \dfrac{\ln^2(3-2x)}{3-2x}dx$;
\task $\Int(1-\sin3x)^2 dx$.
\end{tasks}

\textbf{2.}\hspace{1mm}Integrate the quadratic fractions:
\begin{tasks}[ label-align=left, label-offset={1mm}, label-width={3mm},column-sep={-20pt}, item-indent={0mm} ,, ](2)
\task $\Int\dfrac{-4x+3}{x^2-16} dx$;
\task $\Int\dfrac{dx}{\sqrt{4x^2-12x+8}}$;
\task $\Int\dfrac{(-4x+1)dx}{x^2-2x+17}$;
\task $\Int\dfrac{(2x+3)dx}{\sqrt{15-2x-x^2}}$.
\end{tasks}
\textbf{3.}\hspace{1mm}Integrate by parts or using the \\ suitable substitutions:
\begin{tasks}[ label-align=left, label-offset={1mm}, label-width={3mm}, column-sep={10pt}, item-indent={0mm},  ,after-item-skip=-1mm,](2)
\task $\Int\dfrac{\ln x}{x\sqrt[3]{1-\ln^2x}}dx$;
\task $\Int\dfrac{\sqrt{x^2+1}}{x}dx$;
\task $\Int (3-9x) \sin3xdx$;
\task $\Int (1-8x^2)e^{4x}dx$;
\task $\Int\dfrac{\textrm{arctg}x}{x^2}dx$;
\task $\Int\dfrac{x}{\sin^2x}dx$.
\end{tasks}
\textbf{4.}\hspace{1mm}Integrate the polynomial fractions:
\begin{tasks}[ label-align=left, label-offset={1mm}, label-width={3mm},after-item-skip=-0.5mm, item-indent={0mm} , ](1)
\task $\Int\dfrac{4x^2-8x-20}{(x^2+6x+5)(x-1)} dx$;
\task $\Int\dfrac{x^3+3x^2+4x-1}{x^3+5x^2+8x+4} dx$;
\task $\Int\dfrac{-7x^2-13x-10}{(x^2+4x+5)(x-1)}dx$.
\end{tasks}
\textbf{5.}\hspace{1mm}Integrate trigonometric expressions:
\begin{tasks}[ label-align=left, label-offset={1mm}, label-width={3mm},after-item-skip=-1mm, item-indent={0mm} ,, ](2)
\task $\Int\sin7x\sin4xdx$;
\task  $\Int\cos^2 x\sin^4x dx$;
\task  $\Int\dfrac{2-\cos x}{2+\cos x}dx$.
\end{tasks}
\textbf{6.}\hspace{1mm}Integrate the fractions with radicals:
\begin{tasks}[ label-align=left, label-offset={1mm}, label-width={3mm}, item-indent={0mm},column-sep={0in} , ](2)
\task $\Int\dfrac{dx}{x\sqrt{1+4x-5x^2}}$;
\task $\Int\dfrac{x dx}{(\sqrt{x}+2\sqrt[4]{x})^2}$.
\end{tasks}
\end{minipage}
\vline\;
      \begin{minipage}[t]{100mm}
\textbf{7.}\hspace{1mm}Solve the definite integrals:
\begin{tasks}[ label-align=left, label-offset={1mm}, label-width={3mm},column-sep={20mm}, item-indent={0mm},after-item-skip=-1mm, ,column-sep={0pt} ](2)
\task  $\Int_{0}^{e-1}(x+1)\ln(x+1)dx$;
\task  $\Int_{0}^{\pi/8} \cos^4 4x\sin^3 4x dx$;
\task  $\Int_{4}^{9} \dfrac{x-1}{\sqrt{x}+1}dx$;
\task  $\Int_{0}^{\pi/4} e^{\textrm{tg}x}\dfrac{dx}{\cos^2x}.$
\end{tasks}

\textbf{8.}\hspace{1mm}Find the area of the figure bounded \\
by the curves:
\begin{tasks}[ label-align=left, label-offset={2mm}, label-width={3mm}, item-indent={0mm}, after-item-skip=-1mm,     ](1)
\task $y=(x-1)^3,\ y=4(x-1)$;
\task $\left\lbrace\begin{array}{l} x=3\cos^3t,\\ y=\sin^3t; \end{array}\right. \ y=0 \ (y\geq 0)$;
\task $\rho=2\sin3\phi.$
\end{tasks}

\textbf{9.}\hspace{1mm}Find the arc-length of the curve:
\begin{tasks}[ label-align=left, label-offset={2mm}, label-width={3mm}, item-indent={0mm}, after-item-skip=-1mm,     ](1)
\task  $y=4-\frac{1}{2}\textrm{ch} 2x,\ 0\leq x\leq \ln4$;
\task   $\left\lbrace\begin{array}{l} x=2(t-\sin t),\\ y=2(1-\cos t); \end{array}\right. \ 0\leq t \leq2\pi$;
\task   $\rho=\sqrt{8}e^{-3\phi}, \ 0\leq\phi\leq\pi$.
\end{tasks}

\textbf{10.}\hspace{1mm}Find the area of the surface formed \\ by rotating the curves around the $l$-axis:
\begin{tasks}[ label-align=left, label-offset={2mm}, label-width={3mm}, item-indent={0mm}, after-item-skip=-1mm,    ](1)
\task $y=2x^2, \ 0\leq x\leq \sqrt{3}, \ l=OY$;
\task $\left\lbrace\begin{array}{l} x=1+3\cos t, \\ y=3\sin t; \end{array}\right. \ l=OX$;
\task $\rho=4\sqrt{\cos2\phi}, \ l=o\rho$.
\end{tasks}

\textbf{11.}\hspace{1mm}Find the volume of the body formed \\
by rotating the curves around the $l$-axis:
\begin{tasks}[ label-align=left, label-offset={2mm}, label-width={3mm}, item-indent={0mm}, after-item-skip=-1mm,     ](1)
\task $y=-x^2+4x,\ y=0, \ l=OX$;
\task $y=x^2, \ y=5x^2-16, \ l=OY$;
\task $\rho=2(1-\cos\phi), \ l=o\rho$.
\end{tasks}

\textbf{12.}\hspace{1mm}Solve the improper integrals:
\begin{tasks}[ label-align=left, label-offset={1mm}, label-width={3mm}, item-indent={0mm},   ,   ,column-sep={-40pt}](2)
\task $\Int_{-1}^{\infty}\dfrac{\textrm{arcctg}^5x}{1+x^2}dx$;
\task $\Int_{1}^{2\sqrt{3}}\dfrac{2x+2}{\sqrt{x^2+2x-3}}dx$.
\end{tasks}
\end{minipage}

\section{Personal task 9}

   \begin{minipage}[t]{90mm}
\textbf{1.}\hspace{1mm}Integrate using the table and \\ substitution under differential:
\begin{tasks}[ label-align=left, label-offset={1mm}, label-width={3mm}, column-sep={10pt}, item-indent={0mm} , ,after-item-skip=-1mm](2)
\task! $\Int x^2(1-3\sqrt{x})^3 dx$;
\task $\Int \dfrac{\textrm{tg}^7x}{1-\cos^2x} dx$;
\task $\Int\dfrac{\textrm{arcsin}^3x+x}{\sqrt{1-x^2}}dx$;
\task! $\Int (\sin x+\cos x)^2 dx.$
\end{tasks}
\textbf{2.}\hspace{1mm}Integrate the quadratic fractions:
\begin{tasks}[ label-align=left, label-offset={1mm}, label-width={3mm},column-sep={-20pt}, item-indent={0mm} ,, ](2)
\task $\Int\dfrac{-2x+4}{x^2-25} dx$;
\task $\Int\dfrac{dx}{\sqrt{3-4x-4x^2}}$;
\task $\Int\dfrac{(-2x-7)dx}{x^2+6x+13}$;
\task $\Int\dfrac{(6x+8)dx}{\sqrt{x^2+4x-5}}$.
\end{tasks}

\textbf{3.}\hspace{1mm}Integrate by parts or using the \\ suitable substitutions:
\begin{tasks}[ label-align=left, label-offset={1mm}, label-width={3mm}, column-sep={0pt}, item-indent={0mm},  ,after-item-skip=-1mm,](2)
\task $\Int\dfrac{4x-7}{(x-1)^5}dx$;
\task $\Int\dfrac{dx}{(x-1)\sqrt{2x-x^2}}$;
\task! $\Int x(\sin x-2\cos x) dx$;
\task $\Int \dfrac{x^2-4x+1}{e^{2x}}dx$;
\task $\Int\dfrac{\ln(x-2)}{(x-2)^2}dx$;
\task $\Int(\textrm{arccos}x)^2 dx$.
\end{tasks}

\textbf{4.}\hspace{1mm}Integrate the polynomial fractions:
\begin{tasks}[ label-align=left, label-offset={1mm}, label-width={3mm},after-item-skip=-0.5mm, item-indent={0mm} , ](1)
\task $\Int\dfrac{7x^2+14x-85}{(x^2-6x+5)(x+3)} dx$;
\task $\Int\dfrac{3x^3+15x^2+24x+10}{x^3+5x^2+7x+3} dx$;
\task $\Int\dfrac{4x^2-9x+30}{(x^2-4x+20)(x+2)}dx$.
\end{tasks}

\textbf{5.}\hspace{1mm}Integrate trigonometric expressions:
\begin{tasks}[ label-align=left, label-offset={1mm}, label-width={3mm},after-item-skip=-1mm, item-indent={0mm} ,, ](2)
\task $\Int\cos x\sin 8xdx$;
\task $\Int\dfrac{\cos^3 x}{\sin^6 x} dx$;
\task $\Int\dfrac{2dx}{3+\sin x+5\cos x}$.
\end{tasks}
\textbf{6.}\hspace{1mm}Integrate the fractions with radicals:
\begin{tasks}[ label-align=left, label-offset={1mm}, label-width={3mm}, item-indent={0mm},column-sep={0in} , ](2)
\task $\Int\dfrac{\sqrt{1-x^2}}{x^4}dx$;
\task $\Int\sqrt{\dfrac{1-x}{1+x}}\dfrac{dx}{x}$.
\end{tasks}
\end{minipage}
\vline \;
      \begin{minipage}[t]{100mm}
\textbf{7.}\hspace{1mm}Solve the definite integrals:
\begin{tasks}[ label-align=left, label-offset={1mm}, label-width={3mm},column-sep={-30pt}, item-indent={0mm},after-item-skip=-1mm](2)
\task  $\Int_{6/\pi}^{2/\pi}\cos\Big(\dfrac{1}{x}\Big)\dfrac{dx}{x^2}$;
\task  $\Int_{0}^{\pi/2}x\cos^4x\sin xdx$;
\task  $\Int_{0}^{\sqrt{3}}\dfrac{2x-2}{x^2+9}dx$;
\task  $\Int_{0}^{1}x^2\sqrt{1-x^2}dx.$
\end{tasks}

\textbf{8.}\hspace{1mm}Find the area of the figure bounded\\
by the curves:
\begin{tasks}[ label-align=left, label-offset={2mm}, label-width={3mm}, item-indent={0mm}, after-item-skip=-1mm,     ](1)
\task $y=-x^2+5x+1,\ y=\dfrac{5}{x}$;
\task $\left\lbrace\begin{array}{l} x=1+2\cos t,\\ y=2\sin t;\end{array}\right. \ y=1 \ (y\geq1)$;
\task $\rho=\cos 2\phi.$
\end{tasks}

\textbf{9.}\hspace{1mm}Find the arc-length of the curve:
\begin{tasks}[ label-align=left, label-offset={2mm}, label-width={3mm}, item-indent={0mm}, after-item-skip=-1mm,     ](1)
\task $y=3\ln(9-x^2),\ 0\leq x\leq2$;
\task $\left\lbrace\begin{array}{l} x=3e^t\cos t,\\ y=3e^t\sin t; \end{array}\right. \ 0\leq t\leq\pi$;
\task $\rho=4\cos^2\Big(\dfrac{\phi}{2}\Big)$.
\end{tasks}

\textbf{10.}\hspace{1mm}Find the area of the surface formed \\ by rotating the curves around the $l$-axis:
\begin{tasks}[ label-align=left, label-offset={2mm}, label-width={3mm}, item-indent={0mm}, after-item-skip=-1mm,    ](1)
\task $y=3-\frac{1}{2}\textrm{ch}2x, \ -1\leq x\leq 1, \ l=OX$;
\task $\left\lbrace\begin{array}{l} x=2\cos^3 t, \\ y=2\sin^3 t;\end{array}\right. \ l=OY$;
\task $\rho=6\sin\phi, \ l=o\rho$.
\end{tasks}

\textbf{11.}\hspace{1mm}Find the volume of the body formed \\
by rotating the curves around the $l$-axis:
\begin{tasks}[ label-align=left, label-offset={2mm}, label-width={3mm}, item-indent={0mm}, after-item-skip=-1mm,     ](1)
\task  $y=x^2+6x+9,\ y=-4x, \ l=OX$;
\task  $y=x^2-8,\ y=7x,\ x=0 \ (x\geq0),\ \\ l=OY$;
\task  $\rho=2\phi,\ 0\leq\phi\leq\pi,\ l=o\rho$.
\end{tasks}

\textbf{12.}\hspace{1mm}Solve the improper integrals:
\begin{tasks}[ label-align=left, label-offset={1mm}, label-width={3mm}, item-indent={0mm},   ,   ,column-sep={-40pt}](2)
\task  $\Int_{0}^{\infty}\dfrac{x^2}{e^{x^3}}dx$;
\task $\Int_{0}^{1}\dfrac{1+x-x^2}{\sqrt{1-x^2}}dx$.
\end{tasks}

\end{minipage}

\section{Personal task 10}

   \begin{minipage}[t]{90mm}
\textbf{1.}\hspace{1mm}Integrate using the table and \\ substitution under differential:
\begin{tasks}[ label-align=left, label-offset={1mm}, label-width={3mm}, column-sep={0pt}, item-indent={0mm} , ,after-item-skip=-1mm](2)
\task! $\Int x\sqrt{x}(1-3x\sqrt{x})^2 dx$;
\task $\Int \dfrac{\cos2x}{\cos^2x}dx$;
\task $\Int \dfrac{1-2x}{\sqrt{1-x^2}}dx$;
\task $\Int(1+\cos\frac{x}{2})^2dx$.
\end{tasks}
\textbf{2.}\hspace{1mm}Integrate the quadratic fractions:
\begin{tasks}[ label-align=left, label-offset={1mm}, label-width={3mm},column-sep={0pt}, item-indent={0mm} ,, ](2)
\task $\Int\dfrac{-6x+5}{x^2-9} dx$;
\task $\Int\dfrac{dx}{\sqrt{4x^2-4x-5}}$;
\task $\Int\dfrac{(-4x+2)dx}{x^2-2x+17}$;
\task $\Int\dfrac{-6x+7}{\sqrt{6x-x^2}} dx$.
\end{tasks}
\textbf{3.}\hspace{1mm}Integrate by parts or using the \\ suitable substitutions:
\begin{tasks}[ label-align=left, label-offset={1mm}, label-width={3mm}, column-sep={0pt}, item-indent={0mm},  ,after-item-skip=-1mm,](2)
\task  $\Int\dfrac{1-e^x}{1+e^x}dx$;
\task $\Int\dfrac{dx}{x\sqrt{x^2-9}}$;
\task! $\Int (2-8x^2) \cos 2x dx$;
\task! $\Int (x+1)^4\ln(x+1) dx$;
\task $\Int \sin\sqrt[3]{x}dx$;
\task $\Int \dfrac{x^2\textrm{arctg}x}{1+x^2}dx$.
\end{tasks}

\textbf{4.}\hspace{1mm}Integrate the polynomial fractions:
\begin{tasks}[ label-align=left, label-offset={1mm}, label-width={3mm},after-item-skip=-0.5mm, item-indent={0mm} , ](1)
\task $\Int\dfrac{4x^2+12x+6}{(x^2+3x+2)(x+3)} dx$;
\task $\Int\dfrac{4x^3-7x^2-2x+1}{x^3-x^2-x+1} dx$;
\task $\Int\dfrac{-5x^2-6x+15}{(x^2+4x+5)(x-2)}dx$.
\end{tasks}

\textbf{5.}\hspace{1mm}Integrate trigonometric expressions:
\begin{tasks}[ label-align=left, label-offset={1mm}, label-width={3mm},after-item-skip=-1mm, item-indent={0mm} ,column-sep={20pt}, ](2)
\task! $\Int\cos x\sin x\cos 3xdx$;
\task $\Int\sqrt[4]{\cos^5 x}\sin^5 x dx$;
\task $\Int\dfrac{dx}{\sin x+\cos x}$.
\end{tasks}
\textbf{6.}\hspace{1mm}Integrate the fractions with radicals:
\begin{tasks}[ label-align=left, label-offset={1mm}, label-width={3mm}, item-indent={0mm},column-sep={0in} , ](2)
\task $\Int \dfrac{x^4}{\sqrt{(1-x^2)^3}}dx$;
\task $\Int\dfrac{x dx}{(2+x)\sqrt{1+x}}$.
\end{tasks}
\end{minipage}
\vline\;
\begin{minipage}[t]{100mm}
\textbf{7.}\hspace{1mm}Solve the definite integrals:
\begin{tasks}[ label-align=left, label-offset={1mm}, label-width={3mm},column-sep={-40pt}, item-indent={0mm},after-item-skip=-1mm, ,column-sep={0pt} ](2)
\task $\Int_{-1}^{\ln4-1}(x+1) e^{2x+2}dx$;
\task $\Int_{0}^{\pi}\sin^2x\cos^4xdx$;
\task $\Int_{-1}^{0}\dfrac{2x-4}{\sqrt{3-2x-x^2}}dx$;
\task $\Int_{0}^{\sqrt{3}}\dfrac{\sqrt{x+1}+1}{\sqrt{x+1}-1}dx.$
\end{tasks}

\textbf{8.}\hspace{1mm}Find the area of the figure bounded \\
by the curves:
\begin{tasks}[ label-align=left, label-offset={2mm}, label-width={3mm}, item-indent={0mm}, after-item-skip=-1mm,     ](1)
\task $2x+y=5,\ xy=2$;
\task $\left\lbrace\begin{array}{l} x=3\cos t,\\ y=6\sin t; \end{array}\right.\ x=0 \ (x\geq0)$;
\task $\rho=2(1+\sin\phi)$.
\end{tasks}

\textbf{9.}\hspace{1mm}Find the arc-length of the curve:
\begin{tasks}[ label-align=left, label-offset={2mm}, label-width={3mm}, item-indent={0mm}, after-item-skip=-1mm,     ](1)
\task $y=5-\textrm{ch}x,\ 0\leq x\leq \ln3$;
\task $\left\lbrace\begin{array}{l} x=2\cos^3 t,\\ y=2\sin^3 t; \end{array}\right. \  0\leq t\leq \pi$;
\task $\rho=5e^{\phi}, \ 0\leq\phi\leq \pi$.
\end{tasks}

\textbf{10.}\hspace{1mm}Find the area of the surface formed \\ by rotating the curves around the $l$-axis:
\begin{tasks}[ label-align=left, label-offset={2mm}, label-width={3mm}, item-indent={0mm}, after-item-skip=-1mm,    ](1)
\task $y=2x^2, \ 0\leq x\leq \sqrt{6}, \ l=OY$;
\task $\left\lbrace\begin{array}{l} x=t-\sin t,\\  y=1-\cos t; \end{array}\right. \ 0\leq t\leq2\pi, \ l=OX$;
\task $\rho=3\sqrt{\sin2\big(\phi-\frac{\pi}{4}\big)}, \ l=o\rho$.
\end{tasks}

\textbf{11.}\hspace{1mm}Find the volume of the body formed \\
by rotating the curves around the $l$-axis:
\begin{tasks}[ label-align=left, label-offset={2mm}, label-width={3mm}, item-indent={0mm}, after-item-skip=-1mm,     ](1)
\task  $y=x^2-4x,\ y=0, \ l=OX$;
\task $y=9-x^2,\ y=8x,\ x=0 \ (x\geq0),\ \\ l=OY$;
\task $\rho=2\cos\phi, \ l=o\rho$.
\end{tasks}

\textbf{12.}\hspace{1mm}Solve the improper integrals:
\begin{tasks}[ label-align=left, label-offset={1mm}, label-width={3mm}, item-indent={0mm},   ,   ,column-sep={-40pt}](2)
\task $\Int_{0}^{\infty}\dfrac{2x-3}{x^2+4}dx$;
\task  $\Int_{1}^{e}\dfrac{\ln x-1}{x\sqrt{1-\ln^2x}}dx$.
\end{tasks}
\end{minipage}

\section{Personal task 11}

   \begin{minipage}[t]{90mm}
\textbf{1.}\hspace{1mm}Integrate using the table and \\ substitution under differential:
\begin{tasks}[ label-align=left, label-offset={1mm}, label-width={3mm}, column-sep={0pt}, item-indent={0mm} , ,after-item-skip=-1mm](2)
\task $\Int \dfrac{(-x^2+3)^3 dx}{x^4}$;
\task $\Int (\textrm{tg}x+\textrm{ctg}x)^2 dx$;
\task $\Int x^2e^{-3x^3+1}dx$;
\task $\Int\sin^2(x+2)dx$.
\end{tasks}

\textbf{2.}\hspace{1mm}Integrate the quadratic fractions:
\begin{tasks}[ label-align=left, label-offset={1mm}, label-width={3mm},column-sep={0pt}, item-indent={0mm} ,, ](2)
\task $\Int\dfrac{1-6x}{x^2-9} dx$;
\task $\Int\dfrac{dx}{\sqrt{12x-9x^2}}$;
\task $\Int\dfrac{(x+2)dx}{x^2+8x+7}$;
\task $\Int\dfrac{(-2x-4)dx}{\sqrt{x^2-2x-8}}$.
\end{tasks}

\textbf{3.}\hspace{1mm}Integrate by parts or using the \\ suitable substitutions:
\begin{tasks}[ label-align=left, label-offset={1mm}, label-width={3mm}, column-sep={-30pt}, item-indent={0mm},  ,after-item-skip=-1mm,](2)
\task $\Int\dfrac{4^x-2^x}{\sqrt{4^x+1}}dx$;
\task $\Int\dfrac{\sqrt{4-x^2}}{x}dx$;
\task $\Int 2x^2 \cos4x dx$;
\task $\Int\arcsin3x dx$;
\task $\Int e^{\sqrt[3]x}dx$;
\task $\Int\ln(x+\sqrt{x^2+1})dx$.
\end{tasks}

\textbf{4.}\hspace{1mm}Integrate the polynomial fractions:
\begin{tasks}[ label-align=left, label-offset={1mm}, label-width={3mm},after-item-skip=-0.5mm, item-indent={0mm} , ](1)
\task $\Int\dfrac{3x^2+7x-4}{(x^2+3x+2)(x-1)} dx$;
\task $\Int\dfrac{x^4+2x^3-4x^2-10x-12}{x^3-4x} dx$;
\task $\Int\dfrac{4x^2+14x+16}{(x^2+2x+2)(x+2)}dx$.
\end{tasks}

\textbf{5.}\hspace{1mm}Integrate trigonometric expressions:
\begin{tasks}[ label-align=left, label-offset={1mm}, label-width={3mm},after-item-skip=-1mm, item-indent={0mm} ,column-sep={0pt}, ](2)
\task $\Int\sin5x\cos 2xdx$;
\task $\Int\cos^4 2x\sin^2 2x dx$;
\task $\Int\dfrac{dx}{3-2\sin x+\cos x}$.
\end{tasks}
\textbf{6.}\hspace{1mm}Integrate the fractions with radicals:
\begin{tasks}[ label-align=left, label-offset={1mm}, label-width={3mm}, item-indent={0mm},column-sep={-40pt} , ](1)
\task $\Int\dfrac{\sqrt{x^2-9}}{x^3}dx$;
\task $\Int\dfrac{dx}{(\sqrt[4]{x+3}-1)\sqrt{x+3}}$.
\end{tasks}
\end{minipage}
\vline\;
      \begin{minipage}[t]{100mm}
\textbf{7.}\hspace{1mm}Solve the definite integrals:
\begin{tasks}[ label-align=left, label-offset={1mm}, label-width={3mm},column-sep={-40pt}, item-indent={0mm},after-item-skip=-1mm, ,column-sep={0pt} ](2)
\task $\Int_{1}^{e}(2x-1)\ln xdx$;
\task $\Int_{0}^{\pi/8}\dfrac{\sqrt[5]{\textrm{tg}^3 2x}}{\cos^2 2x}dx$;
\task $\Int_{0}^{2} x^2\sqrt{1+x^3}dx$;
\task $\Int_{0}^{26}\dfrac{\sqrt[3]{x+1}}{5+\sqrt[3]{x+1}}dx.$
\end{tasks}

\textbf{8.}\hspace{1mm}Find the area of the figure bounded \\
by the curves:
\begin{tasks}[ label-align=left, label-offset={2mm}, label-width={3mm}, item-indent={0mm}, after-item-skip=-1mm,     ](1)
\task $y=-x^2+4x+1,\ y=-x+5$;
\task $\left\lbrace\begin{array}{l} x=2\cos^3t,\\ y=2\sin^3t; \end{array}\right. \ x=0 \ (x\geq 0)$;
\task $\rho=3\cos2\phi$.
\end{tasks}

\textbf{9.}\hspace{1mm}Find the arc-length of the curve:
\begin{tasks}[ label-align=left, label-offset={2mm}, label-width={3mm}, item-indent={0mm}, after-item-skip=-1mm,     ](1)
\task $y=4\ln x,\ \sqrt{3}\leq x\leq \sqrt{15}$;
\task $\left\lbrace\begin{array}{l} x=3(t-\sin t),\\ y=3(1-\cos t); \end{array}\right. \  0\leq t\leq 2\pi$;
\task $\rho=3e^{-3\phi}, \ 0\leq\phi\leq\pi$.
\end{tasks}

\textbf{10.}\hspace{1mm}Find the area of the surface formed \\ by rotating the curves around the $l$-axis:
\begin{tasks}[ label-align=left, label-offset={2mm}, label-width={3mm}, item-indent={0mm}, after-item-skip=-1mm,    ](1)
\task $y=x^3, \ -2\leq x\leq 2, \ l=OX$;
\task $\left\lbrace\begin{array}{l} x=2+\sin t,\\ y=2+\cos t; \end{array}\right. \ l=OY$;
\task $\rho=6\sin^2\Big(\dfrac{\phi}{2}\Big), \ l=o\rho$.
\end{tasks}

\textbf{11.}\hspace{1mm}Find the volume of the body formed \\
by rotating the curves around the $l$-axis:
\begin{tasks}[ label-align=left, label-offset={2mm}, label-width={3mm}, item-indent={0mm}, after-item-skip=-1mm,     ](1)
\task $y=x^2+4x+4,\ y=3x+6, \ l=OX$;
\task $x=-y^2+4y,\ x=0, \ l=OY$;
\task $\rho=2\sin\phi, \ l=o\rho$.
\end{tasks}

\textbf{12.}\hspace{1mm}Solve the improper integrals:
\begin{tasks}[ label-align=left, label-offset={1mm}, label-width={3mm}, item-indent={0mm},   ,   ,column-sep={-40pt}](2)
\task $\Int_{-\infty}^{0}\dfrac{x^2}{x^6+4}dx$;
\task $\Int_{0}^{1/2}\dfrac{dx}{\sqrt[3]{(1-2x)^2}}$.
\end{tasks}
\end{minipage}

\section{Personal task 12}

   \begin{minipage}[t]{90mm}
\textbf{1.}\hspace{1mm}Integrate using the table and \\ substitution under differential:
\begin{tasks}[ label-align=left, label-offset={1mm}, label-width={3mm}, column-sep={0pt}, item-indent={0mm} , ,after-item-skip=-1mm](2)
\task $\Int \dfrac{3x^2 2^x-4^x}{2^x}dx$;
\task $\Int (5-\textrm{ctg}^2x) dx$;
\task $\Int x^22^{2x^3-1}dx$;
\task $\Int\cos^2(3x+1)dx$;
\end{tasks}

\textbf{2.}\hspace{1mm}Integrate the quadratic fractions:
\begin{tasks}[ label-align=left, label-offset={1mm}, label-width={3mm},column-sep={0pt}, item-indent={0mm} ,, ](2)
\task $\Int\dfrac{4x-7}{x^2-4} dx$;
\task $\Int\dfrac{dx}{\sqrt{6x-9x^2}}$;
\task $\Int\dfrac{(2x-5)dx}{x^2-4x+13}$
\task $\Int\dfrac{(-2x+3)dx}{\sqrt{x^2+2x-3}}$;
\end{tasks}

\textbf{3.}\hspace{1mm}Integrate by parts or using the \\ suitable substitutions:
\begin{tasks}[ label-align=left, label-offset={1mm}, label-width={3mm}, column-sep={-10pt}, item-indent={0mm},  ,after-item-skip=-1mm,](2)
\task $\Int x^2(x+2)^{10} dx$;
\task $\Int\dfrac{2dx}{x^2\sqrt{2x^2+2x+1}}$;
\task! $\Int (4-12x^2) \sin 4x dx$;
\task! $\Int (x-1)^3\ln (x-1)dx$;
\task $\Int e^{\sqrt{2x-4}}dx$;
\task $\Int\dfrac{x^2\textrm{arcctg}x}{x^2+1}dx$.
\end{tasks}

\textbf{4.}\hspace{1mm}Integrate the polynomial fractions:
\begin{tasks}[ label-align=left, label-offset={1mm}, label-width={3mm},after-item-skip=-0.5mm, item-indent={0mm} , ](1)
\task $\Int\dfrac{5x^2-11x-6}{(x^2-3x+2)(x+2)} dx$;
\task $\Int\dfrac{-2x^3-12x^2-5x+70}{(x+2)(x^2+3x-10)} dx$;
\task $\Int\dfrac{3x^2+6x+10}{x^3+2x^2+2x}dx$;
\end{tasks}

\textbf{5.}\hspace{1mm}Integrate trigonometric expressions:
\begin{tasks}[ label-align=left, label-offset={1mm}, label-width={3mm},after-item-skip=-1mm, item-indent={0mm} ,column-sep={0pt}, ](2)
\task $\Int\sin2x\sin 5xdx$
\task $\Int\cos^5 x\sqrt[3]{\sin x} dx$;
\task $\Int\dfrac{-2dx}{9-7\cos x}$;
\end{tasks}
\textbf{6.}\hspace{1mm}Integrate the fractions with radicals:
\begin{tasks}[ label-align=left, label-offset={1mm}, label-width={3mm}, item-indent={0mm},column-sep={-30pt} , ](2)
\task $\Int\dfrac{\sqrt{1-9x^2}}{x^4}dx$;
\task $\Int\sqrt{\dfrac{x+1}{x-2}}{\dfrac{dx}{(x+1)^2}}$ 
\end{tasks}
\end{minipage}
\vline \;
      \begin{minipage}[t]{100mm}
\textbf{7.}\hspace{1mm}Solve the definite integrals:
\begin{tasks}[ label-align=left, label-offset={1mm}, label-width={3mm},column-sep={-20pt}, item-indent={0mm},after-item-skip=-1mm ](2)
\task $\Int_{-2}^{0}x^2 e^{x+2}dx$;
\task $\Int_{0}^{3\pi/2}\cos^4(\frac{x}{3})\sin^2(\frac{x}{3})dx$;
\task $\Int_{0}^{1/\sqrt{2}}\dfrac{2\textrm{arcsin} x+x}{\sqrt{1-x^2}}dx$;
\task $\Int_{0}^{e^{2}-1}\dfrac{\ln(x+1)dx}{\sqrt{x+1}}.$
\end{tasks}

\textbf{8.}\hspace{1mm}Find the area of the figure bounded \\
by the curves:
\begin{tasks}[ label-align=left, label-offset={2mm}, label-width={3mm}, item-indent={0mm}, after-item-skip=-1mm,     ](1)
\task $y=(x+3)^2,\ y=4x+12$;
\task $\left\lbrace\begin{array}{l} x=3(t-\sin t),\\ y=3(1-\cos t); \end{array}\right. \ y=0, \ 0\leq t\leq 2\pi$;
\task $\rho=2\sin2\phi.$
\end{tasks}

\textbf{9.}\hspace{1mm}Find the arc-length of the curve:
\begin{tasks}[ label-align=left, label-offset={2mm}, label-width={3mm}, item-indent={0mm}, after-item-skip=-1mm,     ](1)
\task $y=2\ln(4-x^2),\ -1\leq x\leq 1$;
\task $\left\lbrace\begin{array}{l} x=2e^t\cos t,\\ y=2e^t\sin t; \end{array}\right. \ 0\leq t\leq \pi$;
\task $\rho=4(1-\cos\phi)$.
\end{tasks}

\textbf{10.}\hspace{1mm}Find the area of the surface formed \\ by rotating the curves around the $l$-axis:
\begin{tasks}[ label-align=left, label-offset={2mm}, label-width={3mm}, item-indent={0mm}, after-item-skip=-1mm,    ](1)
\task $y=\sqrt{1-\dfrac{x^2}{9}}, \ -3\leq x\leq 3, \ l=OX$;
\task $\left\lbrace\begin{array}{l} x=\cos^3 t,\\ y=\sin^3 t; \end{array}\right. \ l=OY$;
\task $\rho=\sqrt{\cos2\phi}, \ l=o\rho$.
\end{tasks}

\textbf{11.}\hspace{1mm}Find the volume of the body formed \\
by rotating the curves around the $l$-axis:
\begin{tasks}[ label-align=left, label-offset={2mm}, label-width={3mm}, item-indent={0mm}, after-item-skip=-1mm,     ](1)
\task $y=-x^2+4x+1,\ y=5-x, \ l=OX$;
\task $x=-y^2+7y,\ x=6, \ l=OY$;
\task $\rho=6\phi, \ 0\leq\phi\pi, \ l=o\rho$.
\end{tasks}

\textbf{12.}\hspace{1mm}Solve the improper integrals:
\begin{tasks}[ label-align=left, label-offset={1mm}, label-width={3mm}, item-indent={0mm},   ,   ,column-sep={-40pt}](2)
\task $\Int_{1}^{\infty}\dfrac{4x-1}{x^2+9}dx$;
\task $\Int_{0}^{1}\dfrac{\textrm{arccos}x}{\sqrt{1-x^2}}dx$.
\end{tasks}

\end{minipage}

\section{Personal task 13}

   \begin{minipage}[t]{90mm}
\textbf{1.}\hspace{1mm}Integrate using the table and \\ substitution under differential:
\begin{tasks}[ label-align=left, label-offset={1mm}, label-width={3mm}, column-sep={0pt}, item-indent={0mm} , ,after-item-skip=-1mm](2)
\task $\Int \dfrac{(2-x\sqrt{x})^3dx}{x^2}$;
\task $\Int (\textrm{tg}x-\textrm{ctg}x)^2 dx$;
\task $\Int\dfrac{\ln^4(2x+3)dx}{2x+3}$;
\task $\Int\cos^2(3-2x)dx$;
\end{tasks}

\textbf{2.}\hspace{1mm}Integrate the quadratic fractions:
\begin{tasks}[ label-align=left, label-offset={1mm}, label-width={3mm},column-sep={-15pt}, item-indent={0mm} ,, ](2)
\task $\Int\dfrac{5x+3}{x^2-16} dx$;
\task $\Int\dfrac{dx}{\sqrt{12+12x-9x^2}}$;
\task $\Int\dfrac{(4x+6)dx}{x^2-4x+8}$;
\task $\Int\dfrac{(2x-7)dx}{\sqrt{x^2-4x-12}}$;
\end{tasks}

\textbf{3.}\hspace{1mm}Integrate by parts or using the \\ suitable substitutions:
\begin{tasks}[ label-align=left, label-offset={1mm}, label-width={3mm}, column-sep={5pt}, item-indent={0mm},  ,after-item-skip=-1mm,](2)
\task $\Int\dfrac{dx}{\sqrt{1-e^{2x}}}$;
\task $\Int\dfrac{\sqrt{4+x^2}}{x}dx$;
\task $\Int 12x \sin 3x dx$;
\task $\Int \dfrac{x^2+4x}{e^{2x}}dx$;
\task $\Int\arccos 2x dx$;
\task $\Int\dfrac{2x+3}{\cos^2x}dx$;
\end{tasks}

\textbf{4.}\hspace{1mm}Integrate the polynomial fractions:
\begin{tasks}[ label-align=left, label-offset={1mm}, label-width={3mm},after-item-skip=-0.5mm, item-indent={0mm} , ](1)
\task $\Int\dfrac{8x^2-31x+25}{(x^2-5x+6)(x-1)} dx$;
\task $\Int\dfrac{2x^3-22x^2+52x+49}{x^3-14x^2+49x} dx$;
\task $\Int\dfrac{9x^2-10x+34}{(x^2-2x+10)(x-2)}dx$;
\end{tasks}

\textbf{5.}\hspace{1mm}Integrate trigonometric expressions:
\begin{tasks}[ label-align=left, label-offset={1mm}, label-width={3mm},after-item-skip=-1mm, item-indent={0mm} ,column-sep={0pt}, ](2)
\task $\Int\cos9x\cos2xdx$
\task $\Int \sqrt[4]{\dfrac{\cos x}{\sin^{9}x}}dx$;
\task $\Int\dfrac{6\sin x dx}{13\cos^2x+4\sin^2x}$;
\end{tasks}
\textbf{6.}\hspace{1mm}Integrate the functions with radicals:
\begin{tasks}[ label-align=left, label-offset={1mm}, label-width={3mm}, item-indent={0mm},column-sep={-40pt} , ](2)
\task $\Int x^2\sqrt{9-x^2}dx$;
\task $\Int\dfrac{(\sqrt[6]{x+2}-1)dx}{(x+2)(1+\sqrt[3]{x+2})}$.
\end{tasks}
\end{minipage}
\vline\;
      \begin{minipage}[t]{100mm}
\textbf{7.}\hspace{1mm}Solve the definite integrals:
\begin{tasks}[ label-align=left, label-offset={1mm}, label-width={3mm},column-sep={-40pt}, item-indent={0mm},after-item-skip=-1mm, ,column-sep={0pt} ](2)
\task $\Int_{0}^{e-1}\ln\sqrt{x+1}dx$;
\task $\Int_{0}^{\pi} \cos^8x\sin^3x dx$;
\task $\Int_{0}^{\sqrt{3}} \dfrac{2\textrm{arctg}x-x}{1+x^2}dx$;
\task $\Int_{1}^{5}\dfrac{dx}{x+\sqrt{2x-1}}.$
\end{tasks}

\textbf{8.}\hspace{1mm}Find the area of the figure bounded \\
by the curves:
\begin{tasks}[ label-align=left, label-offset={2mm}, label-width={3mm}, item-indent={0mm}, after-item-skip=-1mm,     ](1)
\task $x+y=6,\ xy=5$;
\task $\left\lbrace\begin{array}{l} x=3\cos^3t,\\ y=3\sin^3t; \end{array}\right. \ x=0 \ (x\geq 0)$;
\task $\rho=2\sqrt{\cos2\phi}.$
\end{tasks}

\textbf{9.}\hspace{1mm}Find the arc-length of the curve:
\begin{tasks}[ label-align=left, label-offset={2mm}, label-width={3mm}, item-indent={0mm}, after-item-skip=-1mm,     ](1)
\task $y=4-e^x,\ 0\leq x\leq \sqrt{15}$;
\task $\left\lbrace\begin{array}{l} x=2(t-\sin t),\\ y=2(1-\cos t); \end{array}\right. \ 0\leq t\leq 2\pi$;
\task $\rho=2\cos\phi+2\sin\phi$.
\end{tasks}

\textbf{10.}\hspace{1mm}Find the area of the surface formed \\ by rotating the curves around the $l$-axis:
\begin{tasks}[ label-align=left, label-offset={2mm}, label-width={3mm}, item-indent={0mm}, after-item-skip=-1mm,    ](1)
\task $y=3x+6, \ -2\leq x\leq 4, \ l=OX$;
\task $\left\lbrace\begin{array}{l} x=2+4\cos t, \\ y=4\sin t; \end{array}\right. \ (x\geq 0) \ l=OY$;
\task $\rho=3e^{2\phi}, \ 0\leq\phi\leq \pi, \ l=o\rho$.
\end{tasks}

\textbf{11.}\hspace{1mm}Find the volume of the body formed \\
by rotating the curves around the $l$-axis:
\begin{tasks}[ label-align=left, label-offset={2mm}, label-width={3mm}, item-indent={0mm}, after-item-skip=-1mm,     ](1)
\task $y=-x^2+4x+3,\ y=-x+7, \ l=OX$;
\task $y=\sqrt{4-x^2}, \ y=x, \ x=0, \ l=OY$;
\task $\rho=4(1-\cos\phi), \ l=o\rho$.
\end{tasks}

\textbf{12.}\hspace{1mm}Solve the improper integrals:
\begin{tasks}[ label-align=left, label-offset={1mm}, label-width={3mm}, item-indent={0mm},   ,   ,column-sep={-50pt}](2)
\task $\Int_{0}^{\infty}xe^{-x}dx$;
\task $\Int_{0}^{\pi/2}\dfrac{\sin2x}{\sqrt{1-\cos^4x}}dx$.
\end{tasks}

\end{minipage}

\section{Personal task 14}

   \begin{minipage}[t]{90mm}
\textbf{1.}\hspace{1mm}Integrate using the table and \\ substitution under differential:
\begin{tasks}[ label-align=left, label-offset={1mm}, label-width={3mm}, column-sep={5pt}, item-indent={0mm} , ,after-item-skip=-1mm](2)
\task $\Int \sqrt[3]{x}(1+\sqrt[3]{x})^2 dx$;
\task $\Int \dfrac{\sin^3x}{1+\cos x} dx$;
\task $\Int x^32^{x^4+3}dx$;
\task $\Int\dfrac{dx}{\textrm{sh}x+\textrm{ch}x}$;
\end{tasks}

\textbf{2.}\hspace{1mm}Integrate the quadratic fractions:
\begin{tasks}[ label-align=left, label-offset={1mm}, label-width={3mm},column-sep={0pt}, item-indent={0mm} ,, ](2)
\task $\Int\dfrac{-4x-7}{x^2-9} dx$;
\task $\Int\dfrac{8dx}{\sqrt{8+6x-9x^2}}$;
\task $\Int\dfrac{-2x-12}{x^2+4x+8}dx$
\task $\Int\dfrac{4x\; dx}{\sqrt{x^2-2x+15}}$;
\end{tasks}

\textbf{3.}\hspace{1mm}Integrate by parts or using the \\ suitable substitutions:
\begin{tasks}[ label-align=left, label-offset={1mm}, label-width={3mm}, column-sep={-10pt}, item-indent={0mm},  ,after-item-skip=-1mm,](2)
\task $\Int\dfrac{2x-3}{(x-2)^5}dx$;
\task $\Int\dfrac{dx}{x\sqrt{5x^2-4x-1}}$;
\task! $\Int 4x(\cos 2x-\sin 2x)^2 dx$;
\task $\Int\dfrac{\ln(x-1)}{(x-1)^2}dx$;
\task $\Int\textrm{arcctg}x dx$;
\task $\Int e^{\sqrt[3]{x-4}}dx$;
\end{tasks}

\textbf{4.}\hspace{1mm}Integrate the polynomial fractions:
\begin{tasks}[ label-align=left, label-offset={1mm}, label-width={3mm},after-item-skip=-0.5mm, item-indent={0mm} , ](1)
\task $\Int\dfrac{8x^2+7x-30}{(x^2-3x+2)(x+2)} dx$;
\task $\Int\dfrac{x^3-2x^2-14x+21}{x^3-5x^2+3x+9} dx$;
\task $\Int\dfrac{9x+64}{(x^2+6x+10)(x-4)}dx$;
\end{tasks}

\textbf{5.}\hspace{1mm}Integrate trigonometric expressions:
\begin{tasks}[ label-align=left, label-offset={1mm}, label-width={3mm},after-item-skip=-1mm, item-indent={0mm} ,column-sep={0pt}, ](2)
\task $\Int\cos2x\sin 7xdx$
\task $\Int\dfrac{\cos^3 x}{\sin^9 x} dx$;
\task $\Int\dfrac{dx}{5+\sin x+4\cos x}$;
\end{tasks}
\textbf{6.}\hspace{1mm}Integrate the fractions with radicals:
\begin{tasks}[ label-align=left, label-offset={1mm}, label-width={3mm}, item-indent={0mm},column-sep={-30pt} , ](2)
\task $\Int\dfrac{x^2}{\sqrt{1-x^2}}dx$;
\task $\Int\dfrac{\sqrt[6]{x}}{\sqrt{x}-\sqrt[3]{x^2}} dx$.
\end{tasks}
\end{minipage}
\vline \;
      \begin{minipage}[t]{100mm}
\textbf{7.}\hspace{1mm}Solve the definite integrals:
\begin{tasks}[ label-align=left, label-offset={1mm}, label-width={3mm},column-sep={-30pt}, item-indent={0mm},after-item-skip=-1mm, ](2)
\task $\Int_{0}^{\pi/2}e^{\sin^2 x}\sin 2xdx$;
\task  $\Int_{0}^{\pi/4}\cos^4 2x\sin^3 2xdx$;
\task  $\Int_{0}^{\sqrt{3}}\dfrac{-12x+9}{x^2+1}dx$;
\task  $\Int_{0}^{\ln 2/2}\dfrac{3e^{4x}}{\sqrt{e^{2x}-1}}dx.$
\end{tasks}

\textbf{8.}\hspace{1mm}Find the area of the figure bounded \\
by the curves:
\begin{tasks}[ label-align=left, label-offset={2mm}, label-width={3mm}, item-indent={0mm}, after-item-skip=-1mm,     ](1)
\task  $y=(x-4)^2,\ y=3x-12$;
\task  $\left\lbrace\begin{array}{l} x=2+4\cos t,\\ y=2\sin t; \end{array}\right. \ y=1 \ (y\geq1)$;
\task  $\rho=3\sin 3\phi.$
\end{tasks}

\textbf{9.}\hspace{1mm}Find the arc-length of the curve:
\begin{tasks}[ label-align=left, label-offset={2mm}, label-width={3mm}, item-indent={0mm}, after-item-skip=-1mm,     ](1)
\task  $y=2+\frac{1}{2}\ln(\sin 2x),\ \frac{\pi}{6}\leq x\leq \frac{\pi}{4}$;
\task  $\left\lbrace\begin{array}{l} x=2\cos t-\cos2t,\\ y=2\sin t-\sin2t; \end{array}\right. \ 0\leq t\leq 2\pi$;
\task  $\rho=4\phi, \ 0\leq\phi\leq\pi$.
\end{tasks}

\textbf{10.}\hspace{1mm}Find the area of the surface formed \\ by rotating the curves around the $l$-axis:
\begin{tasks}[ label-align=left, label-offset={2mm}, label-width={3mm}, item-indent={0mm}, after-item-skip=-1mm,    ](1)
\task $y=\frac{1}{3}\textrm{ch}3x, \ 0\leq x \leq 1, \ l=OX$;
\task $\left\lbrace\begin{array}{l} x=2(t-\sin t), \\ y=2(1-\cos t); \end{array}\right. \ \begin{array}{l} 0\leq t\leq2\pi, \\ l=OX; \end{array}$
\task $\rho=4\sqrt{\cos2\phi}, \ l=o\rho$.
\end{tasks}

\textbf{11.}\hspace{1mm}Find the volume of the body formed \\
by rotating the curves around the $l$-axis:
\begin{tasks}[ label-align=left, label-offset={2mm}, label-width={3mm}, item-indent={0mm}, after-item-skip=-1mm,     ](1)
\task  $y =-x^2+4x,\ y=-x, \ l=OX$;
\task  $x=6-y^2,\ x=2, \ l=OY$;
\task  $\rho=2\cos\phi, \ l=o\rho$.
\end{tasks}

\textbf{12.}\hspace{1mm}Solve the improper integrals:
\begin{tasks}[ label-align=left, label-offset={1mm}, label-width={3mm}, item-indent={0mm},   ,   ,column-sep={-40pt}](2)
\task  $\Int_{0}^{\infty}\dfrac{(4x-5)dx}{4x^2-4x+2}$;
\task  $\Int_{0}^{2}\dfrac{\textrm{arcsin}\big(\dfrac{x}{2}\big)-2x}{\sqrt{4-x^2}}dx$.
\end{tasks}

\end{minipage}

\section{Personal task 15}

   \begin{minipage}[t]{90mm}
\textbf{1.}\hspace{1mm}Integrate using the table and \\ substitution under differential:
\begin{tasks}[ label-align=left, label-offset={1mm}, label-width={3mm}, column-sep={-18pt}, item-indent={0mm} , ,after-item-skip=-1mm](2)
\task $\Int \dfrac{(3-2x)^3}{x^2} dx$;
\task $\Int (\sin x-\dfrac{1}{\sin x})^2dx$;
\task $\Int \dfrac{x^4}{\sqrt[3]{x^5+9}}dx$;
\task $\Int\dfrac{\textrm{6arcsin}^5x}{\sqrt{1-x^2}}dx$.
\end{tasks}
\textbf{2.}\hspace{1mm}Integrate the quadratic fractions:
\begin{tasks}[ label-align=left, label-offset={1mm}, label-width={3mm},column-sep={-10pt}, item-indent={0mm} ,, ](2)
\task $\Int\dfrac{8x+7}{x^2-4} dx$;
\task $\Int\dfrac{dx}{\sqrt{9x^2-4x-12}}$;
\task $\Int\dfrac{(-2x+11)dx}{x^2-8x+17}$;
\task $\Int\dfrac{6x+3}{\sqrt{5-4x-x^2}} dx$.
\end{tasks}
\textbf{3.}\hspace{1mm}Integrate by parts or using the \\ suitable substitutions:
\begin{tasks}[ label-align=left, label-offset={1mm}, label-width={3mm}, column-sep={-20pt}, item-indent={0mm},  ,after-item-skip=-1mm,](2)
\task $\Int\dfrac{dx}{\sqrt{e^x-4}}$;
\task $\Int\dfrac{-4dx}{x^2\sqrt{5x^2+4x-1}}$;
\task $\Int 12x^2 \sin 3x dx$;
\task $\Int \dfrac{\cos^2x}{e^{2x}}dx$;
\task $\Int \ln(4-x^2)dx$;
\task $\Int \textrm{arctg}\sqrt{x}dx$.
\end{tasks}
\textbf{4.}\hspace{1mm}Integrate the polynomial fractions:
\begin{tasks}[ label-align=left, label-offset={1mm}, label-width={3mm},after-item-skip=-0.5mm, item-indent={0mm} , ](1)
\task $\Int\dfrac{6x^2-28x+28}{(x^2-3x+2)(x-3)} dx$;
\task $\Int\dfrac{2x^3-14x^2+30x-22}{(x^2-4x+3)(x-1)} dx$;
\task $\Int\dfrac{x^2+8x}{(x^2+4x+8)(x+2)}dx$.
\end{tasks}
\textbf{5.}\hspace{1mm}Integrate trigonometric expressions:
\begin{tasks}[ label-align=left, label-offset={1mm}, label-width={3mm},after-item-skip=-1mm, item-indent={0mm} ,, ](1)
\task $\Int\cos x\sin x\cos 4xdx$;
\task $\Int\sin^2 2x\cos^4 2x dx$;
\task $\Int\dfrac{12\cos xdx}{4\sin^2 x-5\cos^2 x+9}$.
\end{tasks}
\textbf{6.}\hspace{1mm}Integrate the fractions with radicals:
\begin{tasks}[ label-align=left, label-offset={1mm}, label-width={3mm}, item-indent={0mm},column-sep={-30pt} ,, ](2)
\task $\Int \dfrac{\sqrt{1-x^2}}{x^4}dx$;
\task $\Int\dfrac{dx}{\sqrt{x+3}+\sqrt[3]{x+3}}$.
\end{tasks}

\end{minipage}
\vline\;
      \begin{minipage}[t]{100mm}
\textbf{7.}\hspace{1mm}Solve the definite integrals:
\begin{tasks}[ label-align=left, label-offset={1mm}, label-width={3mm},column-sep={-40pt}, item-indent={0mm},after-item-skip=-1mm,  ](2)
\task $\Int_{0}^{1}(x-1)^2 e^{-x}dx$;
\task $\Int_{0}^{\pi/4}\dfrac{\sin^6x}{\cos^{10}x}dx$;
\task $\Int_{1}^{5/2}\dfrac{(-2x+5)dx}{\sqrt{8+2x-x^2}}$;
\task $\Int_{2}^{3}\dfrac{\sqrt[3]{(x-2)^2}dx}{1+\sqrt[3]{(x-2)^2}}.$\end{tasks}

\textbf{8.}\hspace{1mm}Find the area of the figure bounded \\
by the curves:
\begin{tasks}[ label-align=left, label-offset={2mm}, label-width={3mm}, item-indent={0mm}, after-item-skip=-1mm,     ](1)
\task $x+2y=7,\ xy=3$;
\task $\left\lbrace\begin{array}{l} x=2\cos t,\\ y=2\sin t; \end{array}\right. \ y=1\ (y\geq1)$;
\task $\rho=4(1-\sin\phi)$.
\end{tasks}

\textbf{9.}\hspace{1mm}Find the arc-length of the curve:
\begin{tasks}[ label-align=left, label-offset={2mm}, label-width={3mm}, item-indent={0mm}, after-item-skip=-1mm,     ](1)
\task $y=4-\frac{1}{2}\textrm{ch}2x,\ 0\leq x\leq \frac{1}{2}\ln3$;
\task $\left\lbrace\begin{array}{l} x=2\cos^3 t,\\ y=2\sin^3 t; \end{array}\right. \ 0\leq t\leq \pi$;
\task $\rho=5\sqrt{2}e^{-2\phi}, \ 0\leq\phi\leq \frac{\pi}{2}$.
\end{tasks}

\textbf{10.}\hspace{1mm}Find the area of the surface formed \\ by rotating the curves around the $l$-axis:
\begin{tasks}[ label-align=left, label-offset={2mm}, label-width={3mm}, item-indent={0mm}, after-item-skip=-1mm,    ](1)
\task $y=\frac{1}{8}x^2, \ 0\leq x\leq 3, \ l=OY$;
\task $\left\lbrace\begin{array}{l} x=2+\cos t,\\ y=3+\sin t; \end{array}\right.  \ l=OX$;
\task $\rho=6\cos \phi, \ l=o\rho$.
\end{tasks}

\textbf{11.}\hspace{1mm}Find the volume of the body formed \\
by rotating the curves around the $l$-axis:
\begin{tasks}[ label-align=left, label-offset={2mm}, label-width={3mm}, item-indent={0mm}, after-item-skip=-1mm,     ](1)
\task $y=x^2-4x+5,\ y=x+5, \ l=OX$;
\task $x=4-y^2,\ x=8-2y^2, \ l=OY$;
\task $\rho=4\sin2\phi, \ l=o\rho$.
\end{tasks}

\textbf{12.}\hspace{1mm}Solve the improper integrals:
\begin{tasks}[ label-align=left, label-offset={1mm}, label-width={3mm}, item-indent={0mm},   ,   ,column-sep={-40pt}](2)
\task $\Int_{1}^{\infty}\dfrac{2x^2+1}{x^4+x^2}dx$;
\task  $\Int_{1}^{e}\dfrac{dx}{x\sqrt{1-\ln x}}$.
\end{tasks}
\end{minipage}

\section{Personal task 16}

   \begin{minipage}[t]{90mm}
\textbf{1.}\hspace{1mm}Integrate using the table and \\ substitution under differential:
\begin{tasks}[ label-align=left, label-offset={1mm}, label-width={3mm}, column-sep={0pt}, item-indent={0mm} , ,after-item-skip=-1mm](2)
\task $\Int\dfrac{4x^3e^x+7e^{2x}}{e^x}dx$;
\task $\Int \dfrac{\sin2xdx}{(\cos x+3)^2-9} $;
\task $\Int \dfrac{\ln^3(3-2x)}{3-2x}dx$;
\task $\Int \dfrac{\sin^2\sqrt{x}}{\sqrt{x}}dx.$
\end{tasks}
\textbf{2.}\hspace{1mm}Integrate the quadratic fractions:
\begin{tasks}[ label-align=left, label-offset={1mm}, label-width={3mm},column-sep={-10pt}, item-indent={0mm} ,, ](2)
\task $\Int\dfrac{12x-3}{x^2-9} dx$;
\task  $\Int\dfrac{-10dx}{\sqrt{4x^2-4x-24}}$;
\task  $\Int\dfrac{(12-6x)dx}{x^2-6x+18}$;
\task  $\Int\dfrac{7-4x}{\sqrt{8-2x-x^2}} dx$.
\end{tasks}
\textbf{3.}\hspace{1mm}Integrate by parts or using the \\ suitable substitutions:
\begin{tasks}[ label-align=left, label-offset={1mm}, label-width={3mm}, column-sep={0pt}, item-indent={0mm},  ,after-item-skip=-1mm,](2)
\task  $\Int\sqrt{e^{2x}+4}dx$;
\task $\Int\dfrac{6 dx}{x\sqrt{9x^2-1}}$;
\task! $\Int 4x(\sin 2x+\cos 2x)^2 dx$;
\task $\Int (2x^2-6) e^{2x}dx$;
\task $\Int(x^2+4)\ln xdx$;
\task $\Int\dfrac{\arcsin x}{\sqrt{x+1}}dx$.
\end{tasks}
\textbf{4.}\hspace{1mm}Integrate the polynomial fractions:
\begin{tasks}[ label-align=left, label-offset={1mm}, label-width={3mm},after-item-skip=-0.5mm, item-indent={0mm} , ](1)
\task $\Int\dfrac{-3x^2-4x-11}{(x^2-5x+4)(x+1)} dx$;
\task$\Int\dfrac{2x^3+10x^2+12x-4}{x^3+4x^2+4x} dx$;
\task $\Int\dfrac{6x^2+4x+36}{(x^2+6x+18)(x-2)}dx$.
\end{tasks}
\textbf{5.}\hspace{1mm}Integrate trigonometric fractions:
\begin{tasks}[ label-align=left, label-offset={1mm}, label-width={3mm},after-item-skip=-1mm, item-indent={0mm} ,, ](2)
\task $\Int\sin5x\sin3xdx$
\task $\Int\cos^2 x\sin^4 x dx$;
\task $\Int\dfrac{4dx}{8+\sin x+7\cos x}$.
\end{tasks}
\textbf{6.}\hspace{1mm}Integrate the fractions with radicals:
\begin{tasks}[ label-align=left, label-offset={1mm}, label-width={3mm}, item-indent={0mm},column-sep={-1in} ,, ](2)
\task $\Int\dfrac{x^5 dx}{\sqrt{4-x^2}}$;
\task $\Int\sqrt[3]{\Big(\dfrac{x+1}{x-2}\Big)^2}\dfrac{dx}{(x+1)^2}$.
\end{tasks}
\end{minipage}
\vline\;
      \begin{minipage}[t]{100mm}
\textbf{7.}\hspace{1mm}Solve the definite integrals:
\begin{tasks}[ label-align=left, label-offset={1mm}, label-width={3mm},column-sep={-20pt}, item-indent={0mm},after-item-skip=-1mm,  ](2)
\task $\Int_{3}^{2e+2}\ln(x-2)dx$;
\task $\Int_{\pi/4}^{\pi/2} \dfrac{6\cos^2x}{\sin^4x} dx$;
\task $\Int_{-1}^{\sqrt[3]{31}} 3 x^2\sqrt[5]{x^3+1}dx$;
\task $\Int_{2}^{4}\dfrac{\sqrt{2x-4}}{5+\sqrt{2x-4}}dx.$
\end{tasks}

\textbf{8.}\hspace{1mm}Find the area of the figure bounded \\
by the curves:
\begin{tasks}[ label-align=left, label-offset={2mm}, label-width={3mm}, item-indent={0mm}, after-item-skip=-1mm,     ](1)
\task $y=2x^2-13x+16, \ x+y=6$;
\task  $\left\lbrace\begin{array}{l} x=2\cos^3t,\\ y=2\sin^3t; \end{array}\right. \ y=0 (y\geq 0)$;
\task  $\rho=2\sqrt{\cos2\phi}.$
\end{tasks}

\textbf{9.}\hspace{1mm}Find the arc-length of the curve:
\begin{tasks}[ label-align=left, label-offset={2mm}, label-width={3mm}, item-indent={0mm}, after-item-skip=-1mm,     ](1)
\task $y=\ln(x+\sqrt{x^2-9}),\ 3\leq x\leq 3\sqrt{5}$;
\task $\left\lbrace\begin{array}{l} x=4(t-\sin t),\\ y=4(1-\cos t); \end{array}\right. \ 0\leq t\leq 2\pi$;
\task $\rho=4\sin\phi+4\cos\phi$.
\end{tasks}

\textbf{10.}\hspace{1mm}Find the area of the surface formed \\ by rotating the curves around the $l$-axis:
\begin{tasks}[ label-align=left, label-offset={2mm}, label-width={3mm}, item-indent={0mm}, after-item-skip=-1mm,    ](1)
\task $x^2+\dfrac{y^2}{4}=1, \ l=OY$;
\task $\left\lbrace\begin{array}{l} x=e^t\cos t,\\ y=e^t\sin t; \end{array}\right. \ 0\leq t\leq\pi, \ l=OX$;
\task $\rho=6\sin\phi, \ l=o\rho$.
\end{tasks}

\textbf{11.}\hspace{1mm}Find the volume of the body formed \\
by rotating the curves around the $l$-axis:
\begin{tasks}[ label-align=left, label-offset={2mm}, label-width={3mm}, item-indent={0mm}, after-item-skip=-1mm,     ](1)
\task $y=-x^2+5x+1,\ y=x+1, \ l=OX$;
\task  $y=\sqrt{9-x^2}, \ y=x, \ x=0, \ l=OY$;
\task  $\rho=5(1-\cos\phi), \ l=o\rho$.
\end{tasks}

\textbf{12.}\hspace{1mm}Solve the improper integrals:
\begin{tasks}[ label-align=left, label-offset={1mm}, label-width={3mm}, item-indent={0mm},   ,   ,column-sep={0pt}](2)
\task  $\Int_{0}^{\infty}\dfrac{\textrm{arctg}\big(\dfrac{x}{3}\big)+4x}{x^2+9}dx$;
\task $\Int_{0}^{1}\dfrac{x+1}{\sqrt{2x-x^2}}dx$.
\end{tasks}
\end{minipage}

\section{Personal task 17}

   \begin{minipage}[t]{90mm}
\textbf{1.}\hspace{1mm}Integrate using the table and \\ substitution under differential:
\begin{tasks}[ label-align=left, label-offset={1mm}, label-width={3mm}, column-sep={-10pt}, item-indent={0mm} , ,after-item-skip=-1mm](2)
\task $\Int x(\sqrt[3]{x}-1)^3 dx$;
\task $\Int \dfrac{dx}{\textrm{sh}x+2\textrm{ch}x}$;
\task $\Int x^2\sqrt{x^3-4}dx$;
\task $\Int\dfrac{\textrm{tg}^5x}{\cos^4x}dx$.
\end{tasks}
\textbf{2.}\hspace{1mm}Integrate the quadratic fractions:
\begin{tasks}[ label-align=left, label-offset={1mm}, label-width={3mm},column-sep={-20pt}, item-indent={0mm} ,, ](2)
\task $\Int\dfrac{-2x+6}{x^2-9} dx$;
\task $\Int\dfrac{dx}{\sqrt{15-4x-4x^2}}$;
\task $\Int\dfrac{(4x-7)dx}{x^2+6x+10}$;
\task $\Int\dfrac{(4x-5)dx}{\sqrt{x^2+2x-3}}$.
\end{tasks}
\textbf{3.}\hspace{1mm}Integrate by parts or using the \\ suitable substitutions:
\begin{tasks}[ label-align=left, label-offset={1mm}, label-width={3mm}, column-sep={0pt}, item-indent={0mm},  ,after-item-skip=-1mm,](2)
\task  $\Int\dfrac{3x+7}{\sqrt{x+4}}dx$;
\task  $\Int\dfrac{\ln x-2}{x(4+\ln^2 x)}dx$;
\task $\Int 6x^2 \sin 2x dx$;
\task  $\Int x^2(2^x-8x)dx$;
\task!  $\Int\ln(x-\sqrt{x^2-4})dx$;
\task  $\Int x\dfrac{\sin x}{\cos^5 x}dx.$
\end{tasks}
\textbf{4.}\hspace{1mm}Integrate the polynomial fractions:
\begin{tasks}[ label-align=left, label-offset={1mm}, label-width={3mm},after-item-skip=-0.5mm, item-indent={0mm} , ](1)
\task $\Int\dfrac{13x+19}{(x^2+x-6)(x+1)} dx$;
\task $\Int\dfrac{2x^3+3x^2-28x-22}{x^3+x^2-8x-12} dx$;
\task $\Int\dfrac{8x^2+11x+11}{(x^2-2x+10)(x+3)}dx$.
\end{tasks}
\textbf{5.}\hspace{1mm}Integrate trigonometric expressions:
\begin{tasks}[ label-align=left, label-offset={1mm}, label-width={3mm},after-item-skip=-1mm, item-indent={0mm} ,, ](2)
\task $\Int\sin 2x\cos5x dx$;
\task $\Int\sqrt[4]{\dfrac{\sin x}{\cos^9 x}} dx$;
\task $\Int\dfrac{dx}{5\sin x-12\cos x}$.
\end{tasks}
\textbf{6.}\hspace{1mm}Integrate the functions with radicals:
\begin{tasks}[ label-align=left, label-offset={1mm}, label-width={3mm}, item-indent={0mm},column-sep={-20pt} ,, ](2)
\task $\Int\dfrac{\sqrt{9+x^2}}{x^4}dx$;
\task $\Int\dfrac{dx}{\sqrt{x}-2\sqrt[3]{x}}$.
\end{tasks}
\end{minipage}
\vline\;
      \begin{minipage}[t]{100mm}
\textbf{7.}\hspace{1mm}Solve the definite integrals:
\begin{tasks}[ label-align=left, label-offset={1mm}, label-width={3mm},column-sep={0pt}, item-indent={0mm},after-item-skip=-1mm,  ](2)
\task $\Int_{-1}^{\ln2}(2x+2) e^{2x}dx$;
\task $\Int_{0}^{2\pi}\sin^4(\frac{x}{2})dx$;
\task $\Int_{0}^{1/\sqrt{2}}\dfrac{3\textrm{acrsin}^2x-x}{\sqrt{1-x^2}}dx$;
\task $\Int_{0}^{\pi^3/8}\cos\sqrt[3]{x}dx.$
\end{tasks}

\textbf{8.}\hspace{1mm}Find the area of the figure bounded \\
by the curves:
\begin{tasks}[ label-align=left, label-offset={2mm}, label-width={3mm}, item-indent={0mm}, after-item-skip=-1mm,     ](1)
\task $y=(x+2)^2,\ y=4x+8$;
\task  $\left\lbrace\begin{array}{l} x=2\cos t,\\ y=5\sin t; \end{array}\right. \ x=1 \ (x\geq1)$;
\task  $\rho=4\sin 3\phi.$
\end{tasks}

\textbf{9.}\hspace{1mm}Find the arc-length of the curve:
\begin{tasks}[ label-align=left, label-offset={2mm}, label-width={3mm}, item-indent={0mm}, after-item-skip=-1mm,     ](1)
\task  $y=\frac{1}{4}\ln(\cos4x),\ -\frac{\pi}{24}\leq x\leq \frac{\pi}{24}$;
\task  $\left\lbrace\begin{array}{l} x=4e^{-t}\cos t,\\ y=4e^{-t}\sin t; \end{array}\right. \ 0\leq t \leq 2\pi$;
\task  $\rho=6\sin^2\big(\dfrac{\phi}{2}\big)$.
\end{tasks}

\textbf{10.}\hspace{1mm}Find the area of the surface formed \\ by rotating the curves around the $l$-axis:
\begin{tasks}[ label-align=left, label-offset={2mm}, label-width={3mm}, item-indent={0mm}, after-item-skip=-1mm,    ](1)
\task $y=1+\frac{1}{4}\textrm{ch}4x, \ -1\leq x\leq 1, \ l=OX$;
\task $\left\lbrace\begin{array}{l} x=3(t-\sin t),\\ y=3(1-\cos t); \end{array}\right. \ 0\leq t\leq2\pi,\\ l=OX$;
\task $\rho=2\sqrt{\cos2\phi}, \ l=o\rho$.
\end{tasks}

\textbf{11.}\hspace{1mm}Find the volume of the body formed \\
by rotating the curves around the $l$-axis:
\begin{tasks}[ label-align=left, label-offset={2mm}, label-width={3mm}, item-indent={0mm}, after-item-skip=-1mm,     ](1)
\task $y=-x^2+5x+1,\ xy=5, \ l=OX$;
\task $y=4-x^2,\ y=8-2x^2, \ l=OY$;
\task $\rho=6\phi,\ 0\leq\phi\leq\pi, \ l=o\rho$.
\end{tasks}

\textbf{12.}\hspace{1mm}Solve the improper integrals:
\begin{tasks}[ label-align=left, label-offset={1mm}, label-width={3mm}, item-indent={0mm},   ,   ,column-sep={-40pt}](2)
\task $\Int_{1}^{\infty}\dfrac{dx}{x^2(x^2+1)}$;
\task  $\Int_{0}^{1}\Big(\dfrac{1}{5}\Big)^{\frac{1}{x}}\dfrac{dx}{x^2}$.
\end{tasks}
\end{minipage}

\section{Personal task 18}

   \begin{minipage}[t]{90mm}
\textbf{1.}\hspace{1mm}Integrate using the table and \\ substitution under differential:
\begin{tasks}[ label-align=left, label-offset={1mm}, label-width={3mm}, column-sep={10pt}, item-indent={0mm} , ,after-item-skip=-1mm](2)
\task $\Int\dfrac{3x^4-x^23^x}{x^2} dx$;
\task $\Int \dfrac{4-3\cos^2x}{1-\cos2x} dx$;
\task $\Int \dfrac{\ln^5(1-3x)}{1-3x}dx$;
\task $\Int(\textrm{cth}x+\textrm{th}x)^2 dx$.
\end{tasks}

\textbf{2.}\hspace{1mm}Integrate the quadratic fractions:
\begin{tasks}[ label-align=left, label-offset={1mm}, label-width={3mm},column-sep={0pt}, item-indent={0mm} ,, ](2)
\task $\Int\dfrac{-4x-8}{x^2-16} dx$;
\task $\Int\dfrac{dx}{\sqrt{9x^2-6x-8}}$;
\task $\Int\dfrac{(-4x-11)dx}{x^2+4x+13}$;
\task $\Int\dfrac{(-2x+10)dx}{\sqrt{15+2x-x^2}}$.
\end{tasks}
\textbf{3.}\hspace{1mm}Integrate by parts or using the \\ suitable substitutions:
\begin{tasks}[ label-align=left, label-offset={1mm}, label-width={3mm}, column-sep={0pt}, item-indent={0mm},  ,after-item-skip=-1mm,](2)
\task $\Int\dfrac{e^{2x}-1+e^x}{\sqrt{1-e^{2x}}}dx$;
\task $\Int\dfrac{\sqrt{x^2+4}}{x}dx$;
\task! $\Int (2x^2-6x+8) \cos 2xdx$;
\task $\Int \dfrac{3-6x^2}{e^{3x}}dx$;
\task $\Int\dfrac{\textrm{arccos}\sqrt{x}}{\sqrt{x}}dx$
\task! $\Int 6(x^2+3)\ln(x^2+9)dx$.
\end{tasks}
\textbf{4.}\hspace{1mm}Integrate the polynomial fractions:
\begin{tasks}[ label-align=left, label-offset={1mm}, label-width={3mm},after-item-skip=-0.5mm, item-indent={0mm} , ](1)
\task $\Int\dfrac{9x^2-31x+20}{(x^2-5x+6)(x+1)} dx$;
\task $\Int\dfrac{2x^3+x^2+5x+2}{x^3-x^2-x+1} dx$;
\task $\Int\dfrac{7x^2-36x+13}{x^3-6x^2+13x}dx$.
\end{tasks}
\textbf{5.}\hspace{1mm}Integrate trigonometric expressions:
\begin{tasks}[ label-align=left, label-offset={1mm}, label-width={3mm},after-item-skip=-1mm, item-indent={0mm} ,, ](2)
\task $\Int\sin5x\sin3xdx$;
\task  $\Int\dfrac{\cos^3 x}{\sin^9x}dx$;
\task  $\Int\dfrac{dx}{5+3\sin x}$.
\end{tasks}
\textbf{6.}\hspace{1mm}Integrate the fractions with radicals:
\begin{tasks}[ label-align=left, label-offset={1mm}, label-width={3mm}, item-indent={0mm},column-sep={0in} , ](2)
\task $\Int\dfrac{4(x+2)^{-2}dx}{\sqrt{8+4x+x^2}}$;
\task $\Int\dfrac{(\sqrt{x}+1)dx}{\sqrt[4]{x^3}+4\sqrt[4]{x}}$.
\end{tasks}
\end{minipage}
\vline\;
      \begin{minipage}[t]{100mm}
\textbf{7.}\hspace{1mm}Solve the definite integrals:
\begin{tasks}[ label-align=left, label-offset={1mm}, label-width={3mm},column-sep={-30pt}, item-indent={0mm},after-item-skip=-1mm, ](2)
\task  $\Int_{0}^{e-1}\dfrac{\ln(x+1)}{(x+1)^2}dx$;
\task  $\Int_{0}^{\pi} \cos^4x\sin^4x dx$;
\task  $\Int_{1}^{9} \dfrac{x^2-16}{\sqrt{x}+2}dx$;
\task  $\Int_{0}^{\pi/4}\dfrac{(e^{-\textrm{tg}x}+\sin x)dx}{\cos^2x}.$
\end{tasks}

\textbf{8.}\hspace{1mm}Find the area of the figure bounded \\
by the curves:
\begin{tasks}[ label-align=left, label-offset={2mm}, label-width={3mm}, item-indent={0mm}, after-item-skip=-1mm,     ](1)
\task $y=(x+1)^3,\ y=5(x+1)$;
\task $\left\lbrace\begin{array}{l} x=3\cos^3t,\\ y=4\sin^3t; \end{array}\right. \ y=2 \ (y\geq 2)$;
\task $\rho=4\sqrt{\sin2\phi}.$
\end{tasks}

\textbf{9.}\hspace{1mm}Find the arc-length of the curve:
\begin{tasks}[ label-align=left, label-offset={2mm}, label-width={3mm}, item-indent={0mm}, after-item-skip=-1mm,     ](1)
\task  $y=2-\frac{1}{3}\textrm{ch} 3x,\ 0\leq x\leq \frac{1}{3}\ln3$;
\task   $\left\lbrace\begin{array}{l} x=3(t-\sin t),\\ y=3(1-\cos t); \end{array}\right. \ 0\leq t\leq 2\pi$;
\task   $\rho=\sqrt{2}e^{-2\phi}, \ 0\leq\phi\leq\pi$.
\end{tasks}

\textbf{10.}\hspace{1mm}Find the area of the surface formed \\ by rotating the curves around the $l$-axis:
\begin{tasks}[ label-align=left, label-offset={2mm}, label-width={3mm}, item-indent={0mm}, after-item-skip=-1mm,    ](1)
\task $y=x^3, \ -3\leq x\leq 3, \ l=OX$;
\task $\left\lbrace\begin{array}{l} x=3+\cos t,\\ y=2+\sin t; \end{array}\right. \ l=OY$;
\task $\rho=6(1+\cos\phi), \ l=o\rho$.
\end{tasks}

\textbf{11.}\hspace{1mm}Find the volume of the body formed \\
by rotating the curves around the $l$-axis:
\begin{tasks}[ label-align=left, label-offset={2mm}, label-width={3mm}, item-indent={0mm}, after-item-skip=-1mm,     ](1)
\task $y=-x^2+6x,\ y=-2x^2+12x, \ l=OX$;
\task $y=x^2, \ y=3x^2-18, \ l=OY$;
\task $\rho=8\cos\phi, \ l=o\rho$.
\end{tasks}

\textbf{12.}\hspace{1mm}Solve the improper integrals:
\begin{tasks}[ label-align=left, label-offset={1mm}, label-width={3mm}, item-indent={0mm},   ,   ,column-sep={-40pt}](2)
\task $\Int_{-\infty}^{0}x^3e^{-x^2}dx$;
\task $\Int_{0}^{1}\dfrac{4x-2}{\sqrt{2x-x^2}}dx$.
\end{tasks}
\end{minipage}

\section{Personal task 19}

   \begin{minipage}[t]{90mm}
\textbf{1.}\hspace{1mm}Integrate using the table and \\ substitution under differential:
\begin{tasks}[ label-align=left, label-offset={1mm}, label-width={3mm}, column-sep={-20pt}, item-indent={0mm} , ,after-item-skip=-1mm](2)
\task $\Int \dfrac{(1-3\sqrt{x})^3dx}{\sqrt{x}}$;
\task $\Int \dfrac{18\textrm{ctg}^5x}{1-\cos^2x} dx$;
\task $\Int\dfrac{x^3-x}{\sqrt{1-x^2}}dx$;
\task $\Int (\cos x+\dfrac{3}{\cos x})^2 dx.$
\end{tasks}
\textbf{2.}\hspace{1mm}Integrate the quadratic fractions:
\begin{tasks}[ label-align=left, label-offset={1mm}, label-width={3mm},column-sep={-20pt}, item-indent={0mm} ,, ](2)
\task $\Int\dfrac{-2x+6}{x^2-9} dx$;
\task $\Int\dfrac{dx}{\sqrt{3-8x-16x^2}}$;
\task $\Int\dfrac{(4-2x)dx}{x^2+4x+20}$;
\task $\Int\dfrac{(6x-9)dx}{\sqrt{x^2-4x+5}}$.
\end{tasks}

\textbf{3.}\hspace{1mm}Integrate by parts or using the \\ suitable substitutions:
\begin{tasks}[ label-align=left, label-offset={1mm}, label-width={3mm}, column-sep={0pt}, item-indent={0mm},  ,after-item-skip=-1mm,](2)
\task $\Int\dfrac{2-5x}{(x+2)^7}dx$;
\task $\Int\dfrac{(2+4\ln x) dx}{x\sqrt{\ln^2x+4}}$;
\task $\Int 16x\sin4x dx$;
\task $\Int \dfrac{x^2-4x}{e^{2x}}dx$;
\task $\Int \ln(\sqrt[5]{x}+1)dx$;
\task $\Int\dfrac{2x+4}{\sin^2x}dx$.
\end{tasks}

\textbf{4.}\hspace{1mm}Integrate the polynomial fractions:
\begin{tasks}[ label-align=left, label-offset={1mm}, label-width={3mm},after-item-skip=-0.5mm, item-indent={0mm} , ](1)
\task $\Int\dfrac{7x^2+24x-15}{x^3+2x^2-3x} dx$;
\task $\Int\dfrac{4x^3+17x^2+27x+16}{(x^2+3x+2)(x+1)} dx$;
\task $\Int\dfrac{x^2+14x-15}{(x^2-2x+5)(x+2)}dx$.
\end{tasks}

\textbf{5.}\hspace{1mm}Integrate trigonometric expressions:
\begin{tasks}[ label-align=left, label-offset={1mm}, label-width={3mm},after-item-skip=-1mm, item-indent={0mm} ,, ](2)
\task $\Int\cos x\sin 9xdx$;
\task $\Int\dfrac{\cos^2 x}{\sin^6 x} dx$;
\task $\Int\dfrac{4dx}{4+4\sin x+5\cos x}$.
\end{tasks}
\textbf{6.}\hspace{1mm}Integrate the fractions with radicals:
\begin{tasks}[ label-align=left, label-offset={1mm}, label-width={3mm}, item-indent={0mm},column-sep={-10pt} , ](2)
\task $\Int\dfrac{dx}{x^3\sqrt{1-x^2}}$;
\task $\Int\sqrt{\dfrac{x+2}{x-1}}\dfrac{dx}{(x+2)^2}$.
\end{tasks}
\end{minipage}
\vline \;
      \begin{minipage}[t]{100mm}
\textbf{7.}\hspace{1mm}Solve the definite integrals:
\begin{tasks}[ label-align=left, label-offset={1mm}, label-width={3mm},column-sep={-40pt}, item-indent={0mm},after-item-skip=-1mm, ,column-sep={0pt} ](2)
\task  $\Int_{\pi^2/9}^{\pi^2/4}(1-\sin\sqrt{x})\dfrac{dx}{\sqrt{x}}$;
\task  $\Int_{0}^{\pi/2}x\cos^2xdx$;
\task  $\Int_{0}^{\sqrt{3}}\dfrac{2x+4}{x^2+1}dx$;
\task  $\Int_{0}^{3}x^5\sqrt{(9-x^2)^3}dx.$
\end{tasks}

\textbf{8.}\hspace{1mm}Find the area of the figure bounded
 \\ by the curves:
\begin{tasks}[ label-align=left, label-offset={2mm}, label-width={3mm}, item-indent={0mm}, after-item-skip=-1mm,     ](1)
\task $y=-x^2+4x+1,\ xy=4$;
\task $\left\lbrace\begin{array}{l} x=2\cos t,\\ y=2+4\sin t; \end{array}\right. \ x=1 \ (x\geq1)$;
\task $\rho=3\cos 3\phi.$
\end{tasks}

\textbf{9.}\hspace{1mm}Find the arc-length of the curve:
\begin{tasks}[ label-align=left, label-offset={2mm}, label-width={3mm}, item-indent={0mm}, after-item-skip=-1mm,     ](1)
\task $y=2-e^{2x},\ \frac{1}{4}\ln2\leq x\leq\frac{1}{4}\ln6$;
\task $\left\lbrace\begin{array}{l} x=3(2\cos t-\cos2t),\\\ y=3(2\sin t-\sin2t); \end{array}\right. \ 0\leq t\leq 2\pi$;
\task $\rho=2\cos\phi+2\sin\phi$.
\end{tasks}

\textbf{10.}\hspace{1mm}Find the area of the surface formed \\ by rotating the curves around the $l$-axis:
\begin{tasks}[ label-align=left, label-offset={2mm}, label-width={3mm}, item-indent={0mm}, after-item-skip=-1mm,    ](1)
\task $y=3x+9, \ -3\leq x\leq 3, \ l=OX$;
\task $\left\lbrace\begin{array}{l} x=4\cos^3t,\\ y=4\sin^3t; \end{array}\right.  \ l=OY$;
\task $\rho=\sqrt{\sin2\big(\phi-\frac{\pi}{4}\big)}, \ l=o\rho$.
\end{tasks}

\textbf{11.}\hspace{1mm}Find the volume of the body formed \\
by rotating the curves around the $l$-axis:
\begin{tasks}[ label-align=left, label-offset={2mm}, label-width={3mm}, item-indent={0mm}, after-item-skip=-1mm,     ](1)
\task  $y=(x-2)^2,\ y=x+4, \ l=OX$;
\task  $y=x^2-8,\ y=2x,\ x=0 \ (x\geq0),\\ l=OY$;
\task  $\rho=\phi,\ 0\leq\phi\leq\pi,\ l=o\rho$.
\end{tasks}

\textbf{12.}\hspace{1mm}Solve the improper integrals:
\begin{tasks}[ label-align=left, label-offset={1mm}, label-width={3mm}, item-indent={0mm},   ,   ,column-sep={-40pt}](2)
\task  $\Int_{-1/2}^{\infty}\dfrac{\textrm{arctg}^3\big(\dfrac{x}{2}\big)}{4+x^2}dx$;
\task $\Int_{1}^{2}\dfrac{(2x+1)dx}{\sqrt{4x-x^2-3}}$.
\end{tasks}

\end{minipage}

\section{Personal task 20}

   \begin{minipage}[t]{90mm}
\textbf{1.}\hspace{1mm}Integrate using the table and \\ substitution under differential:
\begin{tasks}[ label-align=left, label-offset={1mm}, label-width={3mm}, column-sep={5pt}, item-indent={0mm} , ,after-item-skip=-1mm](2)
\task $\Int\sqrt{x}(2-x\sqrt{x})^2 dx$;
\task $\Int \dfrac{1+\cos2x}{1-\cos2x}dx$;
\task $\Int \dfrac{-8x}{\sqrt[3]{1-4x^2}}dx$;
\task $\Int(1-3\sin\dfrac{x}{2})^2dx$.
\end{tasks}
\textbf{2.}\hspace{1mm}Integrate the quadratic fractions:
\begin{tasks}[ label-align=left, label-offset={1mm}, label-width={3mm},column-sep={0pt}, item-indent={0mm} ,, ](2)
\task $\Int\dfrac{7-6x}{x^2-1} dx$;
\task $\Int\dfrac{dx}{\sqrt{12+4x-4x^2}}$;
\task $\Int\dfrac{-4xdx}{x^2-2x+5}$;
\task $\Int\dfrac{(14-6x)dx}{\sqrt{x^2-8x+25}}$.
\end{tasks}
\textbf{3.}\hspace{1mm}Integrate by parts or using the \\ suitable substitutions:
\begin{tasks}[ label-align=left, label-offset={1mm}, label-width={3mm}, column-sep={-20pt}, item-indent={0mm},  ,after-item-skip=-1mm,](2)
\task  $\Int\dfrac{(1-e^x)^2}{1+e^x}dx$;
\task $\Int\dfrac{dx}{(x-3)^2\sqrt{6x-x^2}}$;
\task! $\Int (4-8x^2) \cos 4x dx$;
\task $\Int \ln^2(x+2) dx$;
\task $\Int 2^{\sqrt{x+2}}dx$;
\task! $\Int \arccos(x^2-1)dx$.
\end{tasks}

\textbf{4.}\hspace{1mm}Integrate the polynomial fractions:
\begin{tasks}[ label-align=left, label-offset={1mm}, label-width={3mm},after-item-skip=-0.5mm, item-indent={0mm} , ](1)
\task $\Int\dfrac{4x^2-8x-36}{x^3+x^2-9x-9} dx$;
\task $\Int\dfrac{2x^3-9x^2+x+22}{x^3-2x^2-4x-8} dx$;
\task $\Int\dfrac{4x^2+10x+26}{(x^2+6x+13)(x-1)}dx$.
\end{tasks}

\textbf{5.}\hspace{1mm}Integrate trigonometric expressions:
\begin{tasks}[ label-align=left, label-offset={1mm}, label-width={3mm},after-item-skip=-1mm, item-indent={0mm} ,, ](2)
\task $\Int\cos 2x\cos5xdx$;
\task $\Int\sqrt[3]{\cos^4 x}\sin^3 x dx$;
\task! $\Int\dfrac{12dx}{4\cos^2x-4\sin x\cos x-8\sin^2x}$.
\end{tasks}
\textbf{6.}\hspace{1mm}Integrate the fractions with radicals:
\begin{tasks}[ label-align=left, label-offset={1mm}, label-width={3mm}, item-indent={0mm},column-sep={0in} , ](2)
\task $\Int \dfrac{x^2dx}{\sqrt{(4-x^2)^3}}$;
\task $\Int\dfrac{x dx}{(x+5)\sqrt{9+x}}$.
\end{tasks}
\end{minipage}
\vline\;
      \begin{minipage}[t]{100mm}
\textbf{7.}\hspace{1mm}Solve the definite integrals:
\begin{tasks}[ label-align=left, label-offset={1mm}, label-width={3mm},column-sep={-20pt}, item-indent={0mm},after-item-skip=-1mm](2)
\task $\Int_{1}^{\ln2+1}(x-1) e^{x-1}dx$;
\task $\Int_{0}^{\pi}\sin^4x\cos^2xdx$;
\task $\Int_{1}^{5/2}\dfrac{(2x+7)dx}{\sqrt{9+2x-x^2}}$;
\task $\Int_{-1}^{\sqrt{8}}\dfrac{\sqrt{x+1}-2}{\sqrt{x+1}+2}dx.$
\end{tasks}

\textbf{8.}\hspace{1mm}Find the area of the figure bounded \\
by the curves:
\begin{tasks}[ label-align=left, label-offset={2mm}, label-width={3mm}, item-indent={0mm}, after-item-skip=-1mm,     ](1)
\task $2x+y=7,\ xy=3$;
\task $\left\lbrace\begin{array}{l} x=4\cos t,\\ y=2+4\sin t; \end{array}\right. \ x=2 \ (x\geq2)$;
\task $\rho=4(1-\sin\phi)$.
\end{tasks}

\textbf{9.}\hspace{1mm}Find the arc-length of the curve:
\begin{tasks}[ label-align=left, label-offset={2mm}, label-width={3mm}, item-indent={0mm}, after-item-skip=-1mm,     ](1)
\task $y=4\ln(16-x^2),\ -2\leq x\leq 2$;
\task $\left\lbrace\begin{array}{l} x=3\cos^3 t,\\ y=3\sin^3 t; \end{array}\right. \ 0\leq t\leq \frac{\pi}{2}$;
\task $\rho=\phi, \ 0\leq\phi\leq 2\pi$.
\end{tasks}

\textbf{10.}\hspace{1mm}Find the area of the surface formed \\ by rotating the curves around the $l$-axis:
\begin{tasks}[ label-align=left, label-offset={2mm}, label-width={3mm}, item-indent={0mm}, after-item-skip=-1mm,    ](1)
\task $y=4\textrm{ch}\big(\dfrac{x}{4}\big), \ 0\leq x\leq 4, \ l=OX$;
\task $\left\lbrace\begin{array}{l} x=2(t-\sin t),\\ y=2(1-\cos t); \end{array}\right. \ 0\leq t\leq2\pi,\\ l=OX$;
\task $\rho=8\cos\phi, \ l=o\rho.$
\end{tasks}

\textbf{11.}\hspace{1mm}Find the volume of the body formed \\
by rotating the curves around the $l$-axis:
\begin{tasks}[ label-align=left, label-offset={2mm}, label-width={3mm}, item-indent={0mm}, after-item-skip=-1mm,     ](1)
\task  $y=x^2-6x,\ y=2x^2-12x, \ l=OX$;
\task $y=9-x^2,\ x=8,\ l=OY$;
\task $\rho=4\sin 2\phi, \ l=o\rho$.
\end{tasks}

\textbf{12.}\hspace{1mm}Solve the improper integrals:
\begin{tasks}[ label-align=left, label-offset={1mm}, label-width={3mm}, item-indent={0mm},   ,   ,column-sep={-40pt}](2)
\task $\Int_{0}^{\infty}\dfrac{-2x+5}{x^2+1}dx$;
\task  $\Int_{0}^{\pi/2}\dfrac{3\cos x}{\sqrt[4]{1-\sin x}}dx$.
\end{tasks}
\end{minipage}

\section{Personal task 21}

   \begin{minipage}[t]{90mm}
\textbf{1.}\hspace{1mm}Integrate using the table and \\ substitution under differential:
\begin{tasks}[ label-align=left, label-offset={1mm}, label-width={3mm},column-sep={-30pt}, item-indent={0mm} , ](2)
\task $\Int \dfrac{(1-3x^2)^3}{x^5} dx$;
\task $\Int \dfrac{2-\sin x}{\sin^2x} dx$;
\task $\Int x^3e^{2x^4+5}dx$;
\task $\Int\dfrac{(\sin2x-6\sin x)dx}{\cos^2x+9}$.
\end{tasks}
\textbf{2.}\hspace{1mm}Integrate the quadratic fractions:
\begin{tasks}[ label-align=left, label-offset={1mm}, label-width={3mm},column-sep={-18pt}, item-indent={0mm} , ](2)
\task $\int\dfrac{4-8x}{x^2-4} dx$;
\task $\Int\dfrac{dx}{\sqrt{8+6x-9x^2}}$;
\task $\Int\dfrac{(6x+7)dx}{x^2+4x+5}$
\task $\Int\dfrac{(2x-6)dx}{\sqrt{x^2-2x-3}}$.
\end{tasks}
\textbf{3.}\hspace{1mm}Integrate by parts or using the \\ suitable substitutions:
\begin{tasks}[ label-align=left, label-offset={1mm}, label-width={3mm}, column-sep={0pt}, item-indent={0mm},  ](2)
\task $\Int\sqrt{e^{2x}+16}dx$;
\task $\Int\dfrac{dx}{x^2\sqrt{1-4x^2}}$;
\task $\Int (x^3-6) \sin x dx$;
\task $\Int x\ln(x+1)dx$;
\task $\Int\arcsin\sqrt{x}dx$;
\task $\Int\cos (\frac{1}{2}\ln x)dx$.
\end{tasks}
\textbf{4.}\hspace{1mm}Integrate the polynomial fractions:
\begin{tasks}[ label-align=left, label-offset={1mm}, label-width={3mm}, item-indent={0mm} , ](1)
\task $\Int\dfrac{x^2+5x-20}{(x^2-x-2)(x-3)} dx$;
\task $\Int\dfrac{-2x^3-10x^2+x+20}{x^3+3x^2-4}dx$;
\task $\Int\dfrac{-3x^2-4x+44}{(x^2-8x+20)(x+2)}$.
\end{tasks}
\textbf{5.}\hspace{1mm}Integrate trigonometric expressions:
\begin{tasks}[ label-align=left, label-offset={1mm}, label-width={3mm},after-item-skip=-1mm, item-indent={0mm} , ](2)
\task $\Int\sin5x\cos 2xdx$;
\task $\Int\cos^4 2x\sin^2 2x dx$;
\task $\Int\dfrac{-6\sin x dx}{10-\sin^2x+3\cos^2 x}$.
\end{tasks}
\textbf{6.}\hspace{1mm}Integrate the fractions with radicals:
\begin{tasks}[ label-align=left, label-offset={1mm}, label-width={3mm}, item-indent={0mm},column-sep={0pt} ](2)
\task $\Int\dfrac{\sqrt{x^2+9}}{x^4}dx$;
\task $\Int\dfrac{(4-\sqrt{x}) dx}{(2+\sqrt[3]{x})\sqrt{x}}$.
\end{tasks}
\end{minipage}
\vline\;
      \begin{minipage}[t]{100mm}
\textbf{7.}\hspace{1mm}Solve the definite integrals:
\begin{tasks}[ label-align=left, label-offset={1mm}, label-width={3mm}, item-indent={0mm},after-item-skip=-1mm,  , column-sep={0pt}](2)
\task $\Int_{1}^{\ln2+1}(x-1)^2e^{x-1}dx$;
\task $\Int_{0}^{\pi/4}\dfrac{\sin^3x}{\cos^7x}dx$;
\task $\Int_{1}^{\sqrt[3]{10}} x^2\sqrt{x^3-1}dx$;
\task $\Int_{0}^{3}\dfrac{\sqrt[3]{x+1} \ dx}{4+\sqrt[3]{(x+1)^2}}.$
\end{tasks}
\textbf{8.}\hspace{1mm}Find the area of the figure bounded \\ by the curves:
\begin{tasks}[ label-align=left, label-offset={2mm}, label-width={3mm}, item-indent={0mm}, after-item-skip=-1mm,     ](1)
\task $y=2x^2-9x+10,\ y=-x^2+3x+1$;
\task $\left\lbrace\begin{array}{l} x=2\cos^3t,\\ y=4\sin^3t; \end{array}\right. \ x=1 \ (x\geq 1)$;
\task $\rho=2\cos4\phi.$
\end{tasks}
\textbf{9.}\hspace{1mm}Find the arc-length of the curve:
\begin{tasks}[ label-align=left, label-offset={2mm}, label-width={3mm}, item-indent={0mm}, after-item-skip=-1mm,     ](1)
\task $y=2+\ln x,\ 1\leq x\leq \sqrt{15}$;
\task $\left\lbrace\begin{array}{l} x=3(t-\sin t),\\ y=3(1-\cos t); \end{array}\right. \ 0\leq t\leq2\pi$;
\task $\rho=\sqrt{2}e^{-2\phi}, \ 0\leq\phi\leq\pi$.
\end{tasks}

\textbf{10.}\hspace{1mm}Find the area of the surface formed \\ by rotating the curves around the $l$-axis:
\begin{tasks}[ label-align=left, label-offset={2mm}, label-width={3mm}, item-indent={0mm}, after-item-skip=-1mm,    ](1)
\task $y=\dfrac{\sqrt{x}}{6}(x-12), \ 0\leq x\leq 12, \ l=OX$;
\task $\left\lbrace\begin{array}{l} x=1+2\cos t,\\ y=2\sin t; \end{array}\right. \ (x\geq0) \ l=OY$;
\task $\rho=\sqrt{\cos2\phi}, \ l=o\rho$.
\end{tasks}
\textbf{11.}\hspace{1mm}Find the volume of the body formed \\ by rotating the curves around the $l$-axis:
\begin{tasks}[ label-align=left, label-offset={2mm}, label-width={3mm}, item-indent={0mm}, after-item-skip=-1mm,    ](1)
\task $y=x^2-6x+10,\ y=x+10, \ l=OX$;
\task $x=-y^2+6y-8,\ x=0, \ l=OY$;
\task $\rho=4(1+\cos\phi), \ l=o\rho$.
\end{tasks}
\textbf{12.}\hspace{1mm}Solve the improper integrals:
\begin{tasks}[ label-align=left, label-offset={1mm}, label-width={3mm}, item-indent={0mm},   ,   ,column-sep={0pt}](2)
\task $\Int_{0}^{\infty}\dfrac{x+\textrm{arcctg}x}{x^2+1}dx$;
\task $\Int_{0}^{1/3}\dfrac{-x+3}{\sqrt{1-9x^2}}dx$.
\end{tasks}
\end{minipage}

\section{Personal task 22}

   \begin{minipage}[t]{90mm}
\textbf{1.}\hspace{1mm}Integrate using the table and \\ substitution under differential:
\begin{tasks}[ label-align=left, label-offset={1mm}, label-width={3mm},column-sep={-18pt}, item-indent={0mm} , ](2)
\task $\Int \dfrac{(2+\sqrt[3]{x})^3 dx}{x}$;
\task $\Int\dfrac{dx}{\textrm{sh}x-\textrm{ch}x}$;
\task $\Int x^25^{3x^3-2}dx$;
\task $\Int\dfrac{\sin2x}{\sqrt{1-\cos^4x}}dx$.
\end{tasks}
\textbf{2.}\hspace{1mm}Integrate the quadratic fractions:
\begin{tasks}[ label-align=left, label-offset={1mm}, label-width={3mm},column-sep={-18pt}, item-indent={0mm} ,, ](2)
\task $\Int\dfrac{6x-5}{x^2-1} dx$;
\task $\Int\dfrac{dx}{\sqrt{5+8x-4x^2}}$;
\task $\Int\dfrac{(4x-2)dx}{x^2-4x+13}$;
\task $\Int\dfrac{(-4x-10)dx}{\sqrt{x^2+6x+5}}$.
\end{tasks}
\textbf{3.}\hspace{1mm}Integrate by parts or using the \\ suitable substitutions:
\begin{tasks}[ label-align=left, label-offset={1mm}, label-width={3mm}, column-sep={-20pt}, item-indent={0mm},  ](2)
\task $\Int\dfrac{(3x-7)dx}{(x+1)^5}$;
\task $\Int\dfrac{dx}{(x-3)^2\sqrt{6x-x^2}}$;
\task $\Int 6x^2 \sin 2xdx$;
\task $\Int e^x(\cos4x+2) dx$
\task $\Int x^6\ln xdx$;
\task $\Int\dfrac{\textrm{arctg} \sqrt{x}}{\sqrt{x}}dx$.
\end{tasks}
\textbf{4.}\hspace{1mm}Integrate the polynomial fractions:
\begin{tasks}[ label-align=left, label-offset={1mm}, label-width={3mm},after-item-skip=-0.5mm, item-indent={0mm} ,, ](1)
\task $\Int\dfrac{x^2+20x+24}{(x^2-3x-4)(x+2)} dx$;
\task $\Int\dfrac{-3x^3+3x^2+4x-1}{x^3-2x^2+x} dx$;
\task $\Int\dfrac{4x^2-28x-4}{(x^2+2x+2)(x-4)}dx$.
\end{tasks}
\textbf{5.}\hspace{1mm}Integrate trigonometric expressions:
\begin{tasks}[ label-align=left, label-offset={1mm}, label-width={3mm},after-item-skip=-1mm, item-indent={0mm} ,, ](2)
\task $\Int\sin x\sin 7xdx$;
\task $\Int\dfrac{\cos^3 x}{\sin^9 x}dx$;
\task $\Int\dfrac{dx}{5-2\sin x+4\cos x}$.
\end{tasks}
\textbf{6.}\hspace{1mm}Integrate the fractions with radicals:
\begin{tasks}[ label-align=left, label-offset={1mm}, label-width={3mm}, item-indent={0mm},column-sep={-30pt} ](2)
\task $\Int\dfrac{x^3dx}{\sqrt{(9-x^2)^3}}$;
\task $\Int\sqrt{\dfrac{x+1}{x-2}} \dfrac{dx}{(x+1)^2}$.
\end{tasks}
\end{minipage}
\vline\;
      \begin{minipage}[t]{100mm}
\textbf{7.}\hspace{1mm}Solve the definite integrals:
\begin{tasks}[ label-align=left, label-offset={1mm}, label-width={3mm},column-sep={-50pt}, item-indent={0mm},after-item-skip=-1mm,  ](2)
\task $\Int_{0}^{\log_23}x^2 2^xdx$;
\task $\Int_{0}^{2\pi}\cos^4(\frac{x}{4})\sin^2(\frac{x}{4})dx$;
\task $\Int_{0}^{\sqrt[4]{2e-1}}\dfrac{-12x^3}{x^4+1}dx$;
\task $\Int_{1}^{e}\dfrac{\ln x dx}{x\sqrt{1+\ln x}}.$
\end{tasks}

\textbf{8.}\hspace{1mm}Find the area of the figure bounded \\
by the curves:
\begin{tasks}[ label-align=left, label-offset={2mm}, label-width={3mm}, item-indent={0mm}, after-item-skip=-1mm,     ](1)
\task $y=(x+2)^2,\ y=4x+8$;
\task $\left\lbrace\begin{array}{l} x=4\cos t,\\ y=3\sin t; \end{array}\right. \ x=2 \ (x\geq2)$;
\task  $\rho=3\sqrt{\cos 2\phi}.$
\end{tasks}

\textbf{9.}\hspace{1mm}Find the arc-length of the curve:
\begin{tasks}[ label-align=left, label-offset={2mm}, label-width={3mm}, item-indent={0mm}, after-item-skip=-1mm,     ](1)
\task  $y=\sqrt{4-x^2},\ -\sqrt{2}\leq x\leq 2$;
\task  $\left\lbrace\begin{array}{l} x=3e^{-t}\cos t,\\ y=3e^{-t}\sin t; \end{array}\right. \ 0\leq t\leq 2\pi$;
\task  $\rho=6(1+\cos\phi)$.
\end{tasks}

\textbf{10.}\hspace{1mm}Find the area of the surface formed \\ by rotating the curves around the $l$-axis:
\begin{tasks}[ label-align=left, label-offset={2mm}, label-width={3mm}, item-indent={0mm}, after-item-skip=-1mm,    ](1)
\task $y=\sqrt{4x+12}, \ -3\leq x\leq 1, \ l=OX$;
\task $\left\lbrace\begin{array}{l} x=\cos^3 t,\\ y=\sin^3 t; \end{array}\right. \ l=OY$;
\task $\rho=2\sin\phi, \ l=o\rho$.
\end{tasks}

\textbf{11.}\hspace{1mm}Find the volume of the body formed \\
by rotating the curves around the $l$-axis:
\begin{tasks}[ label-align=left, label-offset={2mm}, label-width={3mm}, item-indent={0mm}, after-item-skip=-1mm,     ](1)
\task  $y=x^2+4x+4,\ y=4-x, \ l=OX$;
\task  $x=-y^2+3y,\ x=-2y^2+6y, \ l=OY$;
\task  $\rho=\sin2\phi, \ l=o\rho$.
\end{tasks}

\textbf{12.}\hspace{1mm}Solve the improper integrals:
\begin{tasks}[ label-align=left, label-offset={1mm}, label-width={3mm}, item-indent={0mm},   ,   ,column-sep={-40pt}](2)
\task $\Int_{1}^{\infty}\dfrac{-18x}{4x^4-9}dx$;
\task $\Int_{-\ln2}^{0}\dfrac{e^x}{\sqrt{1-e^{2x}}}dx$.
\end{tasks}
\end{minipage}

\section{Personal task 23}

   \begin{minipage}[t]{90mm}
\textbf{1.}\hspace{1mm}Integrate using the table and \\ substitution under differential:
\begin{tasks}[ label-align=left, label-offset={1mm}, label-width={3mm}, column-sep={5pt}, item-indent={0mm} ](2)
\task $\Int \dfrac{(3-2x^4)^2}{x^2}dx$;
\task $\Int(3+\textrm{ctg}^2x)dx$;
\task $\Int \dfrac{\ln^7(5-3x)}{5-3x}dx$;
\task $\Int\cos^2(2x+5)dx$.
\end{tasks}
\textbf{2.}\hspace{1mm}Integrate the quadratic fractions:
\begin{tasks}[ label-align=left, label-offset={1mm}, label-width={3mm},column-sep={-10pt}, item-indent={0mm} ,, ](2)
\task $\Int\dfrac{6x+5}{x^2-9} dx$;
\task $\Int\dfrac{dx}{\sqrt{5+4x-4x^2}}$;
\task $\Int\dfrac{6xdx}{x^2+2x+10}$;
\task $\Int\dfrac{(4x+4)dx}{\sqrt{7-6x-x^2}}$.
\end{tasks}
\textbf{3.}\hspace{1mm}Integrate by parts or using the \\ suitable substitutions:
\begin{tasks}[ label-align=left, label-offset={1mm}, label-width={3mm}, column-sep={-20pt}, item-indent={0mm},  ](2)
\task $\Int\dfrac{\sqrt{2x+3}dx}{4+\sqrt{2x+3}}$;
\task $\Int\dfrac{x dx}{x^2\sqrt{9+x^2}}$;
\task $\Int \dfrac{x^2-2x}{e^{2x}}dx$;
\task $\Int (8-4x) \sin 4x dx$;
\task $\Int 2x\dfrac{\cos x dx}{\sin^3x}$.
\task $\Int\arccos(1-x^2)dx$;
\end{tasks}
\textbf{4.}\hspace{1mm}Integrate the polynomial fractions:
\begin{tasks}[ label-align=left, label-offset={1mm}, label-width={3mm},after-item-skip=-0.5mm, item-indent={0mm} ,, ](1)
\task $\Int\dfrac{2x^2+10x-10}{(x^2-7x+10)(x+1)} dx$;
\task  $\Int\dfrac{-2x^3+2x^2+15x-16}{x^3-4x^2+4x} dx$;
\task  $\Int\dfrac{4x^2+19x-29}{(x^2-2x+10)(x+3)}dx$.
\end{tasks}
\textbf{5.}\hspace{1mm}Integrate trigonometric expressions:
\begin{tasks}[ label-align=left, label-offset={1mm}, label-width={3mm},after-item-skip=-1mm, item-indent={0mm} , ](2)
\task $\Int\cos2x\cos6xdx$;
\task $\Int\cos^4 x\sin^5 x dx$;
\task $\Int\dfrac{-6dx}{9\sin^2x-6\sin x\cos x-3\cos^2x}$.
\end{tasks}
\textbf{6.}\hspace{1mm}Integrate the fractions with radicals:
\begin{tasks}[ label-align=left, label-offset={1mm}, label-width={3mm}, item-indent={0mm},column-sep={-30pt} ](2)
\task $\Int\dfrac{dx}{x^3\sqrt{x^2-4}}$;
\task $\Int\dfrac{dx}{\sqrt{x}+4\sqrt[3]{x^2}}$.
\end{tasks}
\end{minipage}
\vline\;
      \begin{minipage}[t]{100mm}
\textbf{7.}\hspace{1mm}Solve the definite integrals:
\begin{tasks}[ label-align=left, label-offset={1mm}, label-width={3mm},column-sep={-40pt}, item-indent={0mm},after-item-skip=-1mm,  ](2)
\task  $\Int_{2}^{e+1}\ln^2(x-1)dx$;
\task $\Int_{\pi/4}^{\pi/2} \dfrac{\cos^2x}{\sin^4x} dx$;
\task $\Int_{0}^{3\sqrt{3}} \dfrac{x^2dx}{9+x^2}$;
\task $\Int_{0}^{4}\dfrac{6dx}{x+\sqrt{3x+4}}.$
\end{tasks}

\textbf{8.}\hspace{1mm}Find the area of the figure bounded \\
by the curves:
\begin{tasks}[ label-align=left, label-offset={2mm}, label-width={3mm}, item-indent={0mm}, after-item-skip=-1mm,     ](1)
\task $x+y=9,\ xy=8$;
\task $\left\lbrace\begin{array}{l} x=2\cos^3t,\\ y=2\sin^3t; \end{array}\right. \ x=1 \;(x\geq 1)$;
\task $\rho=4\sin2\phi.$
\end{tasks}

\textbf{9.}\hspace{1mm}Find the arc-length of the curve:
\begin{tasks}[ label-align=left, label-offset={2mm}, label-width={3mm}, item-indent={0mm}, after-item-skip=-1mm,     ](1)
\task $y=\ln(x+\sqrt{x^2-9}),\ 3\leq x\leq 3\sqrt{6}$;
\task $\left\lbrace\begin{array}{l} x=2(t-\sin t),\\ y=2(1-\cos t); \end{array}\right. \ 0\leq t\leq 2\pi$;
\task $\rho=3\sqrt{2}e^{-2\phi}, \ 0\leq\phi\leq \pi$.
\end{tasks}

\textbf{10.}\hspace{1mm}Find the area of the surface formed \\ by rotating the curves around the $l$-axis:
\begin{tasks}[ label-align=left, label-offset={2mm}, label-width={3mm}, item-indent={0mm}, after-item-skip=-1mm,    ](1)
\task $y=2\textrm{ch}\big(\dfrac{x}{2}\big), \ -2\leq x\leq 2, \ l=OX$;
\task $\left\lbrace\begin{array}{l} x=3+\cos t,\\ y=1+\sin t; \end{array}\right. \ l=OY$;
\task $\rho=4\sqrt{\cos2\phi}, \ l=o\rho$.
\end{tasks}

\textbf{11.}\hspace{1mm}Find the volume of the body formed \\
by rotating the curves around the $l$-axis:
\begin{tasks}[ label-align=left, label-offset={2mm}, label-width={3mm}, item-indent={0mm}, after-item-skip=-1mm,     ](1)
\task $y=x^2-2x+10,\ y=x+10, \ l=OX$;
\task $y=\ln x, \ x=0, \ y=0, \ y=\ln5, \ l=OY$;
\task $\rho=2(1-\cos\phi), \ l=o\rho$.
\end{tasks}
\textbf{12.}\hspace{1mm}Solve the improper integrals:
\begin{tasks}[ label-align=left, label-offset={1mm}, label-width={3mm}, item-indent={0mm},column-sep={-40pt}  ,   ](2)
\task $\Int_{0}^{\infty}\dfrac{2x dx}{\sqrt[7]{(x^2+1)^2}}$;
\task $\Int_{0}^{\pi^2/4}\dfrac{\sin\sqrt{x}-1}{\sqrt{x}}dx$.
\end{tasks}
\end{minipage}

\section{Personal task 24}

   \begin{minipage}[t]{90mm}
\textbf{1.}\hspace{1mm}Integrate using the table and \\ substitution under differential:
\begin{tasks}[ label-align=left, label-offset={1mm}, label-width={3mm}, column-sep={5pt}, item-indent={0mm} , ,after-item-skip=-1mm](2)
\task  $\Int x^2(4-3x^2)^3 dx$;
\task $\Int \dfrac{\sin 2x dx}{4-(\sin x-2)^2}$;
\task $\Int x^53^{5x^6-2}dx$;
\task $\Int\textrm{th}^4xdx$.
\end{tasks}
\textbf{2.}\hspace{1mm}Integrate the quadratic fractions:
\begin{tasks}[ label-align=left, label-offset={1mm}, label-width={3mm},column-sep={-20pt}, item-indent={0mm} ,, ](2)
\task $\Int\dfrac{-4x-8}{x^2-16} dx$;
\task $\Int\dfrac{dx}{\sqrt{12x-9x^2}}$;
\task $\Int\dfrac{(2x+6)dx}{x^2-2x+5}$;
\task $\Int\dfrac{(3x-9)dx}{\sqrt{x^2+4x-5}}$.
\end{tasks}
\textbf{3.}\hspace{1mm}Integrate by parts of using the \\ suitable substitutions:
\begin{tasks}[ label-align=left, label-offset={1mm}, label-width={3mm}, column-sep={-20pt}, item-indent={0mm},  ,after-item-skip=-1mm,](2)
\task $\Int\dfrac{(e^x+1)^2}{\sqrt{e^{2x}+1}}dx$;
\task $\Int\dfrac{dx}{(x-1)\sqrt{2x-x^2}}$;
\task! $\Int (4-2x^2)\cos 2x\sin 2x dx$;
\task $\Int (2-3x)3^{x}dx$;
\task $\Int \dfrac{4\ln x}{\sqrt[3]{x}}dx$;
\task! $\Int (5x^4-1)\textrm{arctg}x\ dx$.
\end{tasks}
\textbf{4.}\hspace{1mm}Integrate the polynomial fractions:
\begin{tasks}[ label-align=left, label-offset={1mm}, label-width={3mm},after-item-skip=-0.5mm, item-indent={0mm} , ](1)
\task $\Int\dfrac{-x^2+3x-6}{(x^2-3x+2)(x-3)} dx$;
\task $\Int\dfrac{3x^3-4x^2-8x-10}{x^3-3x-2} dx$;
\task $\Int\dfrac{4x^2-10x-40}{x^3+4x^2+20x} dx$.
\end{tasks}
\textbf{5.}\hspace{1mm}Integrate trigonometric expressions:
\begin{tasks}[ label-align=left, label-offset={1mm}, label-width={3mm},after-item-skip=-1mm, item-indent={0mm} ,column-sep={5pt}, ](2)
\task $\Int\cos3x\sin 9xdx$;
\task $\Int\dfrac{\cos^3 2x}{\sin^7 2x} dx$;
\task $\Int\dfrac{dx}{9-3\sin x+5\cos x}$.
\end{tasks}
\textbf{6.}\hspace{1mm}Integrate the fractions with radicals:
\begin{tasks}[ label-align=left, label-offset={1mm}, label-width={3mm}, item-indent={0mm},column-sep={-30pt} ,, ](2)
\task $\Int\dfrac{\sqrt{9-x^2}}{x^2}dx$;
\task $\Int\dfrac{\sqrt[4]{x}+\sqrt{x}}{1+\sqrt{x}} dx$.
\end{tasks}
\end{minipage}
\vline\;
      \begin{minipage}[t]{100mm}
\textbf{7.}\hspace{1mm}Solve the definite integrals:
\begin{tasks}[ label-align=left, label-offset={1mm}, label-width={3mm},column-sep={-40pt}, item-indent={0mm},after-item-skip=-1mm,  ](2)
\task $\Int_{0}^{2}(4x-8)e^{-2x}dx$;
\task $\Int_{0}^{\pi/6}\cos^2 3x\sin^5 3xdx$;
\task $\Int_{0}^{2}\dfrac{2x-6}{x^2-9}dx$;
\task $\Int_{-1}^{1}x^2\sqrt{1-x^2}dx.$
\end{tasks}
\textbf{8.}\hspace{1mm}Find the area of the figure bounded \\
by the curves:
\begin{tasks}[ label-align=left, label-offset={2mm}, label-width={3mm}, item-indent={0mm}, after-item-skip=-1mm,     ](1)
\task $y=(x-4)^2,\ y=3x-12$;
\task  $\left\lbrace\begin{array}{l} x=t-\sin t,\\ y=1-\cos t; \end{array}\right. \ y=0,\  0\leq t\leq 2\pi$;
\task  $\rho=3\cos 3\phi.$
\end{tasks}

\textbf{9.}\hspace{1mm}Find the arc-length of the curve:
\begin{tasks}[ label-align=left, label-offset={2mm}, label-width={3mm}, item-indent={0mm}, after-item-skip=-1mm,     ](1)
\task $y=\frac{1}{3}+\frac{1}{3}\ln(\sin 3x),\ \frac{\pi}{9}\leq x\leq \frac{2\pi}{9}$;
\task $\left\lbrace\begin{array}{l} x=2(2\cos t-\cos2t),\\\ y=2(2\sin t-\sin2t); \end{array}\right. \ 0\leq t \leq 2\pi$;
\task $\rho=4\phi, \ 0\leq \phi \leq \pi$.
\end{tasks}

\textbf{10.}\hspace{1mm}Find the area of the surface formed \\ by rotating the curves around the $l$-axis:
\begin{tasks}[ label-align=left, label-offset={2mm}, label-width={3mm}, item-indent={0mm}, after-item-skip=-1mm,    ](1)
\task $x^2+\dfrac{y^2}{16}=1, \ l=OY$;
\task $\left\lbrace\begin{array}{l} x=e^t\cos t,\\ y=e^t\sin t; \end{array}\right. \ 0\leq t\leq\pi, \ l=OX$;
\task $\rho=4\sqrt{\sin2\big(\phi-\frac{\pi}{4}\big)}, \ l=o\rho$.
\end{tasks}

\textbf{11.}\hspace{1mm}Find the volume of the body formed \\
by rotating the curves around the $l$-axis:
\begin{tasks}[ label-align=left, label-offset={2mm}, label-width={3mm}, item-indent={0mm}, after-item-skip=-1mm,     ](1)
\task $y=-x^2+6x+1,\ xy=6, \ l=OX$;
\task $x=4y-y^2,\ x=8y-2y^2, \ l=OY$;
\task $\rho=4\sin\phi, \ l=o\rho$.
\end{tasks}

\textbf{12.}\hspace{1mm}Solve the improper integrals:
\begin{tasks}[ label-align=left, label-offset={1mm}, label-width={3mm}, item-indent={0mm},   ,   ](2)
\task  $\Int_{-1}^{\infty}\dfrac{-8x-10}{4x^2-4x+5}dx$;
\task $\Int_{1/\sqrt{2}}^{1}\dfrac{dx}{x^2\sqrt{1-x^2}}$.
\end{tasks}
\end{minipage}

\section{Personal task 25}

   \begin{minipage}[t]{90mm}
\textbf{1.}\hspace{1mm}Integrate using the table and \\ substitution under differential:
\begin{tasks}[ label-align=left, label-offset={1mm}, label-width={3mm}, column-sep={-5pt}, item-indent={0mm} , ,after-item-skip=-1mm](2)
\task $\Int \dfrac{(3x^2e^x-e^{2x})dx}{e^x}$;
\task $\Int \dfrac{1-\textrm{tg}^2\frac{x}{2}}{1+\textrm{tg}^2\frac{x}{2}}dx$;
\task $\Int \dfrac{x^4 dx}{\sqrt[4]{x^5+9}}$;
\task $\Int\dfrac{\sin2x-2\sin x}{\cos^2x-4}dx$.
\end{tasks}
\textbf{2.}\hspace{1mm}Integrate the quadratic fractions:
\begin{tasks}[ label-align=left, label-offset={1mm}, label-width={3mm},column-sep={-20pt}, item-indent={0mm} ,, ](2)
\task $\Int\dfrac{4x-7}{x^2-1} dx$;
\task $\Int\dfrac{dx}{\sqrt{15-6x-9x^2}}$;
\task $\Int\dfrac{(-2x+4)dx}{x^2-8x+12}$;
\task $\Int\dfrac{(4x-2)dx}{\sqrt{x^2-4x-5}}$.
\end{tasks}
\textbf{3.}\hspace{1mm}Integrate by parts or using the \\ suitable substitutions:
\begin{tasks}[ label-align=left, label-offset={1mm}, label-width={3mm}, column-sep={-20pt}, item-indent={0mm},  ,after-item-skip=-1mm,](2)
\task $\Int\dfrac{\sqrt{x}+2}{x-9}dx$;
\task $\Int\dfrac{dx}{x\sqrt{4x^2+1}}$;
\task! $\Int (2-8x-12x^2) \cos2x dx$;
\task $\Int \dfrac{\sin3x}{e^{4x}} dx$;
\task $\Int \dfrac{\ln(x-2)}{(x-2)^3} dx$;
\task $\Int \dfrac{\arcsin x}{\sqrt{4x+4}} dx$.
\end{tasks}
\textbf{4.}\hspace{1mm}Integrate the polynomial fractions:
\begin{tasks}[ label-align=left, label-offset={1mm}, label-width={3mm},after-item-skip=-0.5mm, item-indent={0mm} , ](1)
\task $\Int\dfrac{10x+14}{(x^2+3x+2)(x-1)} dx$;
\task $\Int\dfrac{2x^3-2x^2-22x+4}{x^3+2x^2-4x-8} dx$;
\task $\Int\dfrac{4x^2-24x-54}{(x^2+6x+18)(x-3)}dx$.
\end{tasks}
\textbf{5.}\hspace{1mm}Integrate trigonometric expressions:
\begin{tasks}[ label-align=left, label-offset={1mm}, label-width={3mm},after-item-skip=-1mm, item-indent={0mm} ,, ](2)
\task $\Int\sin6x\sin 4xdx$;
\task $\Int\cos^4 3x\sin^2 3x dx$;
\task $\Int\dfrac{6\sin x dx}{9+5\sin^2 x+6\cos^2 x}$.
\end{tasks}
\textbf{6.}\hspace{1mm}Integrate the functions with radicals:
\begin{tasks}[ label-align=left, label-offset={1mm}, label-width={3mm}, item-indent={0mm},column-sep={-20pt} ,, ](2)
\task $\Int x^3\sqrt{4-x^2}dx$;
\task $\Int\sqrt{\dfrac{x+4}{x-1}}\dfrac{dx}{(x+4)^2}$.
\end{tasks}
\end{minipage}
\vline\;
      \begin{minipage}[t]{100mm}
\textbf{7.}\hspace{1mm}Solve the definite integrals:
\begin{tasks}[ label-align=left, label-offset={1mm}, label-width={3mm},column-sep={-20pt}, item-indent={0mm},after-item-skip=-1mm,  ](2)
\task $\Int_{0}^{1}(x-1)^2 e^{3x}dx$;
\task $\Int_{0}^{\pi/4}\dfrac{\sin^2x}{\cos^6x}dx$;
\task $\Int_{-2}^{0}\dfrac{2x dx}{\sqrt{12-4x-x^2}}$;
\task $\Int_{2}^{4}\dfrac{x+3}{\sqrt{x-1}}dx.$
\end{tasks}

\textbf{8.}\hspace{1mm}Find the area of the figure bounded \\
by the curves:
\begin{tasks}[ label-align=left, label-offset={2mm}, label-width={3mm}, item-indent={0mm}, after-item-skip=-1mm,     ](1)
\task $y=2x^2-9x+10,\ y=-x^2+3x+1$;
\task $\left\lbrace\begin{array}{l} x=6\cos t,\\ y=2\sin t; \end{array}\right. \ x=0\  (x\geq0)$;
\task $\rho=2(1+\sin\phi)$.
\end{tasks}

\textbf{9.}\hspace{1mm}Find the arc-length of the curve:
\begin{tasks}[ label-align=left, label-offset={2mm}, label-width={3mm}, item-indent={0mm}, after-item-skip=-1mm,     ](1)
\task $y=2-\frac{1}{2}\textrm{ch}2x,\ -\ln2 \leq x\leq \ln2$;
\task $\left\lbrace\begin{array}{l} x=\cos^3 t,\\ y=\sin^3 t; \end{array}\right. \ 0\leq t\leq \pi$;
\task $\rho=\sqrt{2}e^{-\phi}, \ 0\leq\phi\leq2\pi$.
\end{tasks}

\textbf{10.}\hspace{1mm}Find the area of the surface formed \\ by rotating the curves around the $l$-axis:
\begin{tasks}[ label-align=left, label-offset={2mm}, label-width={3mm}, item-indent={0mm}, after-item-skip=-1mm,    ](1)
\task $y=x^2, \ 0\leq x\leq 2\sqrt{5}, \ l=OY$;
\task $\left\lbrace\begin{array}{l} x=2(t-\sin t),\\ y=2(1-\cos t); \end{array}\right. \ \begin{array}{l} 0\leq t\leq2\pi,\\ l=OX;\end{array}$
\task $\rho=8\cos\phi, \ l=o\rho$.
\end{tasks}

\textbf{11.}\hspace{1mm}Find the volume of the body formed \\
by rotating the curves around the $l$-axis:
\begin{tasks}[ label-align=left, label-offset={2mm}, label-width={3mm}, item-indent={0mm}, after-item-skip=-1mm,     ](1)
\task $y=(x-1)^3,\ y=16x-16,\ (x\geq1) \\ l=OX$;
\task $x=8+y^2,\ x=16-y^2, \ l=OY$;
\task $\rho=\sin2\phi, \ l=o\rho$.
\end{tasks}

\textbf{12.}\hspace{1mm}Solve the improper integrals:
\begin{tasks}[ label-align=left, label-offset={1mm}, label-width={3mm}, item-indent={0mm},   ,   ,column-sep={0pt}](2)
\task $\Int_{-\sqrt{2}}^{\infty}\dfrac{4x^3-16x}{x^4+4}dx$;
\task  $\Int_{1}^{e^{\pi}}\dfrac{\sin(\ln x)dx}{x}$.
\end{tasks}
\end{minipage}

\section{Personal task 26}

   \begin{minipage}[t]{90mm}
\textbf{1.}\hspace{1mm}Integrate using the table and \\ substitution under differential:
\begin{tasks}[ label-align=left, label-offset={1mm}, label-width={3mm}, column-sep={0pt}, item-indent={0mm} , ,after-item-skip=-1mm](2)
\task $\Int \dfrac{(x-2)^3}{x\sqrt{x}} dx$
\task $\Int (4-3\textrm{ctg}^2 x)dx $;
\task $\Int \dfrac{x^3}{e^{3x^4-9}}dx$;
\task $\Int 8\sin^2\sqrt[3]{x}\cdot\dfrac{dx}{\sqrt[3]{x^2}}$.
\end{tasks}
\textbf{2.}\hspace{1mm}Integrate the quadratic fractions:
\begin{tasks}[ label-align=left, label-offset={1mm}, label-width={3mm},column-sep={0pt}, item-indent={0mm} ,, ](2)
\task $\Int\dfrac{2x-16}{x^2-4} dx$;
\task  $\Int\dfrac{dx}{\sqrt{15+4x-4x^2}}$;
\task  $\Int\dfrac{(-6x+14)dx}{x^2-6x+13}$;
\task  $\Int\dfrac{4x+6}{\sqrt{x^2+6x}} dx$.
\end{tasks}
\textbf{3.}\hspace{1mm}Integrate by parts or using the \\ suitable substitutions:
\begin{tasks}[ label-align=left, label-offset={1mm}, label-width={3mm}, column-sep={-20pt}, item-indent={0mm},  ,after-item-skip=-1mm,](2)
\task $\Int\dfrac{(2\ln x+8)dx}{x\sqrt{4-\ln^2x}}$;
\task $\Int\dfrac{dx}{(x-2)^2\sqrt{4x-x^2}}$;
\task $\Int 16x^2 \sin4x dx$;
\task $\Int \dfrac{4x^2-8x+12}{e^{2x}}dx$;
\task $\Int(x+3)^3\ln x dx$;
\task $\Int(\arccos 2x)^2 dx$.
\end{tasks}
\textbf{4.}\hspace{1mm}Integrate the polynomial fractions:
\begin{tasks}[ label-align=left, label-offset={1mm}, label-width={3mm},after-item-skip=-0.5mm, item-indent={0mm} , ](1)
\task $\Int\dfrac{2x^2-2x+12}{x^3-x^2-4x+4} dx$;
\task$\Int\dfrac{3x^3-10x^2-22x+6}{x^3-2x^2-7x-4} dx$;
\task $\Int\dfrac{8x^2-36x}{(x^2-4x+8)(x-4)}dx$.
\end{tasks}
\textbf{5.}\hspace{1mm}Integrate trigonometric fractions:
\begin{tasks}[ label-align=left, label-offset={1mm}, label-width={3mm},after-item-skip=-1mm, item-indent={0mm} ,, ](2)
\task $\Int\sin7x\cos4xdx$
\task $\Int\cos^2 3x\sin^4 3x dx$;
\task $\Int\dfrac{-18 dx}{9\sin^2x+6\sin x\cos x-8\cos^2x}$.
\end{tasks}
\textbf{6.}\hspace{1mm}Integrate the fractions with radicals:
\begin{tasks}[ label-align=left, label-offset={1mm}, label-width={3mm}, item-indent={0mm},column-sep={-30pt} ,, ](2)
\task $\Int\dfrac{\sqrt{x^2-1}}{x^4}dx$;
\task $\Int\dfrac{\sqrt[3]{x}+2}{\sqrt{x}+\sqrt[3]{x^2}}dx$. 
\end{tasks}
\end{minipage}
\vline\;
      \begin{minipage}[t]{100mm}
\textbf{7.}\hspace{1mm}Solve the definite integrals:
\begin{tasks}[ label-align=left, label-offset={1mm}, label-width={3mm},column-sep={-40pt}, item-indent={0mm},after-item-skip=-1mm,  ](2)
\task $\Int_{-1}^{e-2}\dfrac{\ln^3(x+2)}{\sqrt{x+2}}dx$;
\task $\Int_{\pi/4}^{\pi/2} \dfrac{\cos^3x}{\sin^5x} dx$;
\task $\Int_{0}^{1}\dfrac{x^2+2x}{x^2+4}dx$;
\task $\Int_{3}^{7}2^{\sqrt{x-3}}dx.$
\end{tasks}

\textbf{8.}\hspace{1mm}Find the area of the figure bounded \\
by the curves:
\begin{tasks}[ label-align=left, label-offset={2mm}, label-width={3mm}, item-indent={0mm}, after-item-skip=-1mm,     ](1)
\task $2x+y=7,\ xy=3$;
\task  $\left\lbrace\begin{array}{l} x=\cos^3t,\\ y=2\sin^3t; \end{array}\right. \ x=0 \ (x\geq 0)$;
\task  $\rho=3\sin4\phi.$
\end{tasks}

\textbf{9.}\hspace{1mm}Find the arc-length of the curve:
\begin{tasks}[ label-align=left, label-offset={2mm}, label-width={3mm}, item-indent={0mm}, after-item-skip=-1mm,     ](1)
\task $y=\ln(x+\sqrt{x^2-9}),\ 3\leq x\leq 6$;
\task $\left\lbrace\begin{array}{l} x=4(t-\sin t),\\ y=4(1-\cos t); \end{array}\right. \ 0\leq t\leq 2\pi$;
\task $\rho=2\cos\phi+2\sin\phi$.
\end{tasks}

\textbf{10.}\hspace{1mm}Find the area of the surface formed \\ by rotating the curves around the $l$-axis:
\begin{tasks}[ label-align=left, label-offset={2mm}, label-width={3mm}, item-indent={0mm}, after-item-skip=-1mm,    ](1)
\task $y=2x+6, \ 0\leq x\leq 2, \ l=OX$;
\task $\left\lbrace\begin{array}{l} x=3+2\cos t,\\ y=1+2\sin t; \end{array}\right. \ l=OY$;
\task $\rho=\sqrt{\cos2\phi}, \ l=o\rho$.
\end{tasks}

\textbf{11.}\hspace{1mm}Find the volume of the body formed \\
by rotating the curves around the $l$-axis:
\begin{tasks}[ label-align=left, label-offset={2mm}, label-width={3mm}, item-indent={0mm}, after-item-skip=-1mm,     ](1)
\task $y=x^2+4x+6,\ y=x+6,\ \ l=OX$;
\task  $y=\sqrt{4-x^2}, \ y=x, \ y=0, \ l=OY$;
\task  $\rho=8\sin^2\Big(\dfrac{\phi}{2}\Big), \ l=o\rho$.
\end{tasks}

\textbf{12.}\hspace{1mm}Solve the improper integrals:
\begin{tasks}[ label-align=left, label-offset={1mm}, label-width={3mm}, item-indent={0mm},   ,   ,column-sep={-40pt}](2)
\task  $\Int_{-\infty}^{0}e^x\cos xdx$;
\task $\Int_{1}^{4}\dfrac{dx}{x\sqrt{x-1}}$.
\end{tasks}
\end{minipage}

\section{Personal task 27}

   \begin{minipage}[t]{90mm}
\textbf{1.}\hspace{1mm}Integrate using the table and \\ substitution under differential:
\begin{tasks}[ label-align=left, label-offset={1mm}, label-width={3mm}, column-sep={-10pt}, item-indent={0mm} , ,after-item-skip=-1mm](2)
\task $\Int \dfrac{(2x-1)^3}{\sqrt[3]{x}} dx$;
\task $\Int\dfrac{\cos2x-1}{\cos2x+1}dx$; 
\task $\Int 5x\sqrt[4]{x^2-9}dx$;
\task $\Int \Big(1-3\textrm{tg}x\Big)^2dx$;
\end{tasks}
\textbf{2.}\hspace{1mm}Integrate the quadratic fractions:
\begin{tasks}[ label-align=left, label-offset={1mm}, label-width={3mm},column-sep={-20pt}, item-indent={0mm} ,, ](2)
\task $\Int\dfrac{4x+12}{x^2-9} dx$;
\task $\Int\dfrac{dx}{\sqrt{8-6x-9x^2}}$;
\task $\Int\dfrac{(-4x+16)dx}{x^2-4x+20}$;
\task $\Int\dfrac{12x dx}{\sqrt{x^2+4x-12}}$.
\end{tasks}
\textbf{3.}\hspace{1mm}Integrate by parts or using the \\ suitable substitutions:
\begin{tasks}[ label-align=left, label-offset={1mm}, label-width={3mm}, column-sep={-20pt}, item-indent={0mm},  ,after-item-skip=-1mm,](2)
\task $\Int\dfrac{8dx}{\sqrt{e^{4x}-16}}$;
\task $\Int\dfrac{\sqrt{9+x^2}}{x}dx$;
\task!  $\Int 6x^2(\sin x+\cos x) dx$;
\task  $\Int x3^xdx$;
\task  $\Int\textrm{arctg}\sqrt{x} \cdot\dfrac{dx}{\sqrt{x}}$;
\task  $\Int\ln(x^2-9)dx$.
\end{tasks}
\textbf{4.}\hspace{1mm}Integrate the polynomial fractions:
\begin{tasks}[ label-align=left, label-offset={1mm}, label-width={3mm},after-item-skip=-0.5mm, item-indent={0mm} , ](1)
\task $\Int\dfrac{4x^2-4x+6}{(x^2-4x+3)(x+2)} dx$;
\task $\Int\dfrac{-2x^3-6x^2-11x-8}{x^3+4x^2+4x} dx$;
\task $\Int\dfrac{x^2+26x-15}{(x^2-2x+5)(x+5)}dx$.
\end{tasks}
\textbf{5.}\hspace{1mm}Integrate trigonometric expressions:
\begin{tasks}[ label-align=left, label-offset={1mm}, label-width={3mm},after-item-skip=-1mm, item-indent={0mm} ,column-sep={-20pt} ](2)
\task $\Int\sin x\cos5x dx$;
\task $\Int\cos^3 x\sqrt[3]{\sin^2 x} dx$;
\task $\Int\dfrac{60 dx}{13+12\sin x}$.
\end{tasks}
\textbf{6.}\hspace{1mm}Integrate the fractions with radicals:
\begin{tasks}[ label-align=left, label-offset={1mm}, label-width={3mm}, item-indent={0mm},column-sep={-30pt} ,, ](2)
\task $\Int\dfrac{\sqrt{x^2-4}}{x^3}dx$;
\task $\Int\sqrt{\dfrac{x-3}{x+1}}\dfrac{dx}{(x-3)^2}$.
\end{tasks}
\end{minipage}
\vline\;
      \begin{minipage}[t]{100mm}
\textbf{7.}\hspace{1mm}Solve the definite integrals:
\begin{tasks}[ label-align=left, label-offset={1mm}, label-width={3mm},column-sep={-40pt}, item-indent={0mm},after-item-skip=-1mm,  ](2)
\task $\Int_{0}^{4}(x-4)e^{x}dx$;
\task $\Int_{0}^{2\pi}\cos^2(\frac{x}{4})\sin^4(\frac{x}{4})dx$;
\task $\Int_{0}^{\sqrt{2}}\dfrac{4x^3+8x}{\sqrt{16-x^4}}dx$;
\task $\Int_{1}^{e^2-1}\dfrac{\ln x}{\sqrt{x+1}}dx.$
\end{tasks}

\textbf{8.}\hspace{1mm}Find the area of the figure bounded \\
by the curves:
\begin{tasks}[ label-align=left, label-offset={2mm}, label-width={3mm}, item-indent={0mm}, after-item-skip=-1mm,     ](1)
\task $y=(x-1)^2,\ y=6x-6$;
\task  $\left\lbrace\begin{array}{l} x=6\cos t,\\ y=5\sin t; \end{array}\right. \ x=3 \ (x\geq3)$;
\task  $\rho=4\sqrt{\cos2\phi}.$
\end{tasks}

\textbf{9.}\hspace{1mm}Find the arc-length of the curve:
\begin{tasks}[ label-align=left, label-offset={2mm}, label-width={3mm}, item-indent={0mm}, after-item-skip=-1mm,     ](1)
\task  $y=\frac{1}{3}\ln(\cos3x),\ -\frac{\pi}{18}\leq x\leq \frac{\pi}{18}$;
\task  $\left\lbrace\begin{array}{l} x=3\cos^3 t,\\ y=3\sin^3 t; \end{array}\right. \ x=0 \ (x\geq 0)$;
\task  $\rho=8\cos^2\Big(\frac{\phi}{2}\Big)$.
\end{tasks}

\textbf{10.}\hspace{1mm}Find the area of the surface formed \\ by rotating the curves around the $l$-axis:
\begin{tasks}[ label-align=left, label-offset={2mm}, label-width={3mm}, item-indent={0mm}, after-item-skip=-1mm,    ](1)
\task $x^2+y^2=4, \ l=OY$;
\task $\left\lbrace\begin{array}{l} x=t-\sin t,\\ y=1-\cos t; \end{array}\right. \ 0\leq t\leq\pi,\ l=OX$;
\task $\rho=6\sin\phi,\ l=o\rho$.
\end{tasks}

\textbf{11.}\hspace{1mm}Find the volume of the body formed \\
by rotating the curves around the $l$-axis:
\begin{tasks}[ label-align=left, label-offset={2mm}, label-width={3mm}, item-indent={0mm}, after-item-skip=-1mm,     ](1)
\task $y=-x^2+4x+1,\ xy=4, \ l=OX$;
\task $y=4-x^2,\ y=12-3x^2, \ l=OY$;
\task $\rho=3e^{-\phi}, 0\leq\phi\leq\pi, \ l=o\rho$.
\end{tasks}

\textbf{12.}\hspace{1mm}Solve the improper integrals:
\begin{tasks}[ label-align=left, label-offset={1mm}, label-width={3mm}, item-indent={0mm},   ,   ,column-sep={-40pt}](2)
\task $\Int_{1}^{\infty}\dfrac{dx}{x+\sqrt[3]{x}}$;
\task  $\Int_{0}^{1}\dfrac{x^2}{\sqrt{1-x^2}}dx$.
\end{tasks}
\end{minipage}

\section{Personal task 28}

   \begin{minipage}[t]{90mm}
\textbf{1.}\hspace{1mm}Integrate using the table and \\ substitution under differential:
\begin{tasks}[ label-align=left, label-offset={1mm}, label-width={3mm}, column-sep={10pt}, item-indent={0mm} , ,after-item-skip=-1mm](2)
\task $\Int\dfrac{x^3-3x^24^x}{x^2} dx$;
\task $\Int (\textrm{tg}x-\textrm{ctg}x)^2 dx$;
\task $\Int \dfrac{3\textrm{arcsin}^2x-x}{\sqrt{1-x^2}}dx$;
\task $\Int\dfrac{\cos x+\sin x}{\sin x-\cos x} dx$.
\end{tasks}

\textbf{2.}\hspace{1mm}Integrate the quadratic fractions:
\begin{tasks}[ label-align=left, label-offset={1mm}, label-width={3mm},column-sep={-20pt}, item-indent={0mm} ,, ](2)
\task $\Int\dfrac{8x-3}{x^2-9} dx$;
\task $\Int\dfrac{dx}{\sqrt{8-4x-4x^2}}$;
\task $\Int\dfrac{(-4x+5)dx}{x^2-2x+2}$;
\task $\Int\dfrac{(-4x-4)dx}{\sqrt{x^2-4x-12}}$.
\end{tasks}
\textbf{3.}\hspace{1mm}Integrate by parts or using the \\ suitable substitutions:
\begin{tasks}[ label-align=left, label-offset={1mm}, label-width={3mm}, column-sep={0pt}, item-indent={0mm},  ,after-item-skip=-1mm,](2)
\task $\Int\dfrac{\ln^2x dx}{x\sqrt[3]{4+\ln x}}$;
\task $\Int\dfrac{\sqrt{2x-5}dx}{1+\sqrt{2x-5}}$;
\task $\Int 12x^2\cos^2 2xdx$;
\task $\Int e^{x}\sin 3x dx$;
\task $\Int \dfrac{x \textrm{arccos}x}{\sqrt{1-x^2}} dx$;
\task $\Int(2x+1)\dfrac{\sin xdx}{\cos^5x}$.
\end{tasks}
\textbf{4.}\hspace{1mm}Integrate the polynomial fractions:
\begin{tasks}[ label-align=left, label-offset={1mm}, label-width={3mm},after-item-skip=-0.5mm, item-indent={0mm} , ](1)
\task $\Int\dfrac{-5x^2+10x-8}{x^3-x^2-4x+4} dx$;
\task $\Int\dfrac{x^3-2x^2-13x+38}{(x^2+4x-5)(x-1)} dx$;
\task $\Int\dfrac{5x^2-35}{(x^2+4x+13)(x+1)}dx$.
\end{tasks}
\textbf{5.}\hspace{1mm}Integrate trigonometric expressions:
\begin{tasks}[ label-align=left, label-offset={1mm}, label-width={3mm},after-item-skip=-1mm, item-indent={0mm} ,, ](2)
\task $\Int\sin8x\sin4x dx$;
\task  $\Int\sqrt[3]{\dfrac{\sin^2x}{\cos^8x}}dx$;
\task  $\Int\dfrac{8\sin x dx}{7+9\sin^2x+13\cos^2x}$.
\end{tasks}
\textbf{6.}\hspace{1mm}Integrate the fractions with radicals:
\begin{tasks}[ label-align=left, label-offset={1mm}, label-width={3mm}, item-indent={0mm},column-sep={-30pt} , ](2)
\task $\Int\dfrac{x^2dx}{\sqrt{(9-x^2)^3}}$;
\task $\Int\dfrac{dx}{\sqrt{x+1}+4\sqrt[4]{x+1}}$.
\end{tasks}
\end{minipage}
\vline\;
      \begin{minipage}[t]{100mm}
\textbf{7.}\hspace{1mm}Solve the definite integrals:
\begin{tasks}[ label-align=left, label-offset={1mm}, label-width={3mm},column-sep={-1in}, item-indent={0mm},after-item-skip=-1mm ](2)
\task  $\Int_{1}^{e}\ln^2x dx$;
\task  $\Int_{0}^{\pi} \cos^4 2x\sin^4 2x dx$;
\task  $\Int_{0}^{4} \dfrac{\sqrt{x}dx}{x+4}$;
\task  $\Int_{\sqrt{2}}^{2} x^3\sqrt{x^2-1} dx.$
\end{tasks}

\textbf{8.}\hspace{1mm}Find the area of the figure bounded \\
by the curves:
\begin{tasks}[ label-align=left, label-offset={2mm}, label-width={3mm}, item-indent={0mm}, after-item-skip=-1mm,     ](1)
\task $y=(x+2)^3,\ y=9(x+2)$;
\task $\left\lbrace\begin{array}{l} x=4\cos^3t,\\ y=2\sin^3t; \end{array}\right. \ y=1 \ (y\geq 1)$;
\task $\rho=6(1-\sin\phi).$
\end{tasks}

\textbf{9.}\hspace{1mm}Find the arc-length of the curve:
\begin{tasks}[ label-align=left, label-offset={2mm}, label-width={3mm}, item-indent={0mm}, after-item-skip=-1mm,     ](1)
\task  $y=2-e^x,\ 0\leq x\leq \frac{1}{2}\ln15$;
\task   $\left\lbrace\begin{array}{l} x=3(t-\sin t),\\ y=3(1-\cos t); \end{array}\right. \ 0\leq t \leq 2\pi$;
\task   $\rho=\sqrt{2}e^{-2\phi}, \ 0\leq\phi\leq4\pi$.
\end{tasks}

\textbf{10.}\hspace{1mm}Find the area of the surface formed \\ by rotating the curves around the $l$-axis:
\begin{tasks}[ label-align=left, label-offset={2mm}, label-width={3mm}, item-indent={0mm}, after-item-skip=-1mm,    ](1)
\task $9y^2=4x^3, \ 0\leq x\leq 1, \ l=OY$;
\task $\left\lbrace\begin{array}{l} x=2+\cos t,\\ y=1+\sin t; \end{array}\right. \ l=OX$;
\task $\rho=3\sqrt{\cos2\phi}, \ l=o\rho$.
\end{tasks}
\textbf{11.}\hspace{1mm}Find the volume of the body formed \\
by rotating the curves around the $l$-axis:
\begin{tasks}[ label-align=left, label-offset={2mm}, label-width={3mm}, item-indent={0mm}, after-item-skip=-1mm,     ](1)
\task $y=-x^2+5x+1,\ y=6-x, \ l=OX$;
\task $y=1-x^2, \ y=3-3x^2, \ l=OY$;
\task $\rho=2\sin 2\phi, \ l=o\rho$.
\end{tasks}

\textbf{12.}\hspace{1mm}Solve the improper integrals:
\begin{tasks}[ label-align=left, label-offset={1mm}, label-width={3mm}, item-indent={0mm},   ,   ,column-sep={-40pt}](2)
\task $\Int_{1}^{\infty}\dfrac{\textrm{arctg}^3x}{1+x^2}dx$;
\task $\Int_{0}^{\pi^2/4}\dfrac{\cos^2\sqrt{x}}{\sqrt{x}}dx$.
\end{tasks}
\end{minipage}

\section{Personal task 29}

   \begin{minipage}[t]{90mm}
\textbf{1.}\hspace{1mm}Integrate using the table and \\ substitution under differential:
\begin{tasks}[ label-align=left, label-offset={1mm}, label-width={3mm},column-sep={-20pt}, item-indent={0mm} ,, ](2)
\task $\Int \dfrac{(2+3\sqrt{x})^3}{x^2} dx$;
\task $\Int \textrm{tg}^4x dx$;
\task $\Int \dfrac{x^3}{\sqrt[3]{(1-x^4)^2}}dx$;
\task $\Int (\textrm{th}x+\textrm{cth}x)^2 dx.$
\end{tasks}
\textbf{2.}\hspace{1mm}Integrate the quadratic fractions:
\begin{tasks}[ label-align=left, label-offset={1mm}, label-width={3mm},column-sep={-20pt}, item-indent={0mm} ,, ](2)
\task $\Int\dfrac{2x-6}{x^2-25} dx$;
\task $\Int\dfrac{dx}{\sqrt{8+4x-4x^2}}$;
\task $\Int\dfrac{(4x-9)dx}{x^2-6x+18}$;
\task $\Int\dfrac{(6x-6)dx}{\sqrt{5-4x+x^2}}$.
\end{tasks}

\textbf{3.}\hspace{1mm}Integrate by parts or using the \\ suitable substitutions:
\begin{tasks}[ label-align=left, label-offset={1mm}, label-width={3mm}, column-sep={-50pt}, item-indent={0mm},  ,after-item-skip=-1mm,](2)
\task $\Int\dfrac{5x+6}{(x+1)^7}dx$;
\task $\Int\dfrac{dx}{(x-2)^2\sqrt{4x-x^2}}$;
\task! $\Int (4x^2-8x)\sin x\cos x dx$;
\task $\Int e^{4x}\sin 3x dx$;
\task $\Int\dfrac{\ln x}{x^3}dx$;
\task! $\Int (3x^2+1)\textrm{arcctg}x dx$.
\end{tasks}

\textbf{4.}\hspace{1mm}Integrate the polynomial fractions:
\begin{tasks}[ label-align=left, label-offset={1mm}, label-width={3mm},after-item-skip=-0.5mm, item-indent={0mm} , ](1)
\task $\Int\dfrac{-x^2+15x-20}{(x^2-5x+6)(x+1)} dx$;
\task $\Int\dfrac{2x^3-3x^2-17x-26}{x^3+x^2-8x-12} dx$;
\task $\Int\dfrac{-14x+80}{(x^2-2x+10)(x+4)}dx$.
\end{tasks}

\textbf{5.}\hspace{1mm}Integrate trigonometric expressions:
\begin{tasks}[ label-align=left, label-offset={1mm}, label-width={3mm},after-item-skip=-1mm, item-indent={0mm} ,, ](2)
\task $\Int\cos 2x\sin 8xdx$;
\task $\Int\dfrac{\cos^2 x}{\sin^6 x} dx$;
\task $\Int\dfrac{4dx}{7-2\sin x+6\cos x}$.
\end{tasks}
\textbf{6.}\hspace{1mm}Integrate the fractions with radicals:
\begin{tasks}[ label-align=left, label-offset={1mm}, label-width={3mm}, item-indent={0mm},column-sep={-40pt} , ](2)
\task $\Int\dfrac{\sqrt{4+x^2}}{x^4}dx$;
\task $\Int\sqrt{\dfrac{x+3}{x-2}}\dfrac{dx}{(x+3)^2}$.
\end{tasks}
\end{minipage}
\vline \;
      \begin{minipage}[t]{100mm}
\textbf{7.}\hspace{1mm}Solve the definite integrals:
\begin{tasks}[ label-align=left, label-offset={1mm}, label-width={3mm},column-sep={-1in}, item-indent={0mm},after-item-skip=-1mm ](2)
\task  $\Int_{1/\ln4}^{1/\ln2}\sqrt[x]{e}\cdot\dfrac{dx}{x^2}$;
\task  $\Int_{0}^{\pi}\cos^4x\sin^2xdx$;
\task  $\Int_{0}^{\sqrt{3}}\dfrac{2x+4}{x^2+1}dx$;
\task  $\Int_{0}^{1}x^3\sqrt{1-x^2}dx.$
\end{tasks}

\textbf{8.}\hspace{1mm}Find the area of the figure \\ bounded
by the curves:
\begin{tasks}[ label-align=left, label-offset={2mm}, label-width={3mm}, item-indent={0mm}, after-item-skip=-1mm,     ](1)
\task $y=2x^2-11x+13,\ y=-x^2+4x+1$;
\task $\left\lbrace\begin{array}{l} x=4\cos t,\\ y=2+4\sin t; \end{array}\right. \ x=2 \ (x\geq2)$;
\task $\rho=3\sqrt{\sin 2\phi}.$
\end{tasks}

\textbf{9.}\hspace{1mm}Find the arc-length of the curve:
\begin{tasks}[ label-align=left, label-offset={2mm}, label-width={3mm}, item-indent={0mm}, after-item-skip=-1mm,     ](1)
\task $y=4\ln(16-x^2),\ -2\leq x\leq2$;
\task $\left\lbrace\begin{array}{l} x=2e^{-t}\cos t,\\ y=2e^{-t}\sin t; \end{array}\right.  \ 0\leq t\leq 2\pi$;
\task $\rho=4\cos^2\Big(\dfrac{\phi}{2}\Big)$.
\end{tasks}

\textbf{10.}\hspace{1mm}Find the area of the surface formed \\ by rotating the curves around the $l$-axis:
\begin{tasks}[ label-align=left, label-offset={2mm}, label-width={3mm}, item-indent={0mm}, after-item-skip=-1mm,    ](1)
\task $y=\frac{1}{5}\textrm{ch}5x, \ 0\leq x\leq 1, \ l=OX$;
\task $\left\lbrace\begin{array}{l} x=3\cos^3 t,\\ y=3\sin^3 t; \end{array}\right. \ l=OY$;
\task $\rho=2\cos\phi, \ l=o\rho$.
\end{tasks}
\textbf{11.}\hspace{1mm}Find the volume of the body formed \\
by rotating the curves around the $l$-axis:
\begin{tasks}[ label-align=left, label-offset={2mm}, label-width={3mm}, item-indent={0mm}, after-item-skip=-1mm,     ](1)
\task  $y=x^2-4x+4,\ y=6x-12, \ l=OX$;
\task  $y=x^2-6,\ y=5x,\ x=0 \ (x\geq0),\ l=OY$;
\task  $\rho=6\phi,\ 0\leq\phi\leq\pi,\ l=o\rho$.
\end{tasks}

\textbf{12.}\hspace{1mm}Solve the improper integrals:
\begin{tasks}[ label-align=left, label-offset={1mm}, label-width={3mm}, item-indent={0mm},   ,   ,column-sep={-40pt}](2)
\task  $\Int_{0}^{\infty}x^3e^{-x^2}dx$;
\task $\Int_{0}^{1}\dfrac{x^3-4x}{\sqrt{1-x^4}}dx$.
\end{tasks}

\end{minipage}

\section{Personal task 30}

   \begin{minipage}[t]{90mm}
\textbf{1.}\hspace{1mm}Integrate using the table and \\ substitution under differential:
\begin{tasks}[ label-align=left, label-offset={1mm}, label-width={3mm}, column-sep={-40pt}, item-indent={0mm} , ,after-item-skip=-1mm](2)
\task! $\Int \sqrt{x}(1+2x\sqrt{x})^2 dx$;
\task $\Int \dfrac{\cos2x+1}{\sin^2x}dx$;
\task $\Int \dfrac{3\textrm{arccos}^2x-2x}{\sqrt{1-x^2}}dx$;
\task! $\Int\Big(\textrm{sh}(x+1)-2\Big)^2dx$.
\end{tasks}
\textbf{2.}\hspace{1mm}Integrate the quadratic fractions:
\begin{tasks}[ label-align=left, label-offset={1mm}, label-width={3mm},column-sep={-20pt}, item-indent={0mm} ,, ](2)
\task $\Int\dfrac{-4x+6}{x^2-9} dx$;
\task $\Int\dfrac{dx}{\sqrt{9x^2-6x-8}}$;
\task $\Int\dfrac{(-4x+4)dx}{x^2+2x+17}$;
\task $\Int\dfrac{-6x+6}{\sqrt{6x-x^2}} dx$.
\end{tasks}
\textbf{3.}\hspace{1mm}Integrate by parts or using the \\ suitable substitutions:
\begin{tasks}[ label-align=left, label-offset={1mm}, label-width={3mm},after-item-skip=-1mm, item-indent={0mm} ,column-sep={-15pt}, ](2)
\task  $\Int\dfrac{\sqrt{x-3}\ dx}{x+4\sqrt{x-3}}$;
\task $\Int\dfrac{dx}{x^2\sqrt{4-x^2}}$;
\task $\Int 9x^2 \cos 3x dx$;
\task $\Int x^3\ln(x-1) dx$;
\task $\Int e^{\sqrt[3]{x+4}}dx$;
\task $\Int x\cdot\textrm{arcsin}(x^2)dx$.
\end{tasks}

\textbf{4.}\hspace{1mm}Integrate the polynomial fractions:
\begin{tasks}[ label-align=left, label-offset={1mm}, label-width={3mm},after-item-skip=-0.5mm, item-indent={0mm} , ](1)
\task $\Int\dfrac{3x^2+5x-4}{(x^2+3x+2)(x+3)} dx$;
\task $\Int\dfrac{-x^3+4x^2+6x-28}{(x^2-6x+8)(x-2)} dx$;
\task $\Int\dfrac{5x^2-16x-45}{(x^2+2x+5)(x-4)}dx$.
\end{tasks}

\textbf{5.}\hspace{1mm}Integrate trigonometric expressions:
\begin{tasks}[ label-align=left, label-offset={1mm}, label-width={3mm},after-item-skip=-1mm, item-indent={0mm} ,, ](2)
\task $\Int \sin 5x\cos 3xdx$;
\task $\Int \sin^2x\cos^4x dx$;
\task $\Int\dfrac{12 dx}{4\sin^2x-12\sin x\cos x+5\cos^2x}$.
\end{tasks}
\textbf{6.}\hspace{1mm}Integrate the fractions with radicals:
\begin{tasks}[ label-align=left, label-offset={1mm}, label-width={3mm}, item-indent={0mm},column-sep={-30pt} , ](2)
\task $\Int \dfrac{x^3}{\sqrt{x^2+4}}dx$;
\task $\Int\dfrac{(\sqrt[3]{x}+1)dx}{\sqrt[3]{x^2}+2\sqrt[3]{x}-3}$.
\end{tasks}
\end{minipage}
\vline\;
      \begin{minipage}[t]{100mm}
\textbf{7.}\hspace{1mm}Solve the definite integrals:
\begin{tasks}[ label-align=left, label-offset={1mm}, label-width={3mm},column-sep={-50pt}, item-indent={0mm},after-item-skip=-1mm ](2)
\task $\Int_{0}^{1}\textrm{arctg} x dx$;
\task $\Int_{0}^{\pi/4}\dfrac{\sin^3x}{\cos^5x}dx$;
\task $\Int_{1}^{2}\dfrac{4x+2}{\sqrt{4x-x^2}}dx$;
\task $\Int_{-2}^{2}x^2\sqrt{4-x^2}dx.$
\end{tasks}

\textbf{8.}\hspace{1mm}Find the area of the figure bounded \\
by the curves:
\begin{tasks}[ label-align=left, label-offset={2mm}, label-width={3mm}, item-indent={0mm}, after-item-skip=-1mm,     ](1)
\task $y=(x+2)^3,\ y=4x+8$;
\task $\left\lbrace\begin{array}{l} x=4\cos t,\\ y=6\sin t; \end{array}\right. \ x=2 \ (x\geq2)$;
\task $\rho=3\cos2\phi$.
\end{tasks}

\textbf{9.}\hspace{1mm}Find the arc-length of the curve:
\begin{tasks}[ label-align=left, label-offset={2mm}, label-width={3mm}, item-indent={0mm}, after-item-skip=-1mm,     ](1)
\task $y=1-\frac{1}{3}\textrm{ch}3x,\ 0\leq x\leq \frac{1}{3}\ln2$;
\task $\left\lbrace\begin{array}{l} x=4\cos^3 t,\\ y=4\sin^3 t; \end{array}\right. \ 0\leq t\leq 2\pi$;
\task $\rho=6\cos\phi+6\sin\phi$.
\end{tasks}

\textbf{10.}\hspace{1mm}Find the area of the surface formed \\ by rotating the curves around the $l$-axis:
\begin{tasks}[ label-align=left, label-offset={2mm}, label-width={3mm}, item-indent={0mm}, after-item-skip=-1mm,    ](1)
\task $y=\sqrt{2x+8}, \ -4\leq x\leq 4, \ l=OX$;
\task $\left\lbrace\begin{array}{l} x=3(t-\sin t),\\ y=3(1-\cos t); \end{array}\right. \ \begin{array}{l} 0\leq t\leq2\pi,\\ l=OX;\end{array}$
\task $\rho=2(1+\cos\phi), \ l=o\rho$.
\end{tasks}
\textbf{11.}\hspace{1mm}Find the volume of the body formed \\
by rotating the curves around the $l$-axis:
\begin{tasks}[ label-align=left, label-offset={2mm}, label-width={3mm}, item-indent={0mm}, after-item-skip=-1mm,     ](1)
\task  $y=x^2-4x,\ y=2x^2-8x, \ l=OX$;
\task $x=9-y^2,\ x+3y=9, \ l=OY$;
\task $\rho=e^{-\phi}, \ 0\leq\phi\leq\pi,\ l=o\rho$.
\end{tasks}

\textbf{12.}\hspace{1mm}Solve the improper integrals:
\begin{tasks}[ label-align=left, label-offset={1mm}, label-width={3mm}, item-indent={0mm},   ,   ,column-sep={-40pt}](2)
\task $\Int_{0}^{\infty}\dfrac{2x+5}{x^2+1}dx$;
\task  $\Int_{1}^{e}\dfrac{dx}{x\sqrt[4]{1-\ln x}}$.
\end{tasks}
\end{minipage}

\section{Sample task}

\begin{minipage}[t]{90mm}
\textbf{1.}\hspace{1mm}Integrate using the table and \\ substitution under differential:
\begin{tasks}[ label-align=left, label-offset={1mm}, label-width={3mm}, item-indent={0mm}](2)
\task $\int x^3(1-x^2)^2dx $;
\task $\int (3+\tg^2x)dx$;
\task $\int x^4 e^{2x^5-1}dx$;
\task $\int \sin^2(1+3x)dx$.
\end{tasks}
\textbf{2.}\hspace{1mm}Integrate the quadratic fractions:
\begin{tasks}[ label-align=left, label-offset={1mm}, label-width={3mm},column-sep={-18pt}, item-indent={0mm}](2)
\task $\int \dfrac{-4x+5}{x^2-4} dx$;
\task $\int \dfrac{dx}{\sqrt{6x-9x^2}}$;
\task $\int \dfrac{(-2x-7)dx}{x^2+6x+10}$;
\task $\int \dfrac{(8x-5)dx}{\sqrt{x^2+4x-5}}$.
\end{tasks}
\textbf{3.}\hspace{1mm}Integrate by parts or using the \\ suitable substitutions:
\begin{tasks}[ label-align=left, label-offset={1mm}, label-width={3mm}, column-sep={0pt}, item-indent={0mm}](2)
\task $\int \dfrac{3x+2}{\sqrt{x+4}}dx$;
\task $\int \dfrac{dx}{x\sqrt{3x^2-4x+1}}$;
\task $\int \dfrac{(4x^2-3)}{e^{2x}}dx$;
\task $\int x^2 \ln x dx$;
\task $\int \sin \sqrt{x+1} dx$;
\task $\int x \dfrac{\sin x}{\cos^3 x} dx$.
\end{tasks}
\textbf{4.}\hspace{1mm}Integrate the polynomial fractions:
\begin{tasks}[ label-align=left, label-offset={1mm}, label-width={3mm}, item-indent={0mm}](1)
\task $\Int\dfrac{-3x^2+2x+13}{x^3+2x^2-x-2} dx$;
\task $\Int\dfrac{3x^3-32x+56}{x^3-2x^2-4x+8}dx$;
\task $\Int\dfrac{x^2-2x-9}{(x^2+4x+5)(x-1)}dx$.
\end{tasks}
\textbf{5.}\hspace{1mm}Integrate trigonometric functions:
\begin{tasks}[ label-align=left, label-offset={1mm}, label-width={3mm}, item-indent={0mm},column-sep={-40pt}](2)
\task! $\int\sin 10x\sin 3xdx$;
\task $\int \sqrt[3]{\dfrac{\sin x}{\cos^{13}x}} dx$;
\task $\int\dfrac{dx}{4\cos x+3\sin x+6}$.
\end{tasks}
\textbf{6.}\hspace{1mm}Integrate the fractions with radicals:
\begin{tasks}[ label-align=left, label-offset={1mm}, label-width={3mm}, item-indent={0mm},column-sep={-10pt}](2)
\task $\Int\dfrac{dx}{x^2\sqrt{4-x^2}}$;
\task $\Int\sqrt[3]{\dfrac{x+1}{x-1}}\dfrac{dx}{(x-1)^{3}}$.
\end{tasks}
\end{minipage}
\vline\;
\begin{minipage}[t]{100mm}
\textbf{7.}\hspace{1mm}Solve the definite integrals:
\begin{tasks}[ label-align=left, label-offset={1mm}, label-width={3mm}, item-indent={0mm}, column-sep={-20pt}](2)
\task $\Int_{0}^{1}\ln(1+x^2)dx$;
\task $\Int_{0}^{\pi/2}\sin^3x\sqrt[4]{\cos x}dx$;
\task $\Int_{1}^{\sqrt[8]{2}} \dfrac{4x^3dx}{\sqrt{4-x^8}}$;
\task $\Int_{1}^{16}\dfrac{1+\sqrt{x}}{\sqrt[4]{x}+\sqrt{x}}dx.$
\end{tasks}

\textbf{8.}\hspace{1mm}Find the area of the figure \\ bounded by the curves:
\begin{tasks}[ label-align=left, label-offset={2mm}, label-width={3mm}, item-indent={0mm}](1)
\task $y=2x^2-10x+6,\ y=x^2-3x$;
\task $\left\lbrace\begin{array}{l} x=4\cos^3t,\\ y=4\sin^3t; \end{array}\right. \ y=\frac{1}{2} \ (y\geq \frac{1}{2})$;
\task $\rho=4\sin 3\varphi.$
\end{tasks}

\textbf{9.}\hspace{1mm}Find the arc length of the curves:
\begin{tasks}[ label-align=left, label-offset={2mm}, label-width={3mm}, item-indent={0mm}](1)
\task $y=\dfrac{1}{3}\ln (\cos 3x),\ 0\leq x\leq \dfrac{\pi}{18}$;
\task $\left\lbrace\begin{array}{l} x=4(t-\sin t),\\ y=4(1-\cos t); \end{array}\right. \ 0\leq t\leq 2\pi$;
\task $\rho=\frac{10}{\sqrt{101}}e^{\frac{\varphi}{10}}, \ 0\leq\varphi\leq2\pi$.
\end{tasks}
\textbf{10.}\hspace{1mm}Find the area of the surface formed\\  by rotating the curves around the $l$-axis:
\begin{tasks}[ label-align=left, label-offset={2mm}, label-width={3mm}, item-indent={0mm}](1)
\task $y=\dfrac{1}{2}\ch 2x, \ -1\leq x\leq 1, \ \ l=OX$;
\task $\left\lbrace\begin{array}{l} x=3+\cos t, \\ y=2+\sin t; \end{array}\right.  \ l=OY$;
\task $\rho=\sqrt{\cos 2\phi}, 0\leq \phi \leq \pi \ l=o\rho$.
\end{tasks}
\textbf{11.}\hspace{1mm}Find the volume of the body formed \\ by rotating the curves around the $l$-axis:
\begin{tasks}[ label-align=left, label-offset={2mm}, label-width={3mm}, item-indent={0mm}](1)
\task $y=x^2+2x+5,\ y=5-x, \ l=OX$;
\task $x=5+4y-y^2,\ x=5, \ l=OY$;
\task $\rho=6(1+\cos\varphi), \ l=o\rho$.
\end{tasks}
\textbf{12.}\hspace{1mm}Solve the improper integrals:
\begin{tasks}[ label-align=left, label-offset={1mm}, label-width={3mm}, item-indent={0mm},column-sep={-40pt}](2)
\task $\Int_{-2}^{\infty}\dfrac{2x-1}{x^2+4}dx$;
\task $\Int_{0}^{\pi/2}e^{-\tg x}\dfrac{dx}{\cos^2 x}$.
\end{tasks}
\end{minipage}

\newpage

\begingroup
\addtolength{\jot}{0.2cm}

\begin{enumerate}[label=\textbf{1.(\alph*)}]
\item $\begin{aligned}[t]
       \dot{I}= & \int x^3(1-x^2)^2dx=
    \left| \begin{array}{c}
     \textit{Applying the Newton binom} \\
     \textit{to the function under the integral}
    \end{array}\right| = \\
    = & \int x^3(1-2x^2+x^4)dx= \int (x^3-2x^5+x^7)dx=  \frac{1}{4}x^4-\frac{2}{6}x^6+\frac{1}{8}x^8+C=  \\
    = & \frac{1}{4}x^4-\frac{1}{3}x^6+\frac{1}{8}x^8+C.
    \end{aligned} $
\vspace{0.3cm}

   $ \textbf{Result: } \dot{I} =  \dfrac{1}{4}x^4-\dfrac{1}{3}x^6+\dfrac{1}{8}x^8+C,\; C\in \mathbb{R}.$

\begin{tabular}{|p{6.0cm}|p{7.5cm}|p{2.0cm}|}
\hline
\vspace{0.05mm}$^*$ Solution of Problem 1(a) guided by Ricard Riba is available on-line: &
\vspace{5.5mm} \url{https://youtu.be/x48CikKlF9c} & \vspace{-3mm} \includegraphics[height=20mm]{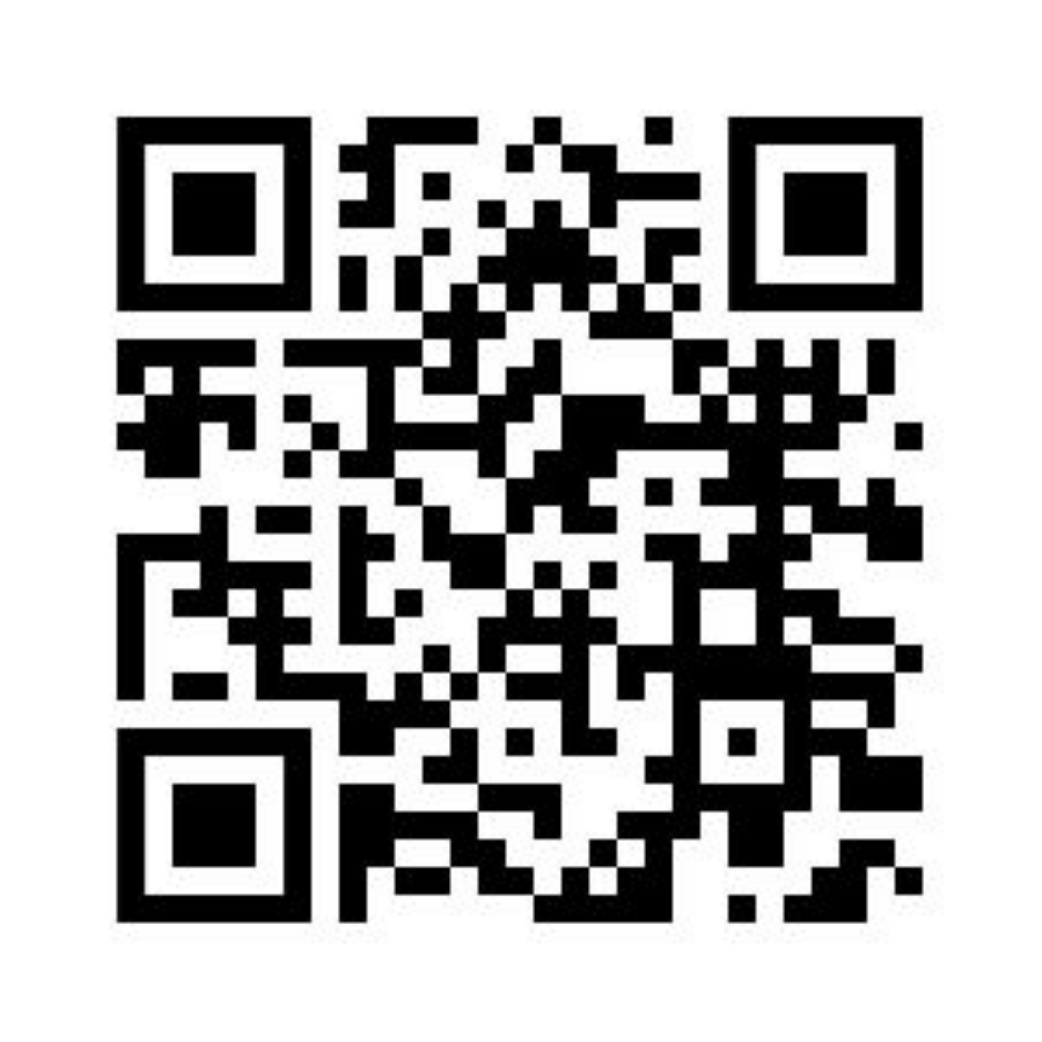}\\
\hline
\end{tabular}

   \item $\begin{aligned}[t]
        \dot{I}= & \int (3+\tg^2x)dx=\int \big(2+(1+\tg^2x)\big)dx=
    \left| \begin{array}{c}
     \textit{Applying the} \\
     \textit{trigonometric identity} \\
     1+\tg^2x=\dfrac{1}{\cos^2x}
    \end{array}\right| = \\
    = & \int \left(2+\frac{1}{\cos^2x}\right)dx= 2x+\int \frac{dx}{\cos^2x}= 2x+\tg x+C.
    \end{aligned}$
    \vspace{0.3cm}

   $ \textbf{Result: } \dot{I} =  2x+\tg x+C,\; C\in \mathbb{R}.$

   \begin{tabular}{|p{6.0cm}|p{7.5cm}|p{2.0cm}|}
\hline
\vspace{0.05mm}$^*$ Solution of Problem 1(b) guided by Irina Blazhievska is available on-line: &
\vspace{5.5mm} \url{https://youtu.be/1MqmZbQ-3qM} & \vspace{-3mm} \includegraphics[height=20mm]{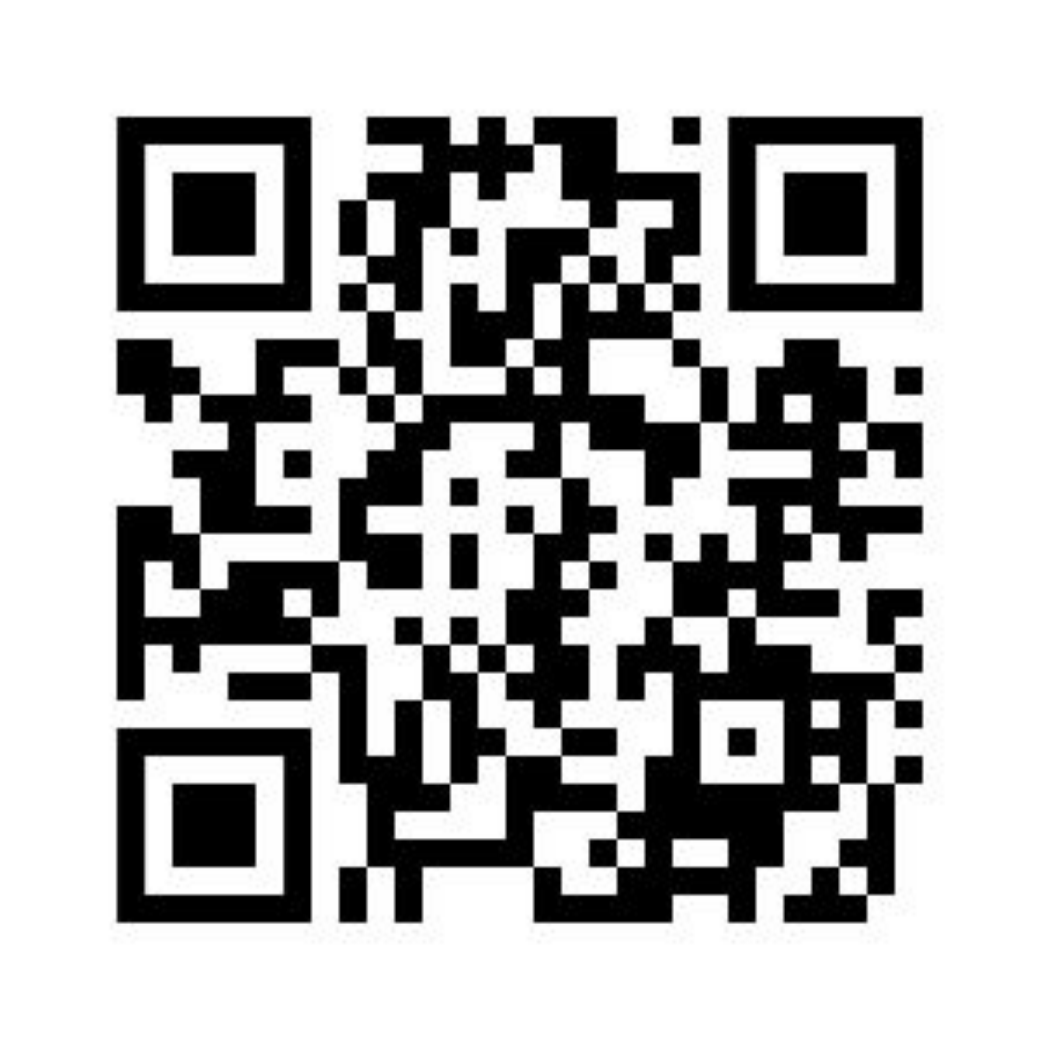}\\
\hline
\end{tabular}

    \item $\begin{aligned}[t]
        \dot{I}= & \int x^4 e^{2x^5-1}dx=
    \left| \begin{array}{c}
     \textit{Checking the exact differential of the power} \\
     d(2x^5-1)=10x^4 dx; \\
     x^4dx=\dfrac{1}{10}d(2x^5-1)
    \end{array}\right| = \\
     = & \frac{1}{10} \int  e^{2x^5-1}d(2x^5-1)=\frac{1}{10} e^{2x^5-1}+C.
     \end{aligned}$
     \vspace{0.3cm}

   $ \textbf{Result: } \dot{I} =  \dfrac{1}{10} e^{2x^5-1}+C,\; C\in \mathbb{R}.$

\begin{tabular}{|p{6.0cm}|p{7.5cm}|p{2.0cm}|}
\hline
\vspace{0.05mm}$^*$ Solution of Problem 1(c) guided by Irina Blazhievska is available on-line: &
\vspace{5.5mm} \url{https://youtu.be/kSn2UvdXWVs} & \vspace{-3mm} \includegraphics[height=20mm]{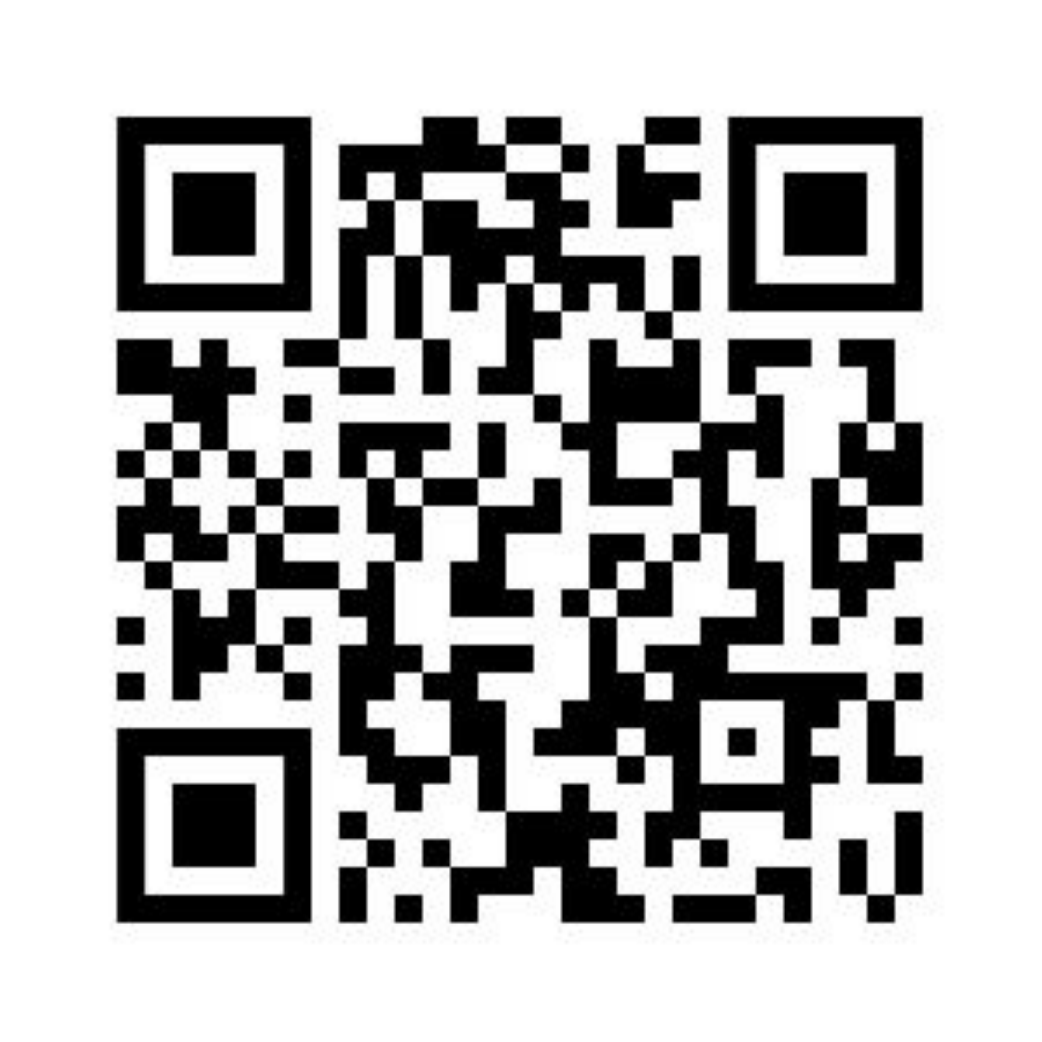}\\
\hline
\end{tabular}

\pagebreak

    \item $\begin{aligned}[t]
        \dot{I}= & \int \sin^2(1+3x)dx=
    \left| \begin{array}{c}
     \textit{Applying the reduction power formula} \\
     \sin^2(1+3x)=\dfrac{1}{2}\big(1-\cos(2+6x)\big) \\
    \end{array}\right| = \\
    = & \dfrac{1}{2}\int \big(1-\cos(2+6x)\big)dx= \dfrac{1}{2}\left(\int dx-\int \cos(2+6x)dx\right)= \\
    = & \left| \begin{array}{c}
     \textit{Substitution under the differential} \\
     dx=\dfrac{1}{6} d(2+6x)
    \end{array}\right| = \\
    = & \dfrac{1}{2}\left(x-\frac{1}{6}\int \cos(2+6x)d(2+6x)\right)+C= \frac{1}{2}x-\frac{1}{12}\sin(2+6x)+C.
     \end{aligned}$
     \vspace{0.3cm}

   $ \textbf{Result: } \dot{I} = \dfrac{1}{2}x-\dfrac{1}{12}\sin(2+6x)+C,\; C\in \mathbb{R}.$

	\end{enumerate}

   \begin{enumerate}[label=\textbf{2.(\alph*)}]
   \item $\begin{aligned}[t]
        \dot{I}= & \int \frac{-4x+5}{x^2-4}dx=
    \left| \begin{array}{c}
     \textit{Checking the differential of the denominator} \\
     d(x^2-4)=2x dx \\[0.2cm]
     \textit{Decomposing the numerator} \\
     -4x+5=-2(2x)+5
    \end{array}\right| = \\
     = & -2\int \frac{2x \;dx}{x^2-4} +5\int \frac{dx}{x^2-4}= -2\int \frac{ d(x^2-4)}{x^2-4} +5\int \frac{dx}{x^2-2^2}= \\
     = & -2\ln|x^2-4|+\frac{5}{4}\ln\left|\dfrac{x-2}{x+2}\right|+C= -\dfrac{3}{4}\ln|x-2|-\dfrac{13}{4}\ln|x+2|+C.
     \end{aligned}$
     \vspace{0.3cm}

   $ \textbf{Result: } \dot{I} =-\dfrac{3}{4}\ln|x-2|-\dfrac{13}{4}\ln|x+2|+C,\; C\in \mathbb{R}.$

   \item $\begin{aligned}[t]
        \dot{I}= & \int \frac{dx}{\sqrt{6x-9x^2}}dx=
    \left| \begin{array}{c}
     \textit{Completing the square inside the radical} \\
     6x-9x^2 =1-(3x-1)^2
    \end{array}\right| = \\
     = & \int \frac{dx}{\sqrt{1-(3x-1)^2}}=  \frac{1}{3}\int \frac{d(3x-1)}{\sqrt{1-(3x-1)^2}}= \frac{1}{3} \arcsin (3x-1)+C.
     \end{aligned}$
     \vspace{0.3cm}

   $ \textbf{Result: } \dot{I} = \dfrac{1}{3} \arcsin (3x-1)+C,\; C\in \mathbb{R}.$

   \begin{tabular}{|p{6.0cm}|p{7.5cm}|p{2.0cm}|}
\hline
\vspace{0.05mm}$^*$ Solution of Problem 2(b) guided by Irina Blazhievska is available on-line: &
\vspace{5.5mm}
   \url{https://youtu.be/KQuEzkh6AqQ} & \vspace{-3mm} \includegraphics[height=20mm]{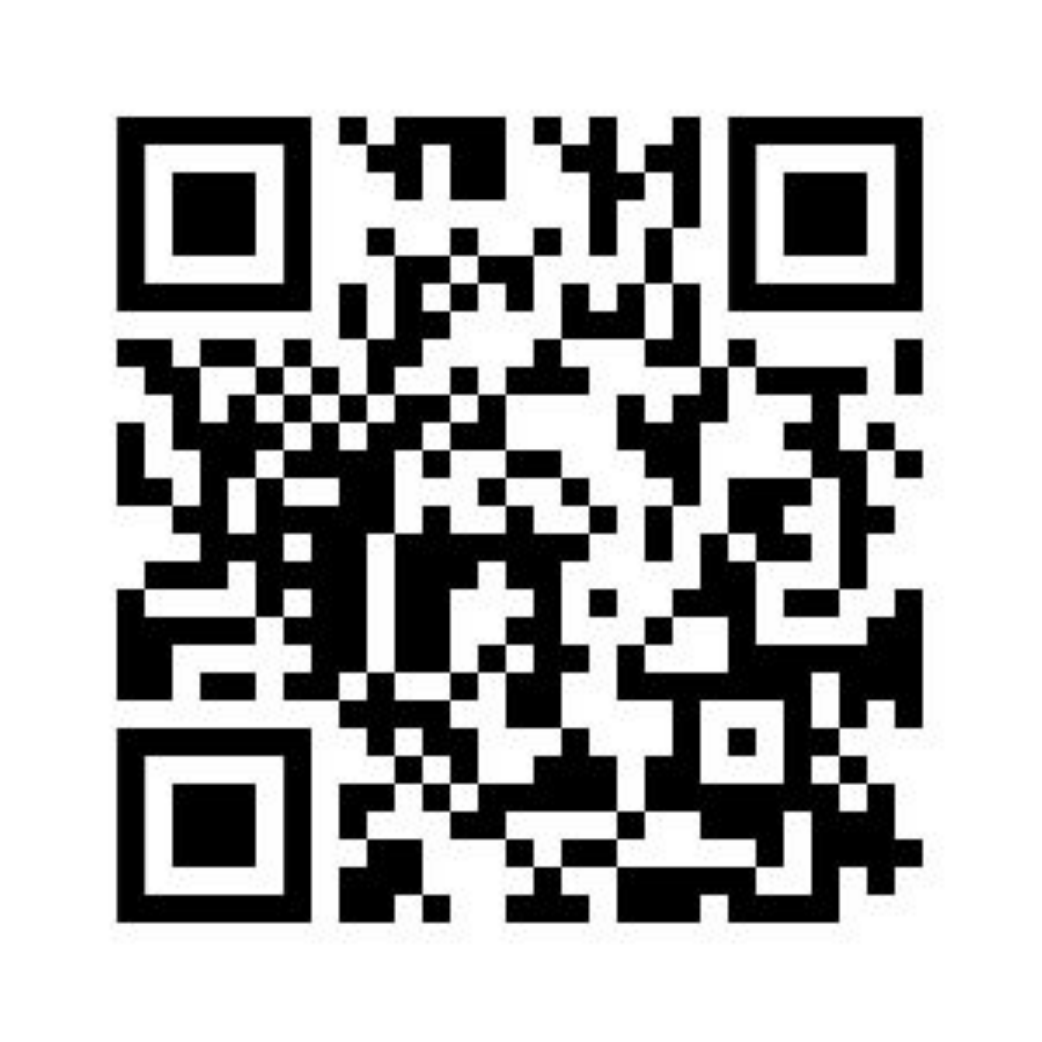}\\
\hline
\end{tabular}

    \item $\begin{aligned}[t]
        \dot{I}= & \int \frac{-2x-7}{x^2+6x+10}dx=
    \left| \begin{array}{c}
     \textit{Checking the differential of the denominator:} \\
     d(x^2+6x+10)=(2x+6) dx; \\[0.2cm]
    \textit{Decomposing the numerator} \\
     -2x-7=-(2x+6)-1
    \end{array}\right| =\\
     = & -\int \frac{2x+6}{x^2+6x+10}dx-\int \frac{dx}{x^2+6x+10} = \left| \begin{array}{c}
     \textit{Completing the square:}\\
     x^2+6x+10=(x+3)^2+1
    \end{array}\right| =\\
     = & -\int \frac{d(x^2+6x+10)}{x^2+6x+10}-\int \frac{d(x+3)}{(x+3)^2+1} = \\
     =&-\ln(x^2+6x+10)- \arctg(x+3)+C.
     \end{aligned}$
     \vspace{0.3cm}

   $ \textbf{Result: } \dot{I} = -\ln(x^2+6x+10)- \arctg(x+3)+C,\; C\in \mathbb{R}.$

      \begin{tabular}{|p{6.0cm}|p{7.5cm}|p{2.0cm}|}
\hline
\vspace{0.05mm}$^*$ Solution of Problem 2(c) guided by Ricard Riba is available on-line: &
\vspace{5.5mm}
   \url{https://youtu.be/RXJcMkNz8Zg} & \vspace{-3mm} \includegraphics[height=20mm]{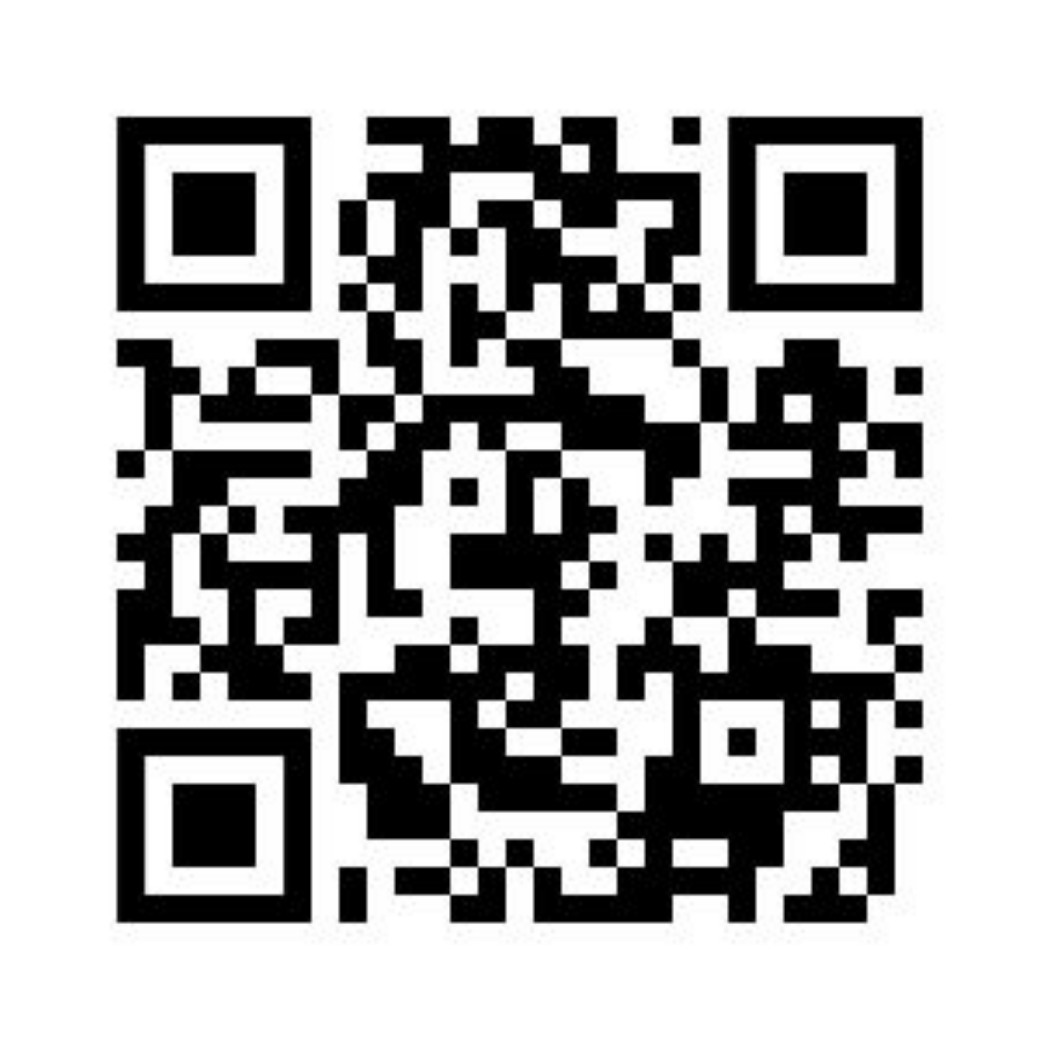}\\
\hline
\end{tabular}

   \item $\begin{aligned}[t]
        \dot{I}= & \int \frac{8x-5}{\sqrt{x^2+4x-5}}dx=
          \left| \begin{array}{c}
          \textit{Checking the differential of the quadratic function:}\\
     d(x^2+4x-5)=(2x+4) dx;\\[0.2cm]
      \textit{Descomposing the numerator} \\
     (8x-5)= 4(2x+4)-21
         \end{array}\right| = \\
         = & 4\int \frac{2x+4}{\sqrt{x^2+4x-5}}dx-21\int \frac{dx}{\sqrt{x^2+4x-5}} = \left| \begin{array}{c}
     \textit{Completing the square:}\\
    x^2+4x-5=\\
    =(x+2)^2-3^2
    \end{array}\right| = \\
     = & 4\int \frac{d(x^2+4x-5)}{\sqrt{x^2+4x-5}}-21\int \frac{d(x+2)}{\sqrt{(x+2)^2-3^2}} = \\
     = & 8\sqrt{x^2+4x-5}-21 \ln \left|x+2+\sqrt{x^2+4x-5}\right|+C.
     \end{aligned}$
     \vspace{0.3cm}

   $ \textbf{Result: } \dot{I} =8\sqrt{x^2+4x-5}-21 \ln \left|x+2+\sqrt{x^2+4x-5}\right|+C,\; C\in \mathbb{R}.$

   \end{enumerate}

    \begin{enumerate}[label=\textbf{3.(\alph*)}]
    \item $\begin{aligned}[t]
        \dot{I}= & \int \frac{3x+2}{\sqrt{x+4}}dx=
    \left| \begin{array}{c}
     \textit{Performing the change of variable:} \\
     t=\sqrt{x+4};
     \\ x=t^2-4, \;dx=2t dt
    \end{array}\right|= \\
    = &\int \frac{3(t^2-4)+2}{t} 2t dt =2\int (3t^2-10)dt= 2(t^3-10t)+C=
     \end{aligned}$

     $\begin{aligned}
    = &    \left| \begin{array}{c}
     \textit{Undoing the change of variable:} \\
     t=(x+4)^{1/2}
    \end{array}\right| = 2(x+4)^{3/2}-20(x+4)^{1/2}+C.
     \end{aligned}$
     \vspace{0.3cm}

   $ \textbf{Result: } \dot{I} =2(x+4)^{3/2}-20(x+4)^{1/2}+C,\; C\in \mathbb{R}.$

         \begin{tabular}{|p{6.0cm}|p{7.5cm}|p{2.0cm}|}
\hline
\vspace{0.05mm}$^*$ Solution of Problem 3(a) guided by Irina Blazhievska is available on-line: &
\vspace{5.5mm}
   \url{https://youtu.be/ZOFSo2oDVzQ} & \vspace{-3mm} \includegraphics[height=20mm]{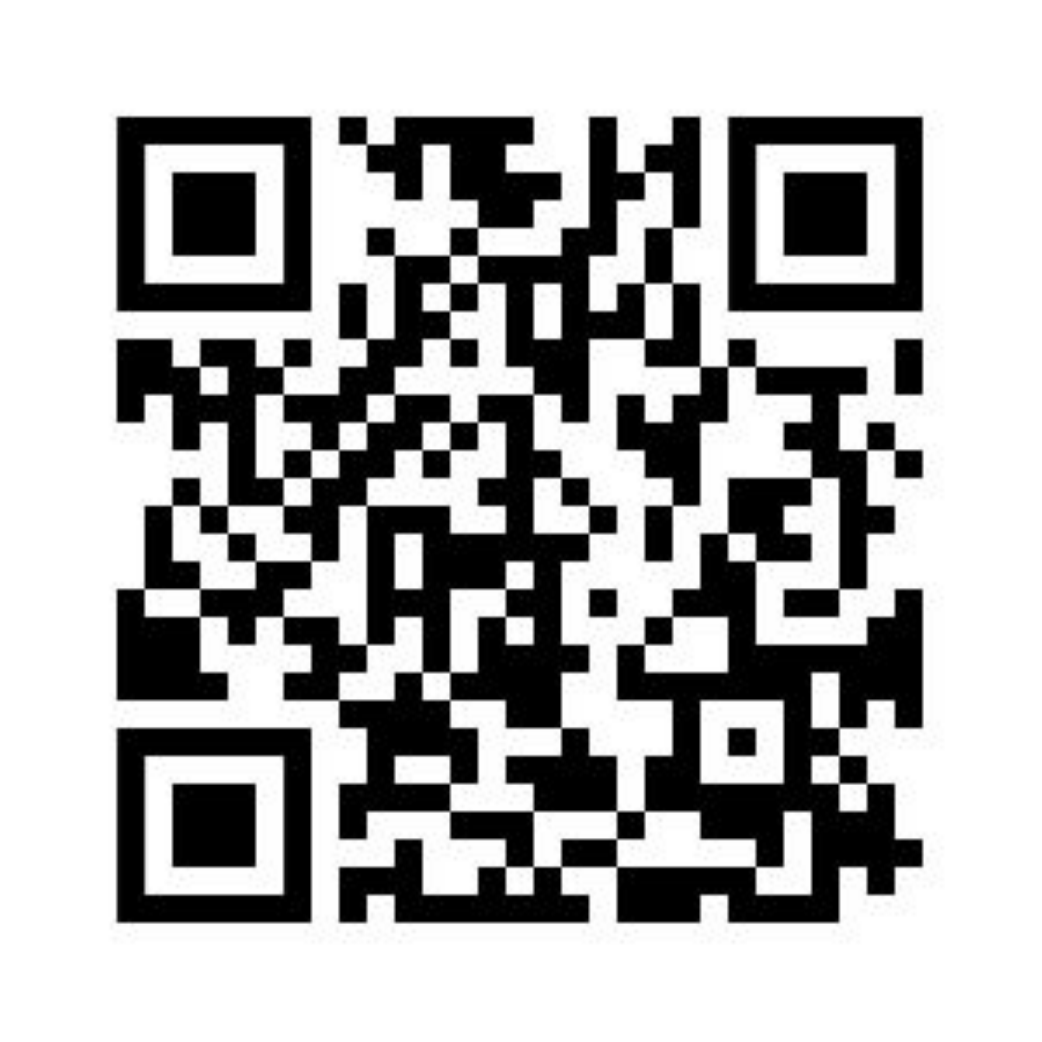}\\
\hline
\end{tabular}

   \item $\begin{aligned}[t]
         \dot{I}= & \int \frac{dx}{x\sqrt{3x^2-4x+1}}=
    \left| \begin{array}{c}
     \textit{Performing the change of variable:}\  t=\dfrac{1}{x};
     \\ x=\dfrac{1}{t}, \; dx=-\dfrac{1}{t^2} dt; \\[0.5cm]
      \sqrt{3x^2-4x+1}=\dfrac{\sqrt{t^2-4t+3}}{t}
    \end{array}\right| = \\
    = & \int \dfrac{-\dfrac{1}{t^2}dt}{\dfrac{1}{t}\cdot \dfrac{\sqrt{t^2-4t+3}}{t}}=  -\int \frac{dt}{\sqrt{t^2-4t+3}}= -\int \dfrac{d(t-2)}{\sqrt{(t-2)^2-1}}= \\
    = & -\ln \left|t-2+\sqrt{t^2-4t+3}\right|+C = \left| \begin{array}{c}
     \textit{Undoing the change of variable:}\\
     t=\dfrac{1}{x}; \\
     \sqrt{t^2-4t+3}=\dfrac{\sqrt{3x^2-4x+1}}{x}
    \end{array}\right| =\\
    = & -\ln \left|\frac{1}{x}-2+\frac{\sqrt{3x^2-4x+1}}{x}\right|+C.
     \end{aligned}$
     \vspace{0.3cm}

   $ \textbf{Result: } \dot{I} =-\ln \left|\dfrac{1}{x}-2+\dfrac{\sqrt{3x^2-4x+1}}{x}\right|+C,\; C\in \mathbb{R}.$

    \item $\begin{aligned}[t]
        \dot{I}= & \int\dfrac{4x^2-3}{e^{2x}}dx=\int (4x^2-3)e^{-2x}dx=
    \left| \begin{array}{c}
     \textit{Using integration by parts:} \\
     \int udv=uv -\int vdu, \\[0.2cm]
     u=4x^2-3; \quad dv=e^{-2x}dx;\\
     du=8x dx; \quad v=-\dfrac{1}{2}e^{-2x}
    \end{array}\right| = \\
    = & -\frac{1}{2}e^{-2x}(4x^2-3)-\int \left(-\frac{1}{2}\right) 8x e^{-2x}dx= \end{aligned}$

     $\begin{aligned}
     = & -\frac{1}{2}e^{-2x}(4x^2-3)+4\int x e^{-2x}dx=
     \left| \begin{array}{c}
     \textit{Using integration by parts:} \\
     \int udv=uv -\int vdu, \\[0.2cm]
     u=x; \quad dv=e^{-2x}dx;\\
     du=dx; \quad v=-\dfrac{1}{2}e^{-2x}
    \end{array}\right| = \\
     = & -\frac{1}{2}e^{-2x}(4x^2-3)+4\left(-\dfrac{1}{2}xe^{-2x}-\int \left(-\frac{1}{2}\right)e^{-2x}dx\right)= \\
     = & -\frac{1}{2}(4x^2-3)e^{-2x}-2xe^{-2x}-e^{-2x}+C= \left(-2x^2-2x+\frac{1}{2}\right)e^{-2x}+C.
     \end{aligned}$
     \vspace{0.3cm}

   $ \textbf{Result: } \dot{I} =\left(-2x^2-2x+\frac{1}{2}\right)e^{-2x}+C,\; C\in \mathbb{R}.$

   \item $\begin{aligned}[t]
        \dot{I}= & \int x^2 \ln x dx=
    \left| \begin{array}{c}
     \textit{Using integration by parts:} \\
     \int udv=uv -\int vdu, \\[0.2cm]
     u=\ln x; \quad dv=x^2 dx;\\
     du= \dfrac{dx}{x}; \quad v=\dfrac{1}{3}x^3
    \end{array}\right| = \frac{1}{3}x^3\ln x-\int \frac{1}{3}x^3\frac{dx}{x} = \\
    = & \frac{1}{3}x^3\ln x-\frac{1}{3}\int x^2dx=  \frac{1}{3}x^3\ln x- \frac{1}{3}\cdot \frac{1}{3}x^3+C= \frac{1}{9}x^3(3\ln x-1)+C.
     \end{aligned}$
     \vspace{0.3cm}

   $ \textbf{Result: } \dot{I} =\dfrac{1}{9}x^3(3\ln x-1)+C,\; C\in \mathbb{R}.$

           \begin{tabular}{|p{6.0cm}|p{7.5cm}|p{2.0cm}|}
\hline
\vspace{0.05mm}$^*$ Solution of Problem 3(d) guided by Ricard Riba is available on-line: &
\vspace{5.5mm}
\url{https://youtu.be/hDYt-m7ZCgM} & \vspace{-3mm} \includegraphics[height=20mm]{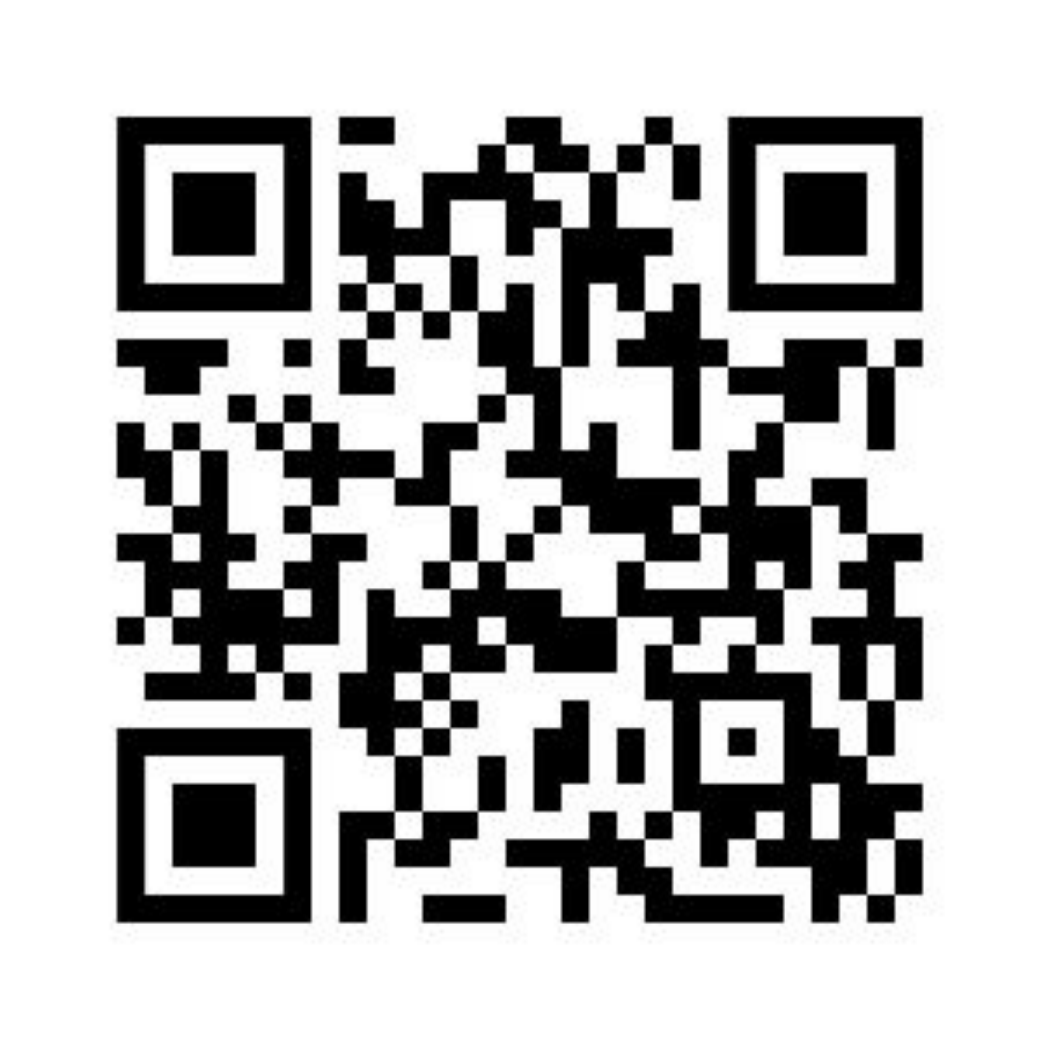}\\
\hline
\end{tabular}

   \item $\begin{aligned}[t]
        \dot{I}= & \int \sin \sqrt{x+1} dx =
    \left| \begin{array}{c}
     \textit{Performing the change of variable:} \  t=\sqrt{x+1};\\
     x=t^2-1,\, dx=2t dt
    \end{array}\right| = \\
    =&\int 2t \sin t dt=    \left| \begin{array}{c}
     \textit{Using integration by parts:} \\
     \int udv=uv -\int vdu, \\[0.2cm]
     u=2t; \quad dv=\sin t dt; \\
     du= 2 dt; \quad v=-\cos t
    \end{array}\right| = -2t \cos t -\int (-\cos t) 2 dt= \\
    &=-2t \cos t +2\sin t +C= \left| \begin{array}{c}
     \textit{Undoing the change of variable:} \ t=\sqrt{x+1}
    \end{array}\right| = \\
    = & -2\sqrt{x+1} \cos \sqrt{x+1}+2\sin \sqrt{x+1}+C.
     \end{aligned}$
     \vspace{0.3cm}

   $ \textbf{Result: } \dot{I} =-2\sqrt{x+1} \cos \sqrt{x+1}+2\sin \sqrt{x+1}+C,\; C\in \mathbb{R}.$

   \item $\begin{aligned}[t]
        \dot{I}= &\int x \frac{\sin x}{\cos^3 x} dx=
    \left| \begin{array}{l}
     \textit{Using integration by parts: } \int udv=uv -\int vdu \\
     u= x; \quad v=\int \dfrac{\sin x}{\cos^3 x} dx=-\int \dfrac{d(\cos x)}{\cos^3 x}=\dfrac{1}{2\cos^{2}x}; \\
     du= dx; \quad dv=\dfrac{\sin x}{\cos^3 x} dx
    \end{array}\right| = \\
    = & \dfrac{x}{2\cos^{2}x}- \dfrac{1}{2}\int\frac{dx}{\cos^2x}=
    \frac{x}{2\cos^2 x}-\frac{1}{2}\tg x +C=  \frac{1}{2}\left(\frac{x}{\cos^2x}-\tg x\right) +C.
     \end{aligned}$
     \vspace{0.3cm}

   $ \textbf{Result: } \dot{I} =\dfrac{1}{2}\left(\dfrac{x}{\cos^2x}-\tg x \right) +C,\; C\in \mathbb{R}.$
\end{enumerate}

\endgroup

\begin{enumerate}[label=\textbf{4.(\alph*)}]
\item $$\dot{I}=\Int\dfrac{-3x^2+2x+13}{x^3+2x^2-x-2} dx$$

Algorithm:

\textbf{Step 1.} The function under the integral is a suitable polynomial fraction:
$$f(x)=\dfrac{Q_2(x)}{P_3(x)}.$$

\textbf{Step 2.} Finding the roots of the denominator $P_3(x)=x^3+2x^2-x-2.$
Possible integer roots may be $\pm 1;$ $\pm 2.$
$$\begin{array}{c}
P_3(1)=1^3+2\cdot 1^2-1-2=0 \Rightarrow x_1=1 \text{ is a root of } P_3(x) \\[0.1cm]
\Rightarrow P_3(x) \text{ is divisible by } (x-1).
\end{array}$$
The decomposition $P_3(x)=(x-1)(x^2+3x+2)=(x-1)(x+2)(x+1),$
implies that all roots are simple and real-valued.

\textbf{Step 3.} Applying the method of unknown coefficients to the suitable fraction:
\begin{align*}
f(x)= & \frac{-3x^2+2x+13}{(x-1)(x+2)(x+1)}=\frac{A}{x-1}+\frac{B}{x+1}+\frac{C}{x+2}= \\
= & \frac{A(x+2)(x+1)+B(x-1)(x+2)+C(x-1)(x+1)}{(x-1)(x+2)(x+1)}.
\end{align*}
Numerator's equality:
$$-3x^2+2x+13=A(x+2)(x+1)+B(x-1)(x+2)+C(x-1)(x+1).$$
Applying equalities by zeros:
\begin{align*}
x=1 & \Rightarrow -3+2+13=A(1+1)(1+2); \ 12=6A \Rightarrow A=2 \\
x=-1 & \Rightarrow -3-2+13=B(-1-1)(-1+2); \ 8=-2B \Rightarrow B=-4 \\
x=-2 & \Rightarrow -3\cdot 4 -4+13=C(-2-1)(-2+1); \ -3=3C\Rightarrow C=-1.
\end{align*}
The resulting decomposition is as follows:
$$f(x)=\frac{-3x^2+2x+13}{(x-1)(x+2)(x+1)}=\frac{2}{x-1}-\frac{4}{x+1}-\frac{1}{x+2}.$$

\textbf{Step 4.} Substitution of the fully decomposed fraction under the integral:
$$ \Int\dfrac{-3x^2+2x+13}{x^3+2x^2-x-2}dx=  \int \left( \frac{2}{x-1}-\frac{4}{x+1}-\frac{1}{x+2} \right) dx= $$
$$= 2\ln|x-1|-4\ln|x+1|-\ln|x+2|+C.$$

\textbf{Result:} $\dot{I}=2\ln|x-1|-4\ln|x+1|-\ln|x+2|+C, C\in \mathbb{R}.$

\item $$\dot{I}=\int \dfrac{3x^3-32x+56}{x^3-2x^2-4x+8}dx$$

Algorithm:

\textbf{Step 1.} The function under the integral is a non-suitable polynomial fraction:
$$f(x)=\dfrac{Q_3(x)}{P_3(x)}.$$
Separation of integer part leads us to the following sum of 0-degree polynomial and a suitable polynomial fraction:
$$f(x)=3+\dfrac{6x^2-20x+32}{x^3-2x^2-4x+8}.$$
\textbf{Step 2.} Finding the roots of the denominator $P_3(x)=x^3-2x^2-4x+8.$
Possible integer roots may be $\pm 1,\pm 2,\pm 4,\pm 8.$
$$\begin{array}{c}
P_3(2)=2^3-2\cdot 2^2-4\cdot2+8 =0\Rightarrow x_1=8 \text{ is a root of } P_3(x)\\[0.1cm]
\Rightarrow P_3(x) \text{ is divisible by } (x-2).
\end{array}$$
The decomposition $P_3(x)=(x-2)(x^2-4)=(x-2)^2(x+2),$
implies that $P_3$ has two real-valued roots, a double root $x=2$ and a simple root $x=-2.$

\textbf{Step 3.} Applying the method of unknown coefficients to the suitable fraction:
\begin{align*}
\dfrac{6x^2-20x+32}{x^3-2x^2-4x+8} = & \frac{6x^2-20x+32}{(x-2)^2(x+2)}=  \frac{A}{x-2}+\frac{B}{(x-2)^2}+\frac{C}{x+2}= \\
= & \frac{A(x-2)(x+2)+B(x+2)+C(x-2)^2}{(x-2)^2(x+2)}.
\end{align*}
Numerator's equality:
\begin{align*}
6x^2-20x+32= & A(x-2)(x+2)+B(x+2)+C(x-2)^2=\\
= & x^2(A+C)+x(B-4C)+(-4A+2B+4C).
\end{align*}
Mix of equalities by zeros with equality of coefficients of some monomial:
\begin{align*}
x=2 & \Rightarrow 6\cdot 2^2-20\cdot 2+ 32=B(2+2); \ 16=4B\Rightarrow B=4 \\
x=-2 & \Rightarrow 6\cdot (-2)^2+32=C(-2-2)^2; \ 96=16C \Rightarrow C=6 \\
x^2 & : \  6=A+C; \ A=6-C=6-6=0\Rightarrow A=0.
\end{align*}
The resulting decomposition is as follows:
$$\dfrac{6x^2-20x+32}{x^3-2x^2-4x+8}=\frac{6x^2-20x+32}{(x-2)^2(x+2)}=\frac{4}{(x-2)^2}+\frac{6}{x+2}.$$

\textbf{Step 4.} Substitution of the fully decomposed fraction under the integral:
$$\int \dfrac{3x^3-32x+56}{x^3-2x^2-4x+8}dx=\int \left( 3+ \frac{4}{(x-2)^2}+\frac{6}{x+2} \right)dx=$$
$$= 3x-\frac{4}{x-2}+6\ln |x+2|+C.$$

\textbf{Result:} $\dot{I}= 3x-\dfrac{4}{x-2}+6\ln |x+2|+C, \ C\in \mathbb{R}.$

\begin{tabular}{|p{6.0cm}|p{7.5cm}|p{2.0cm}|}
\hline
\vspace{0.05mm}$^*$ Solution of Problem 4(b) guided by Irina Blazhievska is available on-line: &
\vspace{5.5mm} \url{https://youtu.be/X_300OO1i8A} & \vspace{-3mm} \includegraphics[height=20mm]{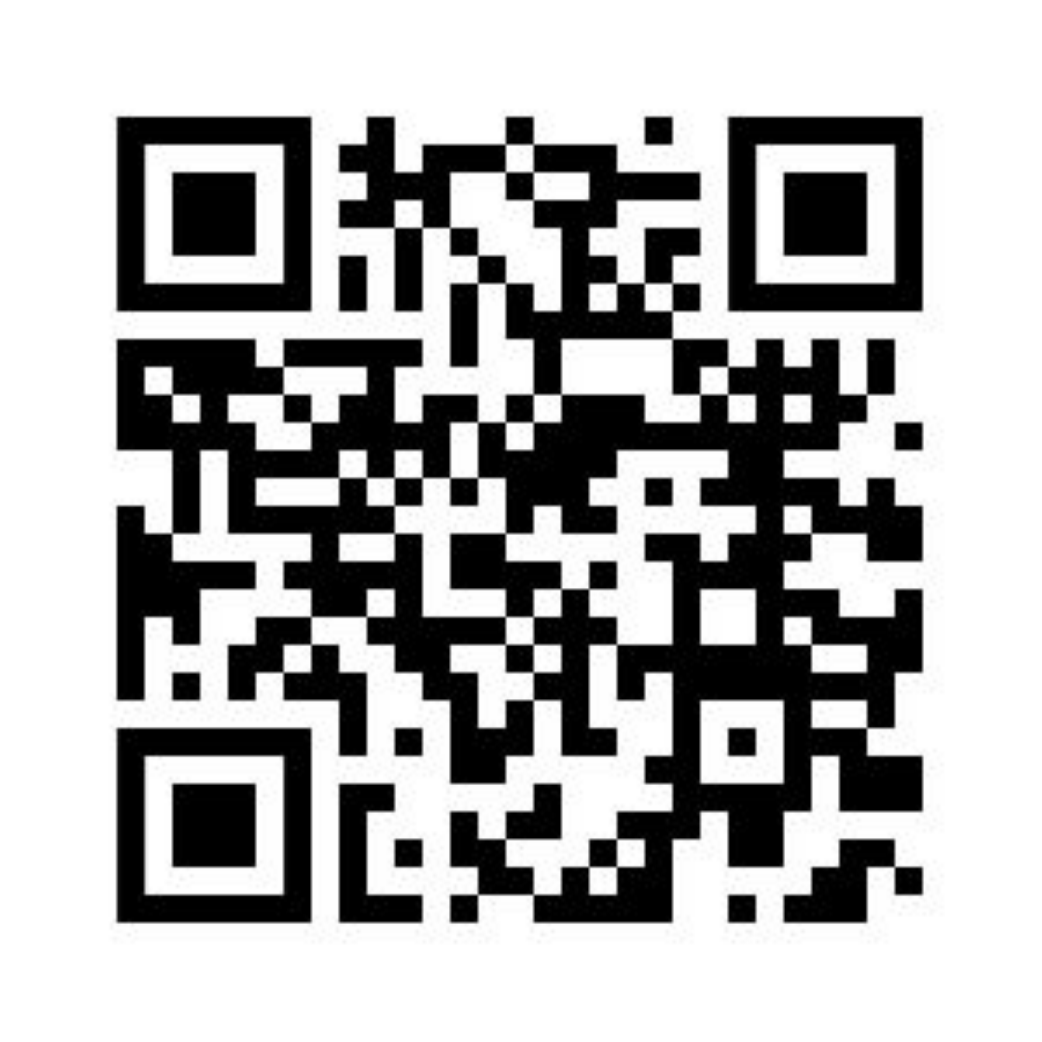}\\
\hline
\end{tabular}

\item
$$\dot{I}=\Int\dfrac{x^2-2x-9}{(x^2+4x+5)(x-1)}dx$$

Algorithm:

\textbf{Step 1.} The function under the integral is a suitable polynomial fraction:
$$f(x)=\dfrac{Q_2(x)}{P_3(x)}.$$
\textbf{Step 2.} Finding the roots of the denominator $P_3(x)=(x^2+4x+5)(x-1).$

The polynomial $x^2+4x+5$ is indecomposable quadratic function since its discriminant is negative.
The decomposition of $P_3(x)$ implies that the unique real-valued root of $P_3(x)$ is $x_1=1.$

\textbf{Step 3.} Applying the method of unknown coefficients to the suitable fraction:
\begin{align*}
f(x)= & \dfrac{x^2-2x-9}{(x^2+4x+5)(x-1)}=\frac{Ax+B}{x^2+4x+5}+\frac{C}{x-1}= \\
= & \dfrac{(Ax+B)(x-1)+C(x^2+4x+5)}{(x^2+4x+5)(x-1)}.
\end{align*}
Numerator's equality:
\begin{align*}
x^2-2x-9= & (Ax+B)(x-1)+C(x^2+4x+5)= \\
= & x^2(A+C)+x(-A+B+4C)+(-B+5C).
\end{align*}
Mix of equalities by zeros with equality of coefficients of some monomial:
\begin{align*}
x=1 & \Rightarrow 1-2-9=C(1+4+5); \ -10=10C \Rightarrow C=-1 \\
x^2 & : \ 1=A+C ; \ A= 1-C=1-(-1)=2 \Rightarrow A=2 \\
x^0 & : \ -9=-B+5C; \ B=9+5C=9+5\cdot (-1)=4\Rightarrow B=4.
\end{align*}
The resulting decomposition is as follows:
$$f(x)=\dfrac{x^2-2x-9}{(x^2+4x+5)(x-1)}=\frac{2x+4}{x^2+4x+5}-\frac{1}{x-1}.$$

\textbf{Step 4.} Substitution of the fully decomposed fraction under the integral:
$$ \Int\dfrac{x^2-2x-9}{(x^2+4x+5)(x-1)}dx=  \int \left( \frac{2x+4}{x^2+4x+5}-\frac{1}{x-1} \right) dx= $$
$$=\int\frac{d(x^2+4x+5)}{x^2+4x+5}-\int\frac{dx}{x-1}= \ln(x^2+4x+5)-\ln|x-1|+C.$$

\textbf{Result:} $\dot{I}=\ln(x^2+4x+5)-\ln|x-1|+C, C\in \mathbb{R}.$

\end{enumerate}

\begingroup
\addtolength{\jot}{0.2cm}

\begin{enumerate}[label=\textbf{5.(\alph*)}]
\item $\begin{aligned}[t]
	\dot{I}= & \int \sin 10 x \sin 3x dx=
    \left| \begin{array}{c}
     \textit{Applying the decomposition of} \\
     \textit{a product of sinus functions:} \\
     \sin \alpha \sin \beta=\dfrac{1}{2}\big(\cos(\alpha-\beta)-\cos(\alpha +\beta)\big)
    \end{array}\right| = \\
    = & \frac{1}{2}\int \big(\cos(10x-3x)-\cos(10x+3x)\big)dx= \frac{1}{2}\int \big(\cos7x-\cos13x\big)dx= \\
    = & \frac{1}{2}\cdot \frac{\sin 7x}{7}-\frac{1}{2}\cdot \frac{\sin 13 x}{13}+C= \frac{1}{14}\sin 7x -\frac{1}{26}\sin 13x +C.
    \end{aligned}$
    \vspace{0.3cm}

$ \textbf{Result: } \dot{I} =  \dfrac{1}{14}\sin 7x -\dfrac{1}{26}\sin 13x +C,\; C\in \mathbb{R}.$

\begin{tabular}{|p{6.0cm}|p{7.5cm}|p{2.0cm}|}
\hline
\vspace{0.05mm}$^*$ Solution of Problem 5(a) guided by Irina Blazhievska is available on-line: &
\vspace{5.5mm}  \url{https://youtu.be/-Zsc5t-YIOk} & \vspace{-3mm} \includegraphics[height=20mm]{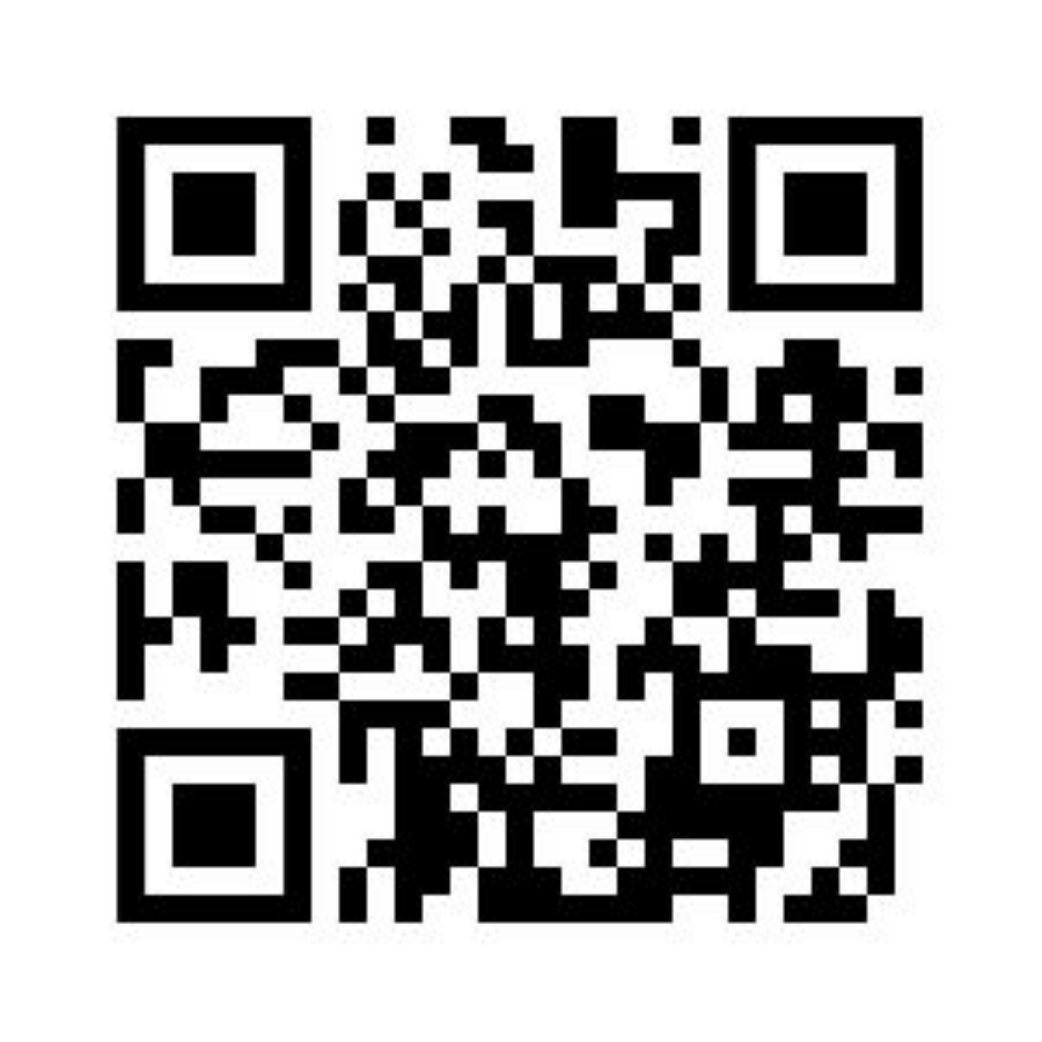}\\
\hline
\end{tabular}

\item $\begin{aligned}[t]
	\dot{I}= & \int \sqrt[3]{\dfrac{\sin x}{\cos^{13} x}}dx=
    \left| \begin{array}{c}
     \textit{Applying the integration of a product} \\
    \textit{$\sin^m x\cos^n x$ with $m=0,$ $n=-4,$ decompose}:
   \\ \cos^{-13/3} x=\cos^{-1/3} x\cos^{-4} x; \\
   \dfrac{1}{\cos^2x}=1+\tg^2x; \
     \dfrac{dx}{\cos^2 x}=d(\tg x)
    \end{array}\right| = \\
    =& \int \sqrt[3]{\dfrac{\sin x}{\cos x}} \frac{dx}{\cos^{4} x}=  \int \tg^{1/3} x \frac{1}{\cos^{2} x}\frac{dx}{\cos^{2} x}=\int \tg^{1/3} x \big(1+\tg^2x\big)d(\tg x)= \\
    = & \int\big(\tg^{1/3} x+\tg^{7/3}x\big)d(\tg x)=\frac{3}{4}\tg^{4/3}x+\frac{3}{10}\tg^{10/3}x+C.
    \end{aligned}$
    \vspace{0.3cm}

$ \textbf{Result: } \dot{I} =  \dfrac{3}{4}\tg^{4/3}x+\dfrac{3}{10}\tg^{10/3}x+C,\; C\in \mathbb{R}.$

\item
$\begin{aligned}[t]
\dot{I}= & \int \frac{dx}{4\cos x+3\sin x+6}=
    \left| \begin{array}{c}
     R(\sin x, \cos x)=\dfrac{1}{4\cos x+3\sin x+6}\\[0.3cm]
     \textit{has a general form.} \\
     \textit{Applying the universal subtitution:} \  t=\tg \dfrac{x}{2}; \\
     dx=\dfrac{2dt}{1+t^2}; \ \sin x=\dfrac{2t}{1+t^2}, \ \cos x=\dfrac{1-t^2}{1+t^2}
    \end{array}\right| = \\
= & \int \frac{\dfrac{2dt}{1+t^2}}{4\dfrac{1-t^2}{1+t^2}+3\dfrac{2t}{1+t^2}+6}=  \int \frac{2dt}{4(1-t^2)+6t+6(1+t^2)}=\int \frac{dt}{t^2+3t+5}=
\end{aligned}$

$\begin{aligned}
=& \left| \begin{array}{c}
   \textit{Completing the square:} \ t^2+3t+5=\\
   =\big(t+\frac{3}{2}\big)^2-\big(\frac{3}{2}\big)^2+5= \big(t+\frac{3}{2}\big)^2+\left(\frac{\sqrt{11}}{2}\right)^2
    \end{array}\right| =\int \dfrac{d\left(t+\frac{3}{2}\right)}{\left(t+\frac{3}{2}\right)^2+\left(\frac{\sqrt{11}}{2}\right)^2}=\\
=&\dfrac{1}{\frac{\sqrt{11}}{2}}\arctg \left(\frac{t+\frac{3}{2}}{\frac{\sqrt{11}}{2}}\right)+C = \left| \begin{array}{c}
   \textit{Undoing}\\
   \textit{the change} \\
   \textit{of variable:} \\
   t=\tg \dfrac{x}{2}
    \end{array}\right| = \dfrac{2}{\sqrt{11}}\arctg \left(\frac{2\tg \dfrac{x}{2}+3}{\sqrt{11}}\right)+C.
\end{aligned}$

$ \textbf{Result: } \dot{I} =  \dfrac{2}{\sqrt{11}}\arctg \left( \dfrac{2\tg \dfrac{x}{2}+3}{\sqrt{11}}\right)+C,\; C\in \mathbb{R}.$

\vspace{-0.2cm}
\begin{tabular}{|p{6.0cm}|p{7.5cm}|p{2.0cm}|}
\hline
\vspace{0.05mm}$^*$ Solution of Problem 5(c) guided by Irina Blazhievska is available on-line: &
\vspace{5.5mm} \url{https://youtu.be/loxi3dwTmho} & \vspace{-3mm} \includegraphics[height=20mm]{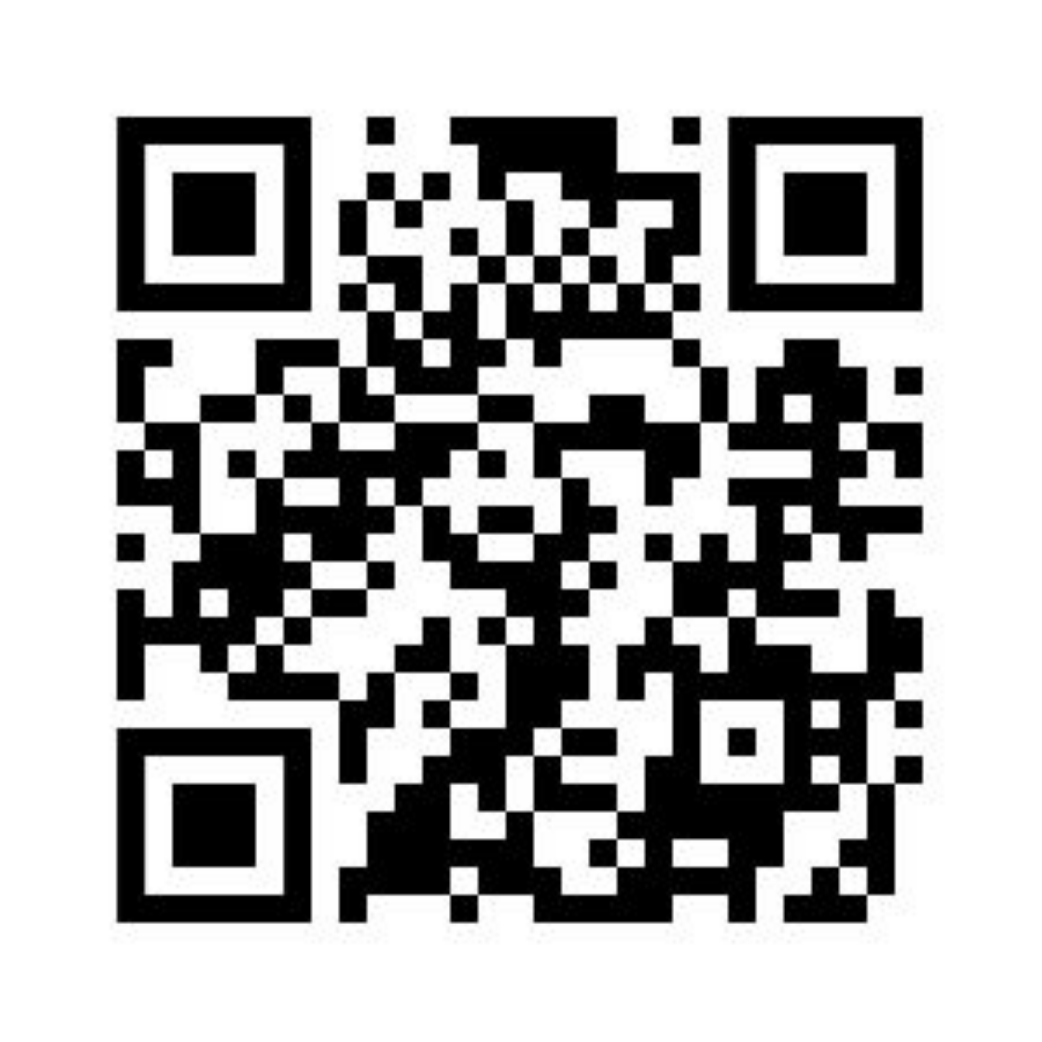}\\
\hline
\end{tabular}
\end{enumerate}

\begin{enumerate}[label=\textbf{6.(\alph*)}]
\item $\begin{aligned}[t]
\dot{I}= & \int \frac{dx}{x^2\sqrt{4-x^2}}=
    \left| \begin{array}{c}
     \textit{Performing the change of variable:} \\
     x=2\sin t;\  dx=2\cos t dt;\\
      \sqrt{4-x^2}=2\cos t
    \end{array}\right| = \int \dfrac{2 \cos t dt}{(2\sin t)^2 2\cos t}=\\
    = & \frac{1}{4}\int \frac{dt}{\sin^2 t}=-\frac{1}{4}\ctg t +C=  \left| \begin{array}{c}
   \textit{Undoing the change of variable:} \\
   \ctg t=\dfrac{2\cos t}{2\sin t}=\dfrac{\sqrt{4-x^2}}{x}
    \end{array}\right| =\\
    =& -\dfrac{1}{4}\dfrac{\sqrt{4-x^2}}{x}=-\dfrac{1}{4}\sqrt{\dfrac{4}{x^2}-1}+C.
    \end{aligned}$
    \vspace{0.3cm}

$ \textbf{Result: } \dot{I} =-\dfrac{1}{4}\sqrt{\dfrac{4}{x^2}-1}+C,\; C\in \mathbb{R}.$

\begin{tabular}{|p{6.0cm}|p{7.5cm}|p{2.0cm}|}
\hline
\vspace{0.05mm}$^*$ Solution of Problem 6(a) guided by Ricard Riba is available on-line: &
\vspace{5.5mm}
\url{https://youtu.be/ba4kyucyLVk} & \vspace{-3mm} \includegraphics[height=20mm]{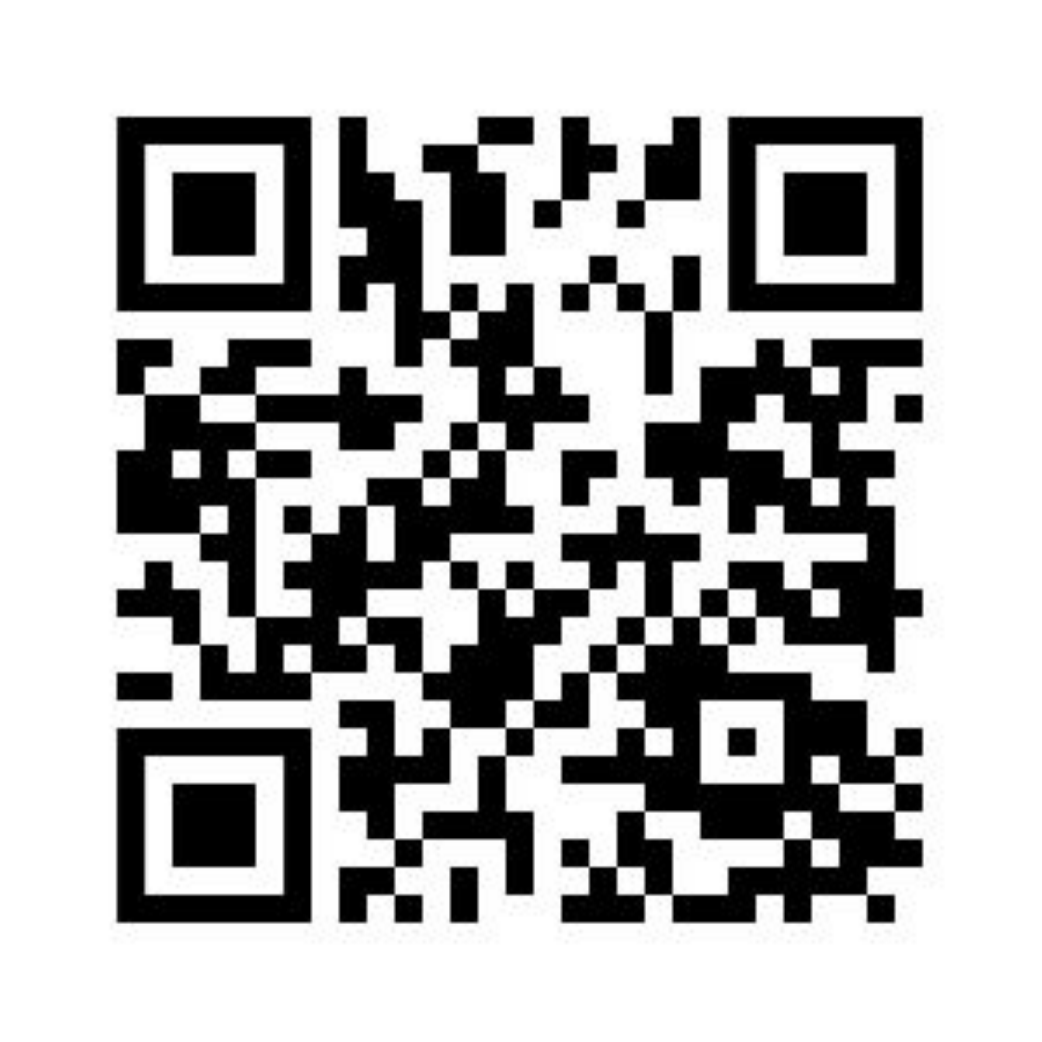}\\
\hline
\end{tabular}

$\bullet$ \textit{Alternative solution 1}
\begin{align*}
\dot{I}= & \int \frac{dx}{x^2\sqrt{4-x^2}}=\int \dfrac{dx}{x^3\sqrt{\dfrac{4}{x^2}-1}}=
\left|\begin{array}{c}
\textit{Substitution under the differential} \\
d\left(\dfrac{4}{x^2}-1\right)=-8\dfrac{dx}{x^3} \\
\dfrac{dx}{x^3}=-\dfrac{1}{8}d\left(\dfrac{4}{x^2}-1\right)
\end{array}\right|= \\
= & -\dfrac{1}{8}\int \left(\dfrac{4}{x^2}-1\right)^{-1/2} d\left(\dfrac{4}{x^2}-1\right)=-\dfrac{2}{8}\left(\dfrac{4}{x^2}-1\right)^{1/2}+C= -\dfrac{1}{4}\sqrt{\dfrac{4}{x^2}-1}+C.
\end{align*}

$\bullet$ \textit{Alternative solution 2}
\begin{align*}
\dot{I}= & \int \frac{dx}{x^2\sqrt{4-x^2}}=
\left|\begin{array}{c}
\textit{Performing the change}\\
\textit{ of variable:} \ t=\dfrac{1}{x}\\
x=\dfrac{1}{t}, \ dx=-\dfrac{1}{t^2}dt
\end{array}\right|= \int \dfrac{-\dfrac{1}{t^2} dt}{\dfrac{1}{t^2} \sqrt{4-\left(\dfrac{1}{t}\right)^2}}=\\
 =&-\int \frac{t}{\sqrt{4t^2-1}}dt = \left|\begin{array}{c}
d(4t^2-1)=8tdt; \\
dt=\dfrac{1}{8}d(4t^2-1)
\end{array}\right|=-\frac{1}{8}\int \frac{d(4t^2-1)}{\sqrt{4t^2-1}}=
\end{align*}
\begin{align*}
=& -\frac{2}{8}\sqrt{4t^2-1}+C =  \left|\begin{array}{c}
\textit{Undoing the change}\\
\textit{of variable:} \ t=\dfrac{1}{x}
\end{array}\right| = -\frac{1}{4}\sqrt{\dfrac{4}{x^2}-1}+C.
\end{align*}

$\bullet$ \textit{Alternative solution 3}
\begin{align*}
\dot{I}= & \int \dfrac{dx}{x^2\sqrt{4-x^2}}=
    \left| \begin{array}{c}
     \textit{Performing the change of variable:} \\
     x=2\tgh t,\  dx=\dfrac{2dt}{\ch^2t};\\
      \sqrt{4-x^2}=\sqrt{4-4\tgh^2 t}=\dfrac{2}{\ch t}
    \end{array}\right| = \int \dfrac{\dfrac{2dt}{\ch^2t}}{(2\tgh t)^2\dfrac{2}{\ch t}}=\\
    = & \frac{1}{4}\int \frac{\ch t dt}{\sh^2 t}=\frac{1}{4}\int \frac{d(\sh t)}{\sh^2t} =-\frac{1}{4}\dfrac{1}{\sh t} +C=   \left| \begin{array}{c}
   \textit{Undoing the change of}\\
   \textit{variable:} \\
   \textrm{sh}t=2\tgh t\dfrac{\ch t}{2}=\dfrac{x}{\sqrt{4-x^2}}\\
    \end{array}\right| =\\
     & =  -\frac{1}{4}\frac{\sqrt{4-x^2}}{x}+C= -\dfrac{1}{4}\sqrt{\dfrac{4}{x^2}-1}+C.
    \end{align*}

\item
$\begin{aligned}[t]
\dot{I}= & \Int\sqrt[3]{\dfrac{x+1}{x-1}}\dfrac{dx}{(x-1)^{3}}=  \left| \begin{array}{c}
     \textit{Performing the change of variable:} \\
     t=\sqrt[3]{\dfrac{x+1}{x-1}}; \ \dfrac{x+1}{x-1}=t^3; \\
    x=\dfrac{t^3+1}{t^3-1}=1+\dfrac{2}{t^3-1}; \  x-1=\dfrac{2}{t^3-1};\\
    dx=-\dfrac{6t^2}{(t^3-1)^2}dt\\
       \end{array}\right| = \\
    = & \int t\cdot \frac{-6t^2}{(t^3-1)^2}\left(\frac{t^3-1}{2}\right)^3dt=  -\frac{6}{8} \int t^3(t^3-1)dt=
     -\frac{3}{4}\int (t^6-t^3)dt= \\
    = & -\frac{3}{4}\left(\frac{t^7}{7}-\frac{t^4}{4}\right)+C=
     -\frac{3}{28}t^7+\frac{3}{16}t^4+C= \\
    = & \left| \begin{array}{c}
   \textit{Undoing the change of variable:} \\
   t=\Big(\dfrac{x+1}{x-1}\Big)^{1/3}
    \end{array}\right| = -\frac{3}{28}\left(\dfrac{x+1}{x-1}\right)^{7/3}+\frac{3}{16}\left(\dfrac{x+1}{x-1}\right)^{4/3}+C.
\end{aligned}$
\vspace{0.3cm}

$ \textbf{Result: } \dot{I} = -\dfrac{3}{28}\left(\dfrac{x+1}{x-1}\right)^{7/3}+\dfrac{3}{16}\left(\dfrac{x+1}{x-1}\right)^{4/3}+C,\; C\in \mathbb{R}.$

\end{enumerate}

\begin{enumerate}[label=\textbf{7.(\alph*)}]
\item
$
\begin{aligned}[t]
\dot{I}= & \int^1_0 \ln (1+x^2) dx=
\left|\begin{array}{c}
\textit{Using integration by parts:}\  \int^b_a udv= uv\mid^b_a-\int^b_a vdu, \\
u=\ln(1+x^2);\quad dv=dx; \\
du=\dfrac{2x dx}{1+x^2};\quad v=x.
\end{array}\right|= \\
= & x\ln(1+x^2)\Big|^1_0-\int^1_0 \frac{2x^2}{1+x^2}dx =  \ln 2- 2\int^1_0 \frac{(x^2+1)-1}{1+x^2}dx = \\
= &\ln 2- 2\int^1_0 \left(1-\frac{1}{1+x^2}\right)dx =  \ln 2-2\Big[x-\arctg x\Big]\Bigg|^1_0 = \\
= & \ln2+ (-2+2\arctg 1)-(0+2\arctg 0)= \ln 2+\frac{\pi}{2}-2.
\end{aligned}
$
\vspace{0.3cm}

$\textbf{Result:}\; \dot{I}=\ln 2+\dfrac{\pi}{2}-2.$

\begin{tabular}{|p{6.0cm}|p{7.5cm}|p{2.0cm}|}
\hline
\vspace{0.05mm}$^*$ Solution of Problem 7(a) guided by Irina Blazhievska is available on-line: &
\vspace{5.5mm} \url{https://youtu.be/jfbI2G23U2M} & \vspace{-3mm} \includegraphics[height=20mm]{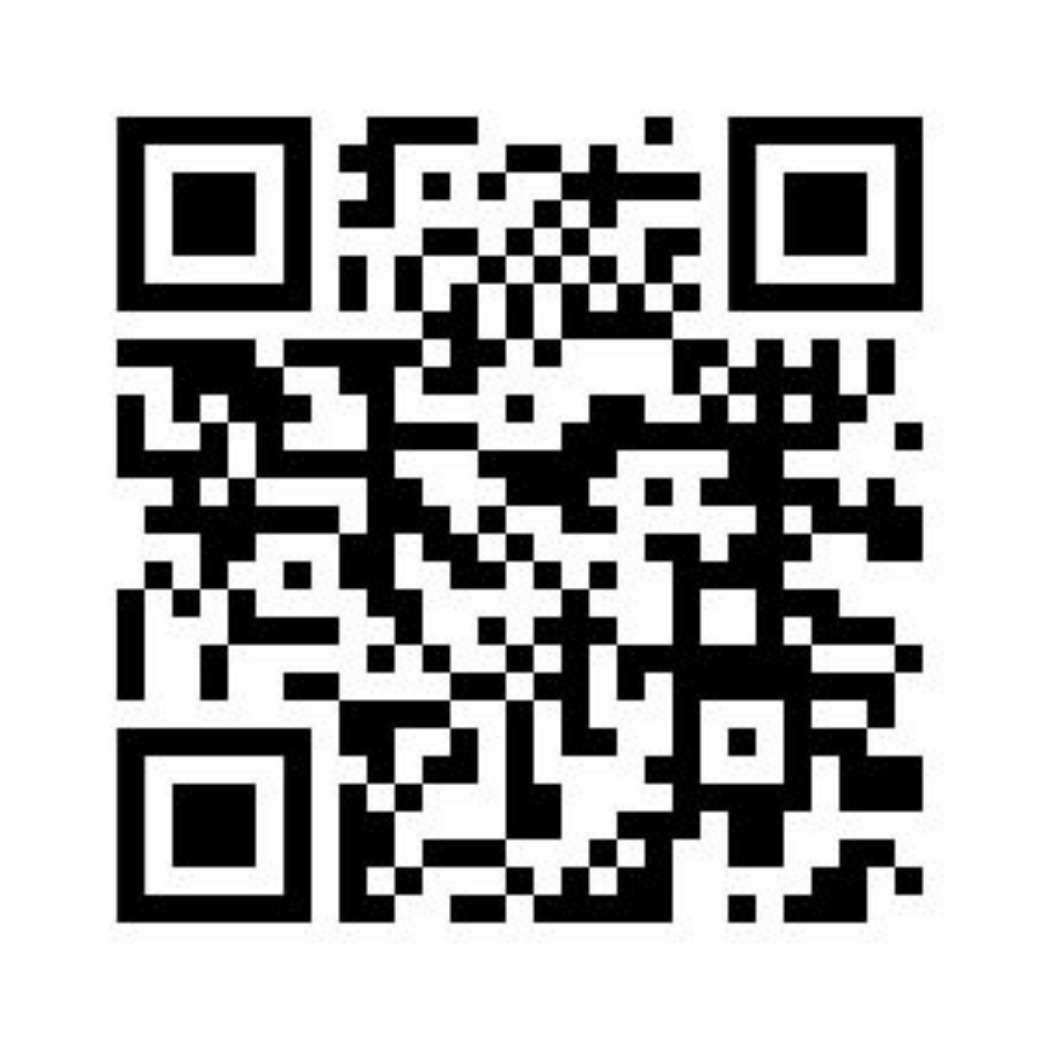}\\
\hline
\end{tabular}

\item
$
\begin{aligned}[t]
\dot{I}= & \int^{\pi/2}_0 \sin^3x\sqrt[4]{\cos x} dx=
\left|\begin{array}{c}
 \textit{Applying the integration of a product} \\
    \textit{$\sin^m x\cos^n x$ with $m=3,$ $n=0$ decompose}:\\
\sin^3 x=\sin^2 x \sin x; \\
\sin^2 x=1-\cos^2 x; \ \sin x dx=-d(\cos x)
\end{array}\right|= \\
= & -\int^{\pi/2}_0 (1-\cos^2 x) \cos^{1/4}x d(\cos x)= \int^{\pi/2}_0 (\cos^{9/4}x-\cos^{1/4}x)d(\cos x)=\\
= & \left[\frac{4}{13}\cos^{13/4}x-\frac{4}{5}\cos^{5/4}x\right]\Bigg|^{\pi/2}_0= 0-\left[\frac{4}{13}-\frac{4}{5} \right]=\frac{32}{65}.
\end{aligned}
$
\vspace{0.3cm}

$\textbf{Result:}\; \dot{I}=\dfrac{32}{65}.$

\begin{tabular}{|p{6.0cm}|p{7.5cm}|p{2.0cm}|}
\hline
\vspace{0.05mm}$^*$ Solution of Problem 7(b) guided by Irina Blazhievska is available on-line: &
\vspace{5.5mm} \url{https://youtu.be/s4VH2LvXh7M} & \vspace{-3mm} \includegraphics[height=20mm]{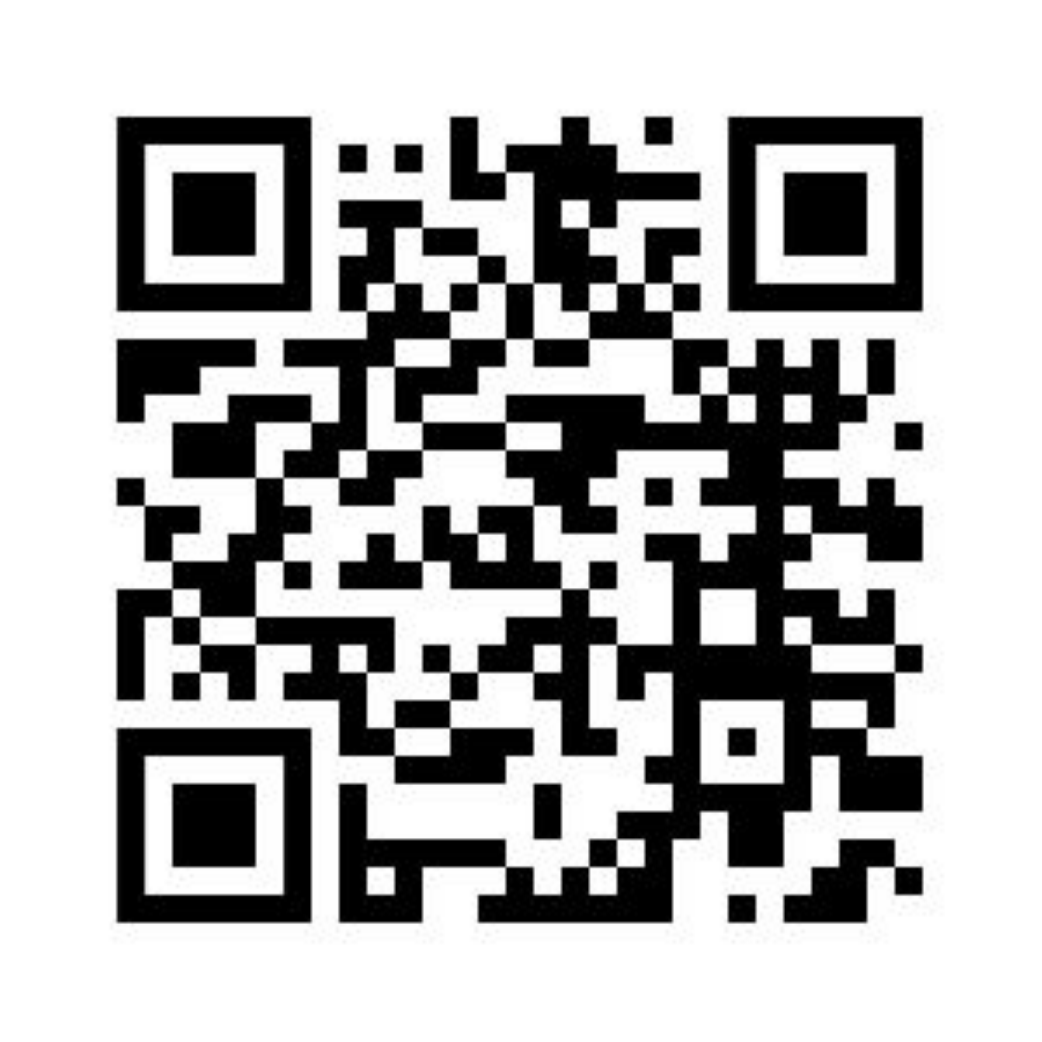}\\
\hline
\end{tabular}

\item
$
\begin{aligned}[t]
\dot{I}= & \int^{\sqrt[8]{2}}_1 \frac{4x^3 dx}{\sqrt{4-x^8}}=\int^{\sqrt[8]{2}}_1 \frac{4x^3 dx}{\sqrt{2^2-(x^4)^2}}=
\left|\begin{array}{c}
\textit{Substitution under the differential}\\
4x^3 dx= d(x^4)
\end{array}\right|= \\
= & \int^{\sqrt[8]{2}}_1 \frac{d(x^4)}{\sqrt{2^2-(x^4)^2}}= \arcsin \frac{x^4}{2} \ \Bigg|^{\sqrt[8]{2}}_1= \arcsin \frac{\sqrt{2}}{2}-\arcsin \frac{1}{2}= \frac{\pi}{4}-\frac{\pi}{6}=\frac{\pi}{12}.
\end{aligned}
$
\vspace{0.3cm}

$\textbf{Result:}\; \dot{I}=\dfrac{\pi}{12}.$

\item
$
\begin{aligned}[t]
\dot{I}= & \int^{16}_1 \frac{1+\sqrt{x}}{\sqrt[4]{x}+\sqrt{x}}dx=
\left|\begin{array}{c}
\textit{Performing the change of variable:}\  t=\sqrt[4]{x}\\
\begin{minipage}[t]{5.5cm}
$x=t^4,\ dx=4t^3dt;$ \\
$x=1\mapsto t=1,$\\
$x=16\mapsto t=2$
\end{minipage}
\begin{minipage}[t]{3cm}
\begin{tabular}[t]{|c|l|l|}
\hline
x & 1& 16 \\ \hline
t & 1 & 2 \\
\hline
\end{tabular}
\end{minipage}
\end{array}\right|= \\
= & \int^{2}_1 \frac{1+t^2}{t+t^2}4t^3dt=  4\int^{2}_1 \frac{t^2+t^4}{1+t}dt=  \left|\begin{array}{c}
\textit{Separating the fraction's integer part}\\
\dfrac{t^2+t^4}{1+t}=t^3-t^2+2t-2+\dfrac{2}{1+t}
\end{array}\right|= \\
= & 4\int^{2}_1 \left(t^3-t^2+2t-2+\dfrac{2}{1+t}\right) dt= \left[t^4-\frac{4}{3}t^3+4t^2-8t+8\ln|1+t|\right]\Bigg|^2_1=
\end{aligned}$

$\begin{aligned}
= & \left(16-\frac{32}{3}+16-16+8\ln3\right)-\left(1-\frac{4}{3}+4-8+8\ln2\right)= \frac{29}{3}+8\ln \frac{3}{2}.
\end{aligned}
$
\vspace{0.3cm}

$\textbf{Result:}\; \dot{I}=\dfrac{29}{3}+8\ln \dfrac{3}{2}.$

\begin{tabular}{|p{6.0cm}|p{7.5cm}|p{2.0cm}|}
\hline
\vspace{0.05mm}$^*$ Solution of Problem 7(d) guided by Irina Blazhievska is available on-line: &
\vspace{5.5mm} \url{https://youtu.be/gCX8MgQrd7A} & \vspace{-3mm} \includegraphics[height=20mm]{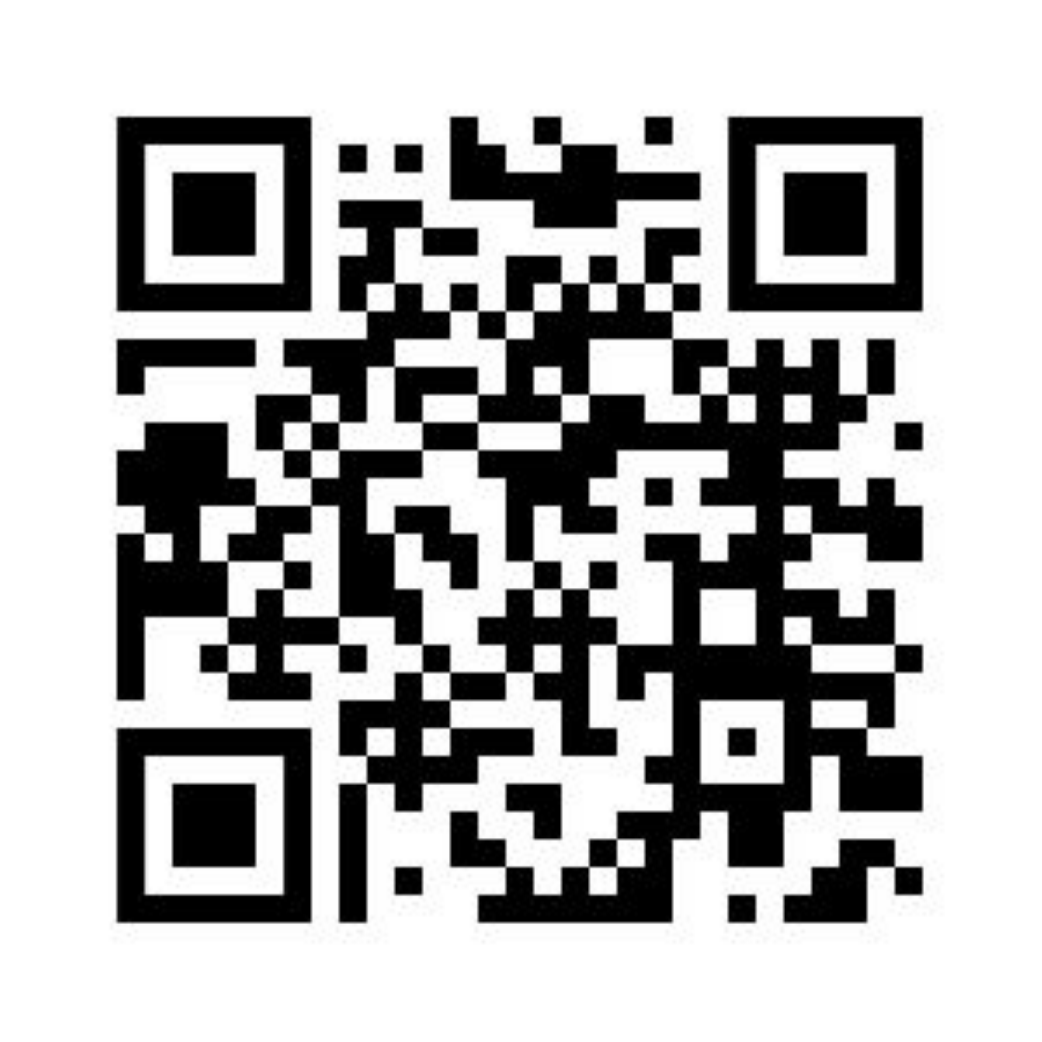}\\
\hline
\end{tabular}

\end{enumerate}

\endgroup
\vspace{0.3cm}

\begin{enumerate}[label=\textbf{8.(\alph*)}]
\item Find the area of the figure bounded by the curves:
$$y=2x^2-10x+6;\quad y=x^2-3x.$$

\begin{tabular}{|p{6.0cm}|p{7.5cm}|p{2.0cm}|}
\hline
\vspace{0.05mm}$^*$ Solution of Problem 8(a) guided by Irina Blazhievska is available on-line: &
\vspace{5.5mm} \url{https://youtu.be/w4wx7Dd547w} & \vspace{-3mm} \includegraphics[height=20mm]{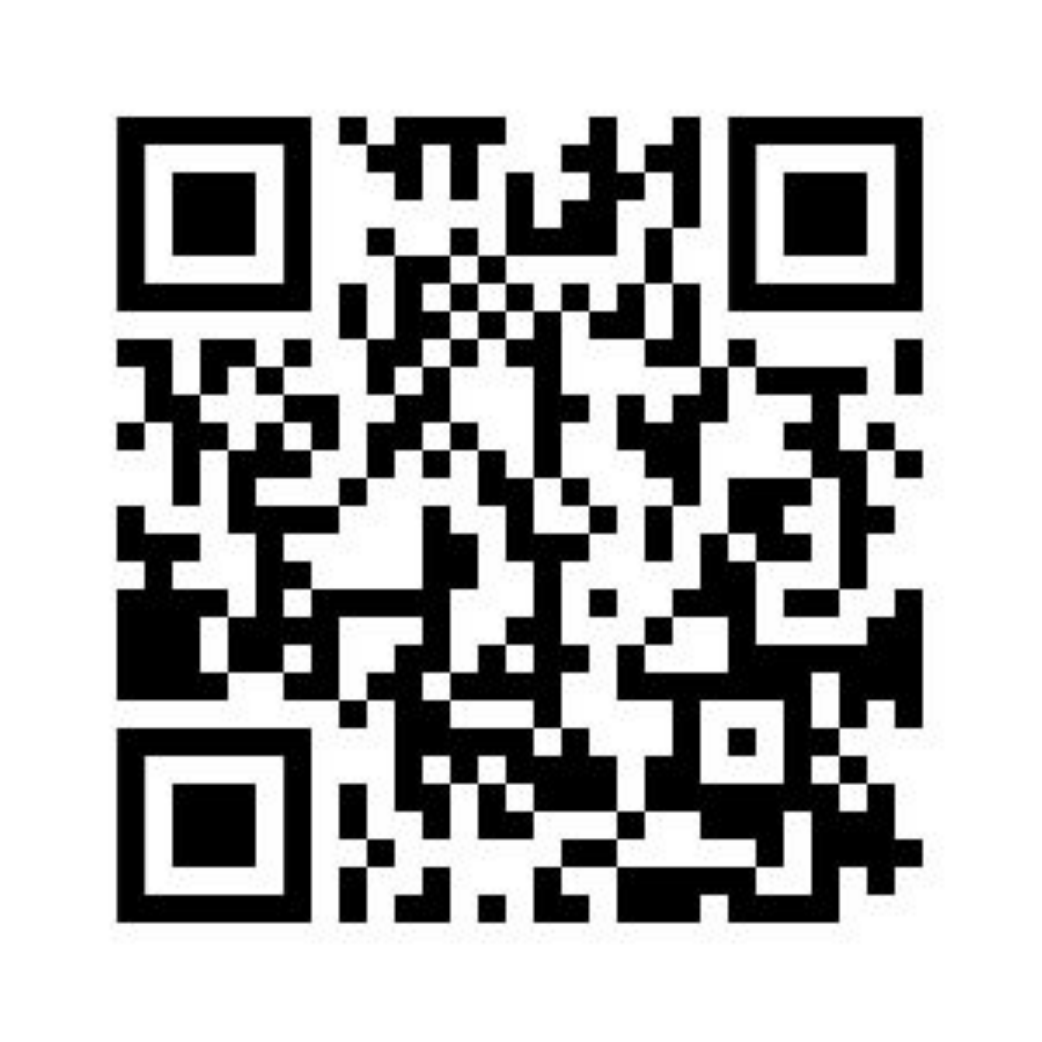}\\
\hline
\end{tabular}

The written version of solution is proposed below. Algorithm:

\textbf{Step 1.} Building the picture.

\parbox[b][7cm][t]{80mm}{\includegraphics[scale=0.8]{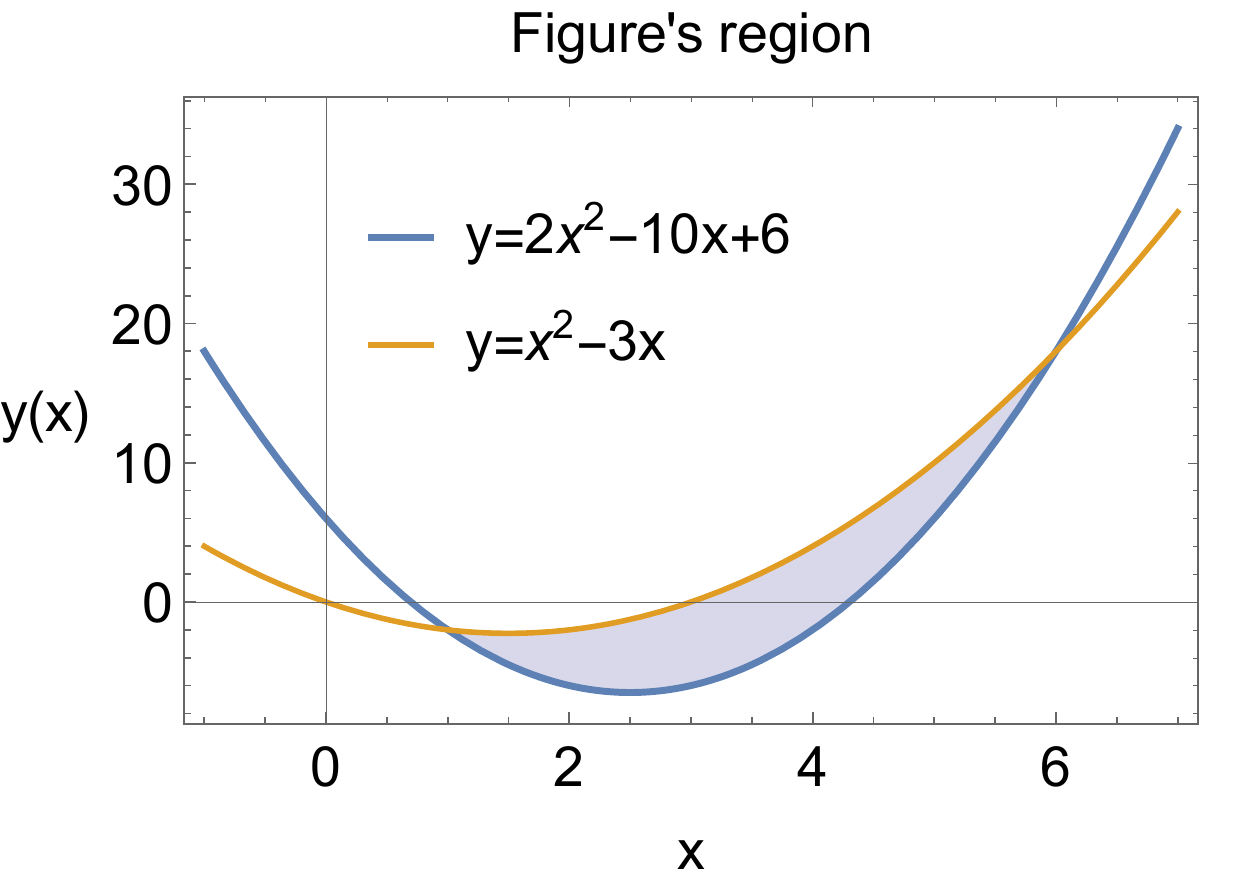}}
\hfill
\parbox[b][7cm][t]{65mm}{
\vspace{4mm}
$y=2x^2-10x+6$ is a $\cup$-shaped parabola with vertex $\left( \frac{5}{2},-\frac{13}{2}\right)$ and points of intersection with OX: $x_1=\frac{5-\sqrt{13}}{2}\approx 0,7$ and  $x_2=\frac{5+\sqrt{13}}{2}\approx 4,3.$

\vspace{4mm}
$y=x^2-3x$ is a $\cup$-shaped parabola with vertex $\left( \frac{3}{2},-\frac{9}{4}\right)$ and points of intersection with OX: $x_1=0$ and  $x_2=3.$}

\textbf{Step 2.} Finding the points of intersection between the curves.
$$\left\lbrace\begin{array}{l}
y=2x^2-10x+6 \\
y=x^2-3x
\end{array}\right.
\qquad\Rightarrow \qquad
\begin{aligned}
2x^2-10x+6= & x^2-3x;\\
x^2-7x+6= & 0; \\
(x-6)(x-1)= & 0;\\
\text{Abscises of points:} & \;
x_1=1,x_2=6.
\end{aligned}
$$

\textbf{Step 3.} Analytical description of the region in Cartesian coordinates:
$$\mathfrak{D}=\left\lbrace 1\leq x \leq 6;\; 2x^2-10x+6\leq y\leq x^2-3x\right\rbrace.$$

\textbf{Step 4.} Applying the suitable formula to find the area.
\begin{align*}
\mathcal{S}(\mathfrak{D})= & \int^b_a\big(y_2(x)-y_1(x)\big)dx= \int^6_1\big(x^2-3x-(2x^2-10x+6)\big)dx=\\
=&\int^6_1 (7x-x^2-6)dx=\left[ \frac{7}{2}x^2-\frac{1}{3}x^3-6x\right]\Bigg|_1^6= \\
=&\frac{7}{2}(36-1)-\frac{1}{3}(216-1)-6(6-1)=5\left(\frac{147-86-36}{6}\right)=\frac{125}{6}.
\end{align*}
$\textbf{Result:}\;\mathcal{S}(\mathfrak{D})=\dfrac{125}{6} \ (\text{square units}).$

\item Find the area of the figure bounded by the curves:
$$x(t)=4\cos^3t, \ y(t)=4\sin^3t,\ y=\frac{1}{2} \ \left(y\geq \frac{1}{2}\right).$$
Algorithm:

\textbf{Step 1.} Building the picture.

\parbox[b][10.5cm][c]{80mm}{\includegraphics[scale=0.8]{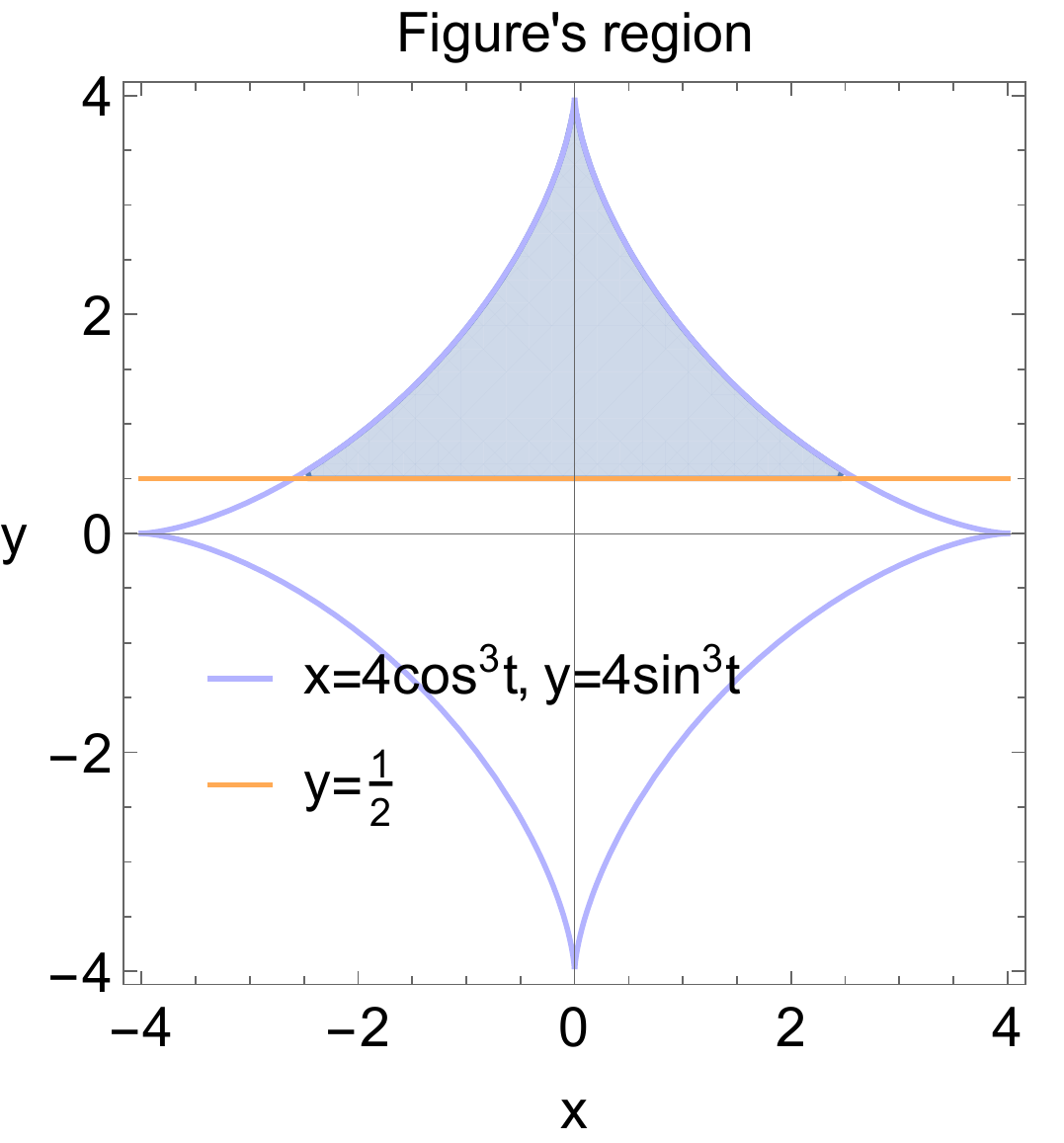}}
\hspace{0.5cm}
\parbox[b][10.5cm][t]{75mm}{
$x(t)=4\cos^3t, \ y(t)=4\sin^3t$ is a parametric representation of an astroid inscribed inside a circle of radius 4.

\vspace{4mm}
$y=\frac{1}{2}$ is the straight line passing through $\big(0,\frac{1}{2}\big)$ and parallel to OX.
\vspace{0.5cm}

\textbf{Step 2.} Finding the points of intersection between the curves.
$$
\begin{array}{c}
\left\lbrace\begin{array}{l}
y=4\sin^3 t \\
y=\dfrac{1}{2}
\end{array}\right.\Rightarrow \
\sin t= \dfrac{1}{2} \\[0.8cm]
\Rightarrow \text{Parameters of points:} \\[0.2cm]
t_1=\dfrac{\pi}{6}, \ t_2=\dfrac{5\pi}{6}.
\end{array}
$$
}

\textbf{Step 3.} Analytical description of the region in Parametric coordinates.

The OY-symmetry of region $\mathfrak{D}$ implies that its area is:
$$\mathcal{S}(\mathfrak{D})=2\mathcal{S}(\mathfrak{D}^+),$$
where $\mathfrak{D}^+$ is the subregion contained in the first quadrant,
$$\mathfrak{D}^+=\left\lbrace 0\leq x\leq 4\cos^3 t; \
\frac{1}{2}\leq y\leq 4\sin^3 t, \ \frac{\pi}{6}\leq t\leq\frac{\pi}{2}
\right\rbrace.$$

\textbf{Step 4.} Applying the suitable formula to find the area.
\begin{align*}
\mathcal{S}(\mathfrak{D})= & 2\mathcal{S}(\mathfrak{D^+})= 2\int^\beta_\alpha x(t) y'(t) dt= \left|\begin{array}{c}
x(t)=4\cos^3t,\\
            y'(t)=(4\sin^3t)'=12\sin^2t\cos t
                 \end{array}\right|=\\
=&2\cdot 4\cdot 12 \int^{\pi/2}_{\pi/6} \cos^4 t \sin^2 t dt=  96 \int^{\pi/2}_{\pi/6} \frac{1+\cos 2 t}{2}\cdot\frac{\sin^2 2t}{4} dt=\\
=&12 \int^{\pi/2}_{\pi/6}  \left(\frac{1-\cos 4t}{2}+\cos 2t \sin^2 2t \right) dt= 12\left[ \frac{t}{2}-\frac{1}{8}\sin 4t +\frac{1}{6}\sin^3 2t\right]\Bigg|^{\pi/2}_{\pi/6}= \\
= &12\left[\frac{1}{2}\left(\frac{\pi}{2}-\frac{\pi}{6}\right)-\frac{1}{8} \left( \sin 2\pi -\sin \frac{4\pi}{6}\right)+\frac{1}{6}\left(\sin^3 \pi -\sin^3 \frac{\pi}{3}\right)\right]= \\
= & 12\left[\frac{\pi}{6}+\frac{\sqrt{3}}{16}-\frac{1}{6}\left(\frac{\sqrt{3}}{2}\right)^3\right]=12\frac{\pi}{6}=2\pi.
\end{align*}
$\textbf{Result:}\;\mathcal{S}(\mathfrak{D})= 2\pi \ (\text{square units}).$

\item Find the area of the figure bounded by the curve:
$$\rho=4\sin 3\phi.$$
Algorithm:

\textbf{Step 1.} Building the picture.

$\rho=4\sin 3\phi$ is a polar representation of a 3-petaled rose inscribed inside a circle of radius 4. Since this curve is constructed from a sin-function and negative radius is not allowed, it has OY-symmetry and it is well-defined when $\frac{2\pi}{3}k\leq \phi \leq \frac{\pi}{3} +\frac{2\pi}{3}k, \ k\in \left\lbrace 0,1,2\right\rbrace.$
\vspace{4mm}

The fact that this curve has period $\frac{2\pi}{3}$ implies that the rose's k-petal is created by rotating 0-petal in counterclockwise direction around the pole an angle $\frac{2\pi}{3}k,$ $k\in \left\lbrace 1,2\right\rbrace.$ Next we show the building of the $0$-petal. For this, we set $\phi \in \left[0,\frac{\pi}{3}\right]$ and consider the additional table.
$$
\begin{array}{|c|c|c|c|c|c|c|c|c|c|}
\hline
&&&&&&&&&\\[-0.4cm]
\phi & 0 & \frac{\pi}{18}  & \frac{\pi}{12} & \frac{\pi}{9} & \frac{\pi}{6} & \frac{2\pi}{9} & \frac{\pi}{4} & \frac{5\pi}{18} & \frac{\pi}{3} \\
&&&&&&&&&\\[-0.4cm]
\hline
&&&&&&&&&\\[-0.4cm]
\rho=4\sin3\phi & 0 & 2 & 2\sqrt{2} & 2\sqrt{3} & 4 & 2\sqrt{3} & 2\sqrt{2} &  2 &  0\\ \hline
\end{array}
$$

\textbf{Step 2.} Analytical description of the region in Polar coordinates.

\parbox[b][9.5cm][t]{80mm}{\includegraphics[scale=0.8]{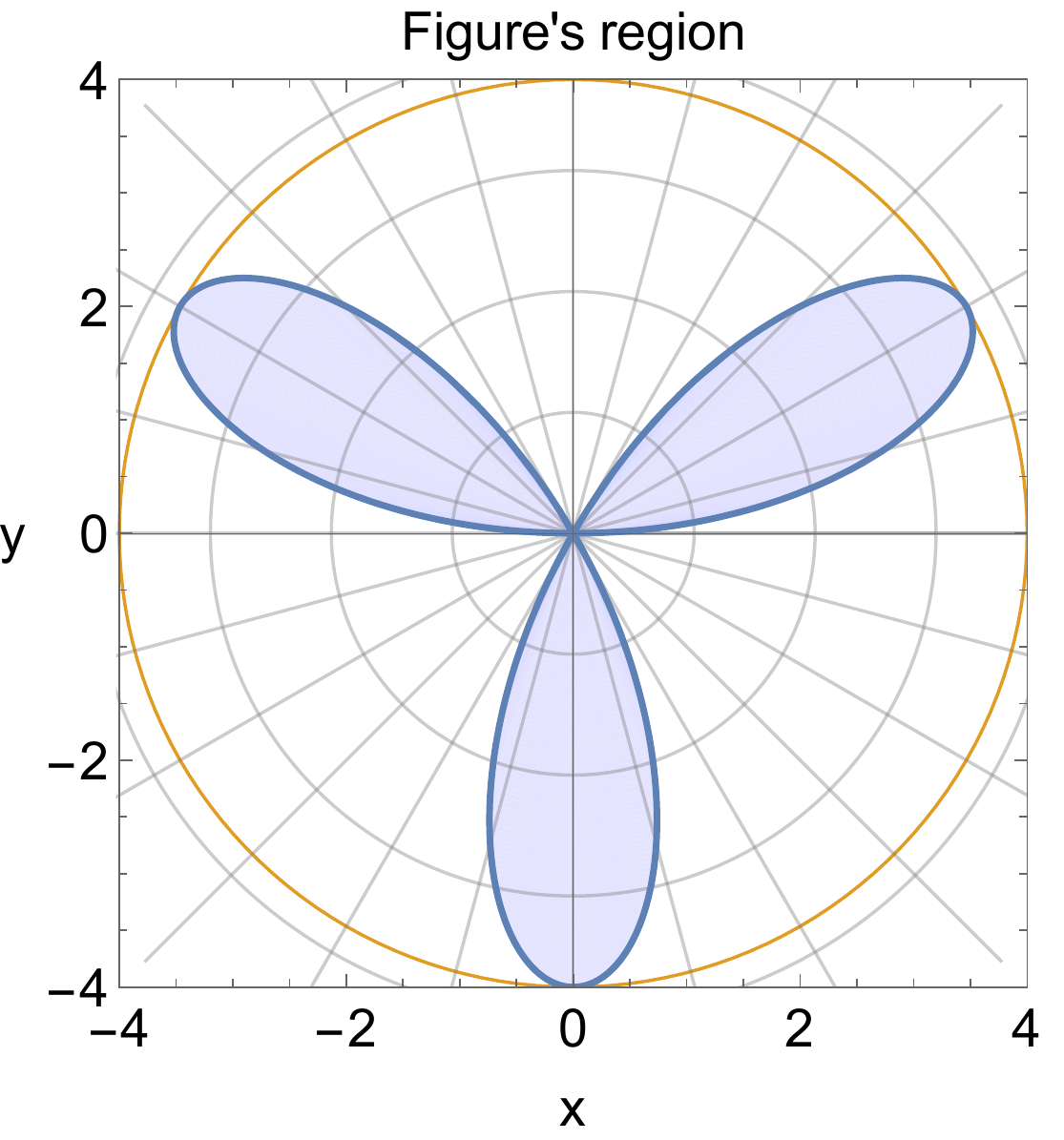}}
\hfill
\parbox[b][9.5cm][t]{75mm}{
\vspace{4mm}

Since the rose is formed by rotation of 0-petal around the pole, all petals have the same metric characteristics. \\ Then the 3 petals have the same area. The multiplication of the $0$-petal's area by $3$ give us the area of the whole region $\mathfrak{D}$:
$$\mathcal{S}(\mathfrak{D})=3\mathcal{S}(\mathfrak{D}_0),$$
where $\mathfrak{D}_0$ is the 0-petal's region,
$$\mathfrak{D}_0=\left\lbrace 0\leq \phi \leq \frac{\pi}{3}; \ 0\leq \rho \leq 4\sin 3\phi \right\rbrace$$}

\textbf{Step 3.} Applying the suitable formula to find the area.
\begin{align*}
\mathcal{S}(\mathfrak{D})= & 3\mathcal{S}(\mathfrak{D}_0)=\frac{3}{2}\int^{\beta}_{\alpha} \rho^2(\phi)d\phi= \frac{3}{2}\int^{\pi/3}_{0}(4\sin 3\phi)^2 d\phi=24\int^{\pi/3}_{0} \sin^2 3\phi d\phi=\\
= & 12\int^{\pi/3}_{0}(1-\cos 6\phi)d\phi = 12\left[\phi -\frac{1}{6}\sin 6\phi\right]\Bigg|^{\pi/3}_0=12\left(\frac{\pi}{3}-\frac{1}{6}\sin 2\pi\right)=4\pi.
\end{align*}
$\textbf{Result:}\;\mathcal{S}(\mathfrak{D})= 4\pi \ (\text{square units}).$

\end{enumerate}

\begin{enumerate}[label=\textbf{9.(\alph*)}]
\item Find the arc length of the curve
$$y=\dfrac{1}{3}\ln (\cos 3x),\  0\leq x\leq \dfrac{\pi}{18}.$$
Algorithm:

\textbf{Step 1.} Building the picture.

\parbox[b][7cm][t]{110mm}{\includegraphics[scale=0.8]{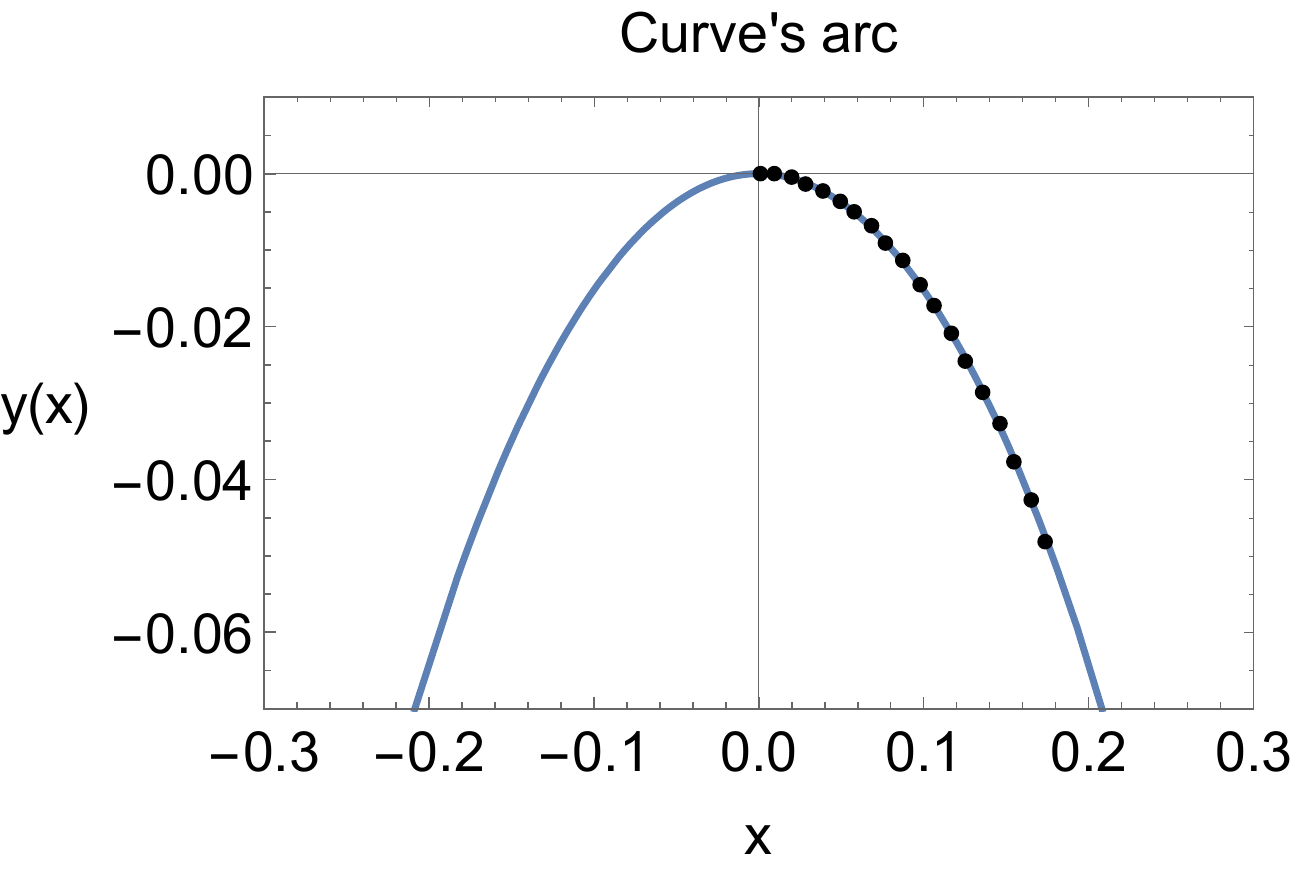}}
\parbox[b][7cm][c]{60mm}{$y=\frac{1}{3}\ln (\cos 3x)$ is a curve constructed from a cos-function. Since the logarithm is only defined for positive values, this curve is only defined on regions $\left[-\frac{\pi}{6}+\frac{2\pi k}{3}, \frac{\pi}{6}+\frac{2\pi k}{3}\right], \ k\in \mathbb{Z}.$

On each region it has a $\bigcap$-like form with zeros-maxima at $x=\frac{2\pi}{3}k,$ and vertical asymptotes: \\ $x=\pm\frac{\pi}{6}+\frac{2\pi k}{3}, k\in \mathbb{Z}.$}

\textbf{Step 2.} Analytical description of the arc segment in Cartesian coordinates:
    $$\Gamma=\left\lbrace y= \dfrac{1}{3}\ln(\cos 3x); \ 0\leq x\leq \dfrac{\pi}{18}\right\rbrace$$
Note that $\left[ 0,\frac{\pi}{18}\right]$ belongs to the region of curve's well-definiteness.

\textbf{Step 3.} Applying of the suitable formula to find the arc length.
    \begin{align*}
        l(\Gamma)= &\int^b_a \sqrt{1+(y_x')^2}dx=
        \left|\begin{array}{c}
            y= \dfrac{1}{3}\ln(\cos 3x); \  y_x'=\dfrac{1}{3}\cdot\dfrac{-3\sin 3x}{\cos 3x}=-\tg 3x; \\[0.5cm]
            1+(y'_x)^2=1+(-\tg 3x)^2=\dfrac{1}{\cos^2 3x}
        \end{array}\right|= \\
        = & \int^{\pi/18}_0\sqrt{1+\tg^2 3x}dx= \int^{\pi/18}_0 \frac{1}{|\cos 3x|}dx=\int^{\pi/18}_0 \frac{1}{\cos 3x}dx= \\
        = & \int^{\pi/18}_0 \frac{\cos 3x}{\cos^2 3x}dx=-\frac{1}{3}\int^{\pi/18}_0 \frac{d(\sin 3x)}{\sin^2 3x-1}dx= -\frac{1}{3}\frac{1}{2}\ln \left|\frac{\sin 3x-1}{\sin 3x+1}\right|\Bigg|^{\pi/18}_0= \\
        = & -\frac{1}{6}\ln \left|\frac{\sin\left( 3\cdot \dfrac{\pi}{18}\right)-1}{\sin \left(3\cdot \dfrac{\pi}{18}\right)+1}\right|+\frac{1}{6}\ln\left|\frac{0-1}{0+1}\right|=-\frac{1}{6}\ln\left|\frac{\dfrac{1}{2}-1}{\dfrac{1}{2}+1}\right|=\frac{1}{6}\ln 3.
    \end{align*}
    \textbf{Result:} $l(\Gamma)=\dfrac{1}{6}\ln 3 \ (\text{units}).$

\newpage

\item Find the arc length of the curve
$$x(t)=4(t-\sin t), \ y(t)=4(1-\cos t), \ t\in [0,2\pi].$$

Algorithm:

\textbf{Step 1.} Building the picture.

\parbox[b][5cm][t]{100mm}{
\includegraphics[scale=0.8]{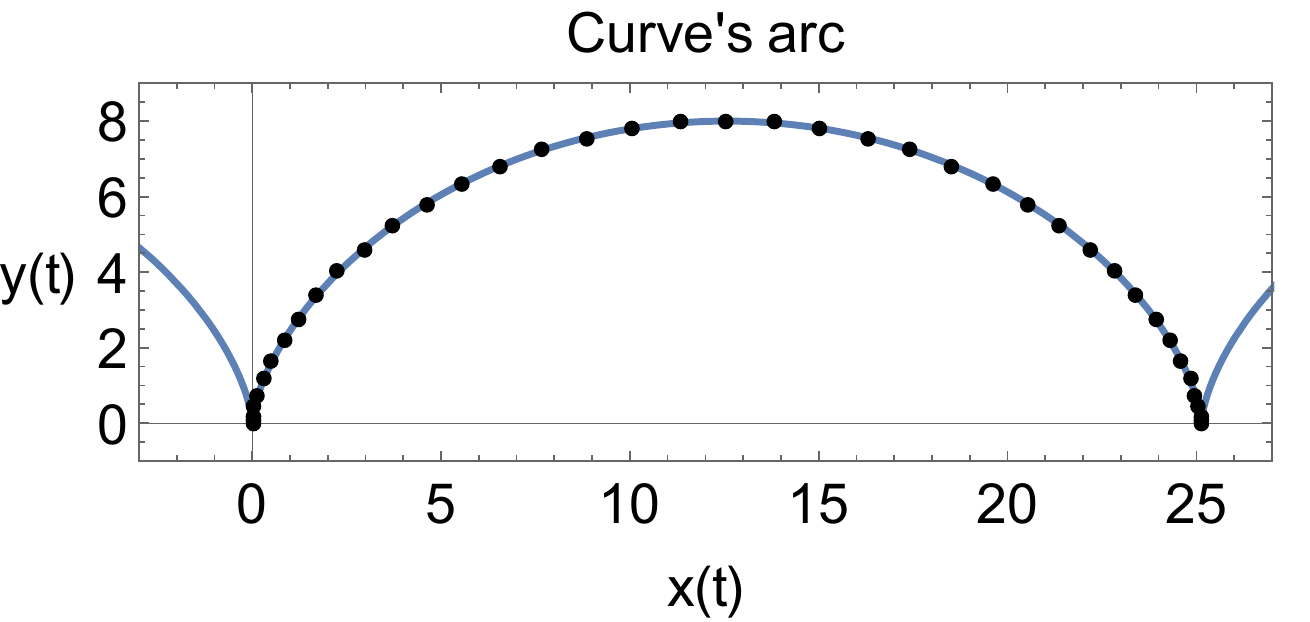}}
\hspace{0.5cm}
\parbox[b][5cm][t]{60mm}{
$$\left\lbrace\begin{array}{l}
x(t)=4(t-\sin t)\\ y(t)=4(1-\cos t)
\end{array}\right.$$
is a parametric representation of a \textit{cycloid} with a horizontal base $y=0.$}

This curve is generated by the trajectory of a marked point in a circle of radius 4 rolling on the "positive" side of the base.
The parameter region $0\leq t\leq2\pi$ creates one full wave of the curve, which intersects OX at points $x=0$ and $x=8\pi$, has maxima at $(4\pi,8)$ and an axis of symmetry at $x=4\pi$.

\textbf{Step 2.} Analytical description of the arc segment in Parametric coordinates:
$$\Gamma=\left\lbrace x=4(t-\sin t), \ y=4(1-\cos t); \ 0\leq t\leq 2\pi \right\rbrace.$$

\textbf{Step 3.} Applying the suitable formula to find the arc length.
\begin{align*}
    l(\Gamma)= & \int^\beta_\alpha \sqrt{(x'_t)^2+(y'_t)^2}dt=
\left|\begin{array}{c}
    x'_t=\Big(4(t-\sin t)\Big)'=4(1-\cos t); \\
    y'_t=\Big(4(1-\cos t)\Big)'=4\sin t; \\
    (x'_t)^2+(y'_t)^2=4^2\big((1-\cos t)^2+\sin^2t\big)
\end{array}\right| = \\
= & 4 \int^{2\pi}_0 \sqrt{(1-\cos t)^2+\sin^2 t}dt= 4 \int^{2\pi}_0 \sqrt{2(1-\cos t)}dt=8\int^{2\pi}_0 \left|\sin\frac{t}{2}\right|dt = \\
= & 8\int^{2\pi}_0 \sin\frac{t}{2}dt =-16\cos \frac{t}{2}\Bigg|^{2\pi}_0=-16(\cos\pi-\cos 0)=32.
\end{align*}
\textbf{Result:} $l(\Gamma)=32 \ (\text{units}).$

\item Find the arc length of the curve $$\rho=\frac{10}{\sqrt{101}}e^{\frac{\phi}{10}}, \ 0\leq \phi \leq 2\pi.$$

\begin{tabular}{|p{6.0cm}|p{7.5cm}|p{2.0cm}|}
\hline
\vspace{0.05mm}$^*$ Solution of Problem 9(c) guided by Irina Blazhievska is available on-line: &
\vspace{5.5mm} \url{https://youtu.be/AJbHsr1eOes} & \vspace{-3mm} \includegraphics[height=20mm]{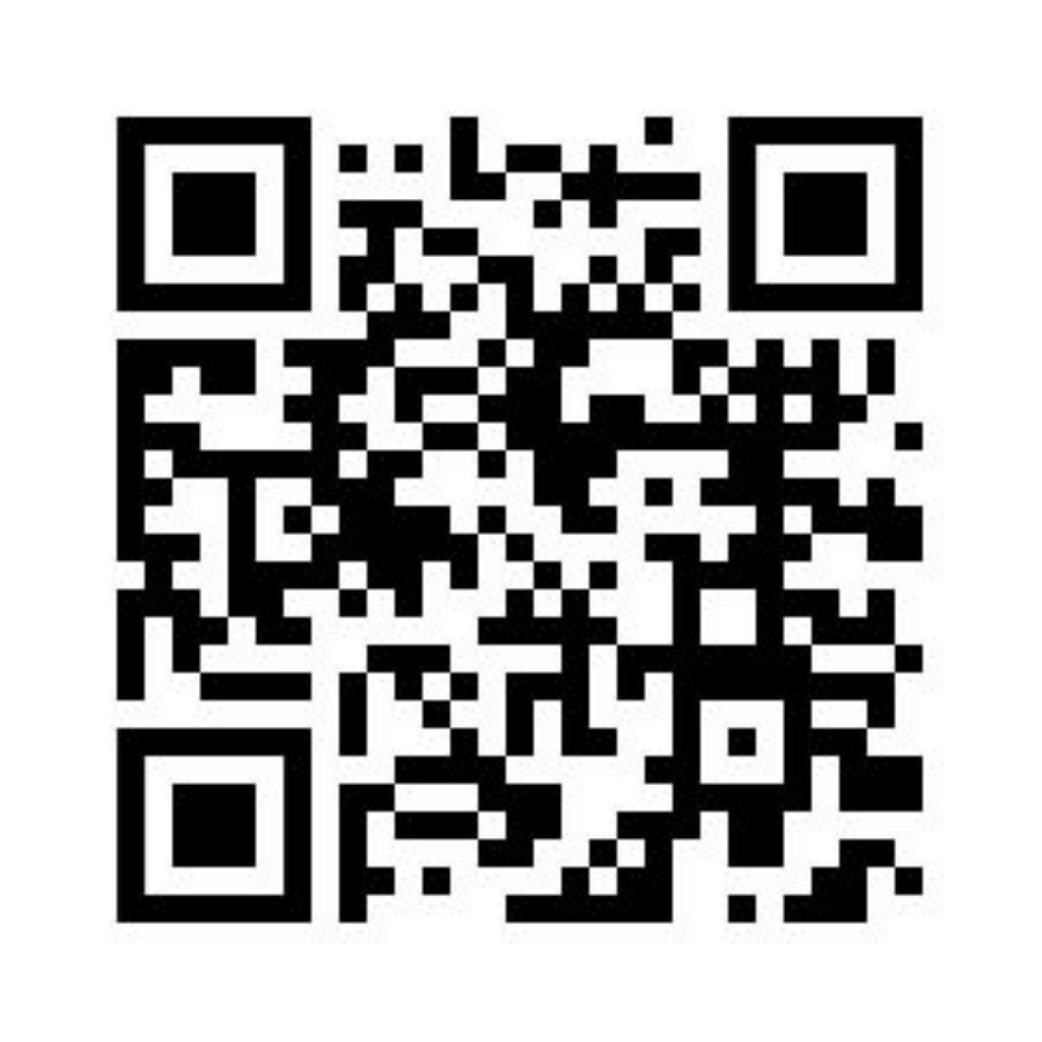}\\
\hline
\end{tabular}

The written version of solution is proposed below. Algorithm:

\textbf{Step 1.} Building the picture.

\parbox[b][9.5cm][t]{85mm}{\includegraphics[scale=0.8]{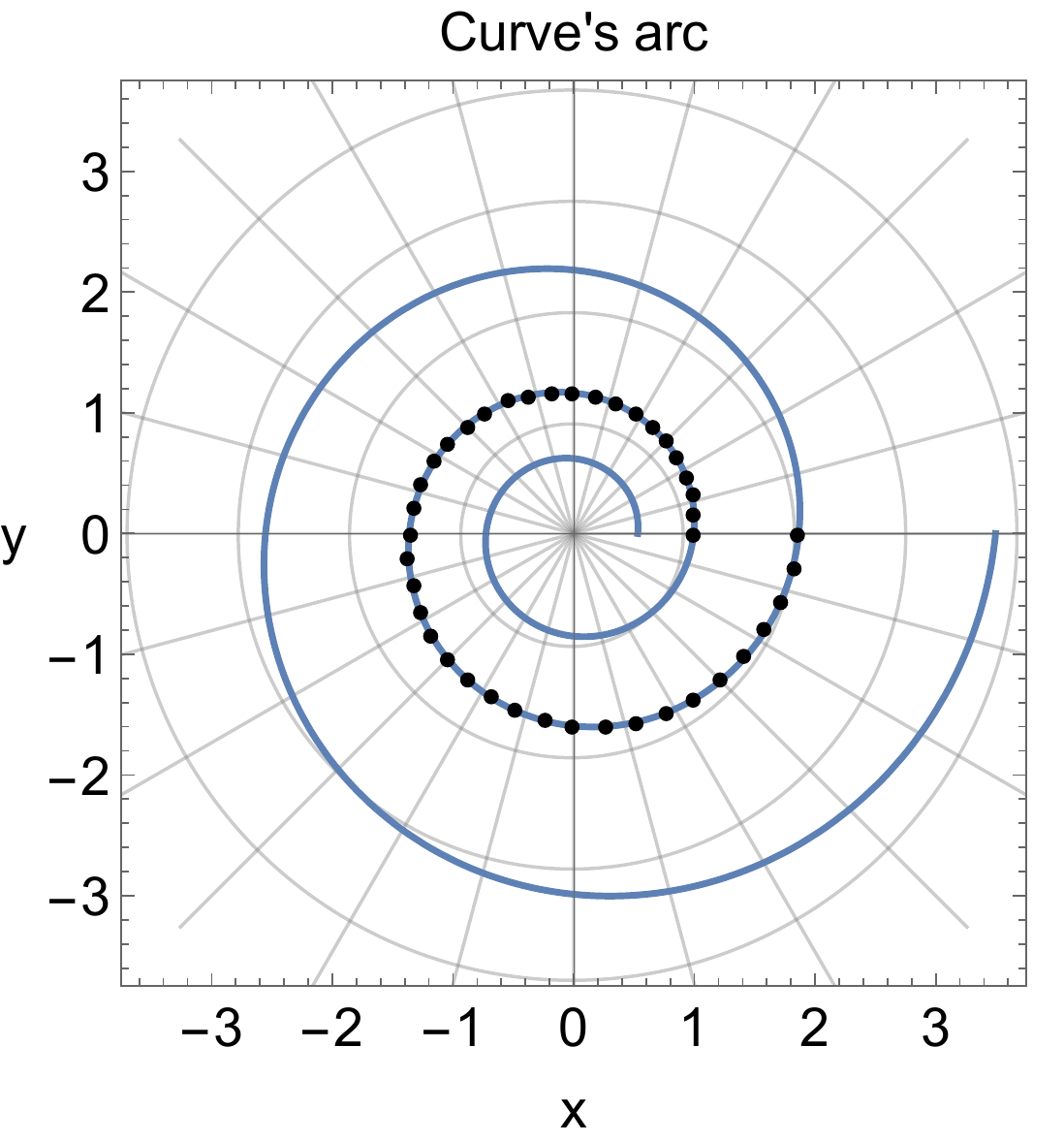}}
\hspace{0.5cm}
\parbox[b][9.5cm][t]{70mm}{
$\rho=\frac{10}{\sqrt{101}}e^{\frac{\phi}{10}}$ is a polar representation of a logarithmic spiral. It has 2 parameters; the initial radius $\frac{10}{\sqrt{101}}$ and the rate of spiral's increasing $k=\frac{1}{10}$ \ ($k=\frac{\rho'_\phi}{\rho}$).
\vspace{0.4cm}

Since $k>0$, the spiral rotates around the pole in counterclockwise direction with increasing polar radius.
\vspace{0.4cm}

The angle domain $0\leq\phi\leq 2\pi$ generates one full turn with continuous increase of radius from $\frac{10}{\sqrt{101}}$ to $\frac{10}{\sqrt{101}}e^{\frac{\pi}{5}}$.}

\textbf{Step 2.} Analytical description of the arc segment in Polar coordinates:
$$\Gamma=\left\lbrace \rho =\frac{10}{\sqrt{101}}e^{\frac{\phi}{10}};
                              \ 0\leq \phi \leq 2\pi \right\rbrace.$$
\vspace{-0.4cm}
\textbf{Step 3.} Applying the suitable formula to find the arc length.
\begin{align*}
    l(\Gamma)= & \int^{\beta}_{\alpha}\sqrt{\rho^2 +(\rho_\phi')^2}d\phi=
        \left|\begin{array}{c}
            \rho=\frac{10}{\sqrt{101}}e^{\frac{\phi}{10}};  \rho_\phi'=\frac{10}{\sqrt{101}}\cdot\frac{1}{10}e^{\frac{\phi}{10}}=\frac{1}{\sqrt{101}}e^{\frac{\phi}{10}}; \\
            \rho^2+(\rho'_\phi)^2=\frac{100}{101}e^{\frac{2\phi}{10}}+\frac{1}{101}e^{\frac{2\phi}{10}}=\big(e^{\frac{\phi}{10}}\big)^2
        \end{array}\right|=\\
       &= \int^{2\pi}_{0}\sqrt{\big(e^{\frac{\phi}{10}}\big)^2} d\phi= \int^{2\pi}_0 e^{\frac{\phi}{10}}d\phi=10e^{\frac{\phi}{10}} \Big|^{2\pi}_0=10(e^{\frac{2\pi}{10}}-e^0)=10(e^{\frac{\pi}{5}}-1).
\end{align*}
\vspace{-0.4cm}
\textbf{Result:} $l(\Gamma)= 10\big(e^{\frac{\pi}{5}}-1\big)$ \ (units).
\end{enumerate}

\begin{enumerate}[label=\textbf{10.(\alph*)}]
\item Find the surface area generated by rotating the curve
$$y=\dfrac{1}{2}\ch 2x, \ -1\leq x\leq 1, \ \text{around the axis}\  l=OX.$$
Algorithm:

\textbf{Step 1.} Building the picture.

\parbox[b][7cm][t]{90mm}{
\includegraphics[scale=0.75]{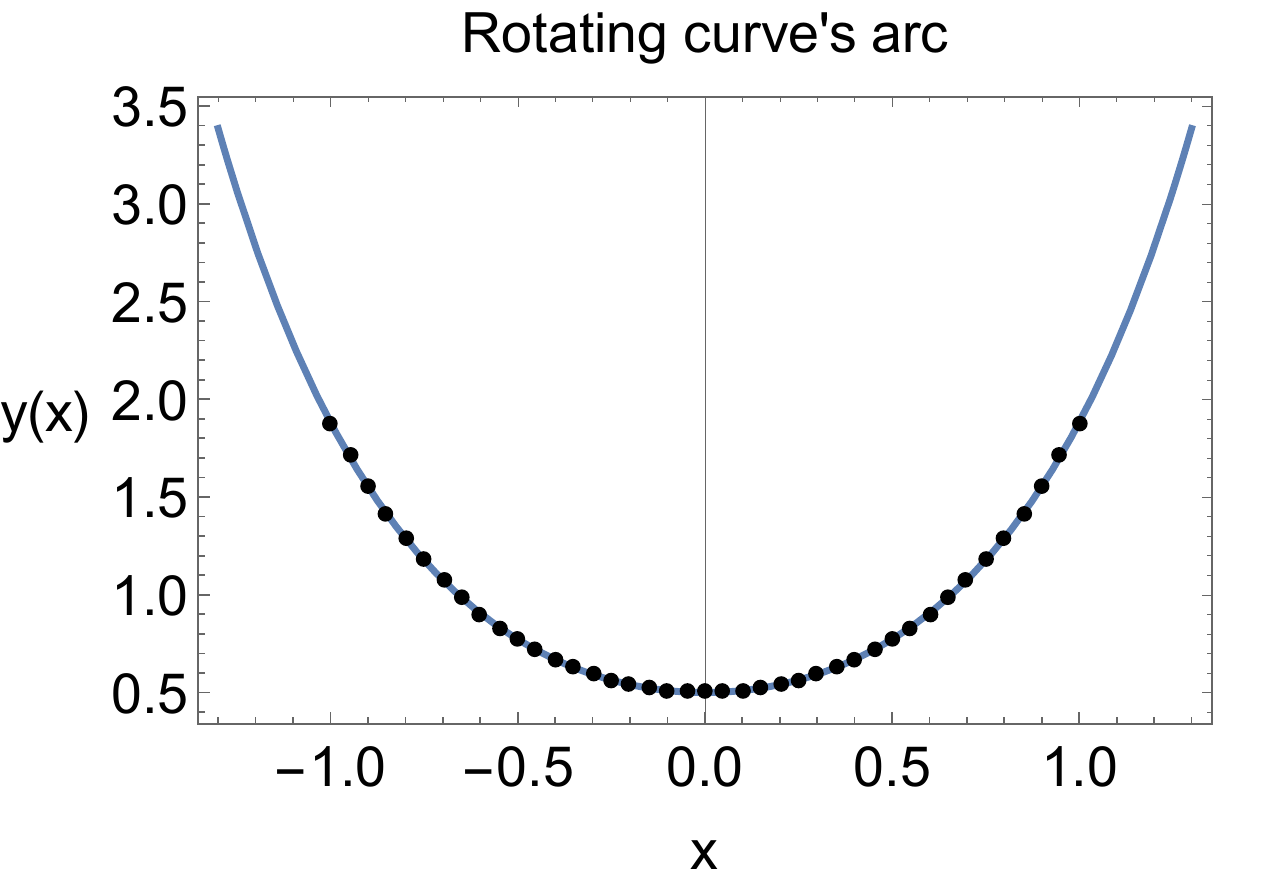}}
\hspace{0.5cm}
\parbox[b][7cm][t]{70mm}{\vspace{-0.5cm}
$y=\frac{1}{2}\ch 2x$ is a catenary with parameter $a=\frac{1}{2},$ minima $(0,\frac{1}{2}).$

This curve is an idealized hanging chain or cable sags under its own weight when it is supported only at its ends.}

\parbox[b][6cm][t]{90mm}{
\vspace{-0.3cm}
The surface of revolution of the catenary curve around OX is a \textit{catenoid.} Note that the catenoid has a minimal surface area.
\vspace{0.2cm}

\textbf{Step 2.} Analytical description of the rotating arc segment (the generatrix) in Cartesian coordinates:
\vspace{-0.2cm}
$$\Gamma=\left\lbrace y=\dfrac{1}{2}\ch 2x ; \ -1\leq x\leq 1 \right\rbrace$$
\textbf{Step 3.} Applying the suitable formula to find the surface area.
}
\hspace{2.5cm}
\parbox[b][4cm][t]{70mm}{
\vspace{-6cm}
\hspace{-2cm}
\includegraphics[scale=0.75]{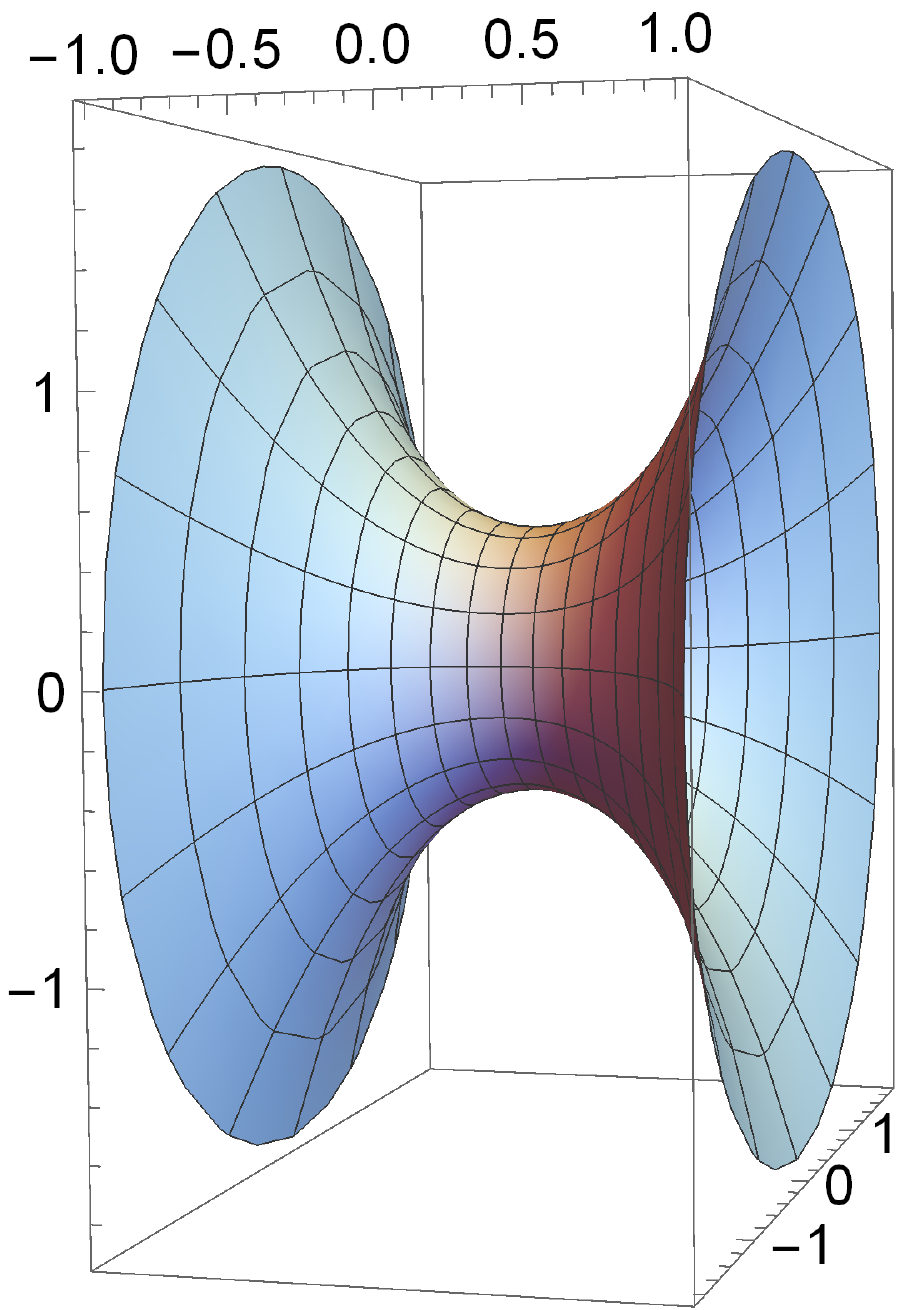}}
\begin{align*}
S_{OX}= &  \ 2\pi \int^b_a y(x)\sqrt{1+(y'_x)^2}dx = \left|\begin{array}{c}
            y= \dfrac{1}{2}\ch 2x; \  y_x'=\dfrac{1}{2}\cdot2\sh 2x=\sh 2x \\
            1+(y'_x)^2=1+(\sh 2x)^2=\ch^2 2x
        \end{array}\right|= \\
= & \ 2\pi \int^1_{-1}\frac{1}{2}\ch 2x \sqrt{1+\left(\sh 2x\right)^2}dx= \pi \int^1_{-1}\ch^2 2x \ dx= \pi \int^1_{-1} \dfrac{1+\ch 4x}{2} dx = \\
= & \frac{\pi}{2}\left[ x+\frac{1}{4}\sh 4x\right]\Bigg|^1_{-1}= \frac{\pi}{2}\left(2+\frac{1}{2}\sh 4\right)=\pi\Big(1+\dfrac{1}{4}\sh 4\Big).
\end{align*}
\vspace{-0.5cm}
\textbf{Result:} $S_{OX}= \pi\Big(1+\dfrac{\pi}{4}\sh 4\Big)$ \ (square units).

\item Find the surface area generated by rotating the curve
$$x=3+\cos t, \ y=2+\sin t, \ \text{around the axis} \ l=OY.$$

\begin{tabular}{|p{6.0cm}|p{7.5cm}|p{2.0cm}|}
\hline
\vspace{0.05mm}$^*$ Solution of Problem 10(b) guided by Irina Blazhievska is available on-line: &
\vspace{5.5mm} \url{https://youtu.be/nXzXyHIw0w8} & \vspace{-3mm} \includegraphics[height=20mm]{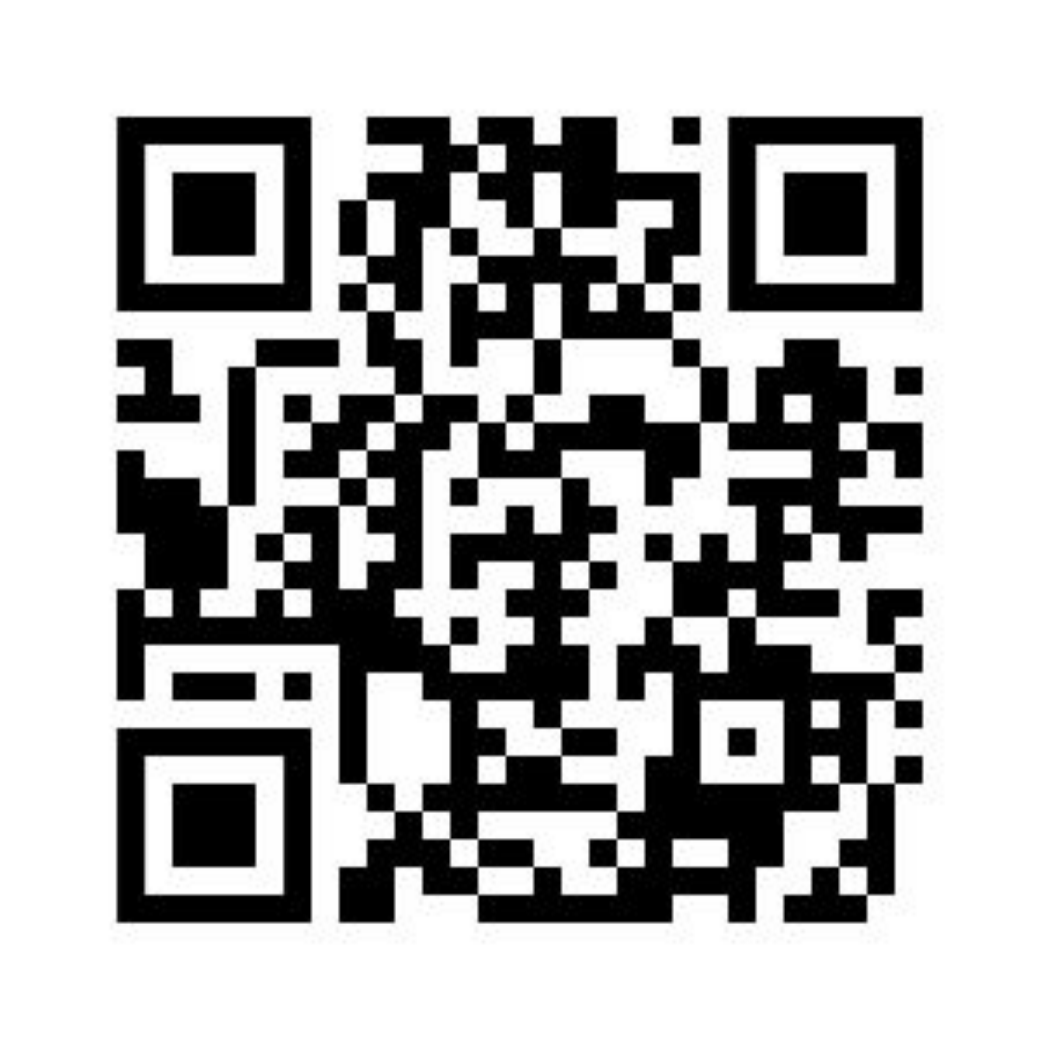}\\
\hline
\end{tabular}

The written version of solution is proposed below. Algorithm:

\textbf{Step 1.} Building the picture.

\parbox[b][6.5cm][t]{80mm}{
\includegraphics[scale=0.7]{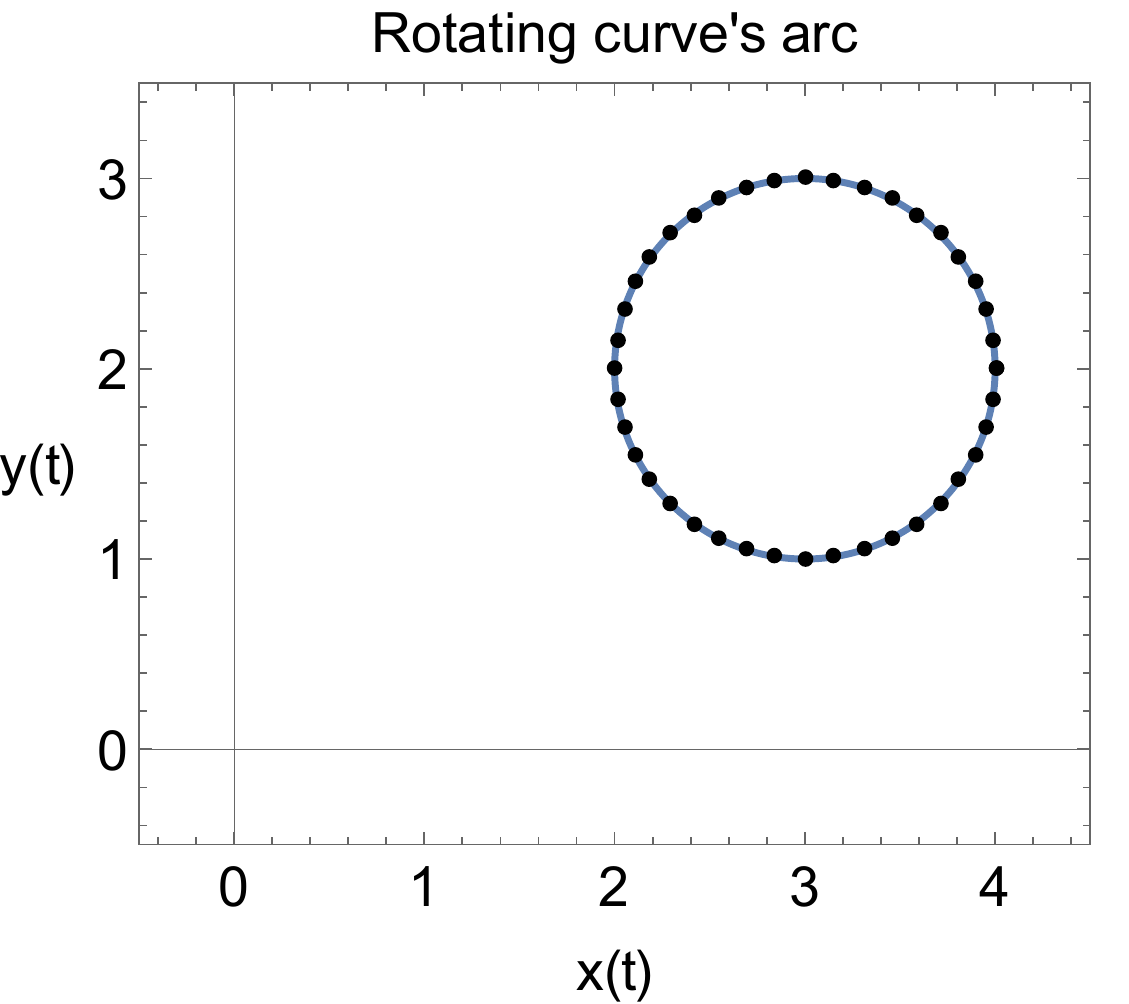}}
\hspace{0.5cm}
\parbox[b][6.5cm][t]{70mm}{
$$\left\lbrace\begin{array}{l}
x(t)=3+\cos t \\
y(t)=2+\sin t
\end{array}\right.$$
is a parametric representation of a circle of radius 1 with center in $(3,2)$: $(x-3)^2+(y-2)^2=1^2.$

The empty intersection of the curve with OY implies that there are no restrictions for parameter's domain: $0\leq t\leq 2\pi.$
The surface of revolution with this rotating curve is a ring-torus.
}

\parbox[b][6cm][b]{60mm}{
\textbf{Step 2.} Analytical description of the rotating arc segment (the generatrix) in Parametric coordinates:
$$\Gamma=\left\lbrace\begin{array}{c} x=3+\cos t, \\ y= 2+\sin t; \\ 0 \leq t \leq 2 \pi \end{array}\right\rbrace$$}
\hspace{0.5cm}
\parbox[b][6cm][b]{80mm}{
\includegraphics[scale=0.75]{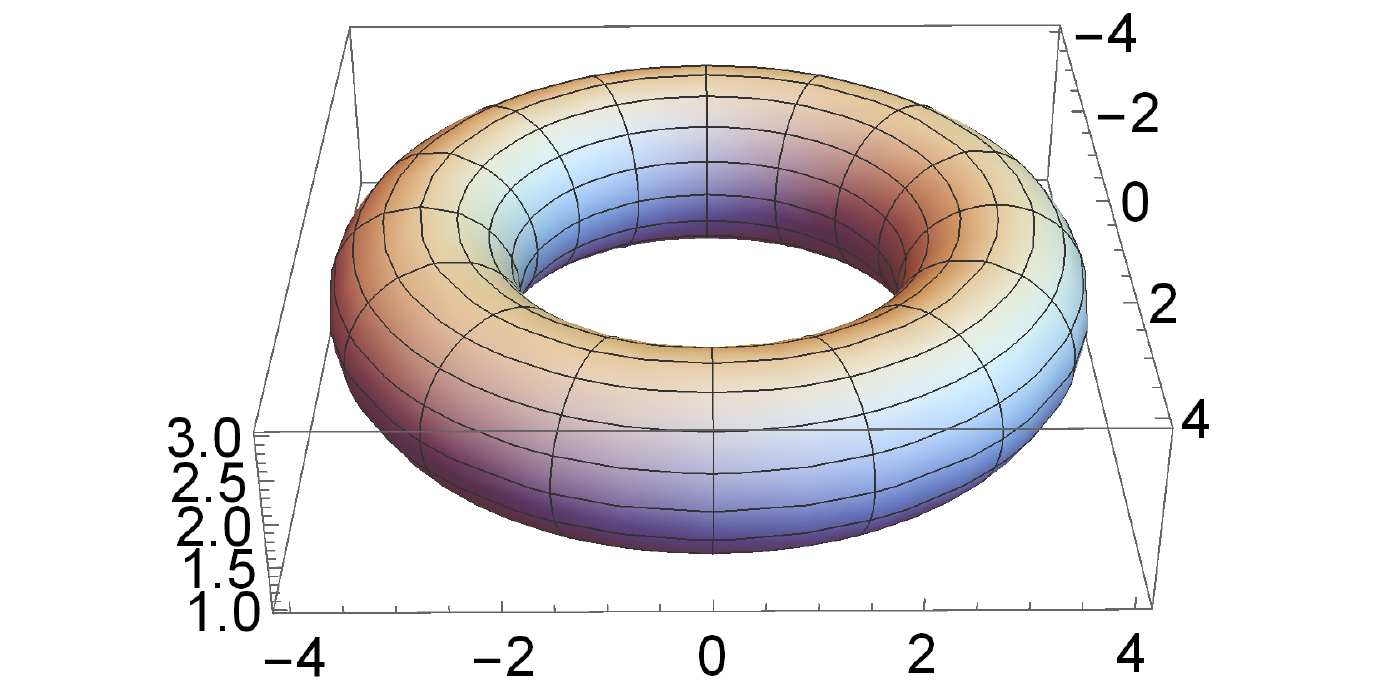}}

\textbf{Step 3.} Applying the suitable formula to find the surface area.
\begin{align*}
S_{OY}= & 2\pi \int^\beta_\alpha x(t)\sqrt{(x'_t)^2+(y'_t)^2} \ dt= \left|\begin{array}{c}
            x(t)=3+\cos t; \ \ y(t)=2+\sin t; \\
            x'_t=(3+\cos t)'=-\sin t; \\
            y'_t=(2+\sin t)'=\cos t; \\
            (x'_t)^2+(y'_t)^2=(-\sin t)^2+(\cos t)^2=1
        \end{array}\right|=
\end{align*}
\begin{align*}
= & 2\pi \int ^{2\pi}_0 (3+\cos t)\sqrt{(-\sin t)^2+(\cos t)^2} \ dt= 2\pi \int^{2\pi}_0 (3+\cos t)\ dt= \\
= & 2\pi \Big[ 3t +\sin t\Big]\Bigg|^{2\pi}_0= 2\pi \left[3(2\pi -0)+(\sin 2\pi -\sin 0)\right]=12\pi^2.
\end{align*}
\textbf{Result:} $S_{OY}= 12\pi^2$  (square units).

\item Find the area surface area generated by rotating the curve
$$\rho=\sqrt{\cos 2\phi} \ \text{around the axis} \ l=o\rho $$
Algorithm:

\textbf{Step 1.} Building the picture.

\parbox[b][70mm][t]{80mm}{
\includegraphics[scale=0.7]{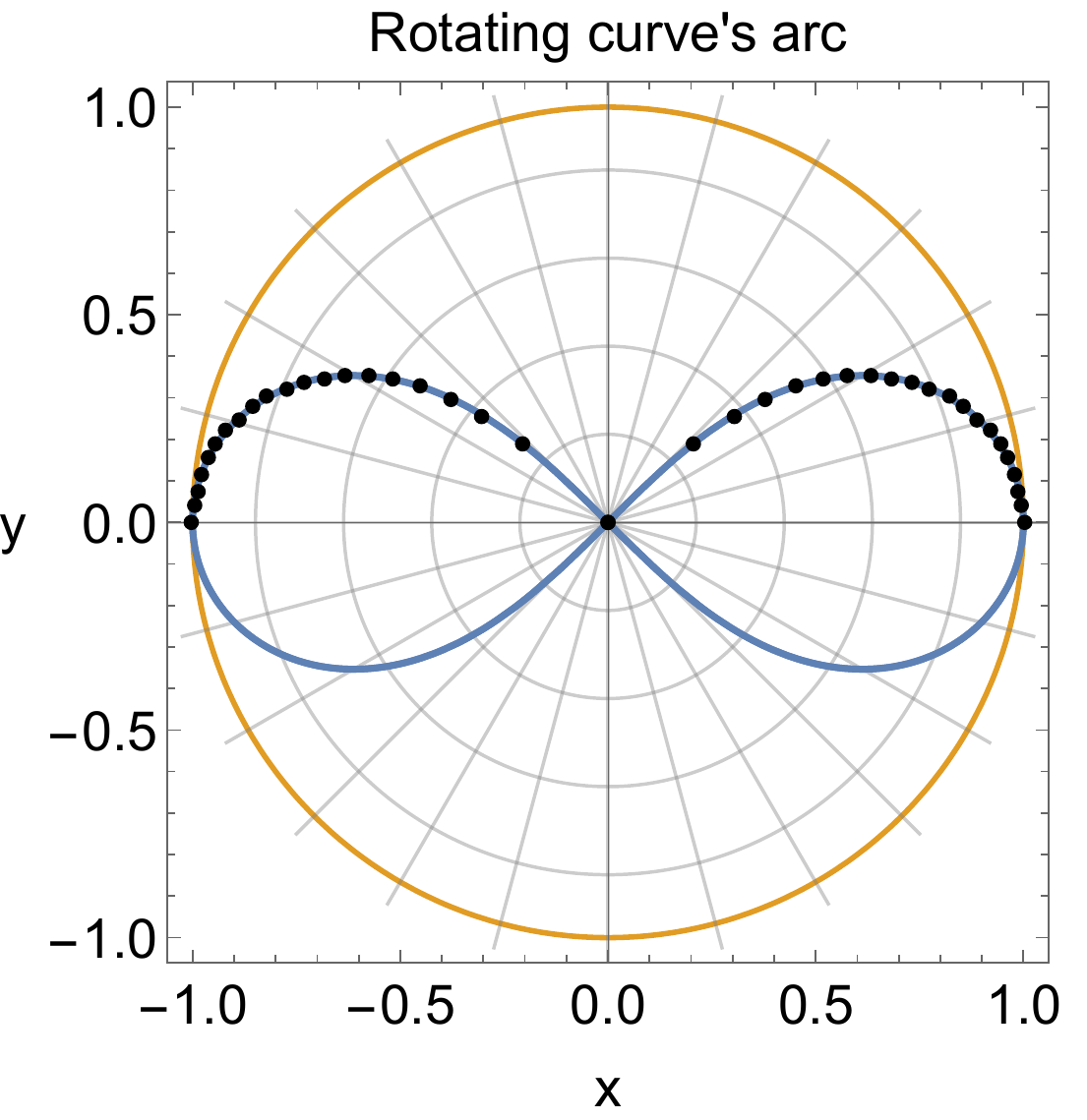}}
\parbox[b][70mm][t]{80mm}{
$\rho=\sqrt{\cos 2\phi}$ is a Polar representation of an $\infty$-shaped lemniscate inscribed inside a circle of radius 1. This curve is known as \textit{Bernoulli's lemniscate.} Since this curve is constructed from a cos-function and negative radius is not allowed, it has OX-symmetry and it is well-defined on $\phi \in [-\frac{\pi}{4},\frac{\pi}{4}]\bigcup[\frac{3\pi}{4},\frac{5\pi}{4}]$.
\vspace{0.5cm}

OX-symmetry of the curve implies the restriction on polar angle domain for rotating part: $\phi \in [0,\frac{\pi}{4}]\cup [\frac{3\pi}{4};\pi].$}

\parbox[b][55mm][c]{60mm}{\vspace{1.5cm}
The surface of revolution generated by the lemniscate is "Hourglass"-shaped and is attached on the right.}
\hspace{1cm}
\parbox[b][55mm][t]{80mm}{
\vspace{1cm}
\includegraphics[scale=0.7]{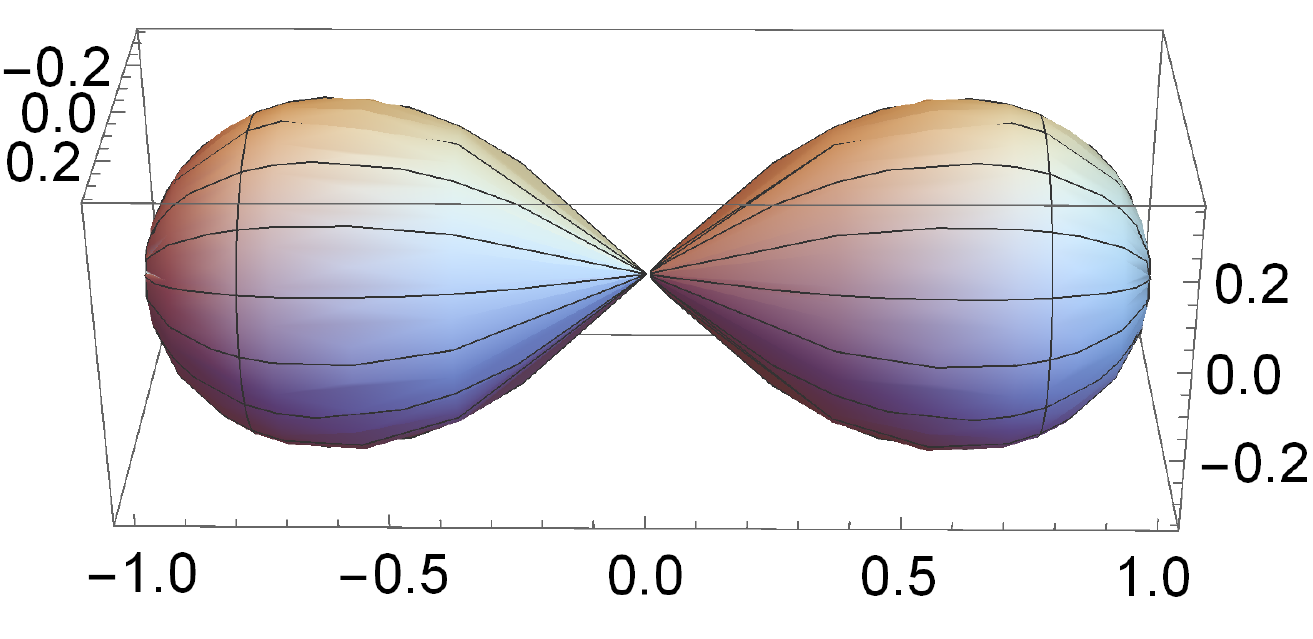}
}

\textbf{Step 2.}
Analytical description of the rotating arc segment (the generatrix) in Polar coordinates.

The mirror symmetry of the surface with respect to YZ-plane implies that both Hourglass' sides have the same metric characteristics.

The multiplication of right-side surface area $S^+_{o\rho}$ by 2 gives us the the surface area of the whole solid of revolution:
 $
 S_{o\rho}=2S^+_{o\rho},
 $
where the generatrix of $S^{+}_{o\rho}$ is located in the first quadrant:

$$\Gamma^+=\left\lbrace \rho=\sqrt{\cos 2\phi}; \ 0\leq \phi \leq \frac{\pi}{4}\right\rbrace.$$

\textbf{Step 3.} Applying the suitable formula to find the surface area.
\begin{align*}
S_{o\rho}= & 2S^+_{o\rho}=2\cdot 2\pi \int^\beta_\alpha \rho(\phi) \sin \phi \sqrt{(\rho)^2+(\rho'_\phi)^2} \ d\phi=\left|\begin{array}{c}
            \rho=\sqrt{\cos2\phi}; \\
            \rho'_\phi=-\dfrac{\sin2\phi}{\sqrt{\cos2\phi}};\\
            (\rho)^2+(\rho'_\phi)^2=\dfrac{1}{\cos2\phi}
        \end{array}\right|=\\
        =& 4\pi \int^{\pi/4}_0 \sqrt{\cos 2\phi} \sin \phi \sqrt{\dfrac{1}{\cos2\phi}}\ d\phi= 4\pi \int^{\pi/4}_0\sin \phi d\phi=-4\pi \cos \phi \Bigg|^{\pi/4}_0 =\\
        =&-4\pi\Big(\cos\bigg(\frac{\pi}{4}\bigg)-\cos0\Big)=-4\pi\Big(\dfrac{\sqrt{2}}{2}-1\Big)=2\pi(2-\sqrt{2}).
\end{align*}
\textbf{Result: } $S_{o\rho}=2\pi\big(2-\sqrt{2}\big)$ \ (square units).
\end{enumerate}

\begin{enumerate}[label=\textbf{11.(\alph*)}]
\item Find the volume of the body formed by rotating region between the curves
$$y=x^2+2x+5, \ y=5-x, \ \text{around the axis} \ l=OX.$$
Algorithm:

\textbf{Step 1.} Building the picture.

\parbox[b][60mm][t]{105mm}{
\includegraphics[scale=0.8]{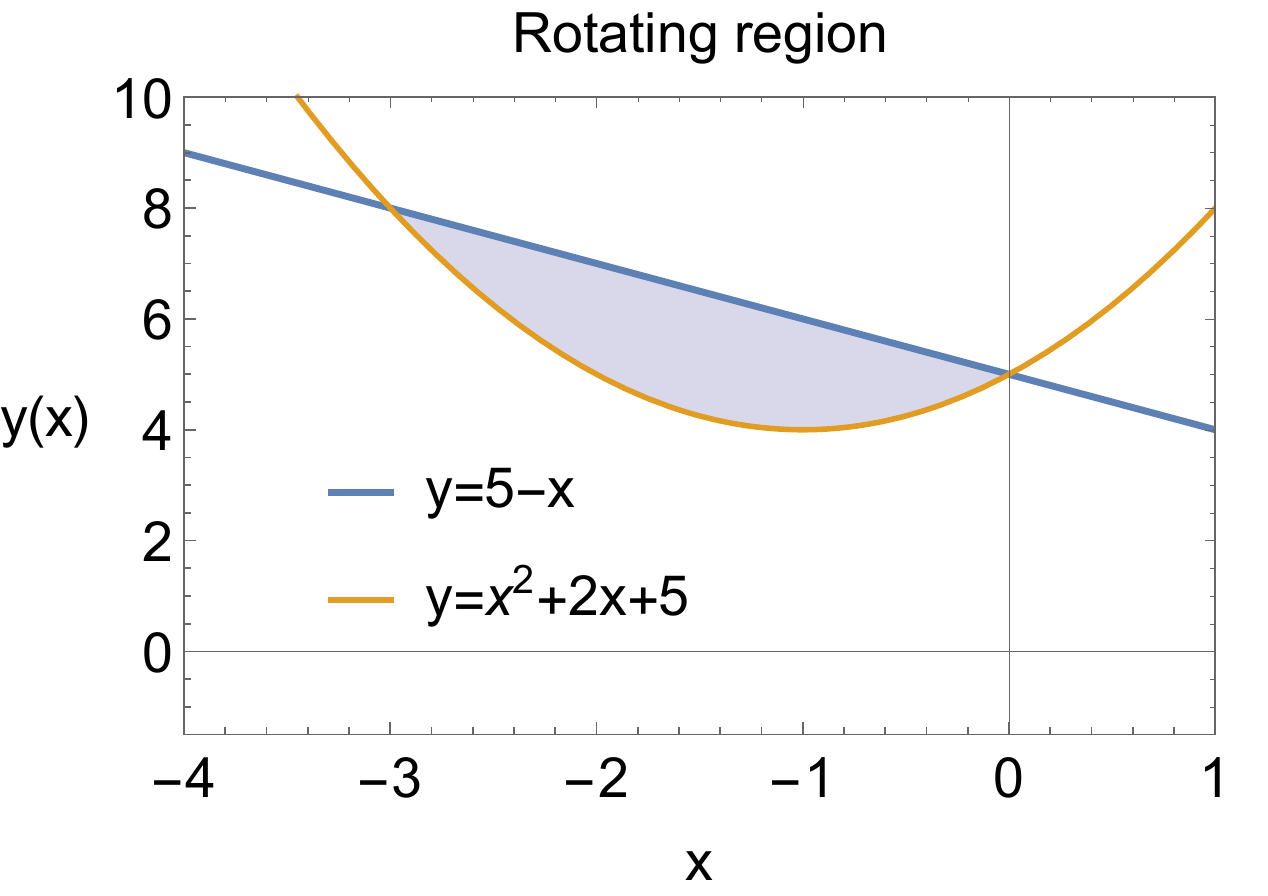}}
\vspace{0.5cm}
\parbox[b][60mm][c]{58mm}{
$y=x^2+2x+5$ is a $\cup$-shaped parabola with vertex $(-1,3)$; it has empty intersection with OX, but intersects $OY$ at $(0,5).$

\vspace{4mm}
$y=5-x$ is a straight line passing through the points $(5,0),
(0,5).$
}

\textbf{Step 2.} Finding the points of intersection between the curves.
\[\left\lbrace\begin{array}{c}
y=x^2+2x+5 \\
y=5-x
\end{array}\right. \quad\Rightarrow\quad
\begin{aligned}
x^2+2x+5= & 5-x; \\
x^2+3x= 0;\ & \ x(x+3)= 0 \\
\text{Abscises of points:} & \ x_1=-3, \ x_2=0.
\end{aligned}\]

\textbf{Step 3.} Analytical description of the rotating region in Cartesian coordinates:
$$\mathfrak{D}=\left\lbrace -3\leq x\leq 0; \ x^2+2x+5 \leq y \leq 5-x \right\rbrace$$

The body generated by revolution of the region $\mathfrak{D}$ around OX is attached below, on the right.

\textbf{Step 4.} Applying the suitable formula to find the volume of the body.

\parbox[b][110mm][t]{75mm}{
\begin{align*}
V_{OX}= & V_2-V_1= \pi\int^b_a \big( y^2_2(x)-y^2_1(x)\big) dx= \\
= & \pi \int^0_{-3}\big((5-x)^2-(x^2+2x+5)^2\big)dx= \\
= & \pi\int^0_{-3}(-x^4 -4x^3-13x^2-30x ) dx= \\
=&\pi\left[-\dfrac{x^5}{5} -x^4-13\dfrac{x^3}{3}-30\dfrac{x^2}{2} \right]\Bigg|^0_{-3}= \\
= & \pi \left[\dfrac{(-3)^5}{5} +(-3)^4+13\dfrac{(-3)^3}{3}+15\cdot(-3)^2 \right]= \\
= & \pi \left[-\dfrac{243}{5}+99\right]=\dfrac{252}{5}\pi=50,4\pi.
\end{align*}
}
\hfill
\parbox[b][110mm][t]{60mm}{\vspace{1.2cm}\includegraphics[scale=0.8]{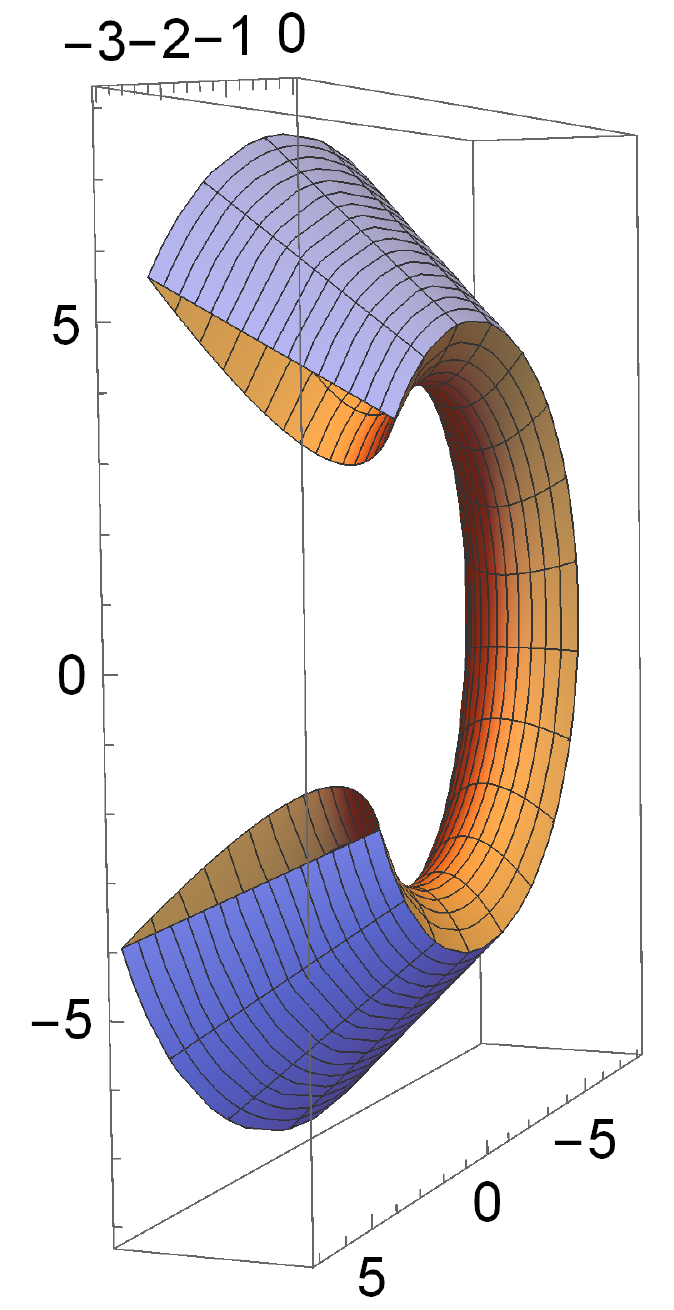}}

\textbf{Result:} $V_{OX}=50,4\pi$ (cubic units).

\item Find the volume of the body formed by rotating the region between the curves
$$x=5+4y-y^2, \ x=5, \ \text{around the axis}\ l=OY.$$
Algorithm:

\textbf{Step 1.} Building the picture.

$x=5+4y-y^2=9-(y-2)^2$ is a $\supset$-shaped horizontal parabola with vertex $(9,2)$ and points of intersection with $OY$: \ $y_1=-1$ and $y_2=5.$

$x=5$ is a straight line passing through $(5,0)$ and parallel to OY.

\parbox[b][10cm][t]{80mm}{\includegraphics[scale=0.8]{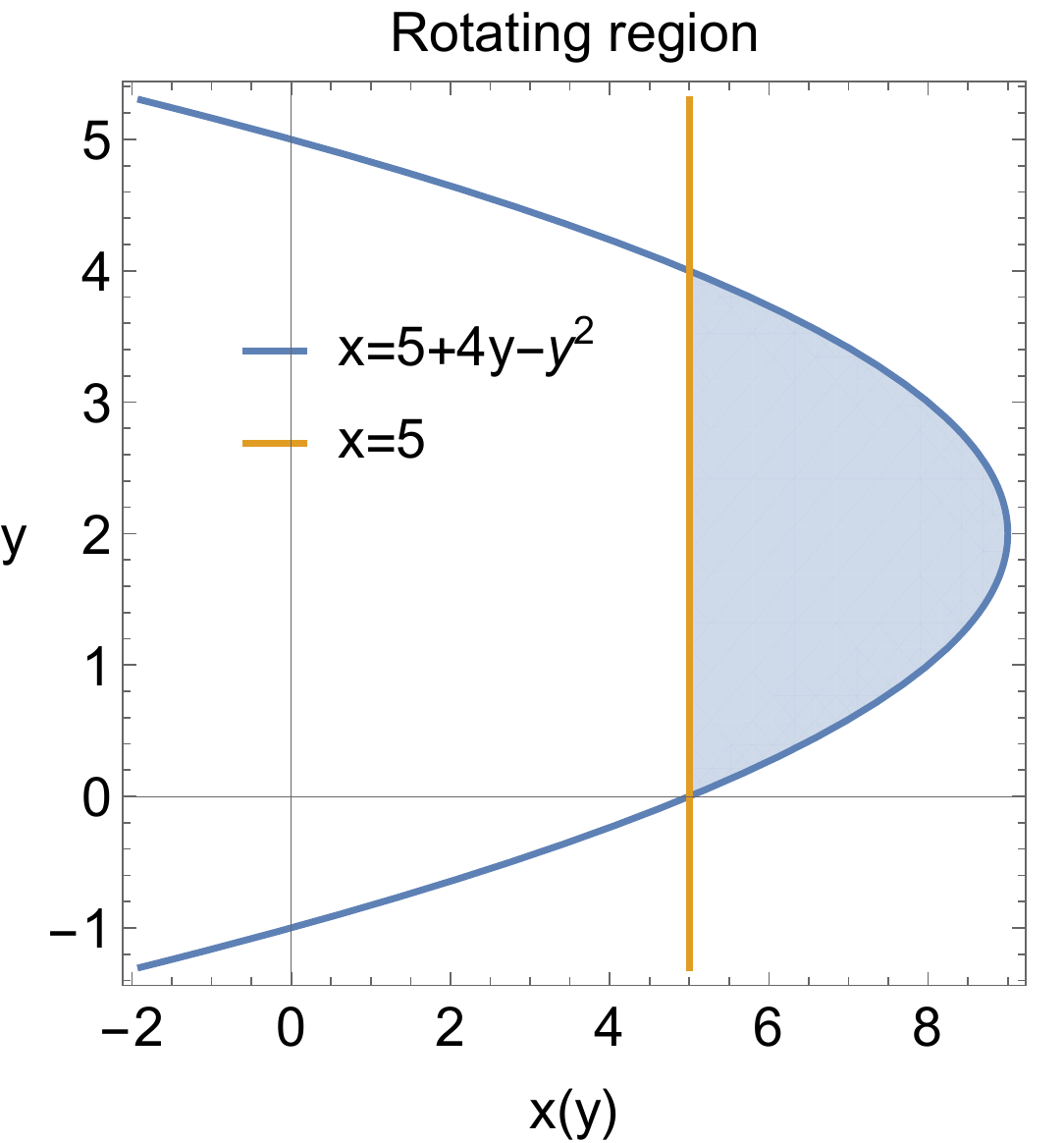}}
\hspace{0.5cm}
\parbox[b][10cm][t]{75mm}{\textbf{Step 2.} Finding the points of intersection between the curves.
$$\left\lbrace\begin{array}{c}
x=5+4y-y^2 \\
x=5
\end{array}\right.\Rightarrow
$$
$$
\begin{aligned}
5+4y-y^2= & 5; -y^2+4y= 0; \\
y(y-4)= & 0 \\
\text{Ordinates of points:} \ & y_1=0, \ y_2=4.
\end{aligned}
$$

\textbf{Step 3.} Analytical description of the rotating region in Cartesian coordinates:
$$\mathfrak{D}=\left\lbrace \begin{array}{c}
                              0\leq y\leq 4; \\
                              5\leq x\leq 5+4y-y^2
\end{array}\right\rbrace.$$}

\parbox[b][4.5cm][t]{30mm}{The body generated by the revolution of the region $\mathfrak{D}$ around OY is attached on the left.}
\hspace{0.5cm}
\parbox[b][4.5cm][t]{120mm}{\vspace{-8mm} \hspace{-0.5cm}
\includegraphics[scale=0.78]{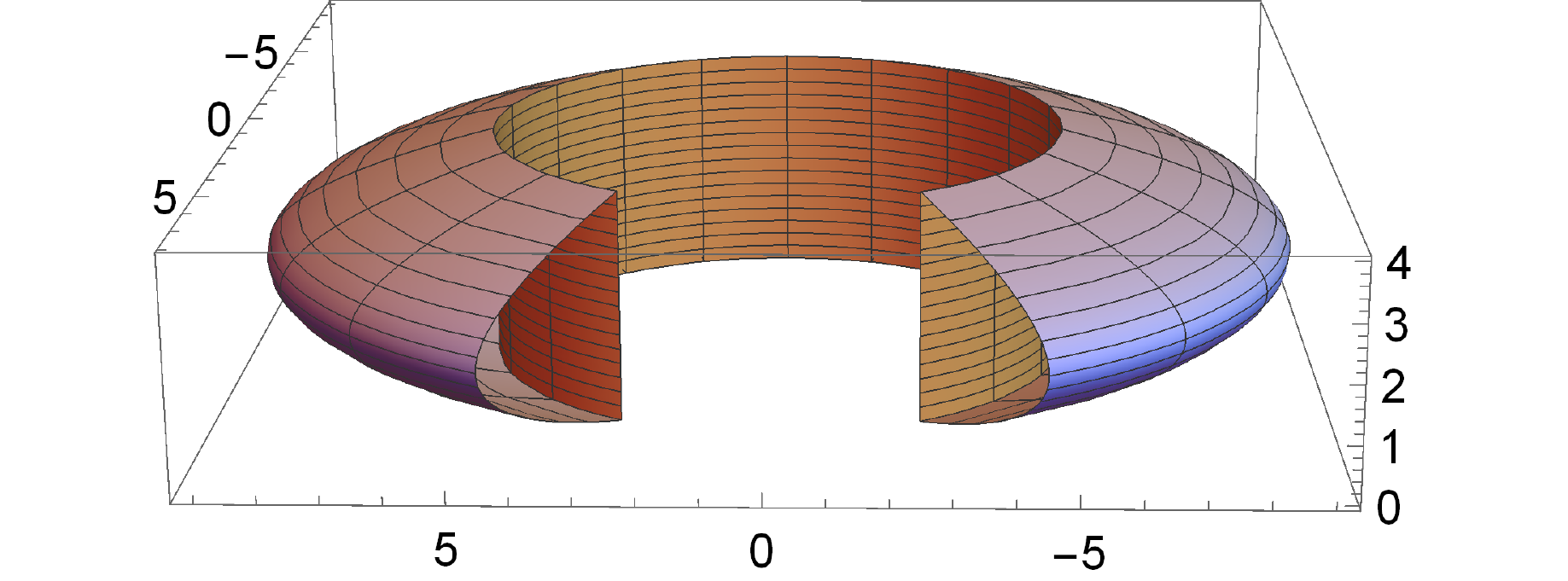}}

\textbf{Step 4.} Applying the suitable formula to find the volume of the body.
\vspace{-0.4cm}
\begin{align*}
V_{OY}= & V_2-V_1= \pi\int^d_c \big(x^2_2(y)-x^2_1(y)\big)dy= \pi \int^4_{0}\big((5+4y-y^2)^2-5^2\big)dy= \\
= & \pi \int^4_{0}(y^4-8y^3+6y^2+40y)dy= \pi \left[\frac{1}{5}y^5-2y^4+2y^3+20y^2\right]\Bigg|^4_0= \\
=&\pi \left[\frac{1}{5}4^5-2\cdot4^4+2\cdot 4^3+20\cdot 4^2\right]=\pi\cdot 4^3\left[\frac{16}{5}-1\right]= \frac{704}{5}\pi=140,8\pi.
\end{align*}
\vspace{-0.4cm}
\textbf{Result:} $V_{OY}=140,8\pi$ (cubic units).

\item Find the volume of the body formed by rotating the curve
$$\rho=6(1+\cos\phi), \ \text{around the axis}\  l=o\rho.$$

\begin{tabular}{|p{6.0cm}|p{7.5cm}|p{2.0cm}|}
\hline
\vspace{0.05mm}$^*$ Solution of Problem 11(c) guided by Irina Blazhievska is available on-line: &
\vspace{5.5mm} \url{https://youtu.be/tr2pXX5GJJE} & \vspace{-3mm} \includegraphics[height=20mm]{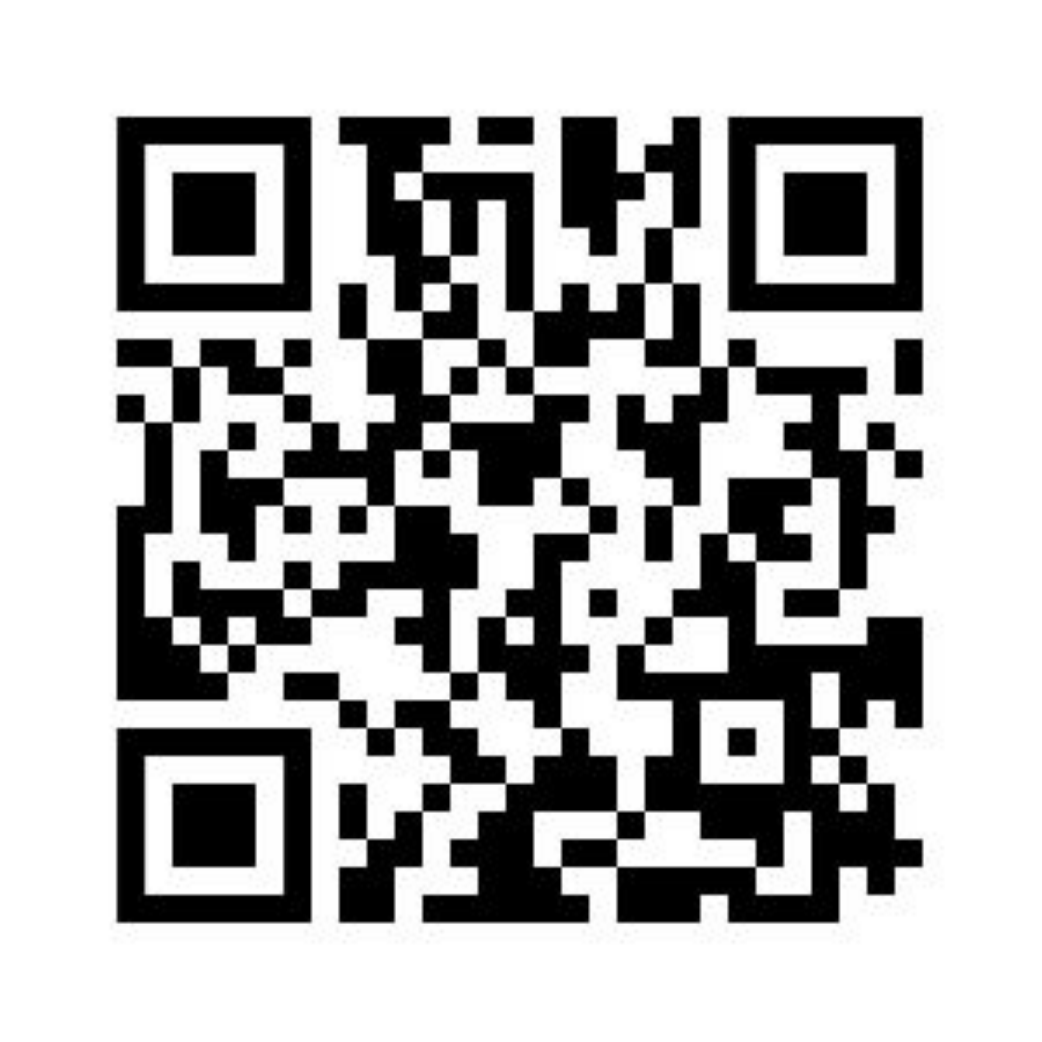}\\
\hline
\end{tabular}

The written version of solution is proposed below. Algorithm:

\textbf{Step 1.} Building the picture.

$\rho=6(1+\cos\phi)$ is a polar representation of a cardioid with a parameter $a=6$. This curve has a heart's shape with a "stalk" at the pole and $o\rho$-symmetry, with angle domain $0\leq\phi\leq2\pi$. Since it is constructed from a cos-function, it has OX-symmetry, which implies the restriction on polar angle domain for rotating part: $\phi \in [0,\pi].$

\textbf{Step 2.} Analytical description of the rotating region in Polar coordinates:
$$\mathfrak{D}=\left\lbrace 0\leq \phi \leq \pi, \ 0\leq \rho\leq 6(1+\cos \phi)\right\rbrace$$
\includegraphics[scale=0.75]{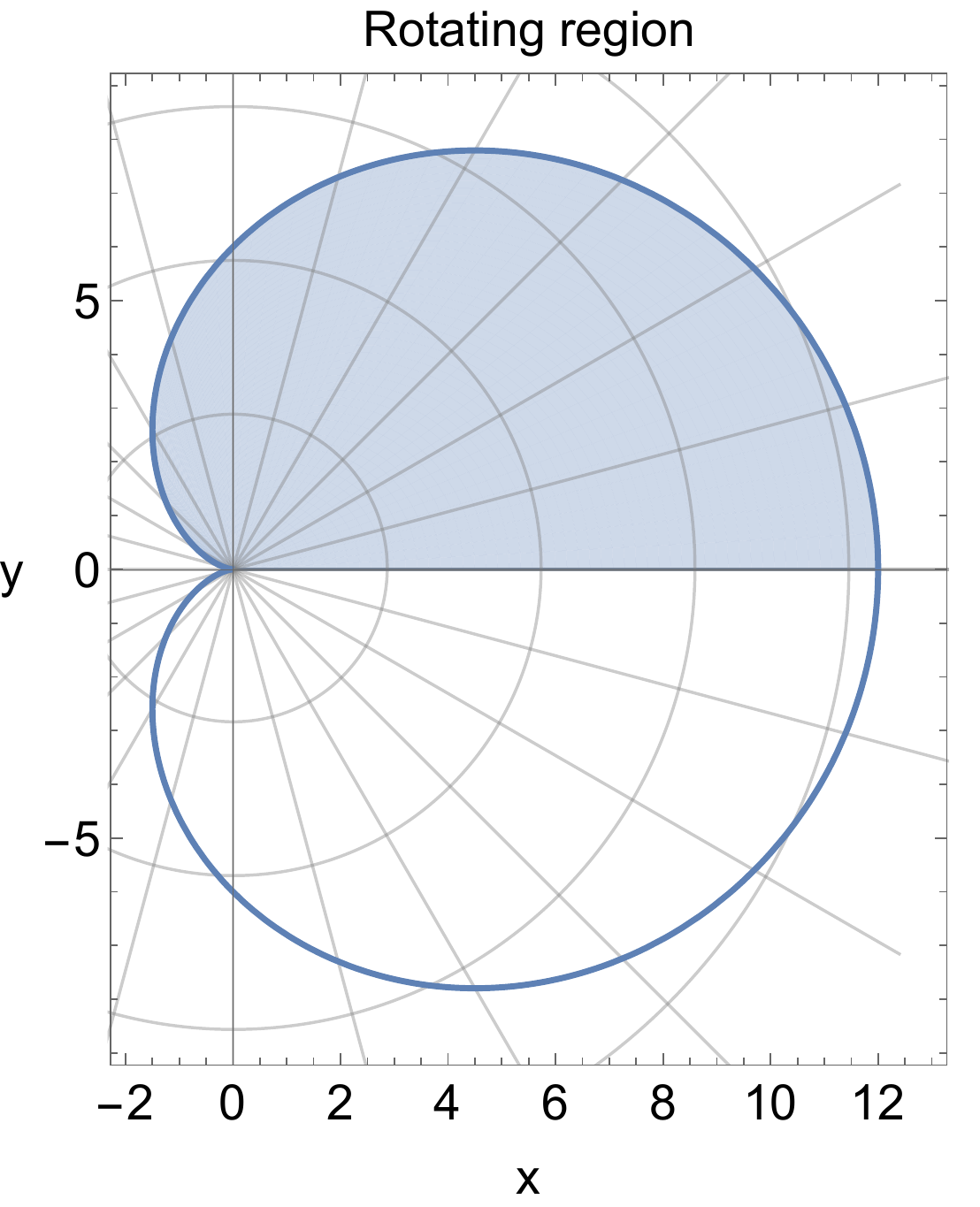}\raisebox{0.1\height}{\includegraphics[scale=0.73]{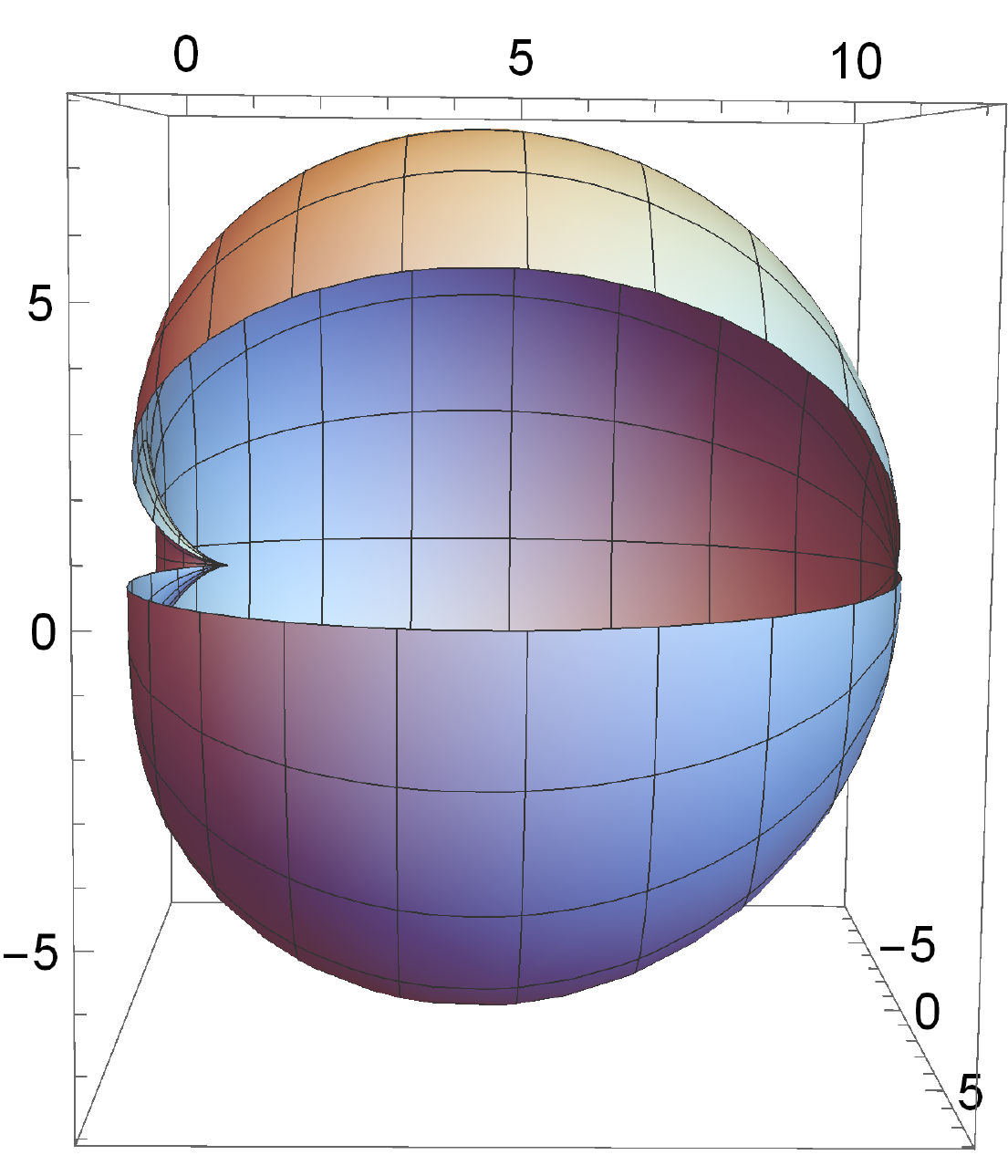}}

The body generated by revolution of the region $\mathfrak{D}$ around the $o\rho$-axis has a shape of an apple with the stalk at the origin; it is attached on the right.

\textbf{Step 3.} Applying the suitable formula to find the volume of the body.
\vspace{-3mm}
\begin{align*}
V_{o\rho}= & \frac{2\pi}{3}\int^\beta_\alpha \rho^3(\phi)\sin \phi \ d\phi= \frac{2\pi}{3}\int^\pi_0 (6(1+\cos \phi))^3\sin \phi \ d\phi=\\
=& -\frac{2\pi}{3}6^3\int^\pi_0 (1+\cos \phi)^3 d(1+\cos \phi)=-6^2\pi \frac{1}{4}\left[(1+\cos \phi)^4\right]\Bigg|^\pi_0= \\
=&-36\pi\left[(1+\cos \pi)^4-(1+\cos 0)^4 \right]=-36\pi \left[0-2^4\right]=476\pi.
\end{align*}
\vspace{-4mm}
\textbf{Result:} $V_{o\rho}=476\pi$ (cubic units).
\end{enumerate}
\vspace{0.1cm}

\begingroup
\addtolength{\jot}{0.15cm}

\begin{enumerate}[label=\textbf{12.(\alph*)}]
\item
$\begin{aligned}[t]
\dot{I}= & \int_{-2}^{\infty}\dfrac{2x-1}{x^2+4}dx=
\left| \begin{array}{c}
f(x)=\dfrac{2x-1}{x^2+4} \textit{ is a bounded continuos} \\
\textit{function for all}\  x\in [-2, +\infty) \\
\Rightarrow \textit{Improper integral of 1st kind}
\end{array} \right|= \\
= & \lim_{N \to \infty}\int^N_{-2}\frac{2x-1}{x^2+4}dx= \lim_{N \to \infty}\int^N_{-2} \left( \frac{2x}{x^2+4}-\frac{1}{x^2+2^2}\right)dx=\\
=&\lim_{N\to \infty}\left(\int^N_{-2}\frac{d(x^2+4)}{x^2+4}-\int^N_{-2}\frac{dx}{x^2+2^2}\right)= \lim_{N\to \infty}\left[ \ln(x^2+4)-\frac{1}{2}\arctg \frac{x}{2}\right]\Bigg|^N_{-2}= \\
= & \lim_{N\to \infty}\left[ \ln(N^2+4)-\frac{1}{2}\arctg\left(\frac{N}{2}\right)-\ln((-2)^2+4)+\frac{1}{2}\arctg\left(\frac{-2}{2}\right)\right]= \\
= & \lim_{N\to \infty}\left[ \ln(N^2+4)-\frac{1}{2}\arctg \left(\frac{N}{2}\right)-\ln8-\frac{\pi}{8}\right]=+\infty-\frac{3\pi}{8}-\ln8=+\infty.
\end{aligned}$

\textbf{Result:} The integral is divergent: $\dot{I}=+\infty.$

\begin{tabular}{|p{6.0cm}|p{7.5cm}|p{2.0cm}|}
\hline
\vspace{0.05mm}$^*$ Solution of Problem 12(a) guided by Irina Blazhievska is available on-line: &
\vspace{5.5mm} \url{https://youtu.be/r1GZaS5Hooo} & \vspace{-3mm} \includegraphics[height=20mm]{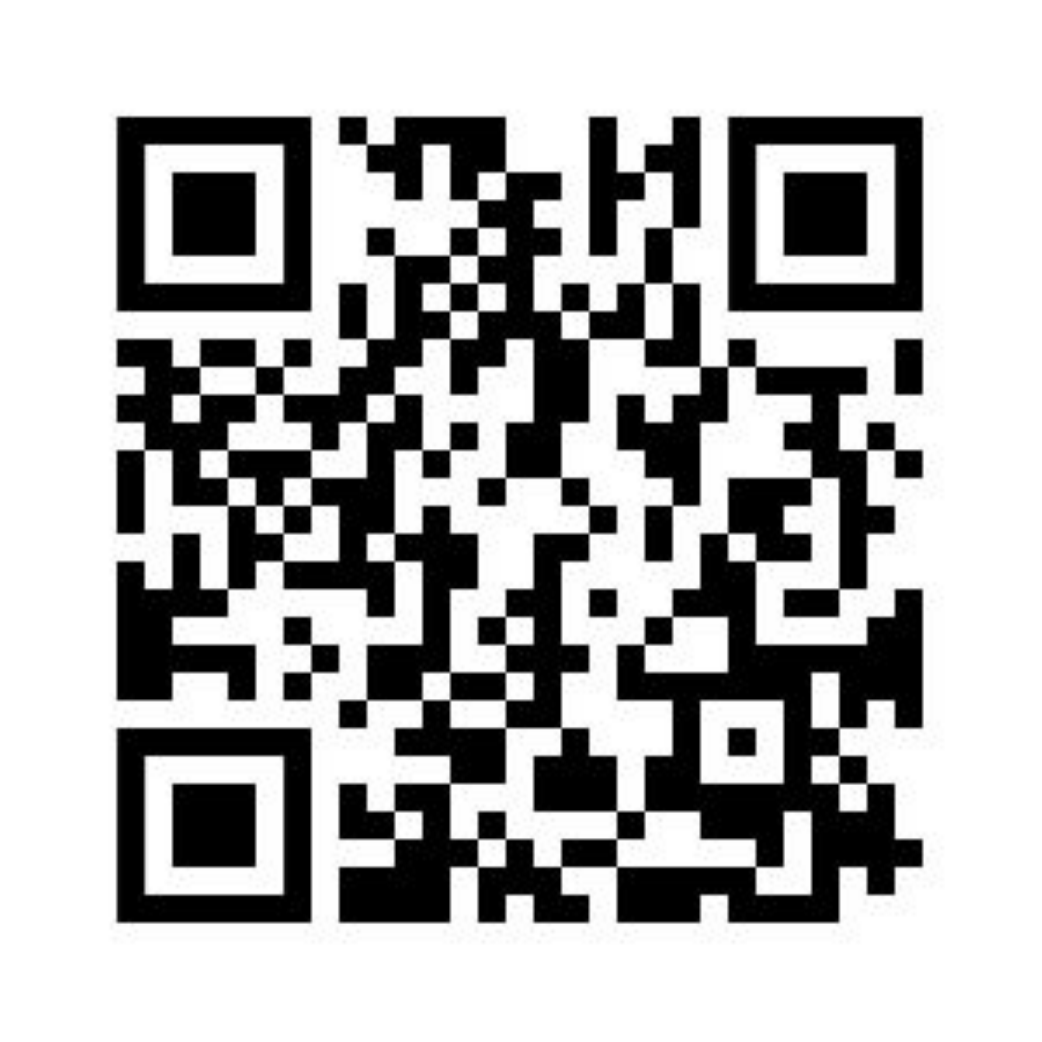}\\
\hline
\end{tabular}

\item
$\begin{aligned}[t]
\dot{I}= & \Int_{0}^{\pi/2}e^{-\tg x}\dfrac{dx}{\cos^2 x}= \left| \begin{array}{c}
f(x)=\dfrac{e^{-\tg x}}{\cos^2 x} \textit{ is a bounded continous} \\
\textit{function for all} \ x\in [0,\frac{\pi}{2});\\
\textit{The endpoint} \ x=\frac{\pi}{2} \textit{ is a singular point of } f(x) \\
\Rightarrow \textit{Improper integral of 2nd kind}
\end{array} \right|= \\
= & \lim_{\varepsilon\to 0^+}\int^{\pi/2-\varepsilon}_0 e^{-\tg x}\frac{dx}{\cos^2 x}=\lim_{\varepsilon\to 0^+}\int^{\pi/2-\varepsilon}_0 e^{-\tg x} d(\tg x)=-\lim_{\varepsilon\to 0^+} e^{-\tg x}\Big|^{\pi/2-\varepsilon}_0= \\
=& -\lim_{\varepsilon\to 0^+} \left[e^{-\tg (\pi/2-\varepsilon)}-e^{\tg 0}\right]= 1-\lim_{\varepsilon\to 0^+} e^{-\ctg \varepsilon}=1-0=1.
\end{aligned}$

\textbf{Result:} The integral is convergent: $\dot{I}=1.$

\begin{tabular}{|p{6.0cm}|p{7.5cm}|p{2.0cm}|}
\hline
\vspace{0.05mm}$^*$ Solution of Problem 12(b) guided by Irina Blazhievska is available on-line: &
\vspace{5.5mm} \url{https://youtu.be/lbBT5mV_pc4} & \vspace{-3mm} \includegraphics[height=20mm]{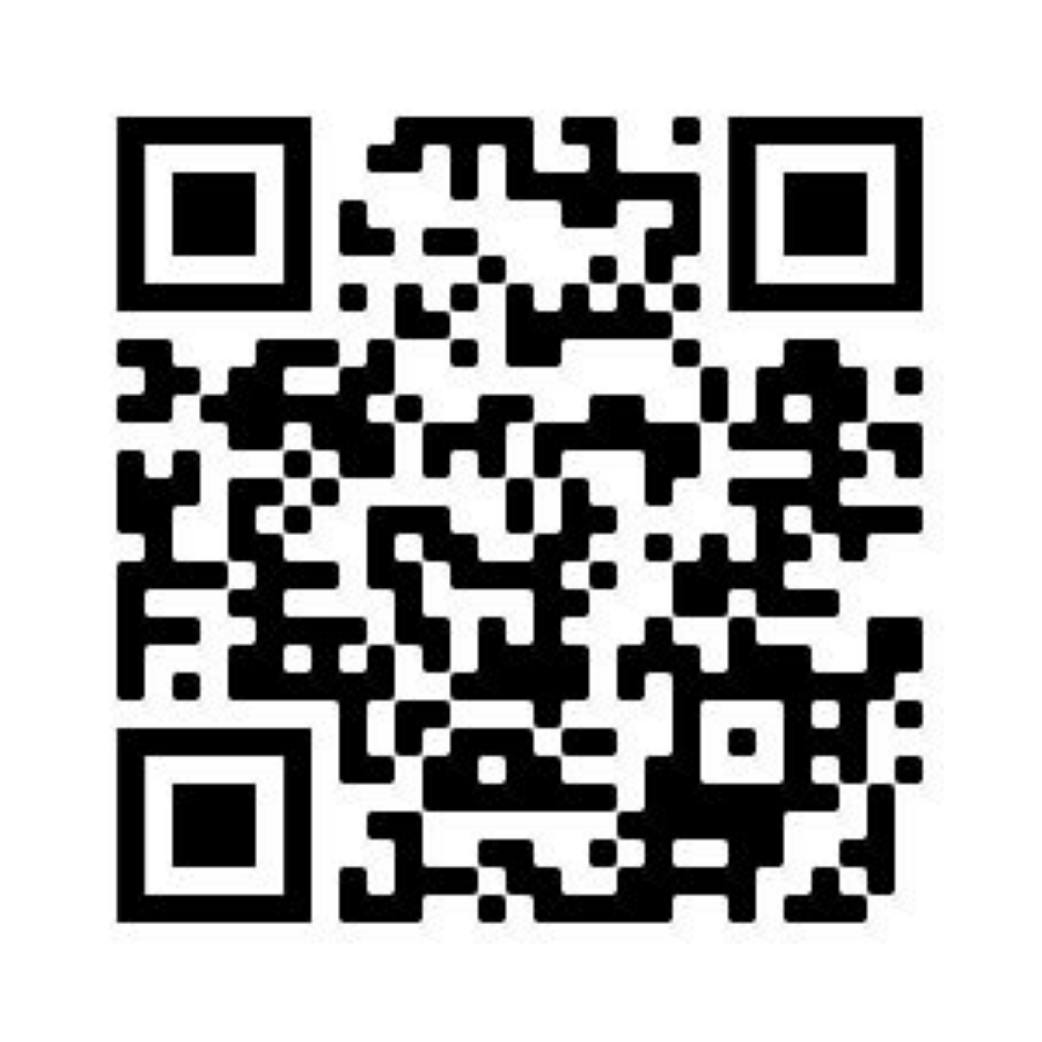}\\
\hline
\end{tabular}
\end{enumerate}
\endgroup

\newpage
\section{Author's Video-Lessons on Techniques of Integration}
\begin{center}
(based on Sample task, page 39)
\end{center}
\vspace{0.5cm}

\begin{tabular}{|p{6.0cm}|p{7.5cm}|p{2.0cm}|}
\hline
\multicolumn{3}{|c|}{\textbf{Indefinite Integration}}\\
\hline
\vspace{0.05mm} Reduction to the table & & \\ of integrals: $$\int x^3(1-x^2)^2 dx$$ Guide: Ricard Riba Garcia &
\vspace{5.5mm} https://youtu.be/x48CikKlF9c & \vspace{-3mm} \includegraphics[height=20mm]{1-a.pdf}\\
\hline
\vspace{0.05mm} Reduction to the table & & \\ of integrals: $$\int\big(3+\textrm{tg}^2x\big)dx$$ Guide: Irina Blazhievska &
\vspace{5.5mm} https://youtu.be/1MqmZbQ-3qM & \vspace{-3mm} \includegraphics[height=20mm]{1-b.pdf}\\
\hline
\vspace{0.05mm} Substitution under the & & \\ differential: $$\int x^4 e^{2x^5-1} dx$$ Guide: Irina Blazhievska &
\vspace{5.5mm} https://youtu.be/kSn2UvdXWVs & \vspace{-3mm} \includegraphics[height=20mm]{1-c.pdf}\\
\hline
\vspace{0.05mm} Integration of fractions & & \\ with quadratic functions: $$\int\dfrac{dx}{\sqrt{6x-9x^2}}$$ Guide: Irina Blazhievska &
\vspace{5.5mm}
   https://youtu.be/KQuEzkh6AqQ & \vspace{-3mm} \includegraphics[height=20mm]{2-b.pdf}\\
\hline
\vspace{0.05mm} Integration of fractions & & \\ with quadratic functions: $$\int\dfrac{(-2x-7)dx}{x^2+6x+10}$$ Guide: Ricard Riba Garcia &
\vspace{5.5mm}
   https://youtu.be/RXJcMkNz8Zg & \vspace{-3mm} \includegraphics[height=20mm]{2-c.pdf}\\
\hline
\end{tabular}

\begin{tabular}{|p{6.0cm}|p{7.5cm}|p{2.0cm}|}
\hline
\multicolumn{3}{|c|}{\textbf{Indefinite Integration}}\\
\hline
\vspace{0.05mm} Change of variable: $$\int\dfrac{(3x+2)dx}{\sqrt{x+4}}$$ Guide: Irina Blazhievska &
\vspace{10mm}
   https://youtu.be/ZOFSo2oDVzQ & \vspace{3mm} \includegraphics[height=20mm]{3-a.pdf}\\
\hline
\vspace{0.05mm} Integration by parts: $$\int x^2\ln x dx$$ Guide: Ricard Riba Garcia &
\vspace{10mm}
https://youtu.be/hDYt-m7ZCgM & \vspace{3mm} \includegraphics[height=20mm]{3-d.pdf}\\
\hline
\vspace{0.05mm} Integration of polynomial & & \\ fractions: $$\int \dfrac{(3x^3-32x+56)dx}{x^3-2x^2-4x+8}$$ Guide: Irina Blazhievska &
\vspace{5.5mm} https://youtu.be/X\_300OO1i8A & \vspace{-3mm} \includegraphics[height=20mm]{4-b.pdf}\\
\hline
\vspace{0.05mm} Integration of trigonometric & & \\ functions (products): $$\int \sin 10x \sin3x dx$$ Guide: Irina Blazhievska &
\vspace{5.5mm}  https://youtu.be/-Zsc5t-YIOk & \vspace{-3mm} \includegraphics[height=20mm]{5-a.pdf}\\
\hline
\vspace{0.05mm} Integration of trigonometric & & \\ rational function (fractions): $$\int \dfrac{dx}{4\cos x+3\sin x+6}$$ Guide: Irina Blazhievska &
\vspace{5.5mm} https://youtu.be/loxi3dwTmho & \vspace{-3mm} \includegraphics[height=20mm]{5-c.pdf}\\
\hline
\vspace{0.05mm} Integration of fractions & & \\ with radicals: $$\int\dfrac{dx}{x^2\sqrt{4-x^2}}$$ Guide: Ricard Riba Garcia &
\vspace{5.5mm}
https://youtu.be/ba4kyucyLVk & \vspace{-3mm} \includegraphics[height=20mm]{6-a.pdf}\\
\hline
\end{tabular}

\begin{tabular}{|p{6.0cm}|p{7.5cm}|p{2.0cm}|}
\hline
\multicolumn{3}{|c|}{\textbf{Definite Integration}}\\
\hline
\vspace{0.05mm} Integration by parts: $$\int_{0}^{1}\ln(1+x^2)dx$$ Guide: Irina Blazhievska &
\vspace{15mm} https://youtu.be/jfbI2G23U2M & \vspace{7mm} \includegraphics[height=20mm]{7-a.pdf}\\
\hline
\vspace{0.05mm} Substitution under the & & \\ differential: $$\int_{0}^{\pi/2}\sin^3x\sqrt[4]{\cos x}dx$$ Guide: Irina Blazhievska &
\vspace{10mm} https://youtu.be/s4VH2LvXh7M & \vspace{3mm} \includegraphics[height=20mm]{7-b.pdf}\\
\hline
\vspace{0.05mm} Change of variable: $$\int_1^{16}\dfrac{(1+\sqrt{x})dx}{\sqrt[4]{x}+\sqrt{x}}$$ Guide: Irina Blazhievska &
\vspace{15mm} https://youtu.be/gCX8MgQrd7A & \vspace{7mm} \includegraphics[height=20mm]{7-d.pdf}\\
\hline
\end{tabular}

\vspace{0.5cm}
\begin{tabular}{|p{6.0cm}|p{7.5cm}|p{2.0cm}|}
\hline
\multicolumn{3}{|c|}{\textbf{Geometric Applications of Definite Integrals}}\\
\hline
\vspace{0.05mm} $\bullet$ Cartesian coordinates: & & \\ Area of the figure bounded & & \\ by the curves: $$y=2x^2-10x+6, \ y=x^2-3x$$ Guide: Irina Blazhievska &
\vspace{1mm} https://youtu.be/w4wx7Dd547w & \vspace{-7mm} \includegraphics[height=20mm]{8-a.pdf}\\
\hline
\vspace{0.05mm} $\bullet$ Polar coordinates: & & \\ Arc length of the curve: $$\rho(\phi)=\dfrac{10}{\sqrt{101}}e^{\frac{\phi}{10}}, \ \phi\in[0,2\pi]$$ Guide: Irina Blazhievska &
\vspace{5.5mm} https://youtu.be/AJbHsr1eOes & \vspace{-3mm} \includegraphics[height=20mm]{9-c.pdf}\\
\hline
\end{tabular}

\begin{tabular}{|p{6.0cm}|p{7.5cm}|p{2.0cm}|}
\hline
\multicolumn{3}{|c|}{\textbf{Geometric Applications of Definite Integrals}}\\
\hline
\vspace{0.05mm} $\bullet$ Parametric coordinates: & & \\ Surface area generated by & & \\ rotating the curve around & & \\ OY-axis: $$x=3+\cos t,\ y=2+\sin t$$ Guide: Irina Blazhievska  &
\vspace{-3mm} https://youtu.be/nXzXyHIw0w8 & \vspace{-10mm} \includegraphics[height=20mm]{10-b.pdf}\\
\hline
\vspace{0.05mm} $\bullet$ Cartesian coordinates: & & \\ Volume of the body & & \\ generated by rotating & & \\ around OX-axis the region & & \\ bounded by the curve: $$\rho(\phi)=6(1+\cos\phi)$$ Guide: Irina Blazhievska  &
\vspace{-5mm} https://youtu.be/tr2pXX5GJJE & \vspace{-12mm} \includegraphics[height=20mm]{11-c.pdf}\\
\hline
\end{tabular}

\vspace{0.5cm}
\begin{tabular}{|p{6.0cm}|p{7.5cm}|p{2.0cm}|}
\hline
\multicolumn{3}{|c|}{\textbf{Improper Integrals}}\\
\hline
\vspace{0.05mm} $\bullet$ Improper integral & & \\ of 1st kind: $$\int_{-2}^{\infty}\dfrac{(2x-1)dx}{x^2+4}$$ Guide: Irina Blazhievska &
\vspace{7mm} https://youtu.be/r1GZaS5Hooo & \vspace{-1mm} \includegraphics[height=20mm]{12-a.pdf}\\
\hline
\vspace{0.05mm} $\bullet$ Improper integral & & \\ of 2nd kind:  $$\int_{0}^{\pi/2}e^{-\textrm{tg} x}\dfrac{dx}{\cos^2x}$$ Guide: Irina Blazhievska&
\vspace{7mm} https://youtu.be/lbBT5mV\_pc4 & \vspace{-1mm} \includegraphics[height=20mm]{12-b.pdf}\\
\hline
\end{tabular}

\newpage

\pagestyle{empty}

\vspace{-0.9cm}
\section*{External links}
Not only well-organized libraries of Math study, but also the interactive graphics and infinity samples are proposed on the following cites:
\begin{itemize}
  \item http://www.mathcentre.ac.uk/students/courses/
  \item https://www.mathcurve.com
  \item http://www.mathematische-basteleien.de/index.htm
  \item http://mathworld.wolfram.com
  \item https://math.stackexchange.com
  \item http://old.nationalcurvebank.org/volrev/volrev.htm
 \end{itemize}

\end{document}